\documentclass[preprint,11pt]{article}

\usepackage{appendix}

\usepackage[nottoc,notlot,notlof]{tocbibind}

\usepackage{subfiles}

\usepackage{indentfirst}
\usepackage{setspace}
\usepackage{enumitem}
\usepackage[left=1.00in,top=1.00in,right=1.00in,bottom=1.00in]{geometry}
\usepackage{mdframed}

\usepackage{float}
\usepackage{graphicx}
\usepackage{caption}
\usepackage{subcaption}
\usepackage{array}

\usepackage{amsmath}
\usepackage{amsfonts}
\usepackage{bm}
\usepackage{bbm}
\usepackage{mathrsfs}
\usepackage[cal=boondox]{mathalfa}
\usepackage{color}

\usepackage{mathtools}
\mathtoolsset{showonlyrefs}

\usepackage[bbgreekl]{mathbbol}
\usepackage{amssymb,amsthm,mathrsfs,tikz,tikzscale}

\title{Mesh Quality Metrics for Isogeometric Bernstein--B\'{e}zier Discretizations}
\author{Luke Engvall, John A. Evans}
\begin{document}


\newtheorem{thm}{Theorem}[section]
\newtheorem{cor}[thm]{Corollary}
\newtheorem{lem}[thm]{Lemma}
\newtheorem{pro}[thm]{Proposition}

\newtheorem{innercustomthm}{Theorem}
\newenvironment{customthm}[1]
  {\renewcommand\theinnercustomthm{#1}\innercustomthm}
  {\endinnercustomthm}
  
\theoremstyle{definition}
\newtheorem{definition}{Definition}[section]
\newtheorem{defi}[thm]{Definition}
  \newtheorem{innerthmproof}{Proof of Theorem}
\newenvironment{thmproof}[1]
  {\renewcommand\theinnerthmproof{#1}\innerthmproof}
  {\endinnerthmproof}

\newcommand{\fgrad}        {\mbox{{$\boldsymbol{\nabla}\mathbf{F}$}}}
\newcommand{\wgrad}      {\mbox{{$\boldsymbol{\nabla}w$}}}
\newcommand{\fgradn}[1]  {\mbox{{$\boldsymbol{\nabla}^{#1}\mathbf{F}$}}}
\newcommand{\wgradn}[1]{\mbox{{$\boldsymbol{\nabla}^{#1}w$}}}
\newcommand{\fdet}          {\mbox{{$\mathrm{det}\boldsymbol{\nabla}\mathbf{F}$}}}
\newcommand{\bnabla}          {\mbox{{$\boldsymbol{\nabla}$}}}

\newcommand{\deriv}{D}
\newcommand{\derivn}[1]{D^{#1}}
\newcommand{\derivx}[1]{D_{#1}}
\newcommand{\derivnx}[2]{D_{#2}^{#1}}

\newcommand{\deldel}[2]{\dfrac{\partial #1 }{\partial #2 } }

\newcommand{\linf}[2]{{\left|\left|#1\right|\right|_{L^\infty\left(#2\right)}}}
\newcommand{\ltwo}[2]{{\left|\left|#1\right|\right|_{L^2\left(#2\right)}}}

\newcommand{\bvec}[1]{\mathbf{#1}}
\newcommand{\bi}{\mathbf{i}}
\newcommand{\bj}{\mathbf{j}}
\newcommand{\bk}{\mathbf{k}}
\newcommand{\br}{\mathbf{r}}
\newcommand{\bs}{\mathbf{s}}
\newcommand{\bn}{\mathbf{n}}
\newcommand{\bp}{\mathbf{p}}
\newcommand{\be}{\mathbf{e}}

\newcommand{\xprojv}{\mathbf{v}}
\newcommand{\bv}{\bm{v}}

\newcommand{\bbi}{\mathbb{i}}
\newcommand{\bbI}{\mathbb{I}}
\newcommand{\bbk}{\mathbb{k}}

\newcommand{\bvi}{\mathbf{I}}

\newcommand{\bxu}{\mathbf{x}}
\newcommand{\bx}{\bm{x}}
\newcommand{\bX}{\widetilde{\mathbf{X}}}
\newcommand{\bP}{\mathbf{P}}
\newcommand{\bF}{\mathbf{F}}
\newcommand{\bN}{\mathbf{N}}
\newcommand{\bT}{\mathbf{T}}
\newcommand{\bJ}{\mathbf{J}}
\newcommand{\bxi}{{{\boldsymbol{\xi}}}}
\newcommand{\bphi}{{{\boldsymbol{\phi}}}}
\newcommand{\blambda}{{\boldsymbol{\lambda}}}
\newcommand{\balpha}{{\boldsymbol{\alpha}}}

\newcommand{\Lnorm}{\left| \left|} 
\newcommand{\Rnorm}{\right| \right|} 
\newcommand{\of}[1]{\left( #1 \right)} 
\newcommand{\normof}[1]{\left| \left| #1 \right| \right|} 
\newcommand{\absof}[1]{ \left| #1 \right|} 
\newcommand{\matof}[1]{ \left[ #1 \right]} 
\newcommand{\setof}[1]{ \left\{ #1 \right\} } 

\newcommand{\Ophys}{\Omega_e} 
\newcommand{\Oref}{\hat{\Omega}}
\newcommand{\Olin}{\overline{\Omega}_e}
\newcommand{\Oproj}{\widetilde{\Omega}_e}

\newcommand{\xphys}{{\bxu_e}}
\newcommand{\xproj}{\widetilde{\bxu}_e}
\newcommand{\xprojd}{\big{[}\widetilde{\bxu}_e\big{]}_d}
\newcommand{\xlin}{\overline{\bxu}_e}

\newcommand{\xvphys}{{\bx}}
\newcommand{\xvproj}{\widetilde{\bx}}
\newcommand{\xvlin}{\overline{\bx}}

\newcommand{\Rproj}{\mathbb{R}^{d+1}}
\newcommand{\Rphys}{\mathbb{R}^{d}}

\newcommand{\BB}{Bernstein--B\'{e}zier }
\newcommand{\boldit}[1]{\textbf{\emph{#1}}}


\date{}
\maketitle

\begin{abstract}
High-order finite element methods harbor the potential to deliver improved accuracy per degree of freedom versus low-order methods.  Their success, however, hinges upon the use of a curvilinear mesh of not only sufficiently high accuracy but also sufficiently high quality.  In this paper, theoretical results are presented quantifying the impact of mesh parameterization on the accuracy of a high-order finite element approximation, and a formal definition of shape regularity is introduced for curvilinear meshes based on these results.  This formal definition of shape regularity in turn inspires a new set of quality metrics for curvilinear finite elements.  Computable bounds are established for these quality metrics using the Bernstein-B\'{e}zier form, and a new curvilinear mesh optimization procedure is proposed based on these bounds.  Numerical results confirming the importance of shape regularity in the context of high-order finite element methods are presented, and numerical results demonstrating the promise of the proposed curvilinear mesh optimization procedure are also provided.  The theoretical results in this paper apply to any piecewise-polynomial or piecewise-rational finite element method posed on a mesh of polynomial or rational mapped simplices and hypercubes.  As such, they apply not only to classical continuous Galerkin finite element methods but also to discontinuous Galerkin finite element methods and even isogeometric methods based on NURBS, T-splines, or hierarchical B-splines.
\end{abstract}

\tableofcontents

\section{Introduction}
High-order finite element methods have risen in popularity in recent years due to their potential to deliver improved accuracy per degree of freedom versus classical low-order methods.  However, in order to deliver optimal convergence rates, high-order finite element methods must be posed on suitable curvilinear meshes of sufficiently high accuracy \cite{sherwin_mesh_2002}.  With this in mind, the goals of this work are (i) to quantify the impact of mesh parameterization on the accuracy of a high-order finite element approximation and (ii) to develop a set of sufficient and computable conditions which guarantee a family of refined but not necessarily nested high-order finite element approximations exhibit optimal convergence rates with respect to the mesh size.  To achieve these goals, we establish a formal definition of shape regularity for curvilinear meshes using a Bramble-Hilbert lemma for polynomial and rational approximations on mapped simplices and hypercubes \cite{bramble_estimation_1970}.  Inspired by this new definition of shape regularity, we introduce a set of quality metrics for curvilinear finite elements, and we establish computable bounds for these quality metrics based on the Bernstein-B\'{e}zier form for curvilinear finite elements defined through a polynomial or rational parametric mapping \cite{prautzsch2013bezier}.

The results presented in this paper apply to any piecewise-polynomial or piecewise-rational finite element method posed on a mesh of polynomial or rational mapped simplices and hypercubes.  As such, the results apply not only to classical continuous Galerkin finite element methods, but also to discontinuous Galerkin finite element methods \cite{arnold2002unified,michoski2016foundations} and even isogeometric analysis (IGA) methods based on Bernstein-B\'{e}zier \cite{engvall_isogeometric_2016,engvall_isogeometric_2017,xia2017isogeometric}, B-spline and Non-Uniform Rational B-spline (NURBS) \cite{hughes2005isogeometric}, Hierarchical B-spline \cite{vuong2011hierarchical}, or T-spline basis functions \cite{bazilevs2010isogeometric}.  It should be mentioned, however, that the results in this paper only pertain to the best approximation properties of a high-order finite element basis and do not take into any method errors that may arise due to the application of a particular finite element method to a problem of interest.

The results presented here are also meant to inform the construction of curvilinear meshes using state-of-the-art mesh generation procedures.  We have taken great care to ensure that the quality metrics presented in this paper may be easily incorporated into existing automated mesh generation and optimization algorithms, and we further propose a simple cost functional for mesh optimization based on these metrics and demonstrate that the functional yields improved curvilinear meshes as compared with classical elasticity-based mesh smoothing techniques.

Naturally, the results here would be of little use if the metrics currently employed in the curvilinear mesh generation community were both cheaper than the ones presented here and also able to be employed to establish sufficient conditions for the generation of high quality curvilinear meshes.  Indeed, there already exist a suite of quality metrics for curvilinear meshes.  The most common of these metrics is the so-called scaled Jacobian metric which is equal to one for a linear finite element mesh and zero for a finite element mesh which contains a element with a non-bijective mapping \cite{dey_curvilinear_1999,persson_curved_2008}.  However, there exist families of refined but non-nested curvilinear meshes with uniformly bounded (from below) scaled Jacobian which exhibit sup-optimal convergence rates with respect to the mesh size.  We present several examples of such families in Section \ref{numerical}.  In fact, as discussed in Subsection \ref{metrics_review}, it is possible to construct a highly skewed finite element with a scaled Jacobian of identically one.  There exist several other quality metrics which expand upon the scaled Jacobian metric \cite{cohen_analysis-aware_2010,escobar_new_2011,gargallo-peiro_distortion_2015,gargallo-peiro_optimization_2015,george_construction_2012,lamata_quality_2013,poya_unified_2016,roca_defining_2011,speleers_optimizing_2015,xie_generation_2013,xu_optimal_2013,xu_high-quality_2014,zhang_solid_2012}, but to the best of our knowledge, none of these metrics are able to fully quantify the impact of mesh parameterization on the accuracy of an arbitrarily high-order finite element approximation.

With the motivation for this paper established, we now present an outline for the remainder of the paper.  We begin by introducing notation  and other preliminaries in Section \ref{preliminaries}. Then, in Section \ref{review}, we provide a review of the relevant literature, and we argue the case that existing mesh quality metrics are not sufficient for quantifying the impact of mesh parameterization on the accuracy of a high-order finite element approximation.  Next, we present interpolation error bounds for rational \BB elements in Section \ref{interp_theory}, and we introduce sufficient conditions for optimal convergence of a high-order finite element approximation in Section \ref{sufficient_conditions}.  With these conditions established,  we present distortion metrics for rational \BB elements in Section \ref{validity_metrics} as well as computable bounds for these metrics.  Finally, we present some numerical results in Section \ref{numerical}, and we provide some concluding remarks and directions for future research in Section \ref{conclusions}.


\section{Notation and Preliminaries}  \label{preliminaries}
Admittedly, the work presented in this paper is notationally intensive.
At the risk of being pedantic, we use this section to briefly introduce notation to be used throughout the remainder of this paper.
We provide a review of multi-index notation (\S \ref{notation}) and derivative notation (\S \ref{deriv_notation}), as well as a review of Bernstein polynomials and B\'{e}zier elements (\S \ref{bern_bez}).
\subsection{Multi-Index Notation} \label{notation}
Throughout this paper, we make heavy use of multi-index notation in order to simplify the equations presented in the sections that follow.
So that the meaning of the equations is unambiguous, we review this notation here.  For a natural number $d$, let $\bn = \{n_1,...,n_d\}$ and $\bk = \{k_1,...,k_d\}$ denote $1 \times d$ multi-indices of non-negative integers and let $\bx$ denote a $1 \times d$ vector of real numbers.  Moreover, let $n$ denote a non-negative integer and let $x$ denote a real number.  We use the following notation to denote common operations on these objects:
\begin{equation*}
|\bn| = \sum\limits_{i = 1}^{d}n_i, \hspace{10pt}  \bn! = \prod\limits_{i=1}^d n_i!, \hspace{10pt}  x^{\bn} =  \prod\limits_{i=1}^d x^{n_i}, \hspace{10pt} \bx^{\bvec{n}} =  \prod\limits_{i=1}^d x_i^{n_i}, \hspace{10pt} {{n}\choose{\bk}} = \dfrac{n!}{\bk!},  \hspace{10pt} {{\bn}\choose{\bk}} = \prod\limits_{i=1}^d{{n_i}\choose{k_i}}.
\end{equation*}
We say that $\bn = \bk$ if all the entries of the multi-indices are equal, $\bn < \bk$ is there exists a $j \in \{1, \ldots, d\}$ such that $n_j < k_j$ and $n_i = k_i$ for $i < j$, and $\bn > \bk$ otherwise.

\subsection{Derivative Notation} \label{deriv_notation}
In order to write derivatives compactly, we use multi-index notation.  In particular, we denote the ${\balpha}^{\textup{th}}$ partial derivative operator with respect to the variables $\bxi$ as: 
\begin{equation*}
  \derivnx{{\balpha}}{\bxi} =
  \dfrac{ \partial^{|{\balpha}|} }
        { \partial \xi_1^{\alpha_1} \partial \xi_2^{\alpha_2} ... \partial \xi_d^{\alpha_d}}.
\end{equation*}
The $\balpha^{\textup{th}}$ partial derivative of a vector--valued function $\bm{f}$ of length $m$ is understood to result in a $m \times 1$ column vector, viz.:
\begin{equation*}
  \derivnx{{\balpha}}{\bxi}\bm{f} =
  \left[ \begin{array}{c}
      \derivnx{{\balpha}}{\bxi}f_1\\
      \vdots\\
      \derivnx{{\balpha}}{\bxi}f_m
    \end{array} \right]
\end{equation*}
and we collect all partial derivatives of order $k = |{\balpha}|$ into the matrix:
\begin{equation*}
  \boldsymbol{\nabla}_{\bxi}^{k}\bm{f} =
  \left[ \begin{array}{cccc}
      \derivnx{{\balpha}_1}{\bxi} \bm{f}   & \derivnx{{\balpha}_2}{\bxi}  \bm{f} & \hdots & \derivnx{{\balpha}_n}{\bxi}  \bm{f} \\
    \end{array} \right]
\end{equation*}
where $\balpha_i < \balpha_j$ for $i < j$.  We note that $ \boldsymbol{\nabla}_{\bxi}^1 = \boldsymbol{\nabla}_{\bxi} $ is the standard gradient operator,
and we take $\boldsymbol{\nabla}_{\bxi}^\bvec{0}$ to be the identity operator.
Additionally, to shorten certain equations, we at times omit the subscript denoting the independent variable.
That is, we write  $ \boldsymbol{\nabla}\bm{f} = \boldsymbol{\nabla}_{\bxi}\bm{f} $ when the choice of the independent variable is unambiguous.

\subsection{\BB Elements} \label{bern_bez}
The element distortion metrics  proposed in this paper make heavy use of rational B\'{e}zier elements and the Bernstein basis polynomials defined on these elements. 
We take this section to briefly review some of the relevant properties of these elements and their basis functions.

Let $\{B_{\bi}(\bxi)\}_{\bi \in I }$ denote the set of Bernstein basis polynomials defined over a reference domain $\hat{\Omega} \subset \mathbb{R}^{d_r}$,
where $I$ is an index set over the degrees of freedom in the element.
Then, a Bernstein polynomial is defined as:
\begin{equation*}
  b(\bxi) = 
  \sum\limits_{\bi\in I} B_\bi \of{ \bxi } \beta_{\bi} \ \ \ \forall \ \bxi \in \Oref 
\end{equation*}
Now, let us define a set of control points $\{\bP_\bi\}_{\bi \in I}$ in $\mathbb{R}^{d_s}$. Then, a B\'{e}zier element is simply defined through the mapping:
\begin{equation*}
  \xphys \of{ \bxi } = 
  \sum\limits_{\bi\in I} B_\bi \of{ \bxi } \bP_{\bi}
\end{equation*}
Thus a B\'{e}zier element is defined via a polynomial pushforward mapping from parametric space, $\mathbb{R}^{d_r}$ to physical space, $\mathbb{R}^{d_s}$.
We note that in general, this mapping holds for any $d_r \leq d_s$.
However, for the purposes of this paper, we consider only the case where $d_r = d_s = d$. 
When $d=1$ we have a curve, when $d=2$ we have triangles or quadrilaterals, and when $d=3$ we have tetrahedra or hexahedra.

\begin{figure}
\centering
\includegraphics[width=0.95\textwidth]{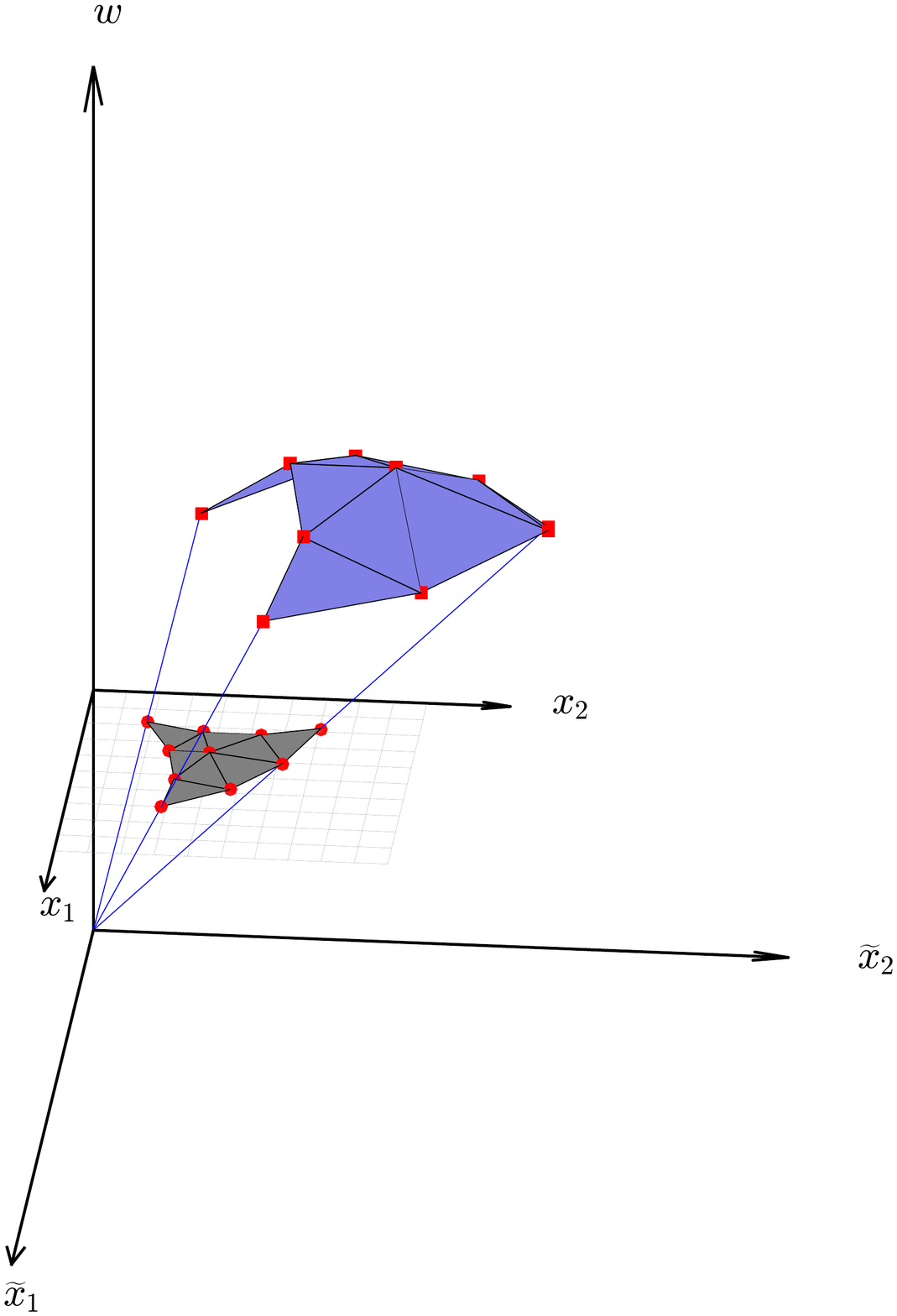}
\caption{Projective transformation from the control net for a projective element $\Oproj \subset \Rproj$ (shown in blue), to the physical element $\Ophys \subset \Rphys$ (shown in grey). 
		Points on the  projective element are given in projective coordinates $\{ \widetilde{x}_1,...,\widetilde{x}_d, w \}$. 
		The physical element is embedded in the $w=1$ plane (gray grid), and points on the physical element are given in terms of the physical coordinates  $\{ x_1,...,x_d\}$.}
\label{proj_CP}
\end{figure}

In addition to the Bernstein basis functions and control points,
let $\{w_\bi\}_{\bi \in I }$ denote a set of \textbf{\textit{control weights}} corresponding to $\{B_{\bi}(\bxi)\}_{\bi \in I }$.
Then we can define a set of corresponding \textbf{\textit{rational}} Bernstein basis functions as:
\begin{equation*}
  R_\bi(\bxi) =
  \dfrac{B_\bi(\bxi) w_\bi}{\sum\limits_{ \bj \in I } B_\bj(\bxi) w_\bj} =
  \dfrac{B_\bi(\bxi) w_\bi}{w(\bxi)}
\end{equation*}
where $w(\bxi) = \sum\limits_{ \bj \in I } B_\bj(\bxi) w_\bj $ denotes the \textbf{\textit{weighting function}} defined over the domain $\hat{\Omega}$.
Then, a rational Bernstein-B\'{e}zier element is defined by the mapping:
\begin{equation*}
  \xphys(\bxi) =
  \sum\limits_{\bi\in I} R_\bi(\bxi)\bP_{\bi} =
  \dfrac{B_\bi(\bxi) \bP_{\bi} w_\bi}{\sum\limits_{ \bj \in I } B_\bj(\bxi) w_\bj}
\end{equation*}
When working with rational B\'{e}zier elements, it is often convenient to consider the corresponding polynomial element in \textbf{\textit{projective space}}, $\Oproj \subset \mathbb{R}^{d+1}$.
The projective element is defined by the mapping:
\begin{equation*}
  \xproj(\bxi) =
  \sum\limits_{\bi\in I} B_\bi(\bxi) \widetilde{\bP}_{\bi} 
\end{equation*}
wherein $\{\widetilde{\bP}_{\bi}\}_{\bi \in I}$ is the set of projective control points, defined as:
\begin{equation*}
  \begin{split}
    \of{ \widetilde{\bP}_\bi }_j \ \ \  =\ & w_\bi \of{ \bP_\bi }_j \ \ \ \  j \in [1, d] \\
    \of{ \widetilde{\bP}_\bi }_{d+1} =\ & w_\bi
  \end{split}
\end{equation*}
This allows us to write a rational B\'{e}zier element in $\mathbb{R}^{d}$ as the projective transformation of a polynomial B\'{e}zier element in $\mathbb{R}^{d+1}$, viz:
\begin{equation*}
  \xphys(\bxi) = \dfrac{\left[ \xproj(\bxi) \right]_{d}}{w(\bxi)}
\end{equation*}
wherein, $\left[ \xproj(\bxi) \right]_{d}$ denotes the first $d$ components of $\xproj(\bxi)$.
We illustrate the control points for a rational \BB element and the corresponding element in projective space in Fig. \ref{proj_CP}. 

Thus, it is readily seen that a rational Bernstein-B\'{e}zier element is defined by:
1) the Bernstein basis polynomials,
2) the B\'{e}zier control points, and
3) the control weights.
In previous work, we have given in-depth descriptions of a variety of Bernstein-B\'{e}zier elements \cite{engvall_isogeometric_2017}.
We briefly review the construction of the Bernstein basis polynomials for various elements here, and refer readers desiring more information to the previous work.

\subsubsection{Simplicial Elements} \label{bez_simplex}

The simplest Bernstein-B\'{e}zier elements are the simplicial elements: curves (1-simplices), triangles (2-simplices), and tetrahedra (3-simplices).
When working with simplicial elements, it is most common to work in barycentric coordinates.
However, as we are particularly interested in the derivatives of Bernstein polynomials with respect to the unit Cartesian reference triangle,
it is also convenient to work explicitly in terms of Cartesian coordinates.
As such, we present both forms, and switch between forms as necessary.

\noindent \textbf{Cartesian Coordinates}

We define the reference domain for a $d$-simplex as:
\begin{equation*}
  \hat{\Omega} = \setof{ \bxi \in \mathbbm{R}^d : 0 \leq \xi_j \leq 1, \sum\limits_{j=1}^d \xi_j \leq 1 }
\end{equation*}
and define the index set for the simplicial Bernstein polynomials of degree $p$ as:
\begin{equation*}
  I^p := \setof{ \bi = \setof{i_1,...,i_d } : 0 \leq i_j \leq p, |\bi| \leq p }
\end{equation*}
Then, we can define the simplicial Bernstein basis polynomials as:
\begin{equation*}
  B_{\bi}^p (\bxi) = p!\prod\limits_{j = 1}^d \dfrac{ \of{ \xi_j }^{i_j} }{i_j} \dfrac{ \of{ 1-|\bxi| }^{i_{d+1}}}{p-|\bi|}
\end{equation*}

\noindent \textbf{Barycentric Coordinates}

For a $d$-simplex, we define the set of $d+1$ barycentric coordinates $\{\lambda_j\}_{j=1}^{d+1}$ corresponding to a point $\bxi$ in the reference domain $\hat{\Omega}$ to be:
\begin{equation*}
  \begin{split}
    \lambda_j = \xi_j \ \ \ j \in [0,d] \\
    \lambda_{d+1} = 1-\sum\limits_{j=1}^d \lambda_j
    \end{split}
\end{equation*}
Then the reference domain is defined in terms of barycentric coordinates as: 
\begin{equation*}
  \hat{\Omega} = \{ \boldsymbol{\bblambda} \in \Rproj : 0 \leq \lambda_j \leq 1, |\mathbbm{\bblambda}| = 1\}
\end{equation*}
and the corresponding index set for the simplicial Bernstein polynomials of degree $p$ is defined as:
\begin{equation*}
  I^p_{\bbi} := \{ \mathbb{i} = \{i_1,...,i_{d+1}\} : 0 \leq i_j \leq p, |\mathbb{i}| = p \}.
\end{equation*}
This allows us to write the simplicial Bernstain polynomials incredibly compactly as:
\begin{equation*}
  B_{\mathbb{i}}^p (\bblambda) = {{p}\choose{\mathbb{i}}}  \bblambda^{\mathbb{i}}
\end{equation*}
We note that for simplicial elements, we differentiate between the Cartesian and barycentric forms by either explicitly specifying the dependence on $\bxi$ or $\bblambda$, or implicitly by the index $\bi$ or $\mathbb{i}$. 
The index $\bi $ is taken to always be a multi-index of length $d$, whereas we understand $\mathbb{i}$ to denote a multi-index of length $d+1$ such that $|\mathbb{i}| = p$.
\subsubsection{Tensor Product Elements}

We define \BB quadrilaterals in $\mathbb{R}^2$ and hexahedra in $\mathbb{R}^3$ using a tensor product construction.
The reference domain is given by:
\begin{equation*}
  \hat{\Omega} = (0,1)^d
\end{equation*}
and the index set for a tensor product construction is given by:
\begin{equation*}
  I^{\bp} := \{ \bi = \{i_1,...,i_{d}\} \ | \ 0 \leq i_j \leq p_j \}.
\end{equation*}
Then, the Bernstein basis polynomials are defined as:
\begin{equation*}
  B_{\bi}^\bp (\bxi) = \prod\limits_{j=1}^dB_{i_j}^{p_j}(\xi_j)
  \end{equation*}
where $B_{i_j}^{p_j}$ are simply the univariate Bernstein polynomials (i.e., one-dimensional simplicial Bernstein polynomials) from Subsubsection \ref{bez_simplex}.

We note that we distinguish between simplicial and tensor product constructions implicitly by using a scalar superscript $p$ for the simplicial basis functions $B^p_\bi$, and a vector super script $\bp$ for the tensor product basis functions $B^\bp_\bi$.

\section{Literature Review} \label{review}
Before proceeding to the novel contributions of this work, it is critical to understand the existing theory regarding error bounds for finite elements, and how these error bounds motivate the need for element quality metrics. 
First, we briefly review the isoparametric concept as it applies to finite elements (\S \ref{isoparametric}).
We then review the fundamental interpolation theory for both linear (\S \ref{theory_review}) and curvilinear (\S \ref{curve_review}) finite elements.
Finally, we present the current state of the art with regards to element quality metrics for curvilinear finite elements  (\S \ref{metrics_review}).
\subsection{The Isoparametric Concept} \label{isoparametric}
\begin{figure}[t]
\centering
  \includegraphics[width=0.85\textwidth]{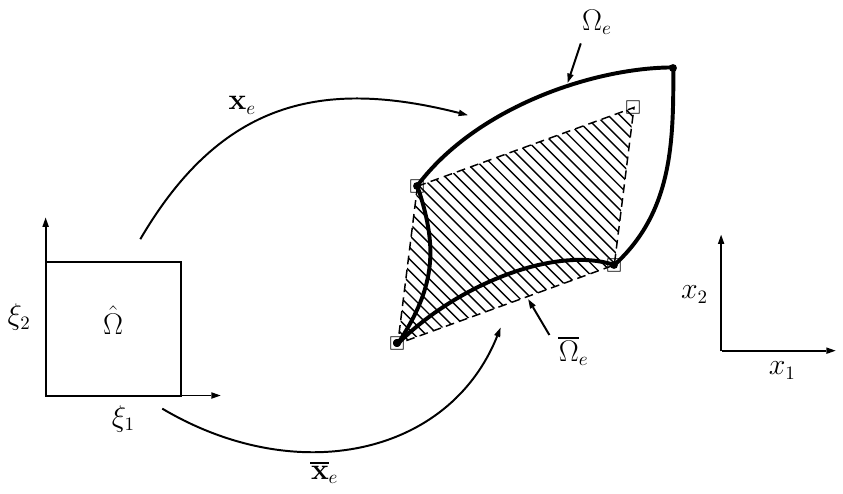}
  \caption{Isoparametric mappings from a reference element  $\Oref$ in parametric space to physical space. 
  The element $\Ophys$ (shown by the bold line) is defined by the higher-order mapping $\xphys$.
  The corresponding linear element $\Olin$ (shown by the dashed line) is defined by a purely affine mapping $\xlin$.}
  \label{fig:mappings}
\end{figure}

Simply put, the isoparamtric concept allows us to define an element in physical space in terms of a mapping from a reference element in parametric space.
Let us denote a unit reference element in parametric space $\Oref$.
Then we denote the element in physical space $\Ophys$, and denote a mapping $\xphys: \Oref \rightarrow \Ophys$ that maps points on the parametric element to points on the physical element. 
In the case of higher-order finite elements, this is a higher-order mapping. Finally, we also consider the element $\Olin$ which is the \textbf{\textit{purely linear}} physical element. That is, $\Olin$ is defined by an affine mapping $\xlin: \Oref \rightarrow \Olin$.
We assume that the mapping $\xphys$ is known. Then we define the affine mapping $\xlin$ as:
\begin{equation*}
  \xlin(\bxi) = \xphys(\textbf{0}) + \overline{\textbf{J}}\bxi
\end{equation*}
wherein the Jacobian of the affine mapping is defined by:
\begin{equation*}
  \overline{\textbf{J}} = \Bigg{[} [\xphys(\textbf{e}_1) - \xphys(\textbf{0})],...,[\xphys(\textbf{e}_d) - \xphys(\textbf{0})] \Bigg{]}
\end{equation*}
These mappings are illustrated in Fig. \ref{fig:mappings}. 
We note that for simplicial elements, the linear physical element is simply defined as the linear interpolant of the corners of the curvilinear physical element. 
However, for tensor product elements (i.e. quadrilaterals or hexahedra), the linear element will not necessarily interpolate every corner of the curvilinear element, as seen in Fig. \ref{fig:mappings}.
This is due to the fact that the tensor product admits bilinear mappings for quadrilaterals and trilinear mappings for hexahedra.

\subsection{Finite Element Interpolation Theory: Linear Elements} \label{theory_review}

The rigorous study of the mathematical foundations of the finite element method began in earnest in the 1960's and 1970's
\cite{bramble_estimation_1970,bramble_bounds_1971,ciarlet_general_1972,ciarlet_interpolation_1972,zlamal_curved_1973,zlamal_curved_1974,babuska_angle_1976}.
Of particular note are the pioneering works of Bramble and Hilbert \cite{bramble_estimation_1970,bramble_bounds_1971} and Ciarlet and Raviart \cite{ciarlet_general_1972,ciarlet_interpolation_1972}
which led to the \textit{\textbf{interpolation theory for finite elements}}.
As this fundamental theory has close bearing on our current work,
we take this opportunity to briefly review the theory here, introducing notation to be used throughout the remainder of this paper.

\begin{figure}[t] 
\centering
\includegraphics[width=0.35\textwidth]{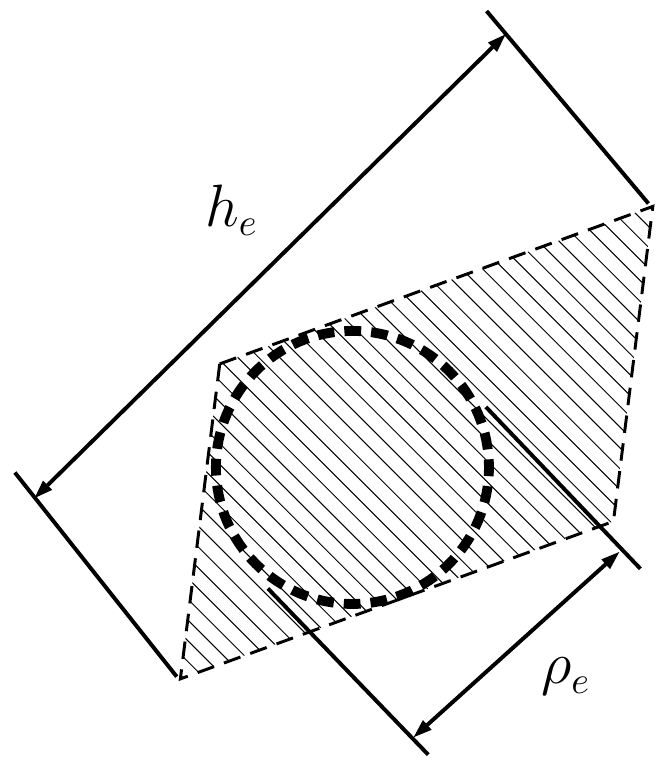}
\caption{Element measures $\rho_e$ and $h_e$ for the linear quadrilateral element $\Olin$ from Fig. \ref{fig:mappings}.}
\label{fig:elem_measures}
\end{figure}

In the finite element method, we approximate a domain $\Omega$ using a set of finite elements, $\{\Ophys\}_{e=1}^E$, where each element $\Ophys \subset \Rphys$ is an open simply connected set, with simply connected boundary. 
Together, this collection of elements forms a finite element discretization or mesh, which we denote as:
\begin{equation*}
\mathcal{M} = \overline{\bigcup\limits_{e=1}^E \Ophys}
\end{equation*}
Traditionally, a finite element mesh is composed of linear triangular or quadrilateral elements in $\mathbb{R}^2$,
and linear tetrahedral or hexahedral elements in $\mathbb{R}^3$.
The study of mesh quality, then, concerns itself with how the shapes and sizes of these elements affect the accuracy of the finite element method.
In order to begin this study, it is first useful to introduce some mesh measures.
Let us denote the diameter of an element as $h_e$, where we measure the diameter as the largest distance between any two vertices of the element.
Next, we denote the diameter of the incircle (in $\mathbb{R}^2$) or insphere (in $\mathbb{R}^3$) of the element as $\rho_e$.
These two metrics are visualized for a quadrilateral element in Fig. \ref{fig:elem_measures}.

Given these element-wise measures, we can then define the corresponding global mesh measures as:
\begin{equation*}
h = \max_{1 \leq e \leq E}{h_e}
\end{equation*}
\begin{equation*}
\rho = \min_{1 \leq e \leq E}{\rho_e}
\end{equation*}
These mesh measures are useful, as they allow us to put certain classifications on our meshes.
Specifically, it allows us to introduce the notion of \textbf{\emph{shape regularity}}.
For a linear element $\Olin$, the \textbf{\emph{element}} shape regularity is given by:
\begin{equation*}
  \sigma_e = \dfrac{h_e}{\rho_e}
\end{equation*}
and for a mesh of linear elements, the \textbf{\emph{global}} shape regularity is given by:
\begin{equation*}
  \sigma = \dfrac{h}{\rho}
\end{equation*}
The notion of shape regularity is important, as it is allows us to compare the shapes of elements independent of their size. 
This is motivated by the fact that we are interested in the effect of element shape on approximation error under mesh refinement, that is as $h \rightarrow 0$.

Let us consider a set of $M$ increasingly refined meshes $\{\mathcal{M}_i\}_{i=1}^M$, with corresponding metrics $\{h_i\}_{i=1}^M$ and $\{\rho_i\}_{i=1}^M$ where $h_1 > h_2 > \hdots > h_M$.
We say a refinement is \boldit{uniform} if each element is simply split into a collection of similarly shaped sub-elements, thereby preserving the existing geometrical structure of the parent mesh. 
Due to the ease of implementation, and constant element shape, regular refinements are frequently used with the finite element method.
However, uniform refinements are constrained by the choice of the initial mesh, and it is not always practical or even possible to perform uniform subdivision on a given mesh.
It is useful then to introduce the concept of \boldit{quasi-uniform} refinement. 
Consider a series of refined meshes $\{\mathcal{M}_i\}_{i=1}^M$ with corresponding global shape regularity metrics $\{ \sigma_i \}_{i=1}^M$.
We say a series of refinements is quasi-uniform if we can bound $\sigma_i$ by some constant $\sigma_0$, viz.:
\begin{equation*}
  \sigma_i \leq \sigma_0 \ \ \ i=1,2,...,M
  \label{reg_cond}
\end{equation*}
In the case of both uniform and quasi-uniform refinements, we say that the set of all elements in the set of meshes belong to a \textbf{\textit{regular family}} of elements.
If a series of refinements is not uniform or quasi-uniform, we say that the refinements are \textit{\textbf{irregular}}, and the elements do not belong to a regular family.
These three classes of refinements are visualized for a simple mesh in Table \ref{mesh_refs}.

     \begin{table}[t]
     \begin{center}
       \caption{Various types of mesh refinement. With uniform refinement, the structure of the original mesh is preserved at each level of refinement. 
      	With quasi-uniform refinement, the structure is not preserved, but all elements belong to a regular family.
	With irregular refinement, the bottom elements become increasingly thin, and therefore do no belong to a regular family. }
     \begin{tabular}{ | >{\centering\arraybackslash} m{2.5cm}  | >{\centering\arraybackslash} m{3.5cm} |  >{\centering\arraybackslash} m{3.5cm} | >{\centering\arraybackslash} m{3.5cm}  | }
     \hline
          & $h$ = 0.5 & $h$ = 0.25 & $h$ = 0.125 \\
               \hline
   	Uniform Refinements &
	\vspace{5pt} \includegraphics[width=3cm,height=3cm,keepaspectratio]{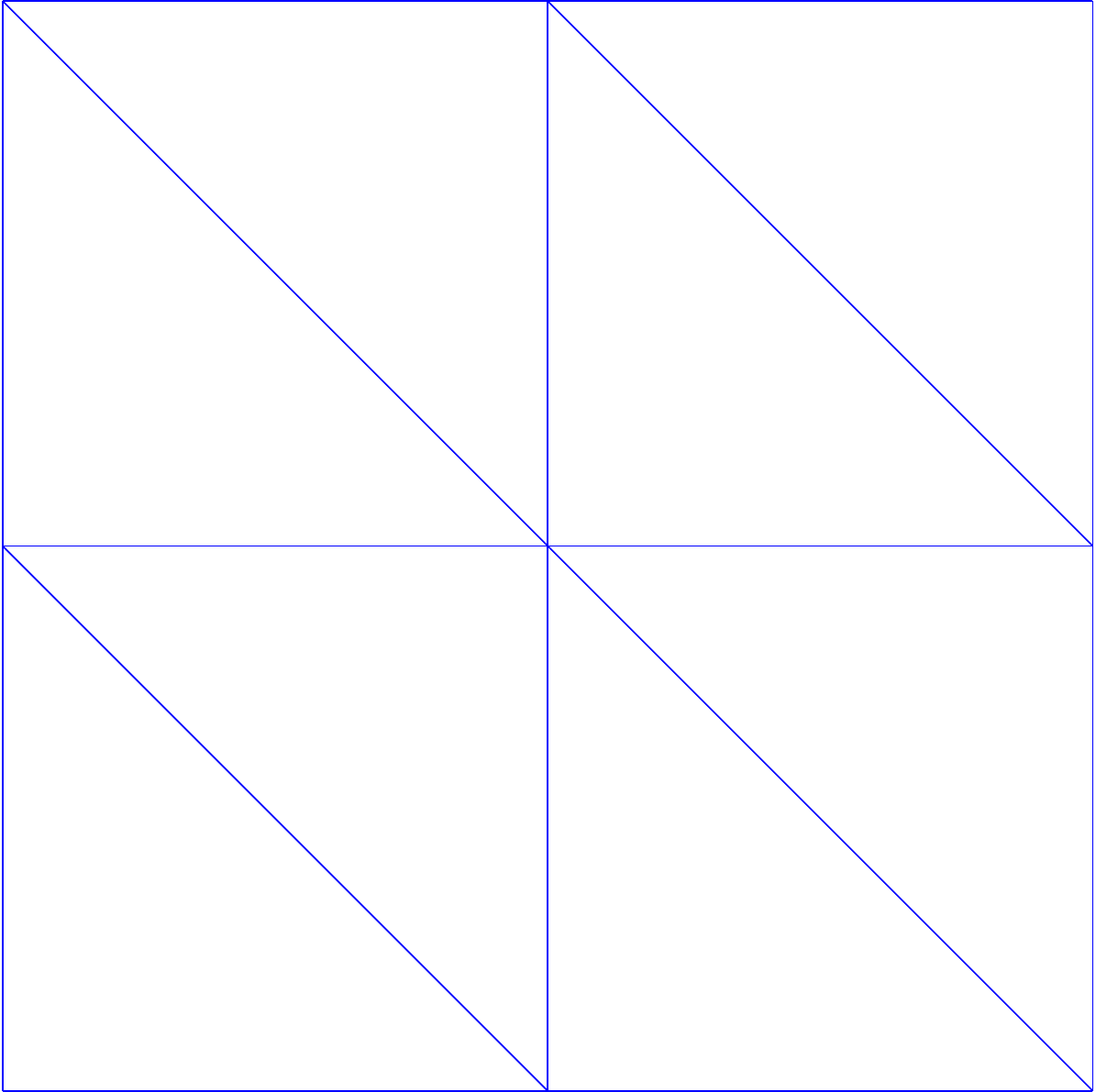} &
	\vspace{5pt} \includegraphics[width=3cm,height=3cm,keepaspectratio]{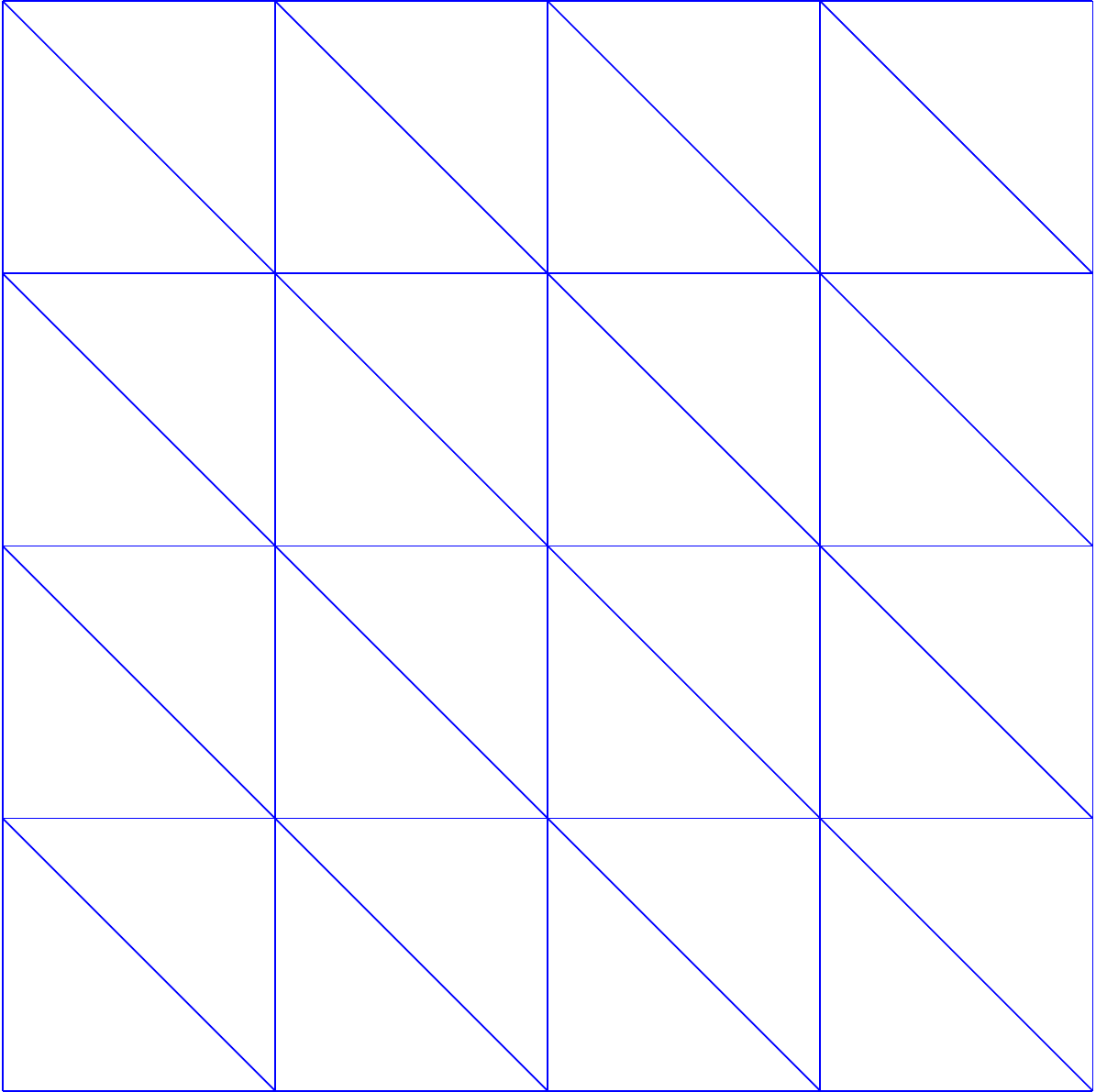} &
	\vspace{5pt} \includegraphics[width=3cm,height=3cm,keepaspectratio]{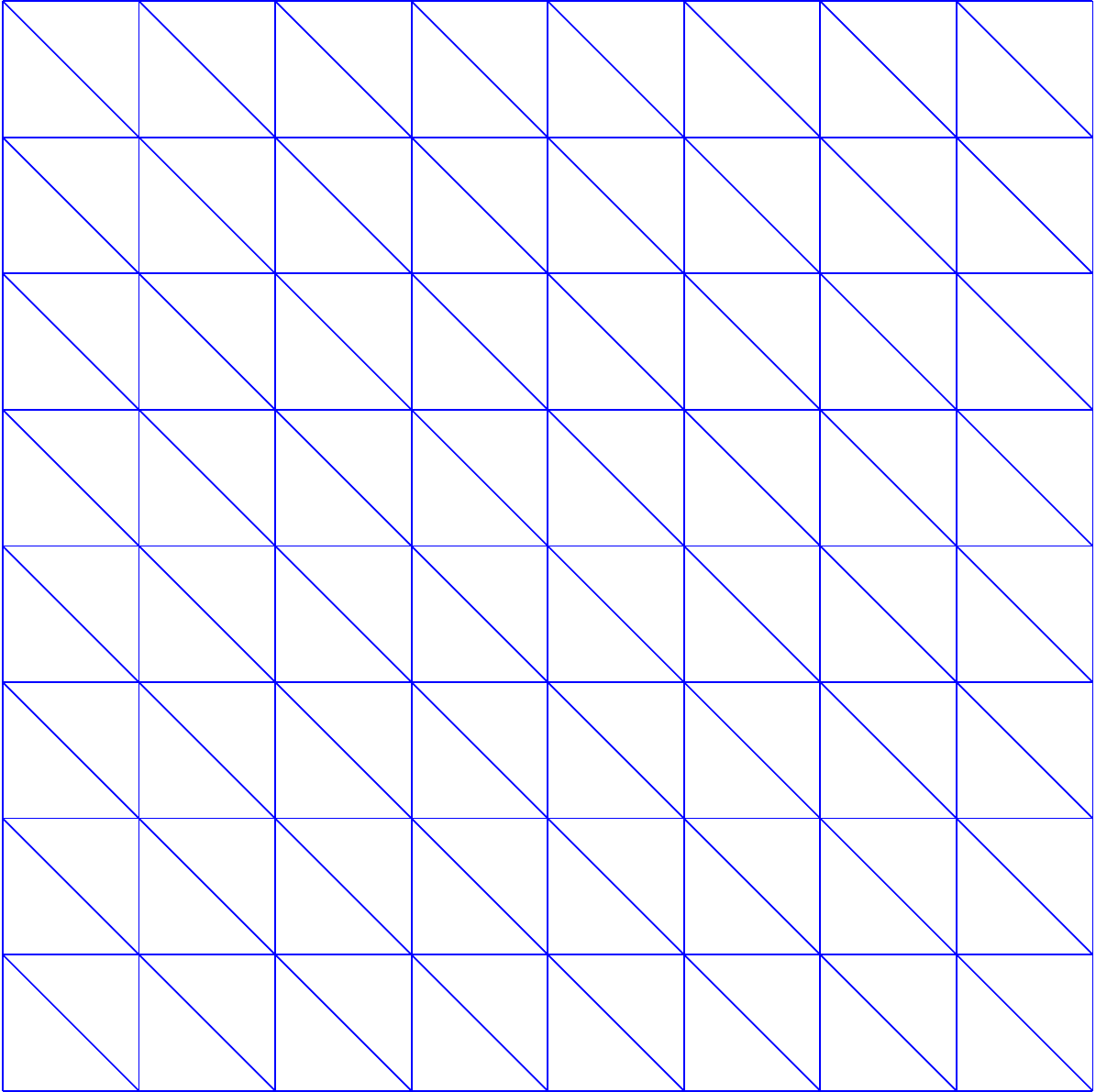}\\
                        \hline
   	Quasi-Uniform Refinements &
	\vspace{5pt}\includegraphics[width=3cm,height=3cm,keepaspectratio]{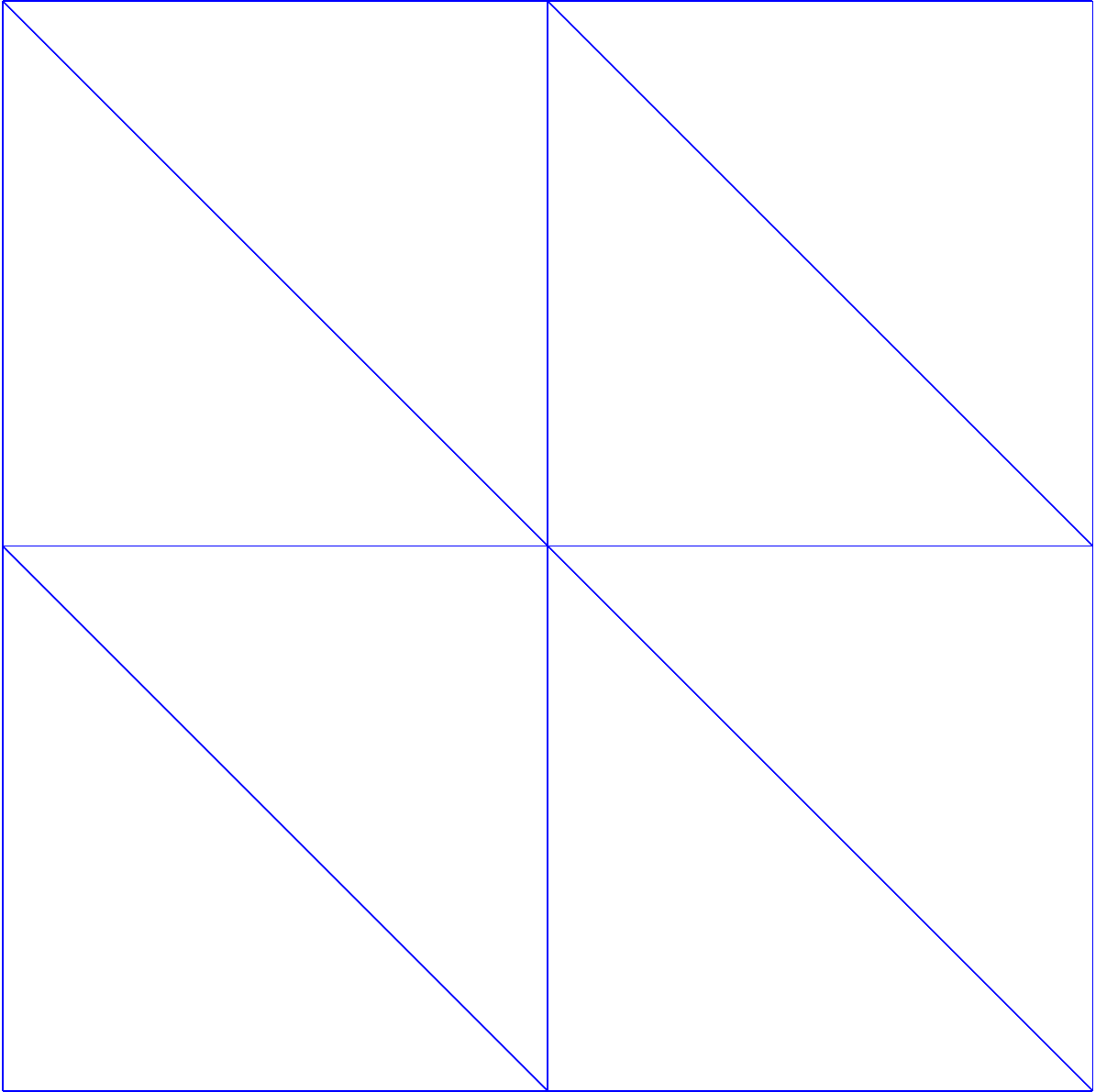} &
	\vspace{5pt}\includegraphics[width=3cm,height=3cm,keepaspectratio]{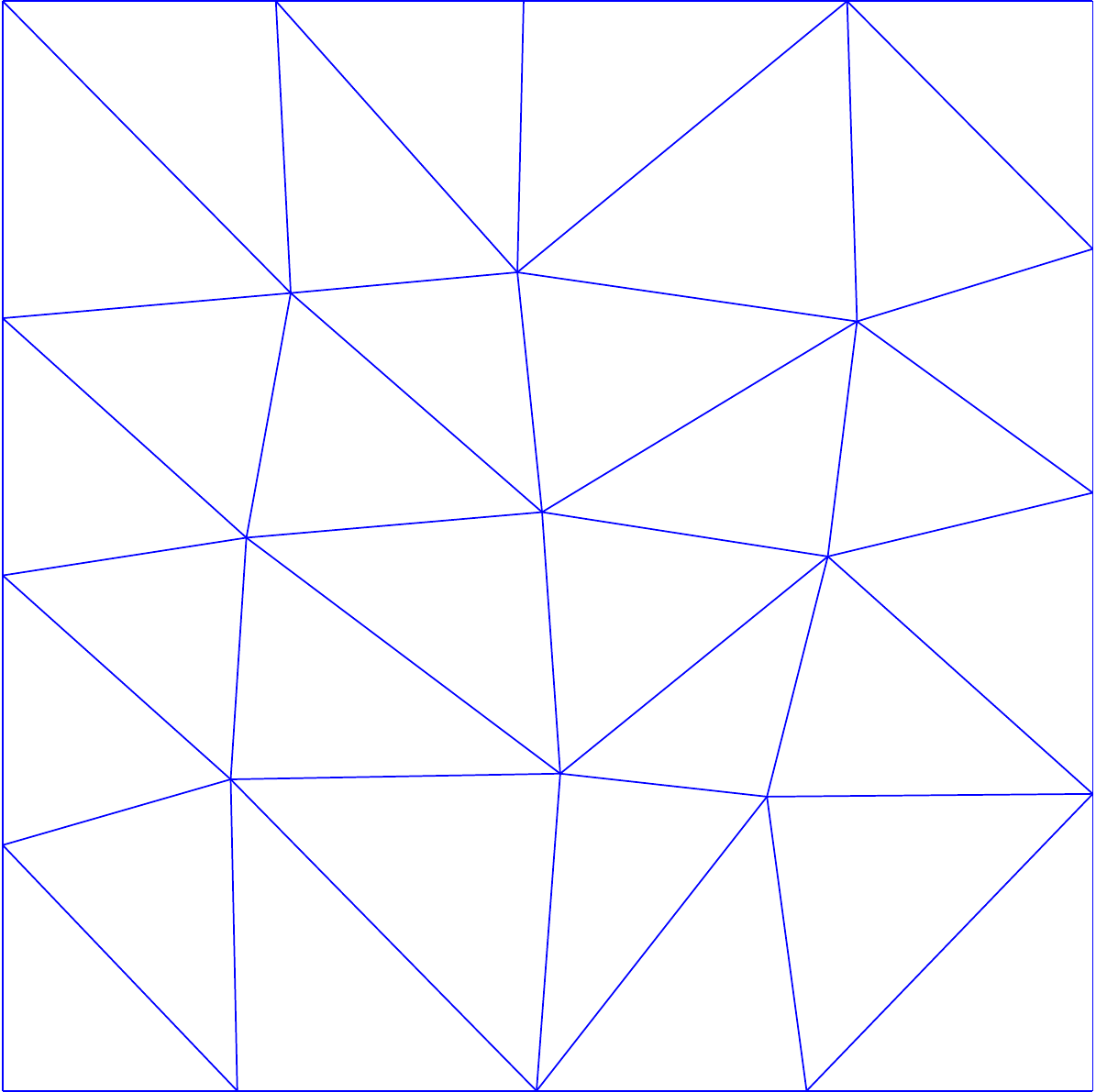} &
	\vspace{5pt}\includegraphics[width=3cm,height=3cm,keepaspectratio]{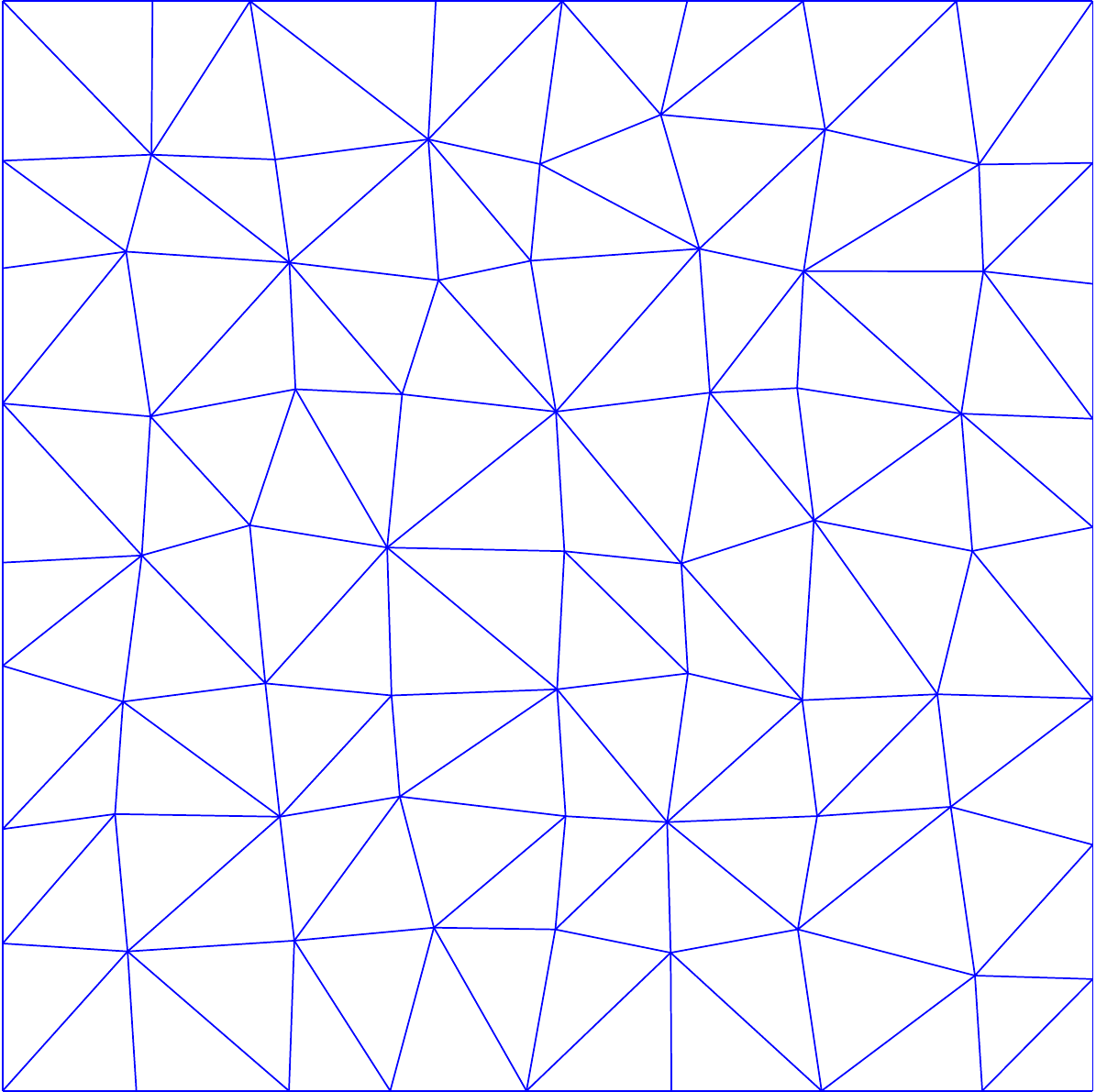}\\
	     \hline
	        	Irregular Refinements &
	\vspace{5pt}\includegraphics[width=3cm,height=3cm,keepaspectratio]{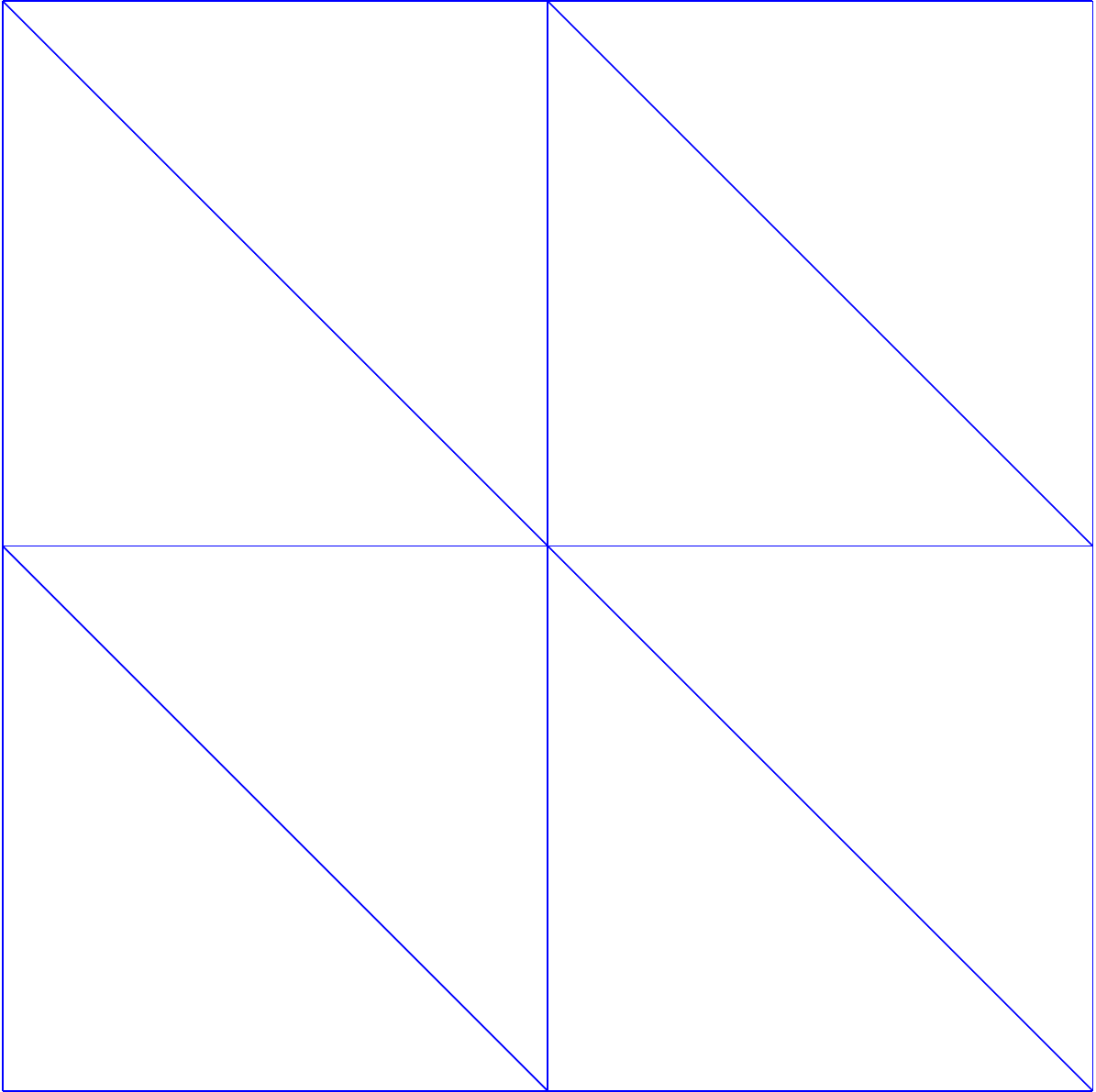} &
	\vspace{5pt}\includegraphics[width=3cm,height=3cm,keepaspectratio]{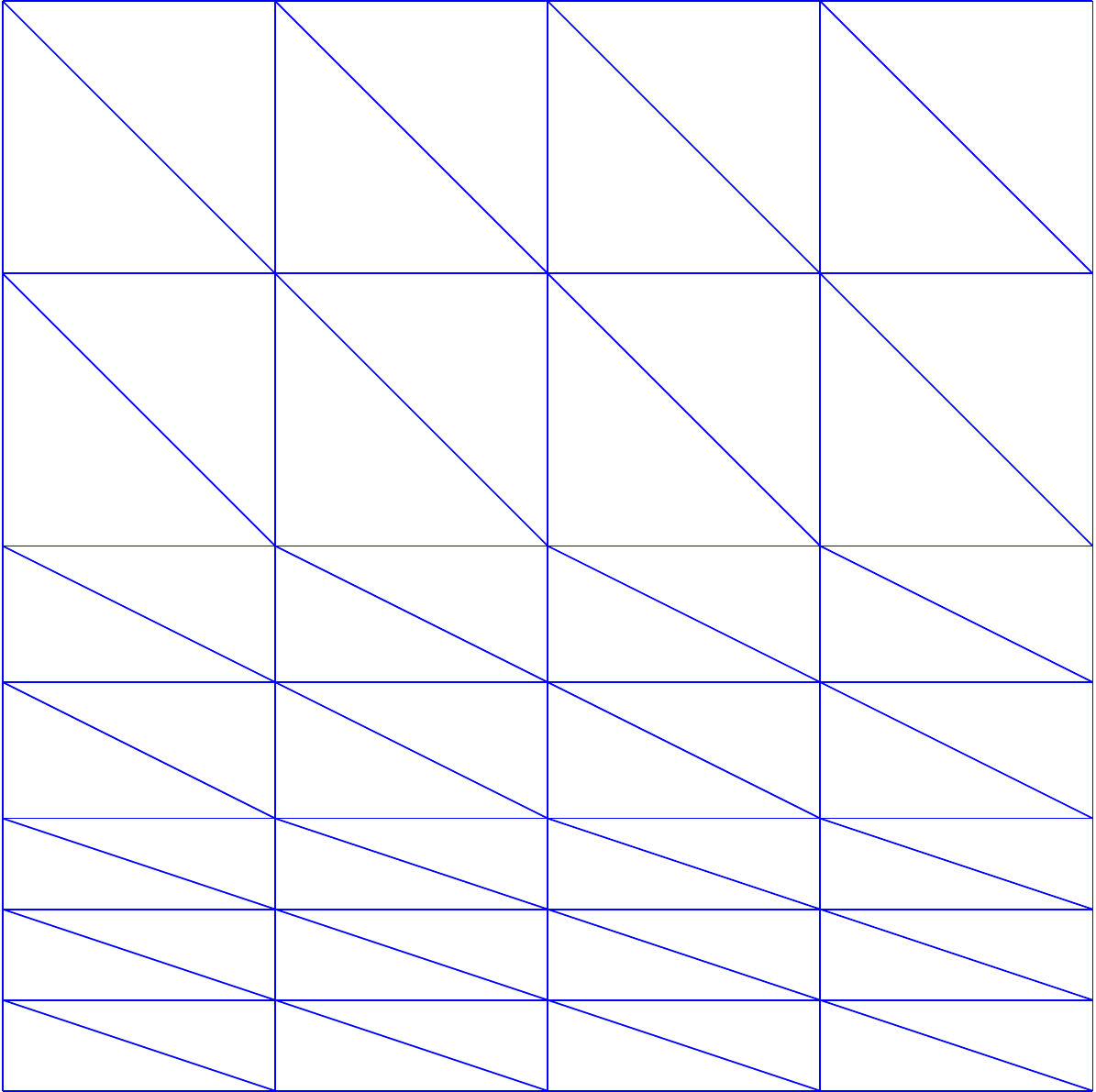} &
	\vspace{5pt}\includegraphics[width=3cm,height=3cm,keepaspectratio]{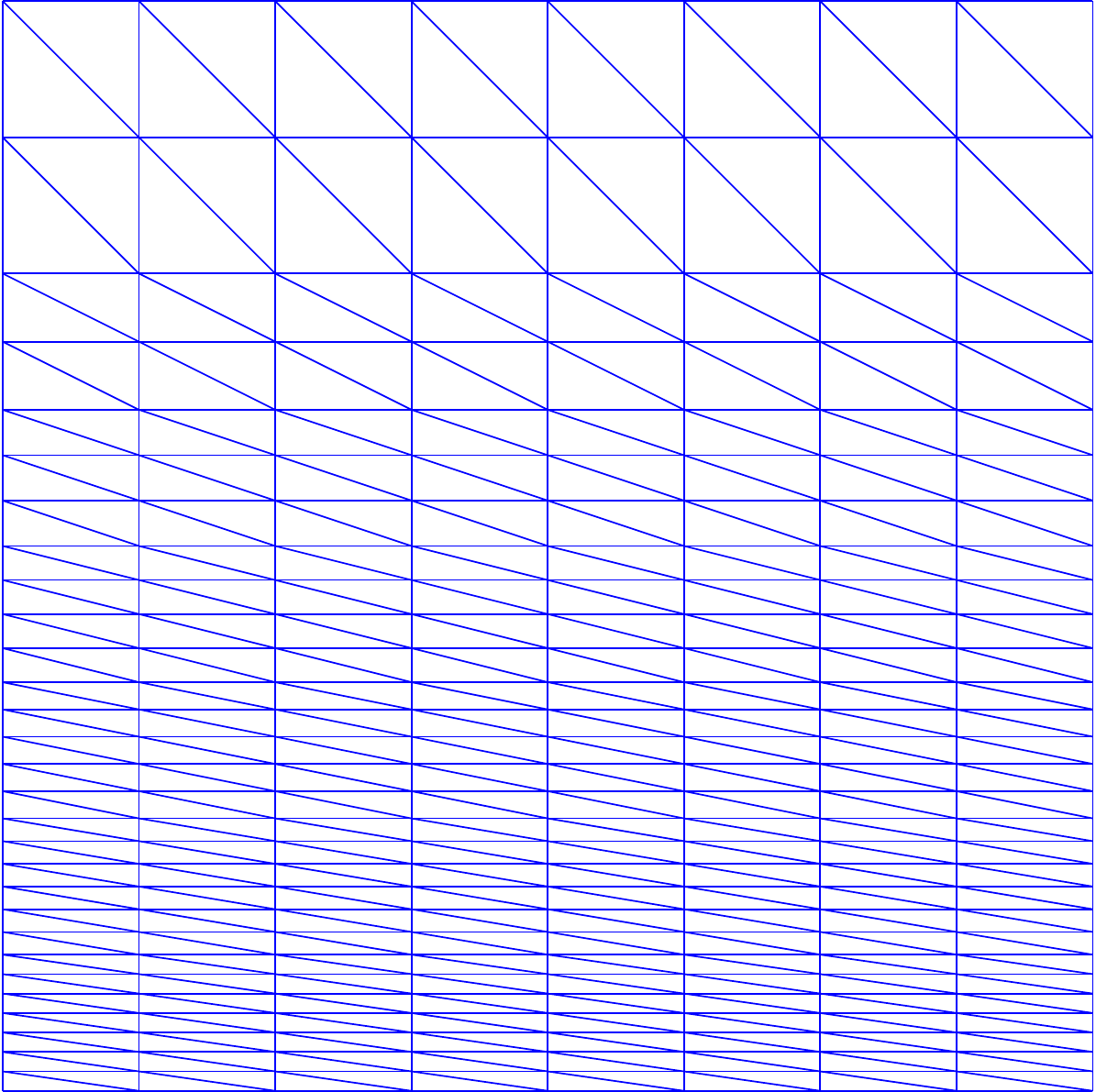}\\
	     \hline
      \end{tabular}
          \label{mesh_refs}
      \end{center}
      \end{table}

With the necessary notation established, we can proceed to state the interpolation theory for finite elements.
Suppose we have a sequence of meshes  $\{\mathcal{M}_i\}_{i=1}^M$.
Given some boundary value problem defined over the domain, we desire to approximate the true solution $u$ using the finite element method. 
That is, we desire to approximate $u$ by $u^h$, a $C^0$-continuous piecewise polynomial of degree $p$, for each mesh in the series.
Then, following Raviart and Ciarlet \cite{ciarlet_general_1972}, it can be shown that the \textbf{\emph{best approximation error}} over a linear mesh with mesh size $h$ is bounded by:
\begin{equation}
  \normof{u - u_h }_{H^m(\Omega)} \leq C\dfrac{h^{p+1}}{\rho^m}|u|_{H^{p+1}(\Omega)}
\end{equation}
wherein $C$ is a constant independent of both the mesh size $h$ and the mesh global shape regularity $\sigma$, and $|| \cdot ||_{H^m(\Omega)}$ denotes the norm:
\begin{equation}
  \normof{ f }_{H^m(\Omega)} = \of{ \sum\limits_{j=0}^m \absof{ f }^2_{H^j(\Omega)} }^{1/2}
\end{equation}
and $| \cdot |_{H^j(\Omega)}$ denotes the seminorm:
\begin{equation}
  |f|_{H^j(\Omega)} = \of {\int\limits_\Omega\sum\limits_{|\balpha|=j} \absof{ D^{\balpha} f }^2 d\Omega } ^{1/2}
\end{equation}
Furthermore, we notice that if every mesh in series belongs to a regular family, this bound simplifies to:
\begin{equation}
  \normof{ u - u_h }_{H^m(\Omega)} \leq Ch^{p+1-m}|u|_{H^{p+1}(\Omega)}
\end{equation}
wherein $C$ is a constant independent of the mesh size $h$ but dependent on the family shape regularity $\sigma_0$.
%
%
%

\subsection{Finite Element Interpolation Theory: Curvilinear Elements}\label{curve_review}
\begin{figure}[t] 
\centering
 \begin{subfigure}[t]{0.24\textwidth}
        \centering
            \includegraphics[width=.9\linewidth]{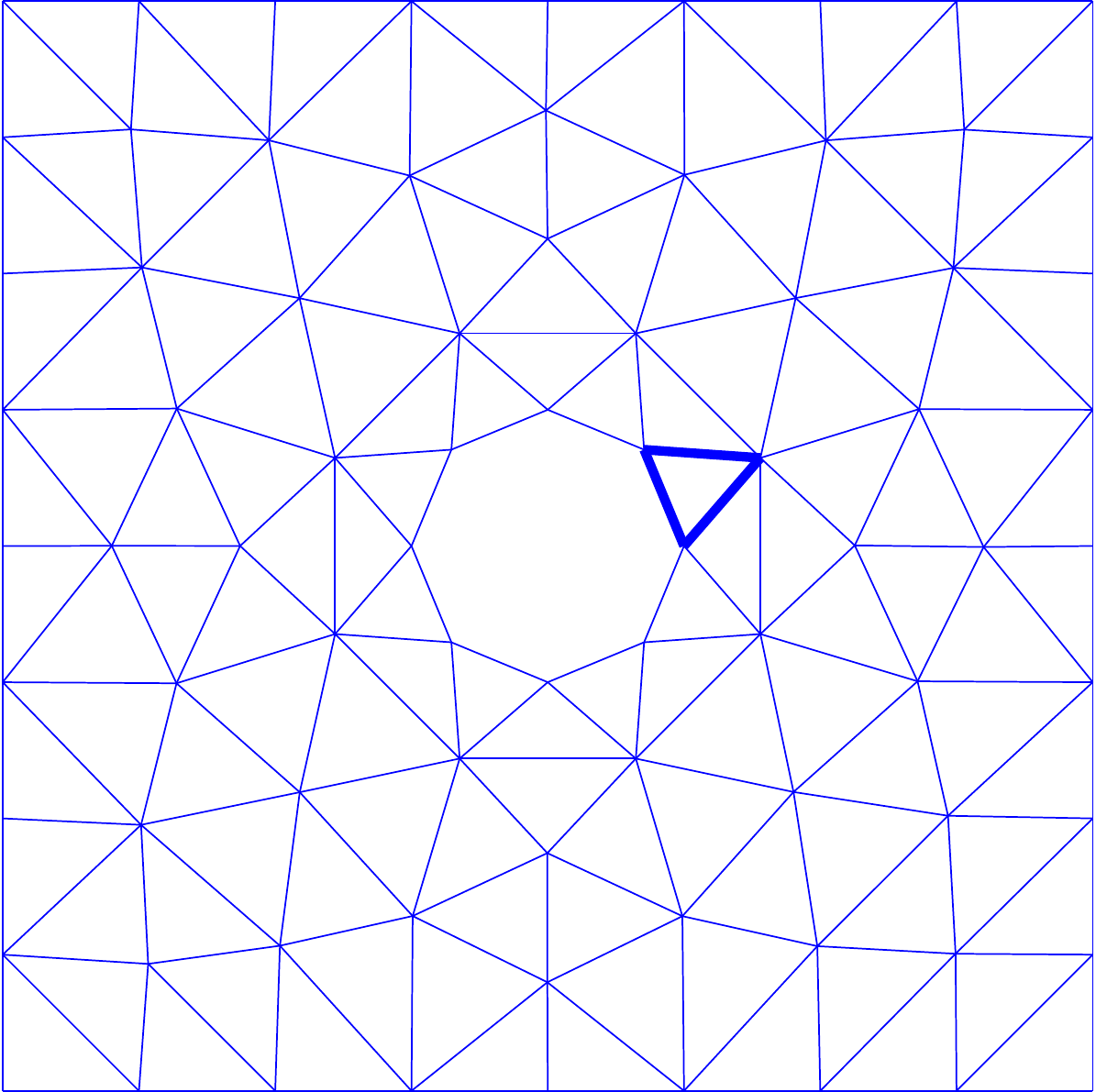}
        \caption{}
    \end{subfigure}
    \begin{subfigure}[t]{0.24\textwidth}
        \centering
            \includegraphics[width=.9\linewidth]{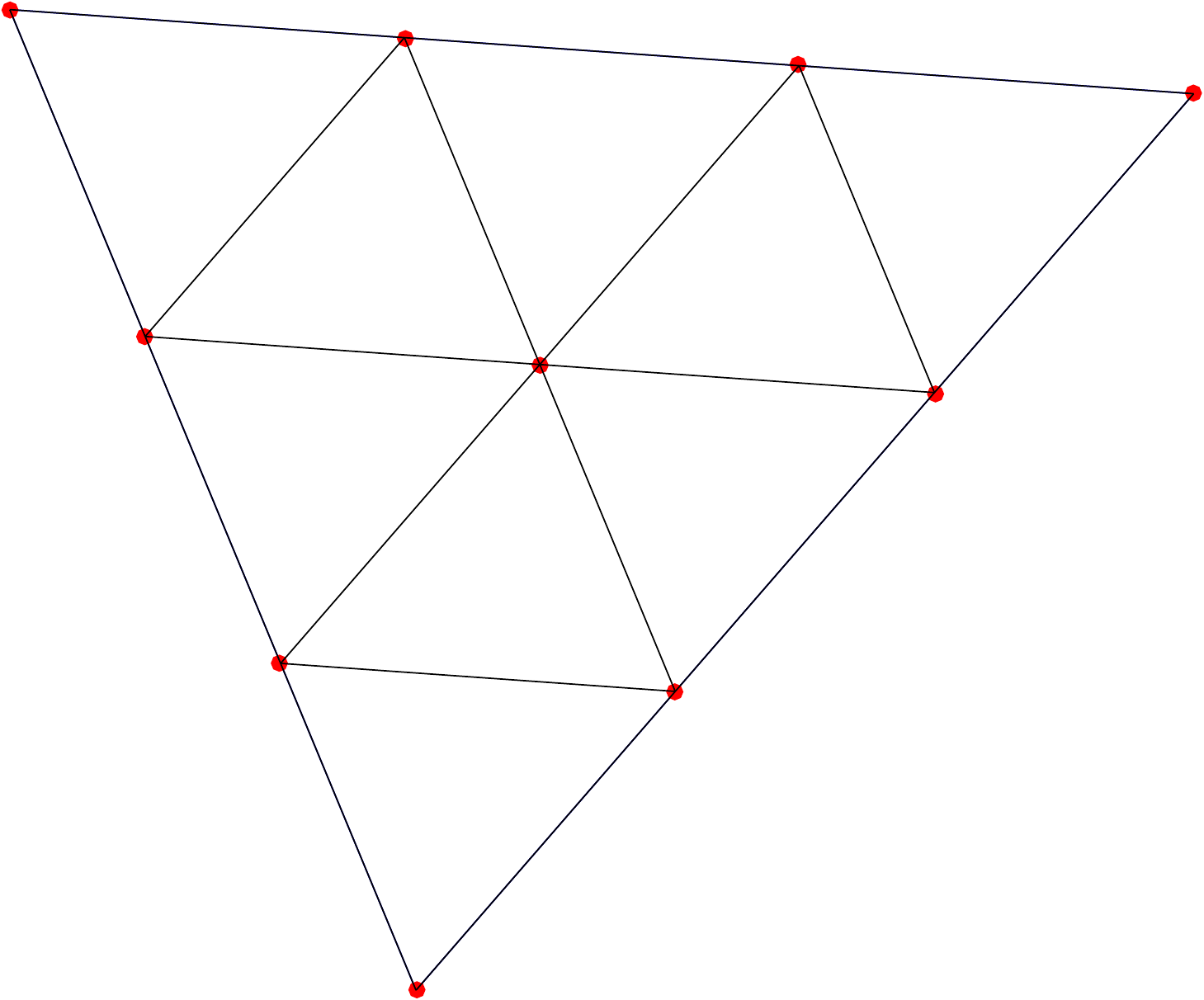}
        \caption{}
    \end{subfigure}
\begin{subfigure}[t]{0.24\textwidth}
        \centering
            \includegraphics[width=.9\linewidth]{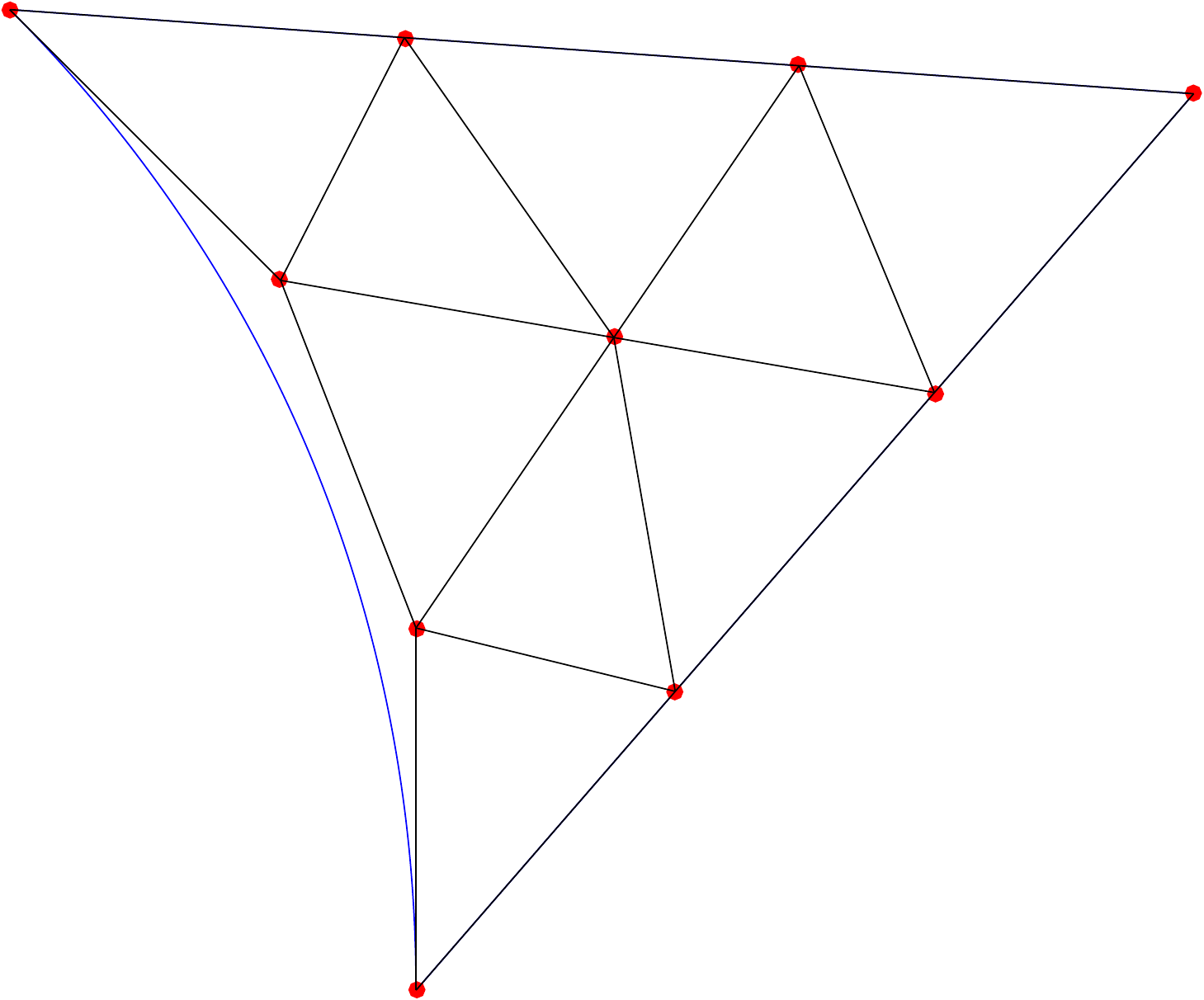}
        \caption{}
    \end{subfigure}
\begin{subfigure}[t]{0.24\textwidth}
        \centering
            \includegraphics[width=.9\linewidth]{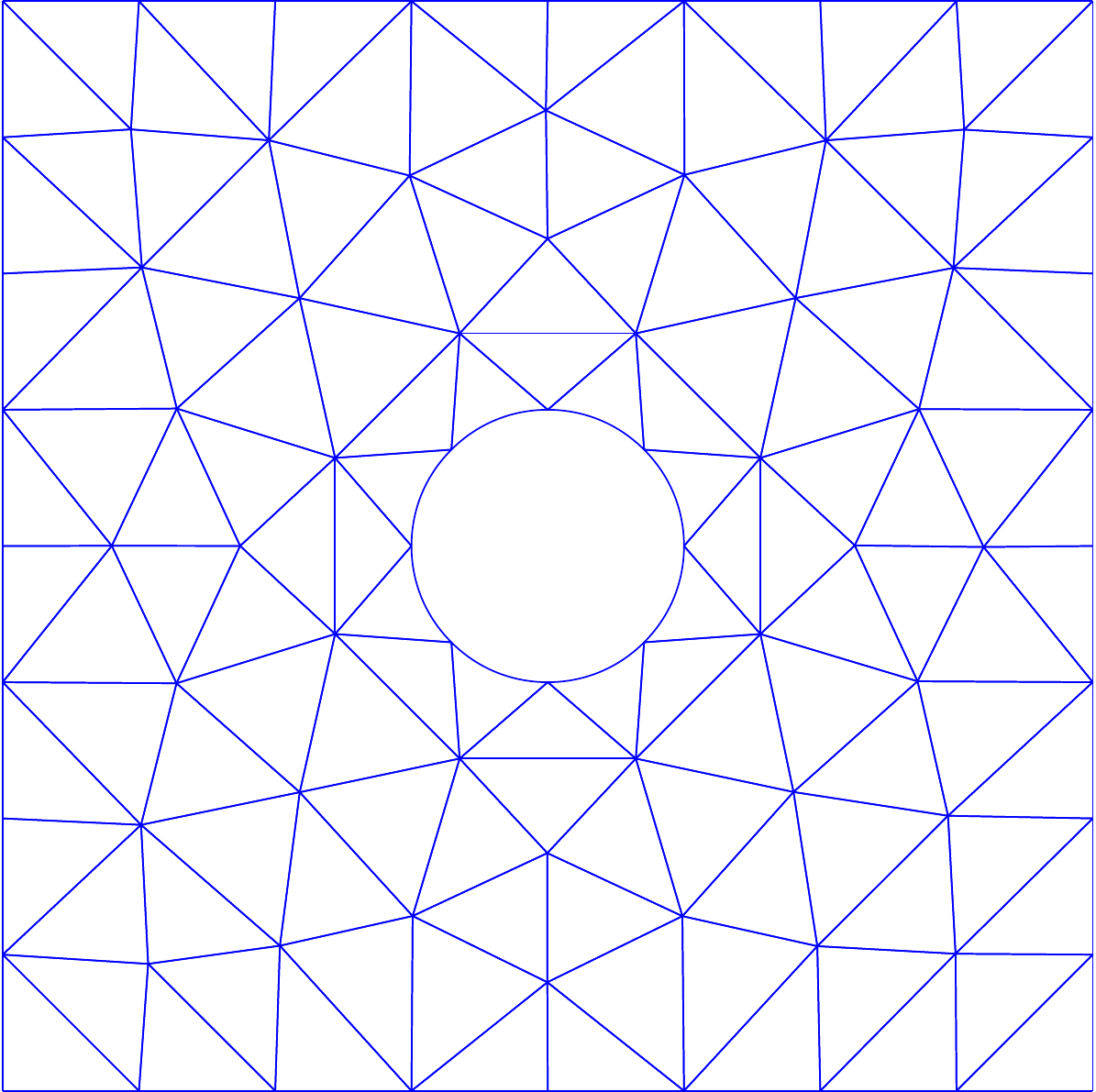}
        \caption{}
    \end{subfigure}
\caption{ Steps for basic curvilinear mesh generation on a plate with a hole. (a) Create an initial linear mesh. A representative boundary element is highlighted in bold. (b) Degree elevate the linear element by inserting higher-order control points. (c) Curve the element to match the boundary by updating control point location. (d) Repeat for each element in the mesh to yield the final curvilinear mesh. }
\label{curved_mesh}
\end{figure}

Similar bounds have also been established for interpolation over curvilinear elements \cite{ciarlet_interpolation_1972,oden_introduction_2012}.
In curvilinear mesh generation, the mesh is most often generated by first constructing a mesh of straight sided elements, and then manipulating element nodes to curve the elements to better match the geometry (see Fig. \ref{curved_mesh}).
Thus with curvilinear finite elements, there is the notion of both the  curvilinear mesh $\mathcal{M}$, as well as the underlying linear mesh $\overline{\mathcal{M}}$.
Now, let us denote a mapping $\xphys: \Oref \rightarrow \Ophys$ that maps a linear master element $\Oref$ to the physical curvilinear element $\Ophys$.
Then, let us assume that the following conditions hold for every element in the mesh.
\begin{enumerate}[label={{Cond. (\ref{curve_review}.\arabic*)}}, align=left] \label{classical_curve_cond}
\item{
The underlying linear elements $\Olin$ belong to a regular family. 
}
\item{ 
The mapping $\xphys$ is invertible with inverse $\xphys^{-1} : \Ophys \rightarrow \Oref$. That is:
\begin{equation*}
\xphys(\bxi) = \bx \Leftrightarrow \xphys^{-1}(\bx)  = \bxi \ \ \ \forall \ \bxi \in \Oref 
\label{invert_cond}
\end{equation*}
}
\item{
The derivatives of the mapping $\xphys$ are bounded as follows:
\begin{equation*}
  \begin{split}
   & \sup\limits_{\bxi \in \Oref}\max\limits_{\absof{\balpha} = k } \absof{ \derivnx{{\balpha}}{\bxi} \xphys } \leq c_kh^k \ \ \ 1 \leq k \leq p + 1 \\
   & \sup\limits_{\bx \in \Ophys}\max\limits_{\absof{\balpha} = 1 } \absof{ \derivnx{{\balpha}}{\bx}\of{ \xphys^{-1}} } \leq c_0h^{-1}
   \label{deriv_cond}
  \end{split}
\end{equation*}
}
\end{enumerate}
Then, for a series of meshes belonging to a regular family, we have the error bound:
\begin{equation}
  \normof{ u - u_h }_{H^m(\Omega)} \leq
    C \dfrac{\sup\limits_{\bxi \in \Oref} \absof{ \det\textbf{J}(\bxi)} }
    {\inf\limits_{\bxi \in \Oref} \absof{ \det\textbf{J}(\bxi) } }
    h^{p+1-m}||u||_{H^{p+1}(\Omega)}
    \label{curve_bound}
\end{equation}
wherein $\textbf{J}(\xi) = \bnabla \bx_e(\xi)$ is the Jacobian matrix and $C > 0$ is a constant independent of the mesh size $h$ but dependent on the family shape regularity $\sigma_0$ as well as the constants $c_i$ for $i = 0, \ldots, p+1$ appearing in Cond. (3.3.3).  Here, we take the mesh size $h_e$ of each curvilinear element $\Ophys$ to be that of the corresponding linear element $\Olin$ (as discussed in \S \ref{isoparametric}) and the global mesh size to be $h = \max_e h_e$.
From this, we see that for a curvilinear mesh to exhibit similar convergence rates to a linear mesh, several criteria must hold. 
First, as before, the underlying linear elements must be shape regular (Cond. \ref{classical_curve_cond}.1). 
However, we must also ensure that the higher-order mapping is invertible (Cond. \ref{classical_curve_cond}.2), and that its derivatives are bounded (Cond. \ref{classical_curve_cond}.3). 
Put simply, we must ensure that the curvilinear elements are not \textbf{\textit{too}} curvilinear.

\subsection{ Element Distortion and Quality Metrics} \label{metrics_review}
The results of the previous section illuminate some important considerations regarding the effect of element shape on the convergence rates of $p$-version finite element methods.
From Cond. \ref{classical_curve_cond}.2 and Cond. \ref{classical_curve_cond}.3, there is a clear need to quantify the magnitude of element distortion, as there is a direct effect on element quality. 
This need has led to the development of so-called \textbf{\emph{element distortion metrics}} and \textbf{\emph{element quality metrics}}. 
The precise definitions of these terms is often quite nebulous, and can often vary from application to application. 
However, for the purposes of this paper we say element distortion metrics quantify the difference between an arbitrary element $\Ophys$, and some ideal element $\Omega_{ideal}$.
Element quality metrics, then, are are simply taken as the inverse of element distortion metrics.
That is, element quality will increase as distortion decreases, and vice-versa.

In general, since element distortion metics and element quality metrics are closely related, we will refer to both with the umbrella term \textbf{\emph{element metrics}}.
To motivate the need for the novel element metrics presented in this paper, we review the existing curvilinear element metrics currently in use in the literature.
We then argue that none of the existing element metrics are sufficient for guaranteeing optimal convergence rates of the $p$-version finite element method over curvilinear elements.

We begin our study of curvilinear element metrics by recognizing that Eq. \eqref{curve_bound} contains the term:
\begin{equation*}
  1 \leq \dfrac{\sup\limits_{\bxi \in \Oref} \absof{ \det\textbf{J} } } 
                     {\inf\limits_{\bxi \in \Oref} \absof{ \det\textbf{J} } } < \infty
\end{equation*}
If the mapping $\xphys$ becomes singular, then ${\inf\limits_{\bxi \in \Oref} \absof{ \det\textbf{J} } } = 0$, and the error bounds in Eq. \eqref{curve_bound} will tend towards infinity.
It is perhaps due to this observation that the overwhelming majority of curvilinear quality metrics are based on some measure of the Jacobian matrix. 
Of these Jacobian based quality metrics, perhaps the most commonly used is the \boldit{scaled Jacobian} \cite{dey_curvilinear_1999,persson_curved_2008}, defined as:
\begin{equation} \label{scaled_jacobian}
  J_S = \dfrac{\inf\limits_{\bxi \in \Oref}|\det\textbf{J}|} 
                     {\sup\limits_{\bxi \in \Oref}|\det\textbf{J}| }
\end{equation}
From Eq. \eqref{scaled_jacobian}, it is readily apparent that $0 \leq J_S \leq 1$, with element quality increasing as $J_S \rightarrow 1$. 
For a linear element, the Jacobian is constant across the element, and the metric is identically unity. 
For a singular or inverted element, $\inf\limits_{\bxi \in \Oref}|\det\textbf{J}| = 0$, and the metric is zero. 

\begin{figure}[t!]
\centering
  \includegraphics[width=0.33\textwidth]{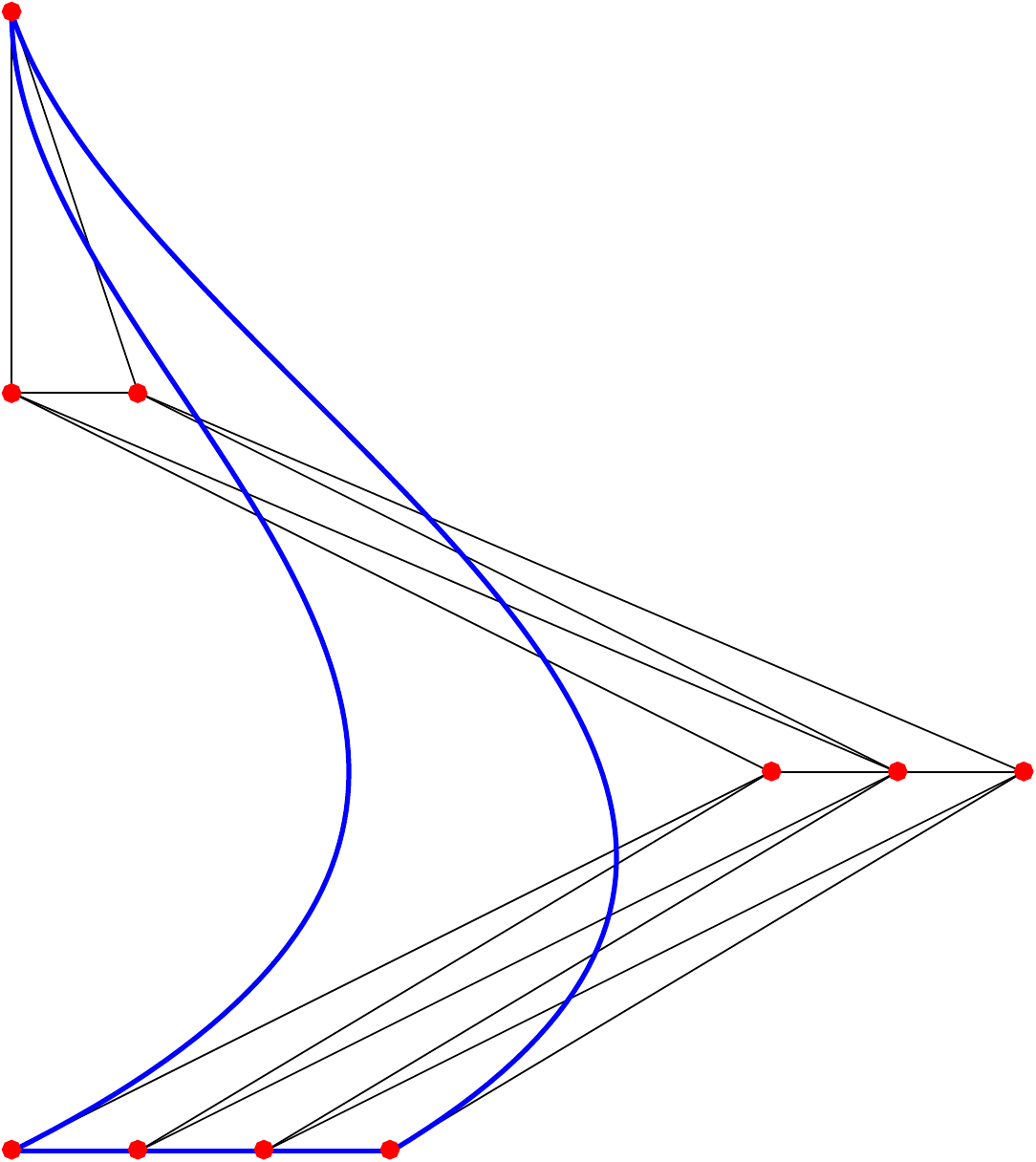}
  \caption{Highly distorted triangular element with a scaled Jacobian of $J_S = 1$.}
  \label{fig:bad_tri}
\end{figure}

Besides the scaled Jacobian, there have been other proposed higher-order quality metrics all based on some measure of the Jacobian matrix,
both for traditional higher-order finite elements \cite{gargallo-peiro_distortion_2015,gargallo-peiro_optimization_2015,george_construction_2012,lamata_quality_2013,poya_unified_2016,roca_defining_2011,xie_generation_2013}, 
and for IGA \cite{cohen_analysis-aware_2010,escobar_new_2011,speleers_optimizing_2015,xu_optimal_2013,xu_high-quality_2014,zhang_solid_2012}.  
Despite the wide array of metrics currently in use, we are not aware of any work relating bounds on these metrics to bounds on higher-order derivatives.
Thus, to our knowledge, none of the existing quality metrics are sufficient for guaranteeing that Cond. \ref{classical_curve_cond}.3 holds.
Because of this, we argue that existing curvilinear element metrics are insufficient for guaranteeing that an arbitrary curvilinear element is well-suited for finite element analysis.
To illustrate a particularly egregious example, consider the cubic Bernstein-B\'{e}zier element shown in Fig. \ref{fig:bad_tri}, which has a scaled Jacobian of $J_S = 1$.
While not all Jacobian based quality metrics  will indicate that this element is of good quality, 
it is troubling that the most commonly used quality metric for curvilinear elements cannot distinguish between this highly skewed element and a purely linear triangle.
We hope this example motivates the need for further study of curvilinear element distortion metrics, and we identify two key challenges to be addressed by the present work in Section \ref{validity_metrics}.

\begin{enumerate}
\item
  {
  \textbf{Bounds on the Jacobian Determinant of Rational Elements}
  
    In general, the Jacobian determinant of a polynomial mapping of degree $p$ is itself a polynomial of degree $p' = d(p-1)$.
    As such, bounding the Jacobian determinant from above and below is difficult for higher-order polynomial elements, and the  task is even more difficult for elements defined by a rational mapping.
    We note that computable bounds have been established for polynomial elements \cite{johnen_geometrical_2013}, but are unaware of any analogous bounds for rational elements.
  }
\item
{
  \textbf{Bounds on the Higher-Order Derivatives of the Parametric Mapping}

We are not aware of any element metrics that quantify the magnitude of higher--order partial 
derivatives of the parametric mapping $\xphys:\Oref \rightarrow \Ophys$.
Furthermore, we are not aware of any attempts to show that bounding any existing metrics implies bounds on higher--order derivatives\footnote{
In their original paper, Ciarlet and Raviart do provide conditions to ensure boundedness of the higher--order derivatives for certain classes of elements \cite{ciarlet_interpolation_1972}. 
However, these conditions are too restrictive to be used effectively with modern automated meshing algorithms.}.
}
\end{enumerate}



\section{Interpolation Theory for Rational \BB Elements} \label{interp_theory}
With the necessary preliminaries established, we now turn our attention to the novel contributions of the present work.
In this section, we present error estimates for rational \BB elements of simplicial or tensor product construction. 
The analysis follows closely the work of Bazilevs et. al. \cite{bazilevs_isogeometric_2006}, wherein interpolation error bounds were derived for IGA using NURBS.  For simplicity, we only analyze the $L^2$ best approximation error associated with one element, though our results extend easily to best approximation error in Sobolev norms (using the techniques outlined in \cite{bazilevs_isogeometric_2006}) as well as meshes comprised of several elements using suitable quasi-interpolation operators (e.g., Cl\'{e}ment-type interpolation operators in the context of standard $C^0$-continuous finite element approximations \cite{clement1975approximation}).

Let us consider a rational \BB element $\Ophys$ with corresponding reference element in parametric space $\Oref$.
For simplicial elements, we denote the space of approximation functions of degree $p$ over the physical curvilinear element as:
\begin{equation*}
  \mathcal{S}^h_p(\Ophys) := \setof{u_h \in L^2(\Ophys):w(u_h \circ \xphys) \in \mathscr{P}^p \of{ \Oref } }
\end{equation*}
where $\mathscr{P}^p \of{ \Oref }$ denotes the space of polynomials of degree $p$ and $w \in \mathscr{P}^p \of{ \Oref }$ denotes the weighting function defined over the reference element.
For tensor product elements, we denote the space of tensor product approximation functions of degree $\bp = \setof{p,...,p}$ over the physical curvilinear element as:
\begin{equation*}
  \mathcal{S}^h_{\bp}(\Ophys) := \setof{u_h \in L^2(\Ophys):w(u_h \circ \xphys) \in \mathscr{Q}^{\bp} \of{ \Oref } }
\end{equation*}
where $\mathscr{Q}^{\bp} \of{ \Oref }$ denotes the space of tensor product polynomials of degree $\bp$ and  $w \in \mathscr{Q}^{\bp} \of{ \Oref }$ denotes the weighting function defined over the reference element.

\begin{figure}[t]
  \centering
  \includegraphics[width=0.4\textwidth]{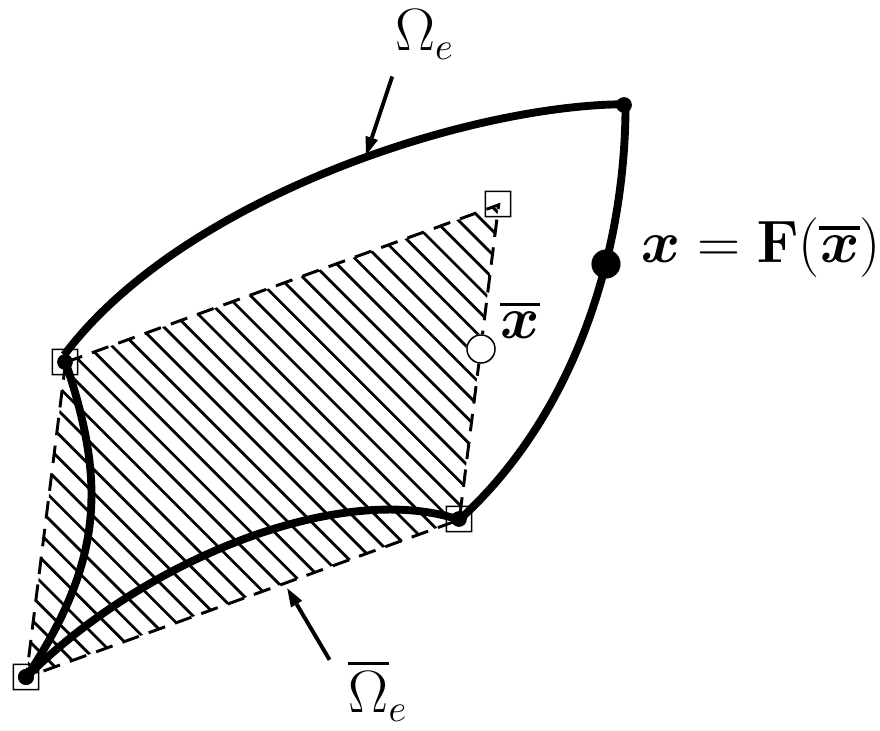}
  \caption{Mapping $\bF$ from the linear element $\Olin$ to the curved element $\Ophys$.}
  \label{mappingF}
\end{figure}

Furthermore, we note that the physical element has corresponding linear element $\Olin$, defined by an affine mapping $\xlin:\Oref \rightarrow \Olin$ as discussed in Subsection \ref{isoparametric}.
We then define a {warping function} $\bF$ that maps the linear element to the curvilinear element, $\bF: \Olin \rightarrow \Ophys$,
and we note that the mapping $\xphys$ is simply the compositions of these two mappings, $\xphys = \bF \circ \xlin $.
See Fig. \ref{mappingF} for an illustration of the mapping $\bF$ for a quadrilateral \BB element.
We can then derive error bounds for the element $\Omega_e$ in terms of the mapping $\bF$. 
%
%
\begin{thm} \label{error_bounds}
  Over the physical element $\Ophys$, there exists a constant $c_s = c_s(\sigma_e,p)$, dependent only on
      the shape regularity of $\Olin$ and the polynomial degree,
  such that for all $u \in H^{p+1}(\Ophys)$, there is an approximation function $u_h \in \mathcal{S}^h_p(\Omega_e)$ if $\Omega_e$ is a simplicial element or $u_h \in \mathcal{S}^h_\bp(\Omega_e)$ if $\Omega_e$ is a tensor product element satisfying the estimate:
  \begin{equation} \label{rational_error_bound}
    || u - u_h||_{L^2(\Ophys)} \leq
    c_sh_e^{p+1}
    c_d(\bF )
    c_v( \bF, w, u)
  \end{equation}
  where:
  \begin{equation*}
    c_d \of{ \bF } = \linf{\det \bnabla_{\xvlin} \bF^{1/2}}{\Olin}  \linf{\det \bnabla_{\xvlin} \bF^{-1/2}}{\Olin} 
  \end{equation*}
  and:
  \begin{equation*}
    c_v(\bF,w,u) =
    \sum\limits_{k=0}^{p+1}\sum\limits_{j=0}^{k}
    c_{\alpha} \of{\bF,j,k}
    \linf{\dfrac{1}{w \circ \xlin^{-1}}}{\Olin}
    \linf{\nabla_{\xvlin}^{p+1-k} \left( w \circ \xlin^{-1} \right)}{\Olin}
    |u|_{H^j(\Ophys)}
  \end{equation*}
  wherein:
  \begin{equation*}
    c_{\alpha} \of{\bF,j,k} =
    \sum\limits_{ \substack{i_1 + i_2 + ... + i_k = j \\ i_1 + 2i_2 + ... + ki_k = k} }
   \of{\linf{\bnabla_{\xvlin} \bF}{\Olin} }^{i_1}
   \of{\linf{\bnabla^2_{\xvlin} \bF}{\Olin} }^{i_2}...
   \of{\linf{\bnabla^k_{\xvlin} \bF}{\Olin} }^{i_k}
  \end{equation*}
  \label{rational_bh}
\end{thm}
\begin{proof}
  By definition, we have:
  \begin{equation}
    \int_{ \Ophys } u^2 d\bx =
    \int_{ \Olin } (u \circ \bF )^2 \det \bnabla_{\xvlin} \bF d\xvlin
  \end{equation}
  Therefore, $\normof{u}_{L^2\of{\Ophys}} = \normof{\det \bnabla_{\xvlin} \bF^{1/2} u \circ \bF}_{L^2 \of{ \Olin } }$, and as a consequence, we can bound:
  \begin{equation}
    \ltwo{u}{\Ophys} \leq
    \linf{ \det \bnabla_{\xvlin} \bF^{1/2} }{ \Olin }
    \ltwo{ u \circ \bF }{ \Olin }
  \end{equation}
  Consequently for each approximation function $u_h \in \mathcal{S}^h_p(\Ophys)$ or $u_h \in \mathcal{S}^h_\bp(\Ophys)$, it follows that:
  \begin{equation}
    \ltwo{u-u_h}{\Ophys} \leq
    \linf{ \det \bnabla_{\xvlin} \bF^{1/2} }{ \Olin }
    \ltwo{ u \circ \bF - u_h \circ \bF }{ \Olin }
  \end{equation}
  Since the weighting function is bounded from above and below, it further follows that
  \begin{equation}
    \ltwo{u-u_h}{\Ophys} \leq
    \linf{ \det \bnabla_{\xvlin} \bF^{1/2} }{ \Olin }
    \linf{\dfrac{1}{w \circ \xlin^{-1}}}{\Olin}
    \ltwo{ \left( w \circ \xlin^{-1} \right) \left( u \circ \bF - u_h \circ \bF \right) }{ \Olin }
  \end{equation}
  Now, we note that by construction of $u_h$, the function $ w \of{ u_h \circ \xphys } $ is an arbitrary polynomial (or tensor product polynomial) of degree $p$.
  Furthermore, we recognize that $ w \of{ u_h \circ \xphys } = w \of{ u_h \circ \bF \circ \xlin } $.
  Thus, because the mapping $\xlin$ is purely affine, the function $ \left( w \circ \xlin^{-1} \right) \of{ u_h \circ \bF }$ is similarly an arbitrary polynomial (or tensor product polynomial) of degree $p$. 
  Therefore, by the classical Bramble--Hilbert lemma, we may select $u_h$ such that the following inquality holds:
  \begin{equation}
    \ltwo{ \left( w \circ \xlin^{-1} \right) \left( u \circ \bF - u_h \circ \bF \right) }{ \Olin } \leq
    c_1 h^{p+1}
    \absof{ \left( w \circ \xlin^{-1} \right) \of { u \circ \bF } }_{H^{p+1} \of{ \Olin } }
  \end{equation}
  where $c_1 = c_1 \of{ p, \sigma_e} $ is a constant that depends only on the polynomial degree $p$ and shape regularity $\sigma_e$ of the linear element $\Olin$.  
  We must now bound the seminorm $\absof{ w \of{ u \circ \bF } }_{ H^{p+1} \of{ \Olin } } $ appearing in the above estimate by an analagous norm over the physical element $ \Ophys $.
  To do so, we first recognize that: 
  \begin{equation}
    \absof{ \left( w \circ \xlin^{-1} \right) \of { u \circ \bF } }_{H^{p+1} \of{ \Olin } }  \leq
    \sum \limits_{k=0}^{p+1}
         {
           \linf{ \bnabla_{\xvlin}^{p+1-k} \left( w \circ \xlin^{-1} \right) }{ \Olin }
           \absof{ u \circ \bF }_{ H^k \of{ \Olin } }
         }
  \end{equation}
  It remains to bound the seminorms $ \absof{ u \circ \bF }_{ H^k \of{ \Olin } } $.
  We may easily obtain control of the $H^1$-seminorm using the estimate:
  \begin{equation}
    \begin{split}
      \absof{ u \circ \bF }_{ H^1 \of{ \Olin } } &=
      \of{ \int_{\Olin} \bnabla_{\xvlin} \of{ u \circ \bF } \cdot \bnabla_{\xvlin} \of{ u \circ \bF } d\xvlin }^{ 1/2 } \\
      &=\of{ \int_{\Ophys} \bnabla_{\xvlin} \bF  \bnabla_{\bx}u \cdot \bnabla_{\xvlin} \bF  \bnabla_{\bx} u \det \bnabla_{\bx} \bF^{-1} d\bx }^{ 1/2 } \\
      &\leq \linf{ \bnabla_{\xvlin} \bF }{\Olin} \linf{ \left(\det \bnabla_{\bx} \bF^{-1}\right)^{1/2 } }{ \Ophys } \absof{ u }_{H^1 \of{ \Ophys } } \\
      &= \linf{ \bnabla_{\xvlin} \bF }{\Olin} \linf{ \det \bnabla_{\xvlin} \bF^{-1/2} }{ \Olin } \absof{ u }_{H^1 \of{ \Ophys } } 
    \end{split}
  \end{equation}
  To obtain control of the higher--order seminorms, we simply recurse on the previous estimate, as is done in \cite{bazilevs_isogeometric_2006}, resulting in the estimate:
  \begin{equation}
    \absof{ u \circ \bF }_{ H^k \of{ \Olin } } \leq
    c_2
    \linf{ \det \bnabla_{\xvlin} \bF^{-1/2} }{ \Olin }
    \sum\limits_{j=0}^k
               {
                 c_{\alpha} \of{\bF,j,k}
                 \absof{ u }_{H^j \of{ \Ophys } }
               }
  \end{equation}
  where $c_2 = c_2 \of{ p, \sigma_e} $ is again a constant that depends only on the polynomial degree $p$ of the basis and the shape regularity of $\Olin$.
  Thus, letting $c_s = c_1 c_2 $, we arrive at the bound presented in Theorem \ref{error_bounds}.
\end{proof}
Theorem \ref{error_bounds} gives valuable theoretical insight into the convergence behavior of rational curvilinear \BB elements, 
as it clearly delineates the effect of the linear shape quality ($c_s$) and curvilinear shape quality ($c_d$ and $c_v$) on the interpolation error bounds.
However, its utility is somewhat limited as the bound is given in terms of the warping function $\bF$, whereas \BB elements are defined by the mapping $\xphys$. 
Since the mapping $\xphys$ is given explicitly by the control points and weights, it is desirable to derive sufficient conditions based on this mapping instead. 
%
We begin by deriving bounds on the gradients $\linf{ \bnabla^k_{\xvlin} \bF}{\Olin}$  in terms of the gradients  $\linf{ \bnabla^k_{\bxi} \xphys}{\Oref} $ (Theorem \ref{gradF_O1}),
as well as bounds on the gradients $\linf{ \bnabla^k_{\xvlin} \left( w \circ \xlin^{-1} \right)}{\Olin}$ in terms of the gradients $\linf{ \bnabla^k_{\bxi} w}{\Oref} $ (Theorem \ref{gradw_O1}).
We then use the results of these theorems to prove Theorem I. 
%
%
%
%
%
%

%
%
%
%
\begin{thm} \label{gradF_O1}
  There exists some constant $c_{r} = c_{r}\of{\sigma_e}$, dependent only on the shape regularity of $\Olin$, such that:
  \begin{equation} \label{k_deriv_bound}
    \linf{\boldsymbol{\nabla}_{\xvlin}^k \bF}{\Olin} \leq 
    c_{r}^k
    \linf{ \boldsymbol{\nabla}_{\bxi}^k \xphys}{ \Oref }
    \of{ \dfrac{1}{h_e} }^k
  \end{equation}

\end{thm}

\begin{proof}
  We first recognize that the warping function $\bF$ can be written as a composition of $\xphys$ and $\xlin^{-1}$, viz:
  \begin{equation*}
    \bF \of{ \xvlin } = \xphys \of { \xlin^{-1} \of{ \xvlin } }
  \end{equation*}
  Thus, the gradient of $\bF$ can be written:
  \begin{equation} \label{gradient_composition}
    \bnabla_{\xvlin}\bF =
    \Bigg{[}  \bnabla_{\bxi} \xphys  \Bigg{]}
    \Bigg{[}  \bnabla_{\bxi} \xlin \Bigg{]}^{-T}
  \end{equation}
  We note that since $\Olin$ is a linear element, the gradient $\bnabla_{\bxi} \xlin $ is constant 
  and $ \normof{  \bnabla_{\bxi} \xlin   }_{L^\infty \of{\Oref} } $ is bounded from below by the radius of the element incircle, $\rho_e$, viz.: 
  \begin{equation} \label{first_deriv_bound}
    \rho_e \leq
    \normof{  \bnabla_{\bxi} \xlin   }_{L^\infty \of{\Oref} } 
  \end{equation}
  and we can similarly bound the norm of the inverse mapping by:
  \begin{equation} \label{inv_deriv_bound}
    \normof{ \bnabla_{\xvlin} \of{ \xlin^{-1} }  }_{L^\infty \of{ \Oref } } \leq
    \dfrac{1}{\rho_e}
  \end{equation}
  Now, let us define a constant $c_{r} = c_{r}(\sigma_e)$, which is independent of mesh size $h_e$ but dependent on the shape regularity of $\Olin$, 
such that:
  \begin{equation*}
    c_{r} \geq \sigma_e = \dfrac{h_e}{\rho_e}
  \end{equation*}
  Then, we can rewrite Eq. \eqref{inv_deriv_bound} as:
  \begin{equation*}
    \normof { \bnabla_{\xvlin} \of{ \xlin^{-1} } }_{L^\infty \of{ \Oref } } \leq
    c_{r}\dfrac{1}{h_e}
  \end{equation*}
  Then, from Eq. \eqref{gradient_composition} and Eq. \eqref{first_deriv_bound}, we can bound the norm of $\bnabla_{\xvlin}\bF$ by:
  \begin{equation*}
    \linf{\bnabla_{\xvlin}\bF}{\Olin} \leq 
    c_{r}
    \normof{
      \bnabla_{\bxi} \xphys 
    }_{L^{\infty} \of{ \Oref } }
    \of{ \dfrac{1} {h_e} }
  \end{equation*}
  Recursing on this process, we can bound the magnitude of the $k^{th}$ total derivative by:
  \begin{equation} \label{k_deriv_bound}
    \linf{\bnabla^k_{\xvlin}\bF}{\Olin} \leq 
    c_{r}^k
    \linf{ \boldsymbol{\nabla}_{\bxi}^k \xphys}{ \Oref }
    \of{ \dfrac{1}{h_e} }^k
  \end{equation}
\end{proof}
%
%
%
%
\begin{thm}\label{gradw_O1}
  The following inequality holds:
  \begin{equation} \label{k_deriv_bound}
    \linf{\bnabla_{\xvlin}^{k} \left( w \circ \xlin^{-1} \right) }{\Olin} \leq 
    c_{r}^k
    \linf{ \bnabla^k_{\bxi} w}{ \Oref }
    \of{ \dfrac{1}{h_e} }^k
  \end{equation}
  where $c_{r} = c_{r}\of{\sigma_e}$ is the constant from Theorem  \ref{gradF_O1}.
\end{thm}
\begin{proof}
  The proof of Theorem \ref{gradw_O1} is identical to the proof for Theorem \ref{gradF_O1}.
\end{proof}



\noindent
\begin{minipage}{1.02\textwidth}
  \begin{mdframed}
    \begin{customthm}{I} \label{error_bounds_x}
      Over the physical element $\Ophys$, there exists a constant $C_{shape} = C_{shape}(\sigma_e,p)$, dependent only on
      the shape regularity of $\Olin$ and the polynomial degree,
      such that for all $u \in H^{p+1}(\Ophys)$, there is an approximation function $u_h \in \mathcal{S}^h_p(\Omega_e)$ if $\Omega_e$ is a simplicial element or $u_h \in \mathcal{S}^h_\bp(\Omega_e)$ if $\Omega_e$ is a tensor product element satisfying the estimate:
      \begin{equation} \label{rational_error_bound}
        || u - u_h||_{L^2(\Ophys)} \leq
        C_{shape}h_e^{p+1}
        C_{det}( \xphys )
        C_{var}( \xphys, w, u)
      \end{equation}
      where:
      \begin{equation*}
        C_{det} \of{ \xphys } = \linf{\det\bnabla_{\bxi}\xphys^{1/2}}{\Oref}  \linf{\det\bnabla_{\bxi}\xphys^{-1/2}}{\Oref} 
      \end{equation*}
      and:
      \begin{equation*}
        C_{var}(\xphys,w,u) =
        \sum\limits_{k=0}^{p+1}\sum\limits_{j=0}^{k}
        \alpha_{j,k}(\xlin)
        \linf{\dfrac{1}{w}}{\Oref}
        \of{ \dfrac{ \linf{\nabla_{\bxi}^{p+1-k}w}{\Oref} }{ h_e^{p+1-k} } }
        |u|_{H^j(\Ophys)}
      \end{equation*}
      wherein:
      \begin{equation*}
        \alpha_{j,k}(\xphys) \leq
        \sum\limits_{ \substack{i_1 + i_2 + ... + i_k = j \\ i_1 + 2i_2 + ... + ki_k = k} }
        \of{ \dfrac{\linf{ \bnabla_{\bxi}      \xphys }{\Oref}}{h_e}     }^{i_1}
        \of{ \dfrac{ \linf{ \bnabla_{\bxi}^2 \xphys }{\Oref} }{ h_e^2 } }^{i_2}...
        \of{ \dfrac{\linf{ \bnabla_{\bxi}^k \xphys }{\Oref} }{ h_e^k }  }^{i_k}
      \end{equation*}
    \end{customthm}
  \end{mdframed}
\end{minipage}
\begin{proof}
  Proving Theorem \ref{error_bounds_x} amounts to simply bounding the constants $c_{d}$ and $c_{v}$ from Theorem \ref{error_bounds} by analogous constants in terms of the mapping $\xphys$. 
  We begin by recognizing that $\xphys = \bF \circ \xlin $, and by extension $\bF = \xphys \circ \xlin^{-1}$. 
  Thus, we can write the gradient of the mapping $\bF$ as $\bnabla_{\xvlin} \bF = \matof{ \bnabla_{\bxi} \xphys } \matof{ \bnabla_{\bxi} \xlin }^{-T} $
  and as a result,  the determinant of $\bnabla_{\xvlin} \bF$ can be written as:
  %
  %
  %
  \begin{equation}
    \det \bnabla_{\xvlin} \bF =  \dfrac{ \det \bnabla_{\bxi} \xphys }{ \det \bnabla_{\bxi} \xlin     } 
  \end{equation}
  allowing us to rewrite $c_{d} \of{ \bF } $ in terms of $\xlin$ and $\xphys$ as:
  \begin{equation}
    c_d \of{ \bF} =
    C_{det}\of{ \xlin, \xphys} =
    \linf{ \of{ \dfrac{ \det \bnabla_{\bxi} \xphys }{ \det \bnabla_{\bxi} \xlin     } }^{1/2} } {\Oref}  
    \linf{ \of{ \dfrac{ \det \bnabla_{\bxi} \xlin     }{ \det \bnabla_{\bxi} \xphys } }^{1/2} } {\Oref} 
  \end{equation}
  Then, recognizing that $\det \bnabla_{\bxi} \xlin$ is a constant, the above equations simplify immediately to: 
  \begin{equation*}
    c_d \of{ \bF} =
    C_{det} \of{ \xphys } =
    \linf{\det\bnabla_{\bxi} \xphys^{1/2}}{\Oref}  \linf{\det\bnabla_{\bxi} \xphys^{-1/2}}{\Oref}
  \end{equation*}
  %
  %
  %
  It remains then to bound the constant $c_{v}$ by a bound in terms of $\xphys$. 
  We begin by substituting the results of Theorem \ref{gradF_O1} into the expression $c_{\alpha}\of{\bF,j,k}$, which yields:
  \begin{equation*}
    \begin{split}
      c_{\alpha}\of{\bF,j,k} \leq
      \sum\limits_{ \substack{i_1 + i_2 + ... + i_k = j \\ i_1 + 2i_2 + ... + ki_k = k} }
       \of{  { c_{r}  \dfrac{\linf{ \bnabla_{\bxi}  \xphys }{\Oref}}{h_e} }    }^{i_1}
      &\of{  { c_{r}^2\dfrac{\linf{ \bnabla_{\bxi}^2 \xphys }{\Oref}}{h_e^2} }  }^{i_2}...\\
   ...&\of{  { c_{r}^k\dfrac{\linf{ \bnabla_{\bxi}^k \xphys }{\Oref}}{h_e^k} }  }^{i_k}
    \end{split}
  \end{equation*}
  Then, recognizing that $c_{r}^{i_1}c_{r}^{2i_2}...c_{r}^{ki_k} = c_{r}^k$, we can factor out the constant $c_{r}$, to arrive at the bound: 
  \begin{equation}\label{alphaF_O1}
    c_{\alpha}\of{ \bF,j,k } \leq c_{r}^k\alpha_{j,k}\of{\xphys}
  \end{equation}
  %
  %
  Finally, employing the results of Eq. \eqref{alphaF_O1} along with Theorem \ref{gradw_O1}, we can bound $c_{v}\of{ \bF, w, u }$ by:
  \begin{equation*}
    c_{v}( \bF,w,u) \leq 
    \sum\limits_{k=0}^{p+1}\sum\limits_{j=0}^{k}
    c_{r}^k\alpha_{j,k}(\xphys)
    \linf{\dfrac{1}{w}}{\Oref}
    c_{r}^{p+1-k}\of{ \dfrac{ \linf{\nabla_{\bxi}^{p+1-k}w}{\Oref} }{ h_e^{p+1-k} } }
    |u|_{H^j(\Ophys)}
  \end{equation*}
  which we can simplify to:
  \begin{equation*}
    c_{v}( \bF,w,u) \leq 
    c_{r}^{p+1} C_{var}(\xphys,w,u)
  \end{equation*}
  Thus, we arrive at the results of Theorem I, where $C_{shape} = c_s c_{r}^{p+1}$.
\end{proof}

The astute reader will recognize that the above results can be simplified even further. 
Indeed, we have left Theorem I in its current form intentionally, as it clearly demonstrates the dependence on $h_e$ in each part of the error bound.
This clear dependence on $h_e$ will prove useful in the next section. 
However, we similarly note that at times it is useful develop size invariant error bounds.
We further simplify the results of Theorem I to achieve error bounds that are independent of the element diameter $h_e$.

\noindent
  \begin{cor} \label{error_bounds_cor}
    Over the physical element $\Ophys$, there exists a constant $C_{shape} = C_{shape}(\sigma_e,p)$, dependent only on
      the shape regularity of $\Olin$ and the polynomial degree,
    such that for all $u \in H^{p+1}(\Ophys)$, there is an approximation function $u_h \in \mathcal{S}^h_p(\Omega_e)$ if $\Omega_e$ is a simplicial element or $u_h \in \mathcal{S}^h_\bp(\Omega_e)$ if $\Omega_e$ is a tensor product element satisfying the estimate:
    \begin{equation} \label{rational_error_bound}
      || u - u_h||_{L^2(\Ophys)} \leq
      C_{shape}
      C_{det}( \xphys )
      C'_{var}( \xphys, w, u)
    \end{equation}
    where:
    \begin{equation*}
      C_{det} \of{ \xphys } = \linf{\det\bnabla_{\bxi}\xphys^{1/2}}{\Oref}  \linf{\det\bnabla_{\bxi}\xphys^{-1/2}}{\Oref} 
    \end{equation*}
    and:
    \begin{equation*}
      C'_{var}(\xphys,w,u) =
      \sum\limits_{k=0}^{p+1}\sum\limits_{j=0}^{k}
      \alpha'_{j,k}(\xphys)
      \linf{\dfrac{1}{w}}{\Oref}
      \linf{\nabla_{\bxi}^{p+1-k}w}{\Oref}
      |u|_{H^j(\Ophys)}
    \end{equation*}
    wherein:
    \begin{equation*}
      \alpha'_{j,k}( \xphys ) =
      \sum\limits_{ \substack{i_1 + i_2 + ... + i_k = j \\ i_1 + 2i_2 + ... + ki_k = k} }
      \of{\linf{ \bnabla_{\bxi}   \xphys }{\Oref} }^{i_1}
      \of{\linf{ \bnabla_{\bxi}^2 \xphys }{\Oref} }^{i_2} ...
      \of{\linf{ \bnabla_{\bxi}^k \xphys }{\Oref} }^{i_k}
    \end{equation*}
  \end{cor}

%


%
%
%
\section{Shape Regular Families of Rational \BB Meshes} \label{sufficient_conditions}
The results of Section \ref{interp_theory} are convenient, as they clearly delineate the effect of element shape on interpolation error bounds.
Namely, $C_{shape}$ quantifies the contribution of \textbf{\textit{linear}} element shape regularity, while $C_{det}$ and $C_{var}$ quantify the contribution of the \textbf{\textit{higher--order mapping}} $\xphys$ and \textbf{\textit{weighting function}} $w$.
By writing the error bounds in terms of the mapping $\xphys$ we may work directly with the Bernstein basis functions defined over the reference element $\Oref$,
and the B\'{e}zier control points in physical space.
Furthermore, the results of Theorem I clearly delineate the dependence of the interpolation error bounds on the higher--order derivatives of the mapping in terms of the element size $h_e$.
In light of this, we desire to develop sufficient conditions for guaranteeing that rational \BB elements will preserve optimal convergence rates under refinement. 
That is, we desire to develop sufficient conditions for guaranteeing that:
\begin{equation}
  \ltwo{ u - u_h }{\Ophys} \leq Ch_e^{p+1}
\end{equation}
where $C$ is some constant independent of element size $h_e$ and the mapping $\xphys$.
Indeed, as discussed in Section \ref{review}, such sufficient conditions have been developed for finite elements based on \textbf{\textit{polynomial}} mappings.
In this section, we develop analogous sufficient conditions for \textbf{\textit{rational}} \BB elements. 
\begin{thm} \label{w_inv_bounded}
  Let us assume there exists a constant $c_w$ such that:
  \begin{equation} \label{w_bounded}
    \normof{ \derivnx{{\balpha}}{\bxi}w }_{L^\infty \of{ \Oref} } \leq
    c_w h_e^{\absof{\balpha}} \ \ \
    \forall \balpha : \absof{\balpha} \leq p+1 
  \end{equation}
  Then, there exists a constant $c_{inv} = c_{inv}(c_w,p)$ such that:
  \begin{equation*} 
    \normof{ \derivnx{{\balpha}}{\bxi} \of{ \dfrac{1}{w} } }_{L^{\infty} \of{ \Oref } } \leq c_{inv}h_e^{\absof{\balpha}}
  \end{equation*}
\end{thm}

\begin{proof}

  By application of the multivariate Fa\`{a} di Bruno's formula \cite{constantine_multivariate_1996}, we write the derivative of the inverse of the weighting function as:
  \begin{equation*}
    \derivnx{{\balpha}}{\bxi} \of{\dfrac{1}{w}} = 
    \sum\limits_{1\leq \absof{\bs} \leq \absof{ \balpha } }
    \of{ -1 }^{|\bs|}\dfrac{|\bs|!}{w^{1+|\bs|}}
    \sum\limits_{p({\balpha},\bs)}
               {\balpha}!
               \prod\limits_{j=1}^{\absof{\balpha}}
               \dfrac{ \of{ D_\bxi^{ \boldsymbol{ \ell }_j } w }^{k_j} }
                     {k_j! \of{ \boldsymbol{\ell}_j! }^{ k_j } }
  \end{equation*}
  where $p({\balpha},\bs)$ denotes the set:
  \begin{equation*}
    p({\balpha},\bs) = \bigcup\limits_{s=1}^{\absof{\balpha}} p_s({\balpha},\bs)
  \end{equation*}
  \begin{equation*}
    \begin{split}
      p_s \of{ {\balpha},\bs } = &
      \bigg{\{} \of{ k_1,...,k_{ \absof{ \balpha } }; \boldsymbol{\ell}_1,...,\boldsymbol{\ell}_{\absof{ \balpha } } } : \\
      & 0<k_i,  \bvec{0} < \boldsymbol{\ell}_1 < ...< \boldsymbol{\ell}_s, \\
      & \sum\limits_{i=1}^sk_i = \bs\ \textup{and} \ \sum\limits_{i=1}^sk_i\boldsymbol{\ell}_i = {\balpha} 
      \bigg{\}}
    \end{split}
  \end{equation*}
  Then, by Eq. \eqref{w_bounded}, we have:
  \begin{equation*}
    \normof{ \derivnx{{\balpha}}{\bxi} \of{ \dfrac{1}{w} } }_{L^\infty(\Omega)} \leq
    \sum\limits_{1\leq |\bs| \leq \absof{\balpha} }
    {
      \dfrac{|\bs|!}{w^{1+|\bs|}}
      \sum\limits_{p({\balpha},\bs)}
      {
        {\balpha}!
        \prod\limits_{j=1}^{\absof{\balpha}}
        \dfrac{ \of{ c_wh_e^{ \absof{\boldsymbol{\ell}_j } } }^{k_j}}
              {k_j! \of{ \boldsymbol{\ell}_j! } ^{k_j}}
      }
    }
  \end{equation*}
  Rearranging, and factoring out terms not dependent on $h_e$ yields:
  \begin{equation} \label{w_leq_Chk}
    \Lnorm \derivnx{{\balpha}}{\bxi}\dfrac{1}{w} \Rnorm_{L^\infty \of{ \Omega } } \leq
    \sum\limits_{1\leq |\bs| \leq \absof{\balpha} }
    \sum\limits_{p \of{ \balpha, \bs } }
    c_{prod}\of{\balpha,\bs}
    \prod\limits_{j=1}^{\absof{\balpha} }
    \of{ h_e^{ \absof{ \boldsymbol{ \ell }_\bj } } }^{k_j}
  \end{equation}
  where:
  \begin{equation*}
    c_{prod} \of{\balpha,\bs} =
    \dfrac{ \absof{ \bs }! \balpha! }{w^{\absof{\bs}+1}}
    \prod\limits_{j=1}^{\absof{\balpha}}
    \dfrac{c_{w}^{k_j} } { k_j! \of{ \boldsymbol{ \ell }_l! }^{ k_j } }
  \end{equation*}
  Finally, noting that:
  \begin{equation*}
    \prod\limits_{j=1}^{\absof{\balpha}}
    \of{ h_e^{ \absof{ \boldsymbol{\ell}_j } } }^{k_j} =
    h_e^{\Big{(} \sum\limits_{j=1}^{\absof{\balpha}}k_j|\boldsymbol{\ell}_j|\Big{)}}
  \end{equation*}
  and that by construction of $p({\balpha},\bs)$:
  \begin{equation*}
    \sum\limits_{j=1}^{\absof{\balpha}}k_j|\boldsymbol{\ell}_j| = \absof{\balpha}
  \end{equation*}
  Inequality \eqref{w_leq_Chk} becomes:
  \begin{equation*}
    \Lnorm \derivnx{{\balpha}}{\bxi}\dfrac{1}{w} \Rnorm_{L^\infty(\Oref)} \leq
    \sum\limits_{1\leq |\bs| \leq \absof{\balpha} }
    \sum\limits_{p({\balpha},\bs)}
    c_{prod} \of{\balpha, \bs}
    h_e^{\absof{\balpha}}
  \end{equation*}
  and the result of  Lemma \ref{w_inv_bounded} follows directly, where $c_{inv} = \sum\limits_{1\leq |\bs| \leq \absof{\balpha} } \sum\limits_{p({\balpha},\bs)}  c_{prod} \of{\balpha, \bs} $. 
\end{proof}
%
%
%
%
\begin{thm} \label{proj_thm}
  Let us assume there exists a constant $C_{proj}$ such that:
  \begin{equation} \label{xproj_bounded}
    \linf{D^{\balpha}_{\bxi}\xproj}{ \Oref} \leq C_{proj}h_e^{\absof{\balpha}} \ \ \ \forall \balpha:\absof{\balpha}=k, \ \ \ k \leq p+1
  \end{equation}
  Then, there exists a constant $C_{grad} = C_{grad}(C_{proj},p)$ such that:
  \begin{equation*}
    \linf{\boldsymbol{\nabla}_{\bxi}^k \xphys}{ \Oref} \leq C_{grad} h_e^k
  \end{equation*}
  \begin{equation*}
    \linf{\boldsymbol{\nabla}_{\bxi}^k w}{ \Oref} \leq C_{grad} h_e^k
  \end{equation*}
\end{thm}
\begin{proof}
  We immediately recognize that since $w = \of{ \xproj }_{d+1}$, the inequality 
  \begin{equation} \label{wgrad_bounded}
     \linf{D^{\balpha}_{\bxi} w}{ \Oref} \leq C_{proj} h_e^k 
  \end{equation}
  holds by definition.
  It remains then to show that the bounds on the derivatives of the projective mapping $\xproj$ imply bounds on the derivatives of the physical mapping $\xphys$.
  We begin by writing the derivative $\derivnx{{\balpha}}{\bxi} \xphys $ by the multi-variate product rule as:
  \begin{equation} \label{prod_rule}
    \derivnx{{\balpha}}{\bxi}\xphys =
    \derivnx{{\balpha}}{\bxi} \of{ \dfrac{\xproj}{w} }=
    \sum\limits_{\bk \in I^{{\balpha}}} {{{\balpha}}\choose{\bvec{k}}}
    \derivnx{{\balpha}-\bk}{\bxi}\of{ \xproj } 
    \derivnx{\bk}{\bxi} \of{ \dfrac{1}{w} }
  \end{equation}
  Then, from Eq. \eqref{wgrad_bounded} and Theorem \ref{w_inv_bounded},  it is readily seen that:
  \begin{equation*}
    \linf{\derivnx{{\balpha}}{\bxi} \of{ \dfrac{1}{w} } }{ \Oref} \leq c_{inv} h_e^{\absof{\balpha}} \ \ \ \forall \ {\balpha} : \absof{\balpha} = k
  \end{equation*}
  Taking  this result, along with Eq. \eqref{xproj_bounded}, and substituting into Eq. \eqref{prod_rule} yields:
  \begin{equation*}
    \linf{\derivnx{{\balpha}}{\bxi} \xphys }{ \Oref} \leq
    \sum\limits_{\bk \in I^{{\balpha}}} {{{\balpha}}\choose{\bk}}
    C_{proj} h_e^{\absof{\balpha}-|\bk|}
    c_{inv} h_e^{|\bk|}
  \end{equation*}
  which reduces immediately to:
  \begin{equation} \label{drdxi_bounded}
    \linf{ \derivnx{{\balpha}}{\bxi} \xphys}{\Oref } \leq C_{proj}c_{inv}h_e^{\absof{\balpha}} \ \ \  \forall \ {\balpha} : \absof{\balpha} = k
  \end{equation}
  Finally, we recognize that if Eq. \eqref{wgrad_bounded} and Eq. \eqref{drdxi_bounded} hold for every derivative of order $k$, then there exists some $C_{grad} = C_{grad} \of{C_{proj},c_{inv}}$ such that:
  \begin{equation} 
    \linf{ \bnabla^{k}_{\bxi} w}{\Oref } \leq C_{grad}h_e^{k}
  \end{equation}
  \begin{equation} 
    \linf{ \bnabla^{k}_{\bxi} \xphys}{\Oref } \leq C_{grad}h_e^{k}
  \end{equation}
  which are exactly the results of Theorem \ref{proj_thm} that we set out to prove.
\end{proof}
%
%
%
%
\noindent
\begin{minipage}{\textwidth}
  \begin{mdframed}
    \begin{customthm}{II} \label{suff_cond_intro}
      For a rational \BB element $\Ophys$ of degree $p$, let the following conditions hold: 
      \begin{enumerate}[label={\textup{ Cond. (\ref{suff_cond_intro}.\arabic*)}},align=left]
      \item{ There exists a constant $C_{max}$ such that:
        \begin{equation*}
          \linf{\det\bnabla_{\bxi} \xphys^{1/2}}{\Oref}  \linf{\det\bnabla_{\bxi} \xphys^{-1/2}}{\Oref} \leq C_{max}
        \end{equation*}
      } 
      \item{There exists a constant $C_{proj}$ such that:
        \begin{equation*}
          \linf{D^{\balpha}_{\bxi} \xproj}{ \Oref} \leq C_{proj}h_e^{\absof{\balpha}} \hspace{10pt} \forall \balpha: \absof{\balpha} \leq p+1
        \end{equation*}
      } 
      \item{There exists a constant $C_{weight}$ such that:
        \begin{equation*}
          \linf{\dfrac{1}{w}}{\Oref} \leq C_{weight}
        \end{equation*}
      } 
      \end{enumerate}
      Then, there exists a constant $C$ only dependent on the element through $C_{max}$, $C_{proj}$, $C_{weight}$, $\sigma_e$, and $p$ such that, for all $u \in H^{p+1}(\Ophys)$, there is an approximation function $u_h \in \mathcal{S}^h_p(\Omega_e)$ if $\Omega_e$ is a simplicial element or $u_h \in \mathcal{S}^h_\bp(\Omega_e)$ if $\Omega_e$ is a tensor product element satisfying:
      \begin{equation*}
        \ltwo{u-u_h}{\Ophys} \leq Ch_e^{p+1} \| u \|_{H^{p+1}(\Ophys)}
      \end{equation*}      
    \end{customthm}
  \end{mdframed}
\end{minipage}

%
%
%
\noindent
  \begin{proof}
    From Theorem I, we immediately see that if Cond. (\ref{suff_cond_intro}.1) holds, then $C_{det} \leq C_{max}$.
    Next, if Cond. (\ref{suff_cond_intro}.2) holds, we can bound $C_{var}$ using the results of Theorem \ref{proj_thm} by:
      \begin{equation*}
        C_{var}(\xphys,w,u) \leq
        \sum\limits_{k=0}^{p+1}\sum\limits_{j=0}^{k}
        \alpha_{j,k}(\xphys)
        \linf{\dfrac{1}{w}}{\Oref}
        \of{ \dfrac{C_{grad}h_e^{p+1-k} }{ h_e^{p+1-k} } }
        |u|_{H^j(\Ophys)}
      \end{equation*}
      wherein:
      \begin{equation*}
        \alpha_{j,k}(\xphys) \leq
        \sum\limits_{ \substack{i_1 + i_2 + ... + i_k = j \\ i_1 + 2i_2 + ... + ki_k = k} }
        \of{ \dfrac{ C_{grad}h_e   }{h_e}     }^{i_1}
        \of{ \dfrac{ C_{grad}h_e^2 }{ h_e^2 } }^{i_2}...
        \of{ \dfrac{ C_{grad}h_e^k }{ h_e^k }  }^{i_k}
      \end{equation*}
      which reduces immediately to:
      \begin{equation*}
        C_{var}(\xphys,w,u) \leq
        C_{grad}
        \linf{\dfrac{1}{w}}{\Oref}
        \sum\limits_{k=0}^{p+1}
        \sum\limits_{j=0}^{k}
        \sum\limits_{ \substack{i_1 + i_2 + ... + i_k = j \\ i_1 + 2i_2 + ... + ki_k = k} }
        C_{grad}^j
        |u|_{H^j(\Ophys)}
      \end{equation*}
      If Cond. (\ref{suff_cond_intro}.3) holds, we have that:
      \begin{equation*}
        C_{var}(\xphys,w,u) \leq
        C_{grad}
        C_{weight}
        \sum\limits_{k=0}^{p+1}
        \sum\limits_{j=0}^{k}
        \sum\limits_{ \substack{i_1 + i_2 + ... + i_k = j \\ i_1 + 2i_2 + ... + ki_k = k} }
        C_{grad}^j
        |u|_{H^j(\Ophys)}
      \end{equation*}      
      and since $|u|_{H^j(\Ophys)} \leq \|u\|_{H^{p+1}(\Ophys)}$ for all $j \leq p+1$, it follows that:
      \begin{equation*}
        C_{var}(\xphys,w,u) \leq
        \left( C_{grad}
        C_{weight}
        \sum\limits_{k=0}^{p+1}
        \sum\limits_{j=0}^{k}
        \sum\limits_{ \substack{i_1 + i_2 + ... + i_k = j \\ i_1 + 2i_2 + ... + ki_k = k} }
        C_{grad}^j \right)
        \|u\|_{H^{p+1}(\Ophys)}
      \end{equation*}
	Thus, by Theorem 1 and the above bounds, we arrive at the results of Theorem II, wherein
	\begin{equation}
	C = C_{shape} C_{max} \left( C_{grad}
        C_{weight}
        \sum\limits_{k=0}^{p+1}
        \sum\limits_{j=0}^{k}
        \sum\limits_{ \substack{i_1 + i_2 + ... + i_k = j \\ i_1 + 2i_2 + ... + ki_k = k} }
        C_{grad}^j \right)
	\end{equation}
	is the desired constant.
  \end{proof}

Now let $M = \{\mathcal{M}_i\}_{i=1}^M$ be a family of rational \BB meshes of the same polynomial degree $p$.  We say that $M$ is a {\bf \emph{shape regular family}} if the following two conditions hold:
\begin{enumerate}[label={\textup{ Cond. (R.\arabic*})},align=left]
\item Each element $\Ophys$ of each mesh $\mathcal{M}_i$ in the family $M$ satisfies Cond. (\ref{suff_cond_intro}.1), Cond. (\ref{suff_cond_intro}.2), and Cond. (\ref{suff_cond_intro}.3) of Theorem II with the same constants $C_{max}$, $C_{proj}$, and $C_{weight}$. 
\item The linear shape regularity $\sigma_e = h_e/\rho_e$ for each element $\Ophys$ of each mesh $\mathcal{M}_i$ in the family $M$ is bounded uniformly from above by a constant $\sigma_0$.
\end{enumerate}
Then, by Theorem II, there exists a \textbf{\textit{universal constant}} $C$ such that for each element $\Ophys$ of each mesh $\mathcal{M}_i$ in the family $M$ and for all $u \in H^{p+1}(\Omega_e)$, there is an approximation function $u_h \in \mathcal{S}^h_p(\Omega_e)$ if $\Omega_e$ is a simplicial element or $u_h \in \mathcal{S}^h_\bp(\Omega_e)$ if $\Omega_e$ is a tensor product element satisfying:
      \begin{equation*}
        \ltwo{u-u_h}{\Ophys} \leq Ch_e^{p+1} \| u \|_{H^{p+1}(\Ophys)}
      \end{equation*}
That is, the constant $C$ is the \textit{\textbf{same}} for \textit{\textbf{every element}} of \textit{\textbf{every mesh}} in the family.


%
\section{Distortion Metrics for Rational \BB Elements }  \label{validity_metrics}
Given the historical role of shape regularity in establishing suitable element metrics for linear finite element meshes, our new definition of shape regularity is an appropriate launching point for constructing element metrics for rational \BB meshes.  In particular, our definition inspires the following new distortion metrics for a rational \BB element:
\begin{enumerate}[label={\textbf{\textup{Distortion Metric \arabic*:}}},align=left,noitemsep]
\item The inverse scaled Jacobian metric: $$\linf{\det\bnabla_{\bxi} \xphys^{1/2}}{\Oref}  \linf{\det\bnabla_{\bxi} \xphys^{-1/2}}{\Oref}$$
\item The scaled derivative metric of order $\balpha$: $$h^{-|\balpha|}_e \linf{D^{\balpha}_{\bxi} \xproj}{ \Oref}$$
\item The inverse weighting metric: $$\linf{\frac{1}{w}}{\Oref}$$
\end{enumerate}
All of the above are truly distortion metrics in that they grow in size as an element is distorted.  Each metric additionally depends on the shape of a curvilinear element but not its size.  However, the above metrics are neither easy nor cheap to compute since they involve finding maxima of quantities such as the determinant Jacobian over the reference element.  One could approximate the three distortion metrics by replacing the $L^{\infty}$-norms appearing in their definition with maxima over a finite set of points, for example:
$$\linf{\det\bnabla_{\bxi} \xphys^{1/2}}{\Oref} \approx \max_{i=1,\ldots,n_{sample}} |\det\bnabla_{\bxi} \xphys (\bxi_i)|^{1/2}$$
but to ensure that such approximations are reasonable, one should also obtain rigorous upper bounds on the size of the distortion metrics.  Consequently, we next turn to the problem of establishing \textbf{\textit{computable upper bounds}} for our three new distortion metrics.  The inverse weighting metric is trivially bounded above by one over the smallest weight, that is:
\begin{equation}
\linf{\frac{1}{w}}{\Oref} \leq \frac{1}{\min_{j \in I} w_j}
\end{equation}
However, establishing computable bounds for the inverse scaled Jacobian metric and the scaled derivative metrics is not a trivial task.
It should be noted before proceeding that our new metrics as well as the upper bounds established in this section not only apply to rational \BB meshes but also standard isoparametric finite element meshes by taking the weighting function to be equal to one.

\subsection{Computable Bounds on the Jacobian Determinant}
We begin by establishing computable upper bounds for the inverse scaled Jacobian metric.
As previously mentioned, efficient algorithms have been proposed for bounding the Jacobian determinant of polynomial elements \cite{johnen_geometrical_2013}, but no such bounds have been proposed for rational elements.
We do note, however, that for every rational element $\Ophys \subset \Rphys$, there is a corresponding projective element $\Oproj\subset \Rproj$ that is defined a polynomial mapping.
Naturally then, we seek a way to compute the Jacobian determinant of the physical element in terms of the projective element.

For an element in $\Rphys$ the differential $d$-form $\omega$ is given by the $d^{th}$ external product of the directional derivatives of $\bnabla_{\bxi} \xphys$. That is:
\begin{equation}
  \omega = \deldel{\xphys}{\xi_1} \wedge...\wedge \deldel{\xphys}{\xi_d}
\end{equation}
where $\wedge$ denotes the wedge product.
For elements in $\mathbb{R}^2$ this yields a 2-form, which is a differential area  $dA$, and for elements in $\mathbb{R}^3$ this yields a 3-form, which is a differential volume element $dV$. 
The Jacobian determinant then, is simply the Hodge dual  of $\omega$, viz:
\begin{equation}
  \det \matof{\bnabla_{\bxi} \xphys} =   *\of{\omega}
\end{equation}
where $*(\cdot)$ is the Hodge star operator. For $\Rphys$, the Hodge star operator denotes the duality between $k$-forms and $(d-k)$-forms. 
For elements in $\Rphys$, $\omega$ is a $d$-form, and as such, $*(\omega)$ is a 0-form, which is a scalar.
Geometrically, the Jacobian determinant gives the area of 2-forms in $\mathbb{R}^2$, and the volume of 3-forms in $\mathbb{R}^3$.

For projective elements in $\Rproj$, we can  define the differential $d$-form $\widetilde{\omega}$ as:
\begin{equation}
  \widetilde{\omega} = \deldel{\xproj}{\xi_1} \wedge...\wedge \deldel{\xproj}{\xi_d}
\end{equation}
As with the physical element, this yields area elements $dA$ when $d=2$ and volume elements $dV$ when $d=3$. 
For projective elements however,  the $d$-form is defined in the $d+1$ dimensional vector space $\mathbb{R}^{d+1}$.
As such, the Hodge dual of $d$-forms in $\mathbb{R}^{d+1}$ are 1-forms, which are simply vectors.

\begin{figure}[t!]
  \centering
  \includegraphics[width=0.85\textwidth]{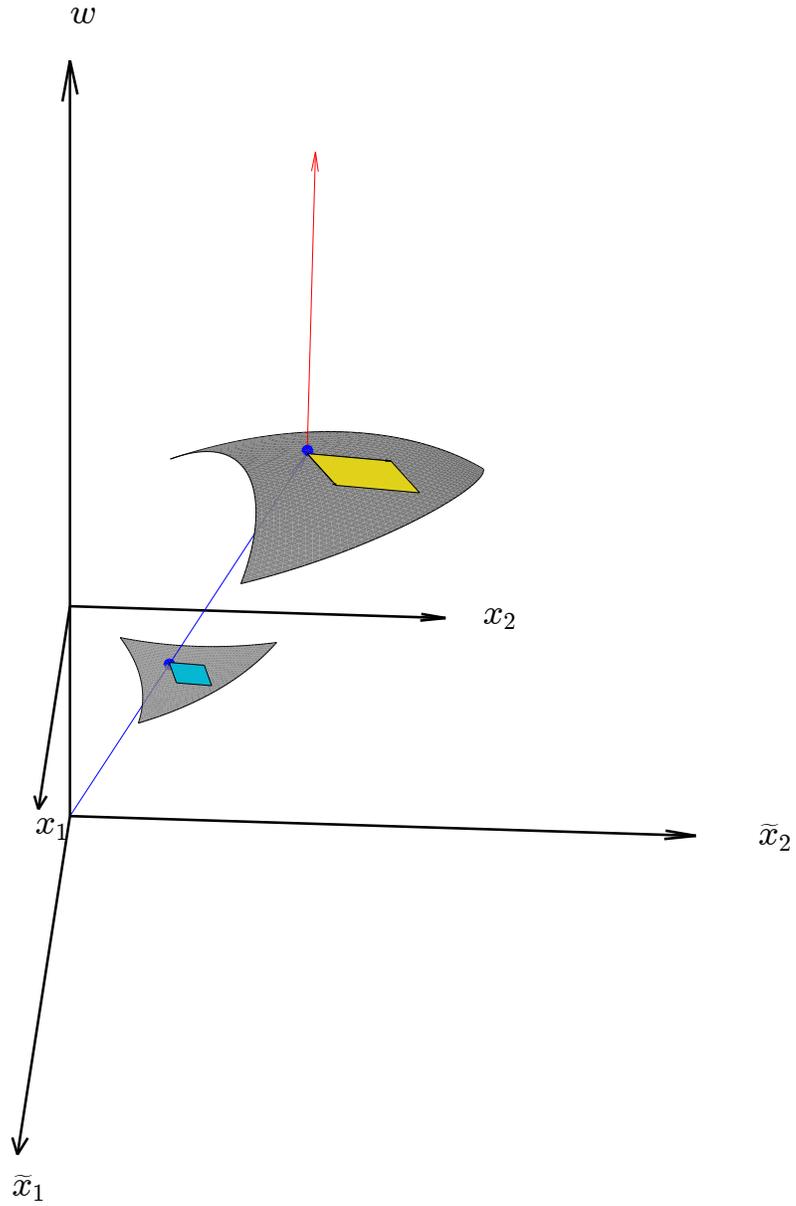}
  \caption{Differential 2-forms on a physical element  and a projective element.}
  \label{PROJ_DET}
\end{figure}

We illustrate these concepts for a rational \BB triangle, shown in Fig. \ref{PROJ_DET}. 
The 2-form for the physical element is visualized by the blue  parallelogram, and the Hodge dual of the 2-form gives the area of the parallelogram.
The 2-form for the projective element is shown by the yellow parallelogram, and the Hodge dual of the 2-form is the corresponding normal vector, and this relation is denoted:
\begin{equation}
  \bN  = *\of{  \deldel{\xproj}{\xi_1}  \wedge \deldel{\xproj}{\xi_2} } =   \deldel{\xproj}{\xi_1} \times \deldel{\xproj}{\xi_2}
\end{equation}
We take care to note that this vector $\bN$ is \textbf{\emph{not}} a unit normal. Rather, the magnitude of $\bN$ is equal to the area of the corresponding 2-form.

Visualizing the differential forms for volumetric elements becomes untenable, as the projective elements are embedded in $\mathbb{R}^4$.
However, for an arbitrary element in $\mathbb{R}^{d+1}$, we can write the vector $\bN$ as:
\begin{equation}\label{sum_minors_d}
  \bN \of{\bxi} = \sum\limits_{i=1}^{d+1} \of{-1}^{d+1+i} M_i \be_i
\end{equation}
wherein the minor $M_i$ is the  determinant of the $d \times d$ matrix formed by deleting the $i^{th}$ row of $\bnabla_{\bxi} \xproj $.
Additionally, let  $\bx = \xphys\of{\bxi}$ denote a coordinate on the physical element, and let $\widetilde{\bx} = \xproj\of{\bxi}$ denote the corresponding coordinate on the projective element.
With this nomenclature established, we can present bounds on $\det \matof{\bnabla_{\bxi} \xphys}$ in terms of the polynomial mapping $\xproj$.
\begin{thm} \label{jac_dot}
  For any rational \BB element, the Jacobian determinant can be calculated by:
  \begin{equation}
    \det \matof{\bnabla_{\bxi} \xphys} = \dfrac{{\bN\of{\bxi}} \cdot \xproj \of{\bxi}}{w^{d+1}} 
  \end{equation}
\end{thm}
\begin{proof}
  First, let us denote the first $d$ components of the mapping $\xproj$ as $\xprojd$. 
  We then recognize that we desire to calculate the determinant of the matrix:
  \begin{equation}
    \det \matof{\bnabla_{\bxi} \xphys} = 
    \det \matof{  \bnabla_{\bxi} \of{ \dfrac{\xprojd}{w} } }  
  \end{equation}
  By the quotient rule we have:
  \begin{equation}
    \det \matof{\bnabla_{\bxi} \xphys} =
    \det \Bigg{[} \dfrac{1}{w} \bigg{[}\bnabla_{\bxi} \xprojd - \bx [\bnabla_{\bxi}w]^T  \bigg{]} \Bigg{]}
  \end{equation}
  and, because for any given $\bxi \in \Oref$, $w$ is some positive constant, we can factor out the weighting function to  write:
  \begin{equation} \label{1_wd_det}
    \det \matof{\bnabla_{\bxi} \xphys} =
    \dfrac{1}{w^d}  \det \Bigg{[} \bnabla_{\bxi} \xprojd - \bx [\bnabla_{\bxi}w]^T  \Bigg{]}
  \end{equation}
  Then, by Eq. \ref{sum_minors_d} we can write the vector $\bN$ as:
  \begin{equation}\label{sum_minors_2}
    \bN \of{\bxi} = \sum\limits_{i=1}^{d+1} \of{-1}^{d+1+i} M_i \be_i
  \end{equation}
  and as a result,  we can write the dot product $\bN\cdot \bx$ as:
  %
  %
  \begin{equation}
    \bN \cdot {\bx}=
    \sum\limits_{i=1}^{d+1} \of{-1}^{d+1+i} M_i x_i =
    \det \left[ \begin{array}{cc}
        \bnabla_{\bxi} \xprojd & \bx  \\
        & \\
        \bnabla_{\bxi} w & 1
      \end{array} \right] 
  \end{equation}
  which in turn can be written as:
  \begin{equation} \label{block_mat_det}
    \bN \cdot {\bx}=
    \det \left[ \begin{array}{cc}
        \bnabla_{\bxi} \xprojd & \bx  \\
        & \\
        \bnabla_{\bxi} w & 1
      \end{array} \right]
    = 
    \det \Bigg{[} \bnabla_{\bxi} \xprojd - \bx [\bnabla_{\bxi}w]^T  \Bigg{]}
  \end{equation}
  Then, recognizing that $  \bN \cdot {\bx} = \dfrac{\bN \cdot \widetilde{\bx}}{ w}$ , and substituting the results of Eq. \eqref{block_mat_det} into Eq. \eqref{1_wd_det}, we get
  \begin{equation}
    \det \matof{\bnabla_{\bxi} \xphys} = \dfrac{1}{w^{d}} \dfrac{\bN \cdot \widetilde{\bx}}{ w}
  \end{equation}
  from which the results of Theorem \ref{jac_dot} follow immediately.
\end{proof}

Conceptually, Theorem \ref{jac_dot} can be thought of as projecting the vector $\bN$ onto the $w=1$ plane along the vector $\widetilde{\bx}$, scaled by $w^{d+1}$.
Alternatively, $\det \matof{\bnabla_{\bxi}\xphys}$ this can be thought of as the apparent magnitude of $\widetilde{\omega}$ as seen by an oberver at the origin.
Either way, we recognize that we can compute the Jacobian determinant as the dot product of two vectors, normalized by the weighting function.
Note that for a rational \BB element of degree $p$ in $\Rphys$, both $\widetilde{\bx}$ and $\bN$ will be vectors in $\Rproj$. 
It is readily seen that the vector $\widetilde{\bx}$ can be written in \BB form, as the control points $\widetilde{\bP}_\bi$ are known.
However, we note that the equation for the surface normal $ \bN = \bN(\bxi)$ can also be written in  \BB form. 
It remains to present a method for calculating the B\'{e}zier coefficients for the surface normal. 
We present formulas for these coefficients for simplicial \BB elements in Theorem \ref{thm:Nk_simplex} and for tensor product elements in Theorem \ref{Nk_TP}.

We begin by considering simplicial elements.
Let $\{ B_\bk^{p'}\}_{\bk \in I^{p'}}$ denote the set of simplicial Bernstein polynomials of degree $p' = d\of{p-1}$, and let $\{\bN_{\bk}\}_{\bk \in I^{p'}}$ denote the set of B\'{e}zier coefficients for the vector $\bN$.
Then, for a simplicial element in projective space  we can write  $\bN \of{ \bxi }$ in \BB form as:
\begin{equation}
  \bN \of{ \bxi } = \sum\limits_{\bk \in I^{p'}} B_\bk^{p'} \bN_\bk \ \  \forall\ \bxi \in \Oref 
\end{equation}
wherein $I^{p'}$ denotes the index set over the simplicial Bernstein polynomials of degree $p'$. 
Then, let us denote a $d$-tuple of multi-indices as $\bvi = \{\bi_1,...,\bi_d\}$, and let us we define the set $\mathcal{I}_\bk^{p'}$ as:
\begin{equation}
  \mathcal{I}^{p'}_\bk := 
  \setof{ \bvi = \{ \bi_j\}_{j=1}^d \ : \ 
    \bi_j \in I^{p-1}, \ 
    \sum\limits_{j=1}^d \bi_j = \bk }
\end{equation}
Finally, let us denote the set of difference vectors in the $\xi_j$ direction as $\{ \Delta \widetilde{\bP}_{\bi_j} \}_{\bi_j \in I^{\bp-\be_j}}$, where we define $\Delta \widetilde{\bP}_{\bi_j}$ as:
\begin{equation}
\Delta \widetilde{\bP}_{\bi_j} = \widetilde{\bP}_{\bi+\be_j} - \widetilde{\bP}_{\bi}
\end{equation}
%
%
%
%
%
%
\begin{thm} \label{thm:Nk_simplex}
  For a simplicial \BB element of degree $p$ in projective space, with projective control points $\{\widetilde{\bP}_{\bi}\}_{\bi \in I^p}$,
  the B\'{e}zier coefficients $\{ \bN_\bk \}_{\bk \in I'}$ for the vector $\bN$ can be calculated as:
  \begin{equation} \label{Nk_simplex}
    \bN_\bk = *\of{
      \sum\limits_{\bvi \in \mathcal{I}_\bk^{p'}}
      \eta_{\bbk}\of{\bvi} 
      \of{ \Delta \widetilde{\bP}_{\bi_{1}} \wedge ... \wedge   \Delta \widetilde{\bP}_{\bi_{d}} }  
    }
  \end{equation}
  wherein the coefficient $\eta_{\bbk}\of{\bvi}$ is defined to be:
  \begin{equation}
    \eta_{\bk}\of{\bvi} =
    \dfrac{ p^d {{p-1}\choose{\bbi_1}}...{{p-1}\choose{\bbi_d}}}{{{p'}\choose{\bbk}}}  
  \end{equation}
  %
  %
  %
\end{thm}

\begin{proof}
  We first recognize that a projective B\'{e}zier element is a $d$-manifold with codimension 1, and as such the surface normal can be found as the Hodge dual of the wedge product of the parametric derivatives, viz:
  \begin{equation} \label{N_wedge}
    \bN = *\of{ \deldel{\xproj}{\xi_1} \wedge ... \wedge \deldel{\xproj}{\xi_d} }
  \end{equation}
  Note that when $d=2$, this is simply the cross product, but the above notation holds for arbitrary $d$.
  Now, recognizing that the partial derivative with respect to the $j^{th}$ parametric coordinate can be found as:
  \begin{equation}
    \deldel{\xproj}{\xi_j} = p\sum\limits_{\bi_j \in I^{p-1}} B_{\bi_j}^{p-1}  \Delta \widetilde{\bP}_{\bi_{j}}
  \end{equation}
  we rewrite Eq. \eqref{N_wedge} as:
  \begin{equation} 
    \bN = *\of{p\sum\limits_{\bi_1 \in I^{p-1}} B_{\bi_1}^{p-1}  \Delta \widetilde{\bP}_{\bi_{1}} \wedge ... \wedge p\sum\limits_{\bi_d \in I^{p-1}} B_{\bi_d}^{p-1}   \Delta \widetilde{\bP}_{\bi_{d}} }
  \end{equation}
  and since the distributive property holds, we can rearrange to yield:
  \begin{equation}
    \bN = *\of{ p^d \sum\limits_{\bi_1 \in I^{p-1}} ... \sum\limits_{\bi_d \in I^{p-1}} B_{\bi_1}^{p-1}...B_{\bi_d}^{p-1} \Delta \widetilde{\bP}_{\bi_{1}} \wedge ... \wedge \Delta \widetilde{\bP}_{\bi_{d}} }
  \end{equation}
  Now, recognizing that we can write the product of the Bernstein basis functions as:
  \begin{equation}
    B_{\bi_1}^{p-1}...B_{\bi_d}^{p-1} = \dfrac{ {{p-1}\choose{\bbi_1}}...{{p-1}\choose{\bbi_d}}}{{{p'}\choose{\bbk}}} B_{\bk}^{p'}
  \end{equation}
  we arrive at:
  \begin{equation}
    \bN = *\of{ p^d \sum\limits_{\bi_1 \in I^{p-1}} ... \sum\limits_{\bi_d \in I^{p-1}}  \dfrac{ {{p-1}\choose{\bbi_1}}...{{p-1}\choose{\bbi_d}}}{{{p'}\choose{\bbk}}} B_{\bk}^{p'} \Delta \widetilde{\bP}_{\bi_{1}} \wedge ... \wedge \Delta \widetilde{\bP}_{\bi_{d}} }
  \end{equation}
  Finally, rearranging the order of summation, we get:
  \begin{equation}
    \bN = *\of{  \sum\limits_{\bk \in I^{p'}} B_{\bk}^{p'} 
      \sum\limits_{\bvi \in \mathcal{I}^{p'}_{\bk}}  
      \dfrac{ p^d {{p-1}\choose{\bbi_1}}...{{p-1}\choose{\bbi_d}}}{{{p'}\choose{\bbk}}}  
      \Delta \widetilde{\bP}_{\bi_{1}} \wedge ... \wedge \Delta \widetilde{\bP}_{\bi_{d}} }
  \end{equation}
  From which Eq. \eqref{Nk_simplex} immediately follows.
\end{proof}
%
%
%
%
%
%
We can now derive similar results for tensor product elements. Let $\{ B_{\bk}^{\bp'}\}_{\bk \in I^{\bp '} }$ denote the set of tensor product Bernstein polynomials of degree $\bp' = d \bp - \boldsymbol{1} $, and let $\{ \bN_{\bk} \}_{\bk \in I^{\bp '}}$ denote the set of B\'{e}zier coefficients for the vector $\bN$  .
Then, for a tensor product element in projective space  we can write  $\bN \of{ \bxi }$ in \BB form as:
\begin{equation}
  \bN(\bxi) = \sum\limits_{\bk \in I^{\bp'}} B_\bk^{\bp'} \bN_{\bk} \ \  \forall\ \bxi \in \Oref 
\end{equation}
Now, let us denote a $d$-tuple of multi-indices as $\bvi = \{\bi_1,...,\bi_d\}$. Then, we define the set of $d$-tuples $\mathcal{I}_\bk^{\bp'}$ as:
\begin{equation}
  \mathcal{I}^{\bp'}_\bk := 
  \setof{ \bvi = \{ \bi_j\}_{j=1}^d \ : \ 
    \bi_j \in I^{\bp-\textbf{e}_j}, \ 
    \sum\limits_{j=1}^d \bi_j = \bk 
  }
\end{equation}
wherein $I^{\bp'}$ denotes the index set over the tensor product Bernstein polynomials of degree $\bp'$. 
As before, let $\{ \Delta \widetilde{\bP}_{\bi_j} \}_{\bi_j \in I^{\bp-\be_j}}$ denote the set of difference vectors in the $\xi_j$ direction.

\begin{thm} \label{Nk_TP}
  For a tensor product \BB element of degree $\bp$ in projective space, with projective control points $\{\widetilde{\bP}_{\bi}\}_{\bi \in I^{\bp}}$,
  the B\'{e}zier coefficients $\{ \bN_\bk \}_{\bk \in I'}$ for the vector $\bN$ can be calculated as:
  \begin{equation} \label{Nk_quad}
    \bN_\bk = *\of{
      \sum\limits_{\bvi \in \mathcal{I}_\bk^{\bp'}}
      \eta_{\bk}\of{\bvi} 
      \of{ \Delta \widetilde{\bP}_{\bi_{1}} \wedge ... \wedge   \Delta \widetilde{\bP}_{\bi_{d}} }  
    }
  \end{equation}
  wherein the coefficient $\eta_{\bk}\of{\bvi}$ is defined to be:
  \begin{equation}
    \eta_{\bk}\of{\bvi} =
    \dfrac{ p_1{{\bp - \textbf{e}_1}\choose{\bi_1}} ... p_d{{\bp - \textbf{e}_d}\choose{\bi_d}} } { {{\bp'}\choose{\bk}} }
  \end{equation}
\end{thm}
\begin{proof}
  The proof follows the proof for simplicial elements almost exactly. 
  We simply recognize that $j^{th}$ directional derivative for a tensor product element is:
  \begin{equation}
    \deldel{\xproj}{\xi_j} = p_j\sum\limits_{\bi_j \in I^{\bp-\be_j}} B_{\bi_j}^{\bp-\be_j}  \Delta \widetilde{\bP}_{\bi_{j}}
  \end{equation}
  and that the product of tensor product Bernstein polynomials  can be written as:
  \begin{equation}
    B_{\bi_1}^{\bp-\be_1}...B_{\bi_d}^{\bp-\be_d} = \dfrac{ {{\bp-\be_1}\choose{\bi_1}}...{{\bp-\be_d}\choose{\bi_d}}}{{{\bp'}\choose{\bk}}} B_{\bk}^{\bp'}
  \end{equation}
  Then, using these identities along with Eq. \eqref{N_wedge}, Eq. \eqref{Nk_quad} can be readily obtained.
\end{proof}
%
%
%
%
%
%
%
\noindent
\begin{minipage}{\textwidth}
  \begin{mdframed}
    \begin{customthm}{III}  \label{J_bound_NP_intro}
      Let $\Ophys \subset \Rphys$ be  a  rational \BB element in physical space with corresponding rational element in projective space $\Oproj \subset \Rproj$.  
      Then, letting $\{\widetilde{\bP}_\bi \}_{\bi \in I} $ denote the projective control points, 
      and letting $\{\bN_{\bk}\}_{\bk \in I'}$ denote the B\'{e}zier coefficients for the normal vector $\bN$, 
      the inverse scaled Jacobian metric is bounded from above by:
      \begin{equation}
        \linf{\det\bnabla_{\bxi} \xphys^{1/2}}{\Oref}  \linf{\det\bnabla_{\bxi} \xphys^{-1/2}}{\Oref} \leq
        \of{
          \dfrac{\max\limits_{j\in I}w_j}
                {\min\limits_{j\in I}w_j}
        }^{d+1}
        \of{
          \dfrac{\max\limits_{\substack{\bi \in I \\ \bk \in I'}}\bN_\bk \cdot \widetilde{\bP}_\bi}
                {\min\limits_{\substack{\bi \in I \\ \bk \in I'}}\bN_\bk \cdot \widetilde{\bP}_\bi}
        }
      \end{equation}
    \end{customthm}
  \end{mdframed}
\end{minipage}

\noindent
  \begin{proof}
    From Theorem \ref{jac_dot}, we recognize that $\det \matof{\bnabla_{\bxi}\xphys}$ is given by:
    \begin{equation}
      \det\matof{\bnabla_{\bxi}\xphys} = \dfrac{\bN \cdot \widetilde{\bx}}{w^{d+1}} 
    \end{equation}
    From this, we can rewrite $\det \matof{\bnabla_{\bxi}\xphys}$ explicitly in terms of the Bernstein basis polynomials as:
    \begin{equation}
      \det \matof{\bnabla_{\bxi}\xphys} =  \dfrac{1}{ \of{ \sum\limits_{\bj \in I} B_{\bj}w_{\bj} }^{d+1}} 
      \sum\limits_{\bi \in I} \sum\limits_{\bk \in I'} B_\bk B_\bi \bN_\bk \cdot \widetilde{\bP}_\bi
    \end{equation}
    Then, because the Bernstein basis polynomials satisfy positivity and partition of unity, we can bound the magnitude of the Jacobian determinant by:
    \begin{equation} \label{J_bound_NP}
      \dfrac
          {\min\limits_{\substack{\bi \in I \\ \bk \in I'}}  
            \bN_\bk \cdot \widetilde{\bP}_\bi}
          { \of{ \max\limits_{\bi \in I} w_{\bi} }^{d+1}} 
          \leq \absof{ \det \matof{\bnabla\xphys} } \leq
          \dfrac
              {\max\limits_{\substack{\bi \in I \\ \bk \in I'}}  
                \bN_\bk \cdot \widetilde{\bP}_\bi}
              { \of{ \min\limits_{\bi \in I} w_{\bi} }^{d+1}} 
    \end{equation}
    Finally, recognizing that the equation for the inverse scaled Jacobian metric can be equivalently written:
    \begin{equation}
      \linf{\det\bnabla_{\bxi} \xphys^{1/2}}{\Oref}  \linf{\det\bnabla_{\bxi} \xphys^{-1/2}}{\Oref} =
      \of{
        \dfrac{\sup\limits_{\bxi \in \Oref} \absof{ \det \matof{\bnabla_{\bxi}\xphys} } }
              {\inf\limits_{\bxi \in \Oref} \absof{ \det \matof{\bnabla_{\bxi}\xphys} } }
      }^{1/2}
    \end{equation}
    we use the results of Eq. \eqref{J_bound_NP} to arrive at the results of Theorem III.
  \end{proof}

\subsection{Computable Bounds on Derivatives of the Mapping $\xphys$}
With a method for calculating bounds on the Jacobian determinant established, we turn our attention to computing bounds for higher-order derivatives so that we can bound the scaled derivative metrics. 
Compared to bounds on the Jacobian determinant, bounds on the higher-order derivatives are relatively easy to derive. 
These bounds are presented below in the proof for Theorems IVa and b. 
With these bounds established, we have succeeded in establishing a set of computable bounds on the three element distortion metrics proposed at the beginning of this section for rational \BB elements. 

\noindent
\begin{minipage}{\textwidth}
  \begin{mdframed}
    \begin{customthm}{IVa} \label{simplex_deriv_bound_intro}
      Let us denote the projective control points of a simplicial \BB element of degree $p$ as $\{\widetilde{\bP}_{\bi}\}_{\bi \in I^p}$.
      Then, the ${\balpha}^{th}$ partial derivative of the mapping $\xproj$ is bounded by:
      \begin{equation}
        \linf{D_{\bxi}^{\balpha} \xproj }{ \Oref } \leq
        \dfrac{p!}{(p - \absof{\balpha})!}
        \max\limits_{\bi \in I^{p-\absof{\balpha}}}
        \absof{
          \sum_{\bj \in I^{\balpha}}
          (-1)^{{\balpha}+\bj}{{{\balpha}}\choose{\bj}}
          \widetilde{\bP}_{\bi + \bj}
        }
      \end{equation}
      Moreover, if the ${\balpha}^{th}$ partial derivative of the mapping $\xproj$ is zero, then:
      \begin{equation}        
        \max\limits_{\bi \in I^{p-\absof{\balpha}}}
        \absof{
          \sum_{\bj \in I^{\balpha}}
          (-1)^{{\balpha}+\bj}{{{\balpha}}\choose{\bj}}
          \widetilde{\bP}_{\bi + \bj}
        } = 0
        \end{equation}
    \end{customthm}
  \end{mdframed}
\end{minipage}

\begin{proof} \label{simplex_deriv_bound}
  Consider a simplicial \BB element in projective space, $\Oproj \subset \Rproj$ defined by control points $\{\widetilde{\bP}_{\bi}\}_{\bi \in I^p}$. 
  We recognize that the derivatives of Bernstein polynomials are themselves Bernstein polynomials of a lower degree \cite{prautzsch2013bezier}. 
  Thus, we can recursively take the derivative of the mapping $\xproj$, which yields the following equation  for the ${\balpha}^{th}$ partial derivative:
  \begin{equation}
    D_{\bxi}^{\balpha} \xproj = 
    \dfrac{p!}{(p - |{\balpha}|)!}\sum\limits_{\bi \in I^{p-|{\balpha}|}}
    \matof{ 
      B_{\bi}^{p-|{\balpha}|}(\bxi) \sum_{\bj \in I^{{\balpha}}}
      (-1)^{{\balpha}+\bj}{{{\balpha}}\choose{\bj}}
      \widetilde{\bP}_{\bi + \bj }
    }
  \end{equation}
  Note, we take care to emphasize that the above sum over  $\bj \in I^{\balpha}$ is a sum over a tensor product index set. 
  This is  a consequence of the fact that the partial derivative $\nabla^{\balpha}$ has an inherently tensor product nature.
  Then, because the Bernstein polynomials satisfy positivity and partition of unity, the desired bound is obtained.
  Now suppose that the ${\balpha}^{th}$ partial derivative of $\xproj$ is zero.  As the basis functions $B_{\bi}^{p-|{\balpha}|}(\bxi)$ are linearly independent, all the coefficients in the above expansion for $D_{\bxi}^{\balpha} \xproj$ must be zero, and hence so is the coefficient of maximum magnitude.
  \end{proof}

\noindent
\begin{minipage}{\textwidth}
  \begin{mdframed}
    \begin{customthm}{IVb} \label{TP_deriv_bound_intro}
      Let us denote the projective control points of a tensor product \BB element of degree $\bp$ as $\{\widetilde{\bP}_{\bi}\}_{\bi \in I^\bp}$.
      Then, the ${\balpha}^{th}$ partial derivative of the mapping $\xproj$ is bounded by:
      \begin{equation}
        \linf{D_{\bxi}^{\balpha} \xproj}{ \Oref } \leq
        \dfrac{\bp!}{(\bp - {\balpha})!}
        \max\limits_{\bi \in I^{\bp-{\balpha}}}
        \absof{
          \sum_{\bj \in I^{\balpha}}
          (-1)^{{\balpha}+\bj}{{{\balpha}}\choose{\bj}}
          \widetilde{\bP}_{\bi + \bj}
        }
      \end{equation}
      Moreover, if the ${\balpha}^{th}$ partial derivative of the mapping $\xproj$ is zero, then:
      \begin{equation}        
        \max\limits_{\bi \in I^{\bp-\balpha}}
        \absof{
          \sum_{\bj \in I^{\balpha}}
          (-1)^{{\balpha}+\bj}{{{\balpha}}\choose{\bj}}
          \widetilde{\bP}_{\bi + \bj}
        } = 0
        \end{equation}
    \end{customthm}
  \end{mdframed}
\end{minipage}

\begin{proof}
  Consider a tensor product \BB element in projective space, $\Oproj \subset \Rproj$ defined by control points $\{\widetilde{\bP}_{\bi}\}_{\bi \in I^\bp}$. 
  As before, we write the derivatives of the Bernstein polynomials as Bernstein polynomials of lower degree.
  This yields the following equation  for the ${\balpha}^{th}$ partial derivative of the mapping $\xproj$:
  \begin{equation}
    D_{\bxi}^{\balpha} \xproj = 
    \dfrac{\bp!}{(\bp - {\balpha})!}\sum\limits_{\bi \in I^{\bp-{\balpha}}}
    \matof{
      B_{\bi}^{\bp-{\balpha}}(\bxi) \sum_{\bj \in I^{{\balpha}}}
      (-1)^{{\balpha}+\bj}{{{\balpha}}\choose{\bj}}
      \widetilde{\bP}_{\bi + \bj }
    }
  \end{equation}
  Then, because the Bernstein polynomials satisfy positivity and partition of unity, the desired bound is obtained.  Now suppose that the ${\balpha}^{th}$ partial derivative of $\xproj$ is zero.  As the basis functions $B_{\bi}^{\bp-\balpha}(\bxi)$ are linearly independent, all the coefficients in the above expansion for $D_{\bxi}^{\balpha} \xproj$ must be zero, and hence so is the coefficient of maximum magnitude.
\end{proof}
%
%
%
%

With the relevant theory established, we now demonstrate a particularly convenient property of the bounds presented in Theorem IVa and Theorem IVb. 
First, let us consider the case of finding the higher-order derivatives of a cubic \BB triangle. 
Table \ref{bound_expressions} shows the bounding expressions for several derivatives of the mapping $\xphys$.
We see  that the bounds on the first derivative can be found by evaluating the expression:
\begin{equation}
  \absof{ 3 \widetilde{\bP}_{\bi + \setof{1,0}} - 3 \widetilde{\bP}_{\bi}} = 3 \absof{ \widetilde{\bP}_{\bi + \setof{1,0}} - \widetilde{\bP}_{\bi}} 
\end{equation}
at the points
$ \widetilde{\bP}_{\bi} \in \left\{\widetilde{\bP}_{\setof{0,0}},\widetilde{\bP}_{\setof{1,0}},\widetilde{\bP}_{\setof{2,0}},\widetilde{\bP}_{\setof{0,1}},\widetilde{\bP}_{\setof{1,1}},\widetilde{\bP}_{\setof{0,2}}\right\} $.
Similarly, the bound on the second derivative is found by evaluating:
\begin{equation}
  \absof{
    6\widetilde{\bP}_{\bi +  \setof{2,0} } - 12\widetilde{\bP}_{\bi +  \setof{1,0}} + 6\widetilde{\bP}_{\bi + \setof{0,0}}
  } =
  6 \absof{
    \left( \widetilde{\bP}_{\bi +  \setof{2,0} } - \widetilde{\bP}_{\bi +  \setof{1,0}} \right) - \left( \widetilde{\bP}_{\bi +  \setof{1,0}} - \widetilde{\bP}_{\bi + \setof{0,0}} \right)
  }
\end{equation}
at the points $ \widetilde{\bP}_{\bi} \in \left\{\widetilde{\bP}_{\setof{1,0}}, \widetilde{\bP}_{\setof{1,0}},\widetilde{\bP}_{\setof{0,1}} \right\}$.
Then, we note that for a simplicial element of degree $p$, all partial derivatives of order $|{\balpha}| = p$ will be constant across the element. 
As such, the third derivative for a cubic \BB triangle can be calculated analytically by the expression:
\begin{eqnarray}
  &\absof{
    6 \widetilde{\bP}_{\setof{3,0}} - 18 \widetilde{\bP}_{\setof{2,0}} + 
    18 \widetilde{\bP}_{\setof{1,0}} - 6 \widetilde{\bP}_{\setof{0,0}}
  } = \\  
  & 6 \absof{
    \left( \left( \widetilde{\bP}_{\bi +  \setof{3,0} } - \widetilde{\bP}_{\bi +  \setof{2,0}} \right) - \left( \widetilde{\bP}_{\bi +  \setof{2,0}} - \widetilde{\bP}_{\bi + \setof{1,0}} \right) \right)
    - \left( \left( \widetilde{\bP}_{\bi +  \setof{2,0} } - \widetilde{\bP}_{\bi +  \setof{1,0}} \right) - \left( \widetilde{\bP}_{\bi +  \setof{1,0}} - \widetilde{\bP}_{\bi + \setof{0,0}} \right) \right)
  }
\end{eqnarray}

From the above, it is apparent that bounds on the derivatives of \BB elements can be calculated using a finite difference stencil involving the projective control points.  
To illustrate this notion, Table \ref{deriv_examples} shows finite difference stencils for each of the  derivatives shown in Table \ref{bound_expressions}.
We then show this stencil applied  to the control points of the element, as well as the resulting B\'{e}zier coefficients for the derivative, shown as vectors on the reference element. 

For clarity, we have shown a non-rational cubic B\'{e}zier triangle in $\mathbb{R}^2$, but the concepts extends readily to elements in projective space. 
We have also included the explicitly calculated stencils for a variety of elements in Appendix \ref{deriv_appendix}. 
This is not an exhaustive list, but we note that the results of Section \ref{validity_metrics} can be used to calculate stencils for \textbf{\emph{any}} simplicial or tensor product \BB element.

\begin{table}[t] 
  \begin{center}
    \caption{Bounds on the derivatives of a cubic \BB triangle. }
    \label{bound_expressions}
    \begin{tabular}{ | >{\centering\arraybackslash} m{3cm} |  >{\centering\arraybackslash} m{9cm} |  }
      \hline
      Derivative & Bound  \\ 
      \hline
      $\deldel{\xproj}{\xi_1} $& \begin{equation} 
        \leq
        \max\limits_{\bi \in I^{p-1}}
        \absof{
          3 \widetilde{\bP}_{\bi + \setof{1,0}} - 3 \widetilde{\bP}_{\bi} 
        }
      \end{equation}\label{first_deriv_ex}
      \\
      \hline
      $\dfrac{\partial^2 \xproj}{\partial \xi_1^2} $ & 
      \begin{equation}  
        \leq 
        \max\limits_{\bi \in I^{p-2}}
        \absof{
          6\widetilde{\bP}_{\bi +  \setof{2,0} } - 12\widetilde{\bP}_{\bi +  \setof{1,0}} + 6\widetilde{\bP}_{\bi + \setof{0,0}}
        }
      \end{equation} \label{second_deriv_ex}
      \\
      \hline
      $\dfrac{\partial^3 \xproj}{\partial \xi_1^3} $& 
      \begin{equation}     \label{third_deriv_ex}
        =
        \absof{
          6 \widetilde{\bP}_{\setof{3,0}} - 18 \widetilde{\bP}_{\setof{2,0}} + 
          18 \widetilde{\bP}_{\setof{1,0}} - 6 \widetilde{\bP}_{\setof{0,0}}
        }
      \end{equation}  \\
      \hline
    \end{tabular}
    \label{deriv_examples}
  \end{center}
\end{table}

\begin{table}[t]
  \begin{center}
    \caption{Example of finding the B\'{e}zier coefficients for several different derivatives using a stencil. }
    \begin{tabular}{ | >{\centering\arraybackslash} m{2cm}  | >{\centering\arraybackslash} m{3cm} |  >{\centering\arraybackslash} m{5cm} |  >{\centering\arraybackslash} m{3cm} |  }
      \hline
      Derivative & Stencil & Stencils Applied to the Physical Triangle & B\'{e}zier Coefficients of the Derivative  \\ 
      \hline
      $\deldel{\xproj}{\xi_1} $&
      \includegraphics[width=0.18\textwidth]{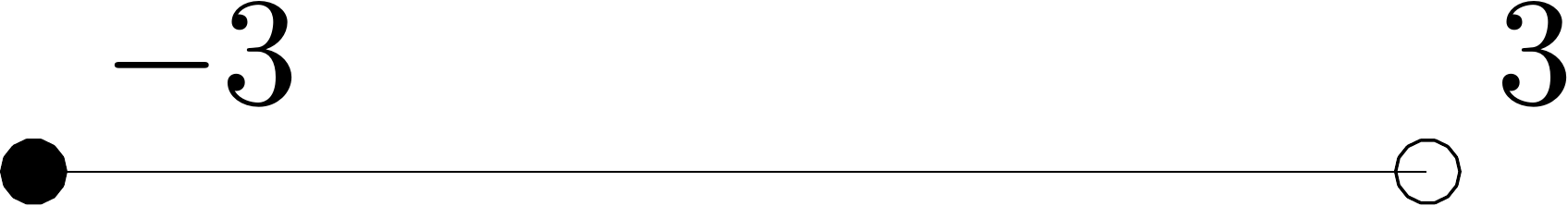} &
      \vspace{5pt} \includegraphics[width=4.5cm,height=4.5cm,keepaspectratio]{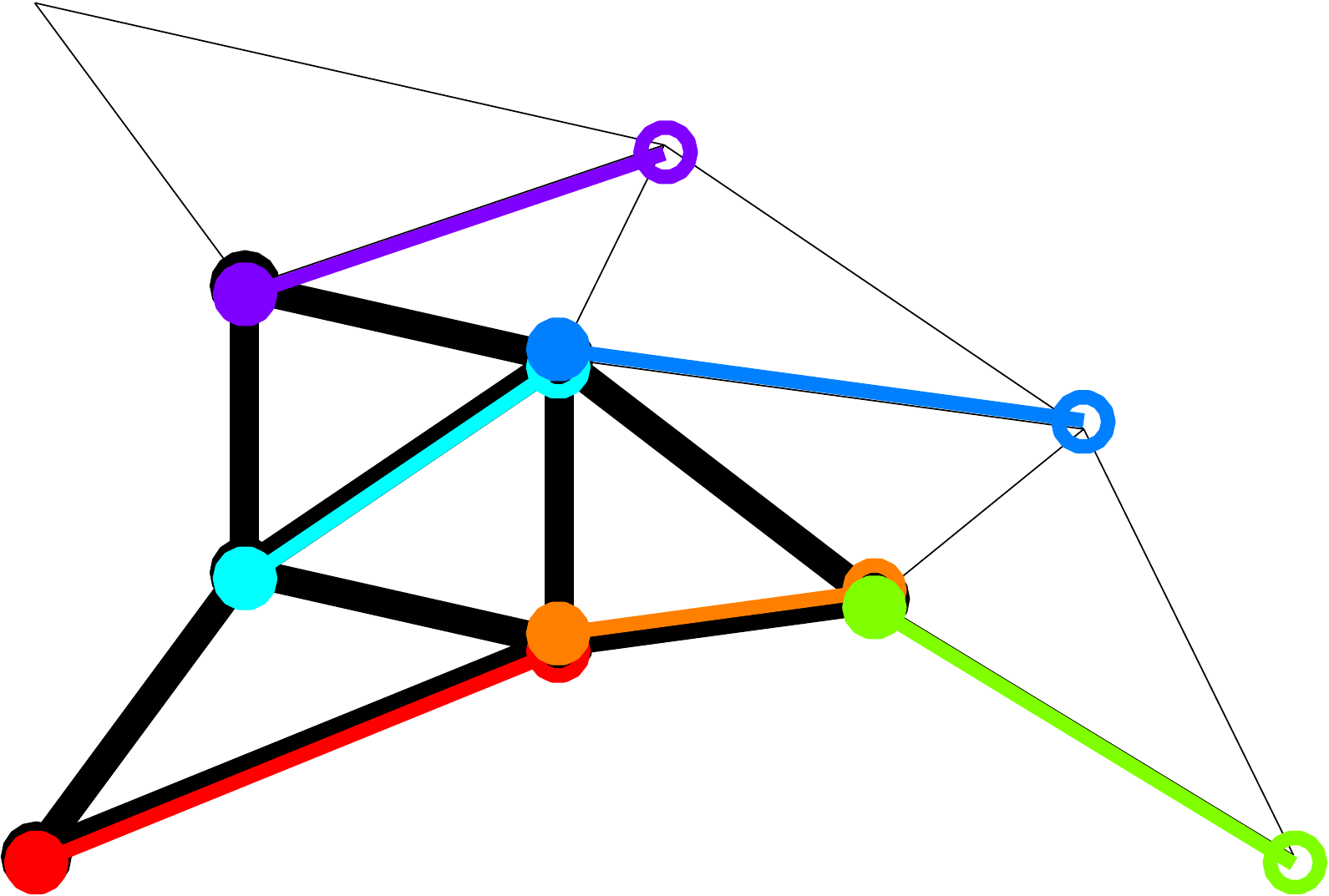} &
      \vspace{5pt} \includegraphics[width=2.5cm,height=2.5cm,keepaspectratio]{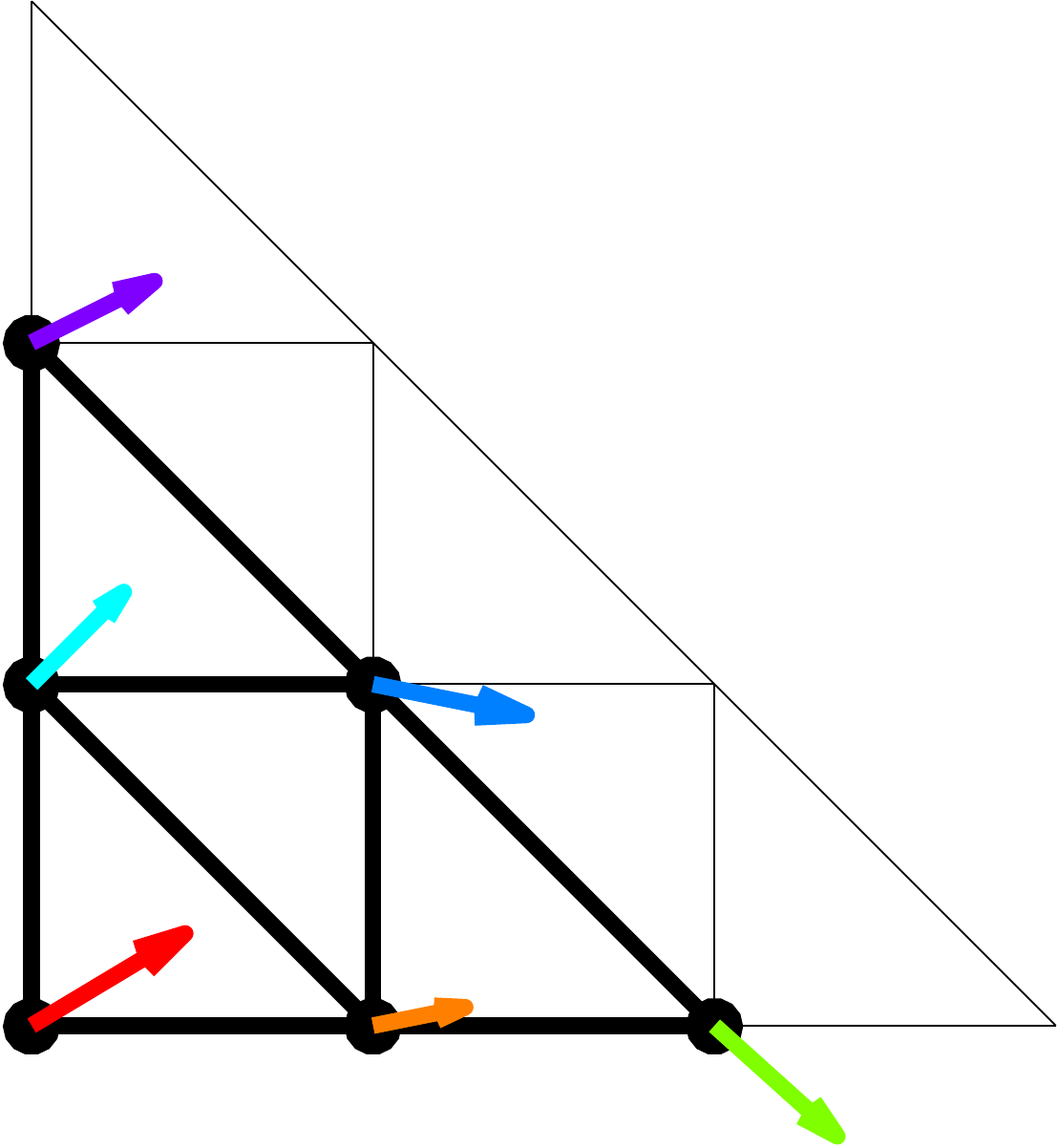}  \hspace{5pt} \\
      \hline
      $\dfrac{\partial^2 \xproj}{\partial \xi_1^2} $&
      \includegraphics[width=0.18\textwidth]{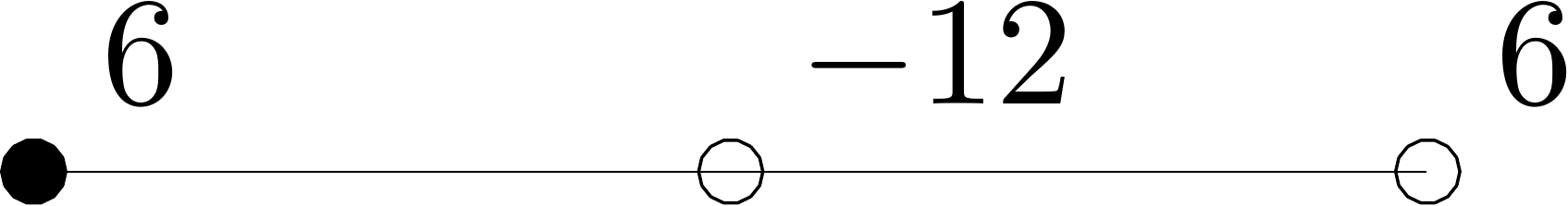} &
      \vspace{5pt} \includegraphics[width=4.5cm,height=4.5cm,keepaspectratio]{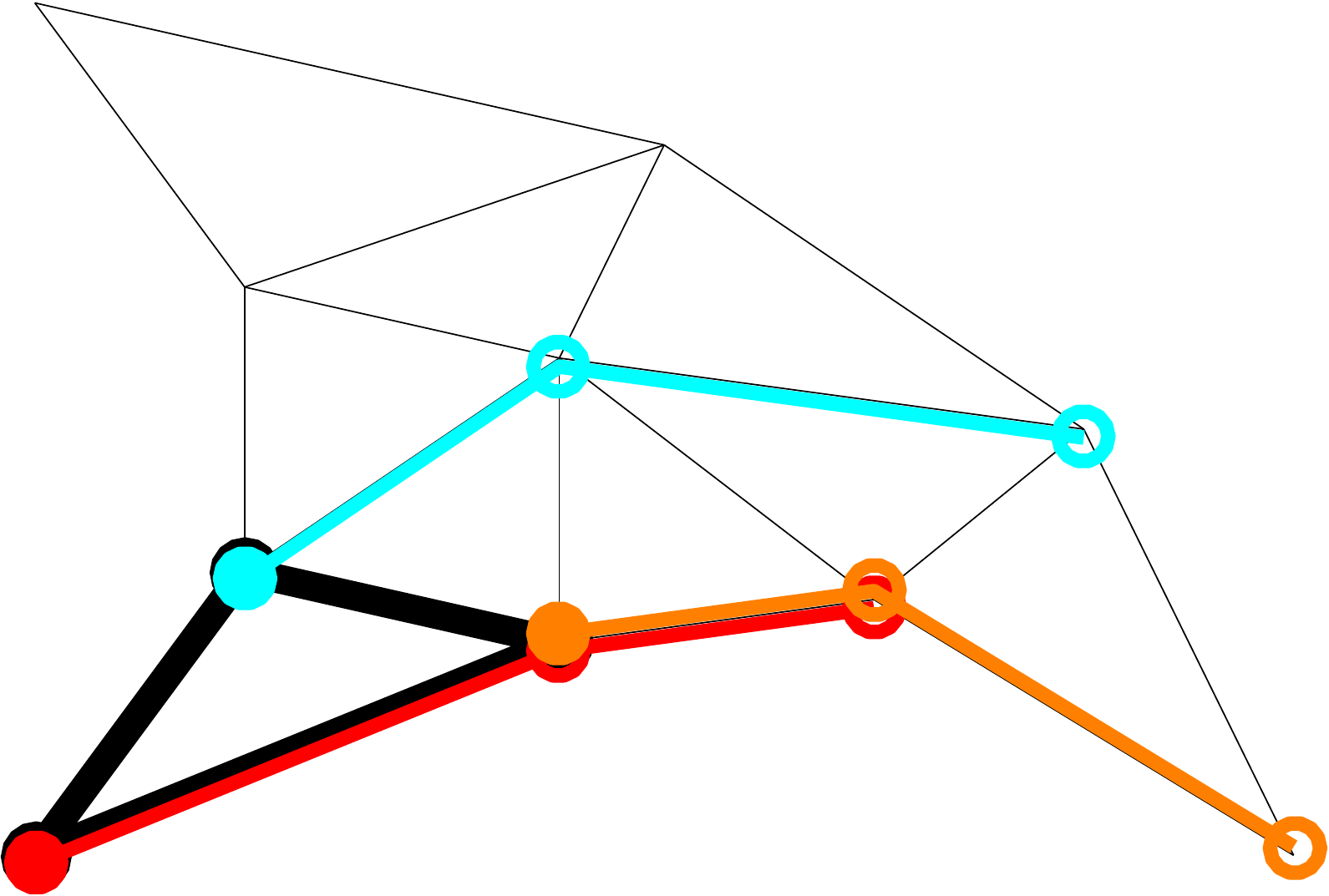} &
      \vspace{5pt} \includegraphics[width=3.0cm,height=3.0cm,keepaspectratio]{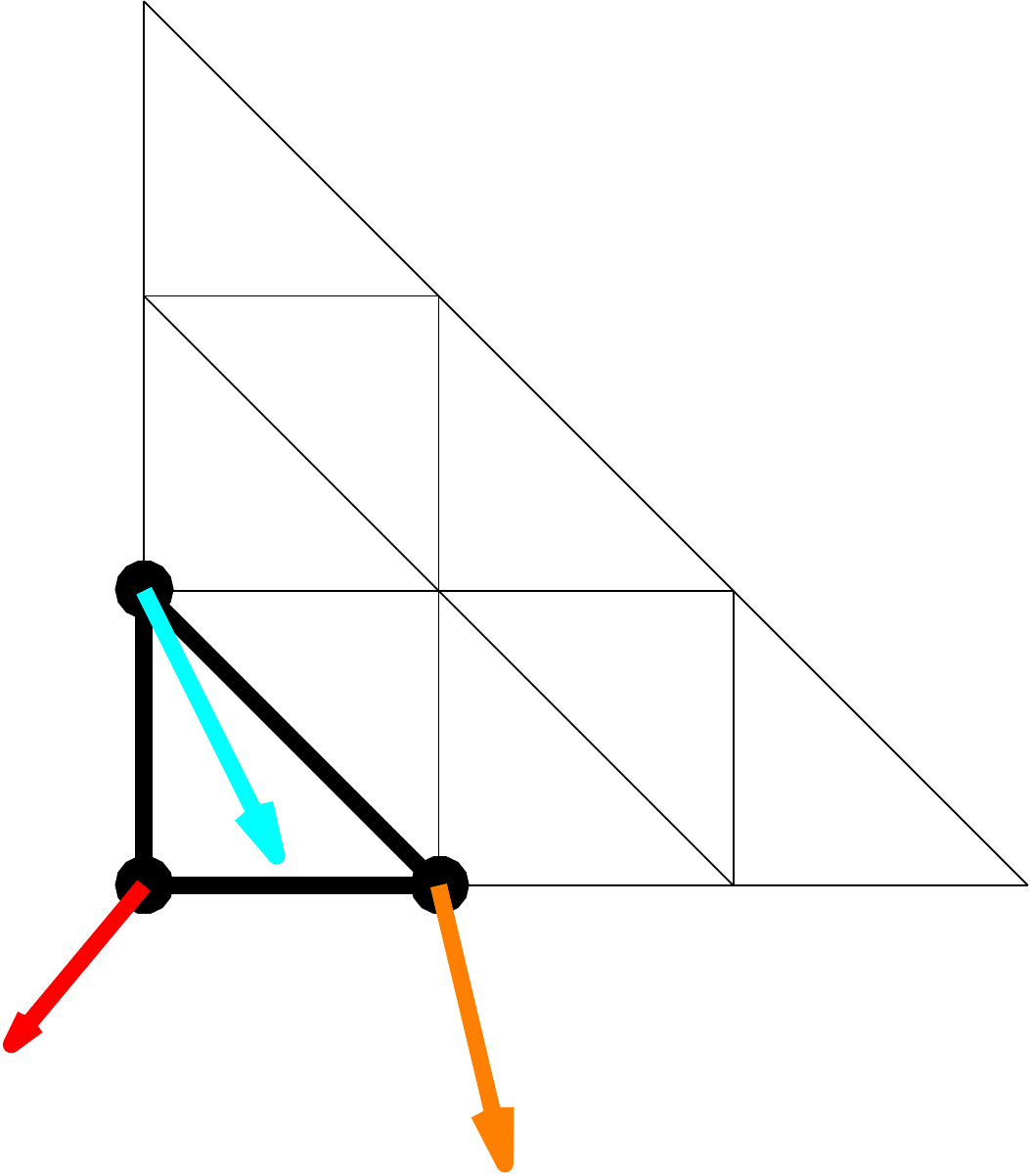} \hspace{10pt}\\
      \hline
      $\dfrac{\partial^3 \xproj}{\partial \xi_1^3} $&
      \includegraphics[width=0.18\textwidth]{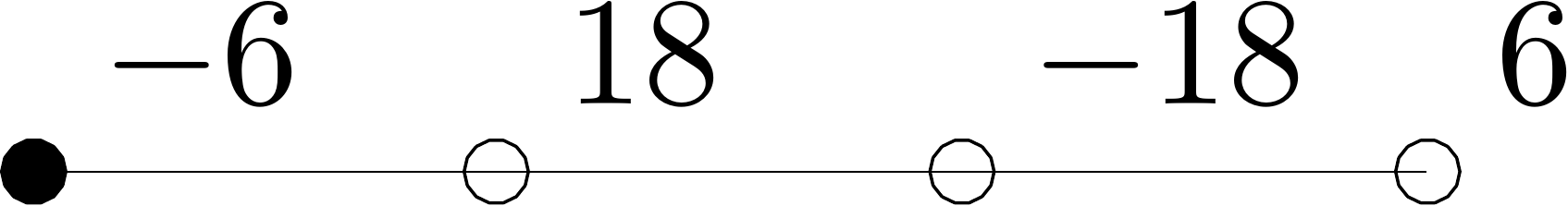} &
      \vspace{5pt} \includegraphics[width=4.5cm,height=4.5cm,keepaspectratio]{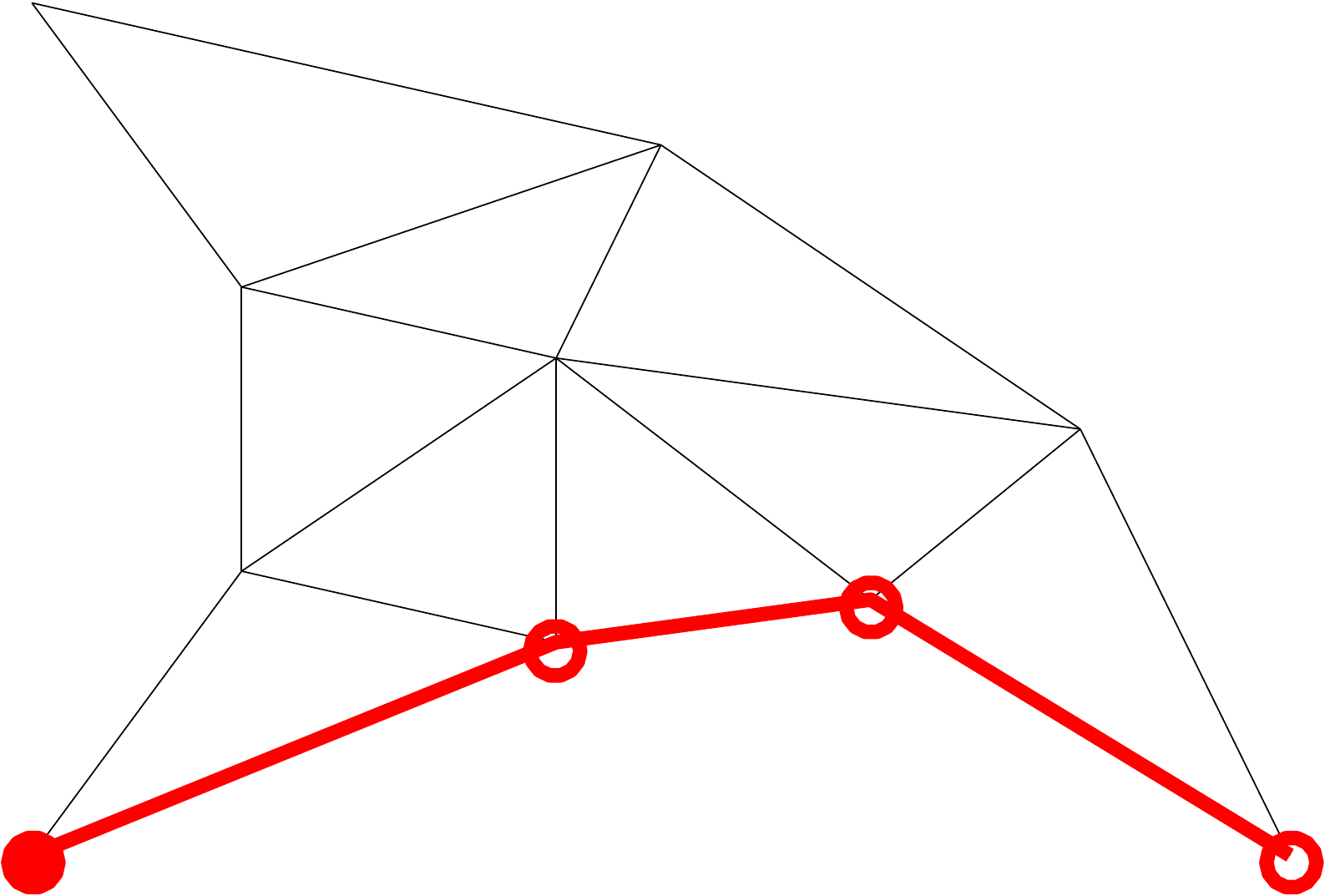} &
      \vspace{5pt} \includegraphics[width=3.5cm,height=3.5cm,keepaspectratio]{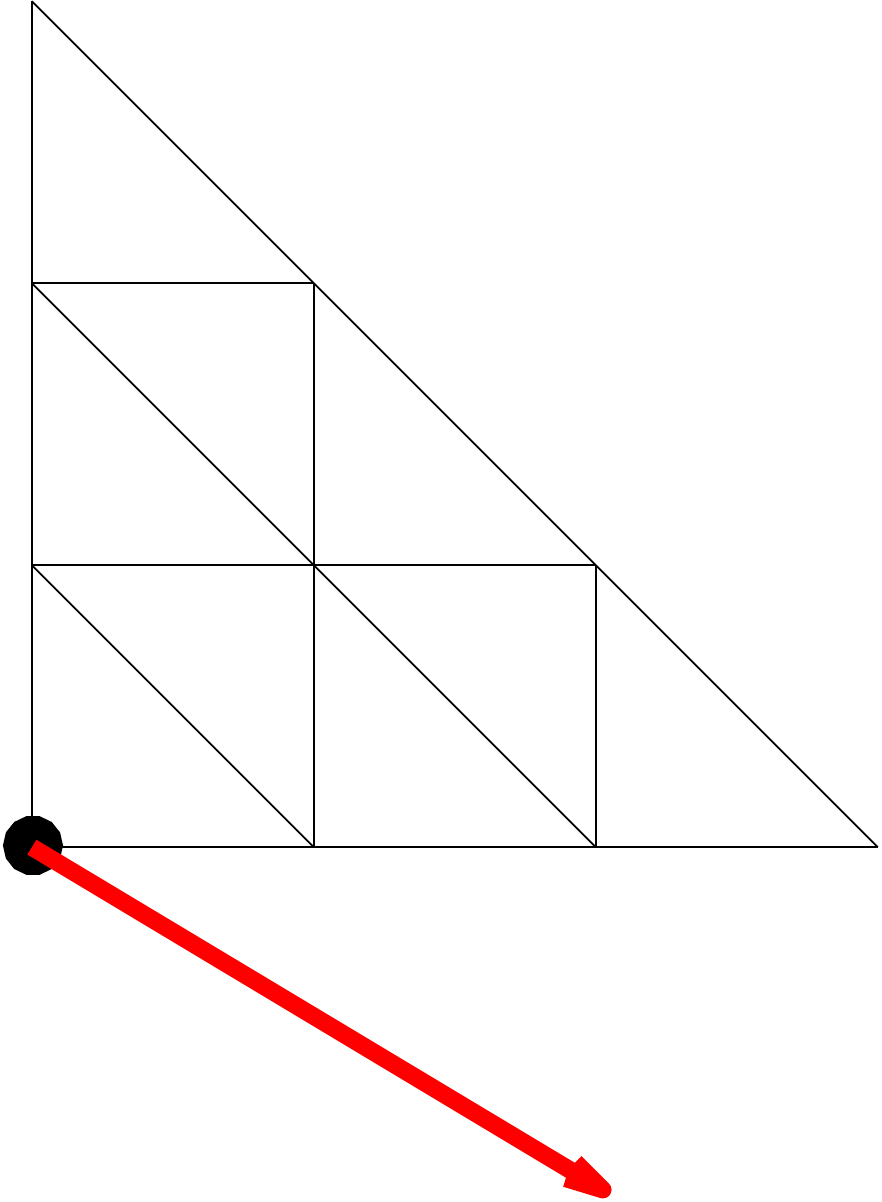}\\
      \hline
    \end{tabular}
    \label{deriv_examples}
  \end{center}
\end{table}


\section{Numerical Examples}\label{numerical}
In this section, we present several numerical examples to demonstrate how our element distortion metrics may be used in practice. Our goals are twofold. First, we desire to confirm our approximation results for shape regular refinements. That is, we want to demonstrate that a shape regular family of \BB meshes exhibits optimal convergence rates.
Second, we wish to demonstrate how the element distortion metrics presented here can be used for mesh optimization.

We provide four examples to benchmark our methods. 
First, we consider a simple rectangular plate, meshed with distorted \textbf{\emph{polynomial}} elements, to study the effect of control point distortion under $h$-refinement. 
Next, we consider a plate with a hole, meshed with distorted \textbf{\emph{rational}} elements, to examine the effect of weighting function distortion under $h$-refinement.
We then consider a quarter annulus, meshed with rational elements, under $p$-refinement, and conclude with an example of how our metrics may be used for mesh optimization.

The examples considered here are relatively simple, but they still demonstrate that poorly shaped elements can have appreciable impacts on solution accuracy. 
We also note that the examples shown here are constrained to the two dimensional case, as this allows for clear and easy visualization of element shape.
However, the implications of these 2D results extend immediately to elements in three dimensions. 

\subsection{Manufactured Solution on a Rectangular Plate}
To demonstrate the use of our validity metrics, we begin by considering several different meshes of a rectangular plate.
To account for both tensor product and simplicial elements, we consider both quadrilateral and triangular meshes.
The triangular meshes are formed by simply bisecting each element in the quadrilateral mesh. 
For both types of elements, we construct an initial mesh, and then create three families of refined meshes.

The initial curvilinear mesh is created by first creating a linear quadrilateral mesh, and degree elevating to \textbf{\emph{non-rational}} bi-cubic B\'{e}zier quadrilaterals.
Then, for each element, we horizontally perturb the middle two rows of control points, $\setof{ \bP_{\setof{i_1,i_2}}}_{i_1=\setof{0,1,2,3},i_2 = \setof{1,2} }$  by some distance:
\begin{equation}
d\bP_{\bi} =
\of{-1}^{i_2}\dfrac{2a}{\of{4m}^{7/4} } \dfrac{a - \absof{\of{\bP_\bi}_{x_1}}}{a}
\label{pert1}
\end{equation}
wherein $m$ denotes the $m^{th}$ mesh in the family, with $m=1$ being the first mesh.

\begin{figure}[t] 
\centering
    \begin{subfigure}[t]{0.49\textwidth}
        \centering
            \includegraphics[width=.9\linewidth]{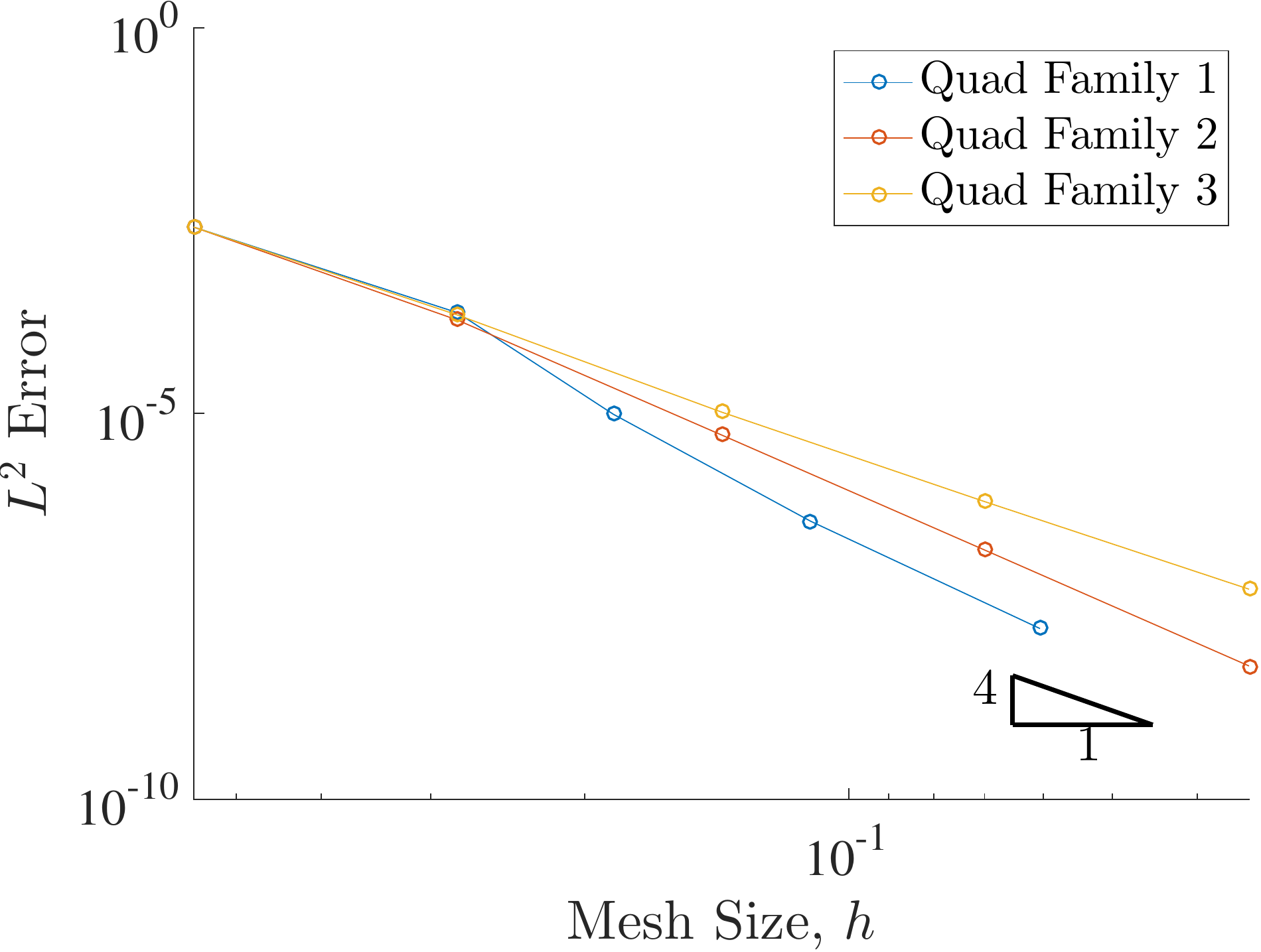}
        \caption{}
    \end{subfigure}
     \begin{subfigure}[t]{0.49\textwidth}
        \centering
            \includegraphics[width=.9\linewidth]{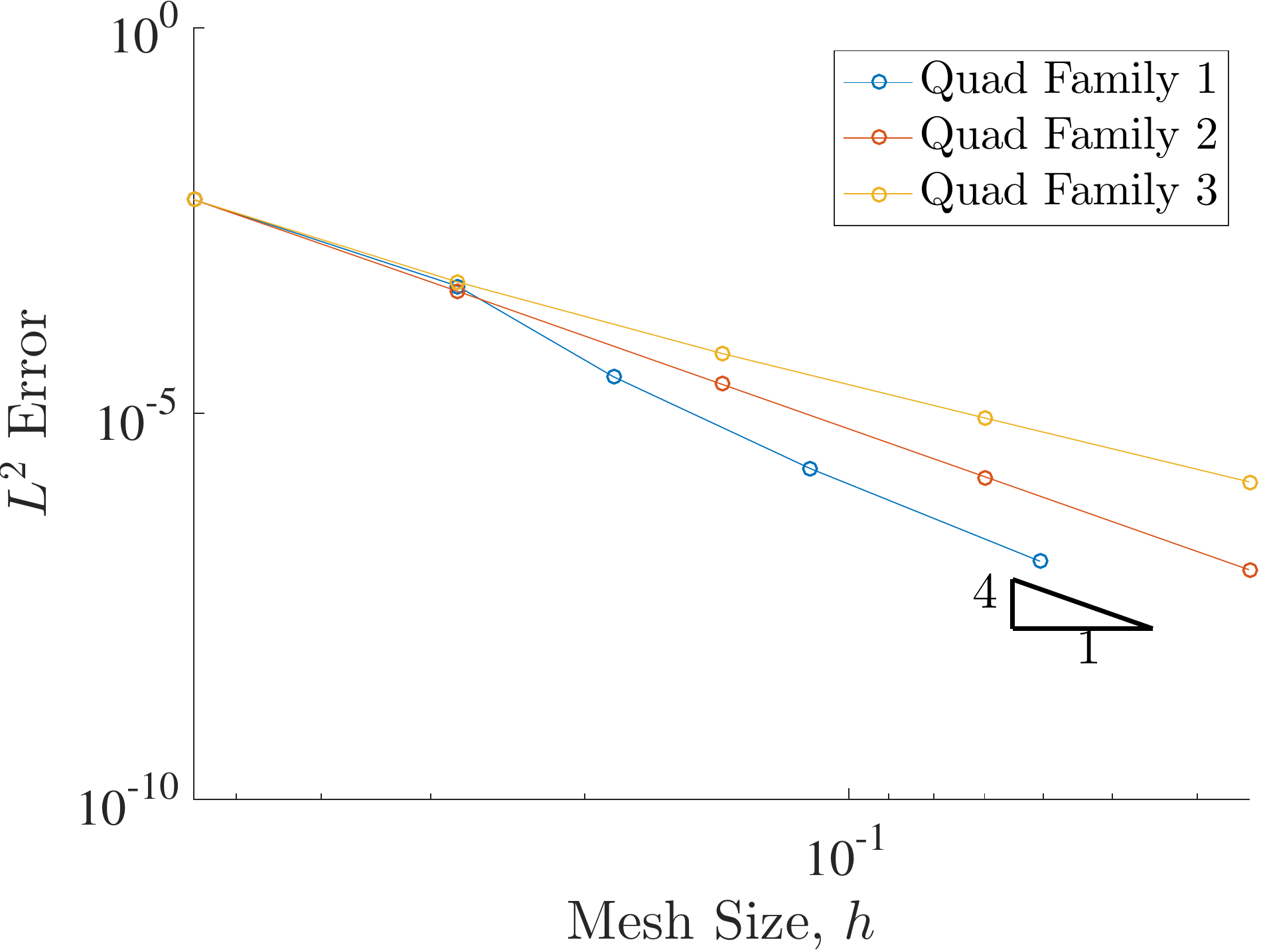}
        \caption{}
    \end{subfigure}
\caption{ Convergence plots for the quadrilateral meshes for the rectangular plate manufactured solution. (a) $L^2$ norm of the error for the $L^2$ projection problem. (b) $L^2$ norm of the error for Poisson's problem. }
\label{plate_quad_norm}
\end{figure}

\begin{figure}[t!] 
\centering
    \begin{subfigure}[t]{0.49\textwidth}
        \centering
            \includegraphics[width=.9\linewidth]{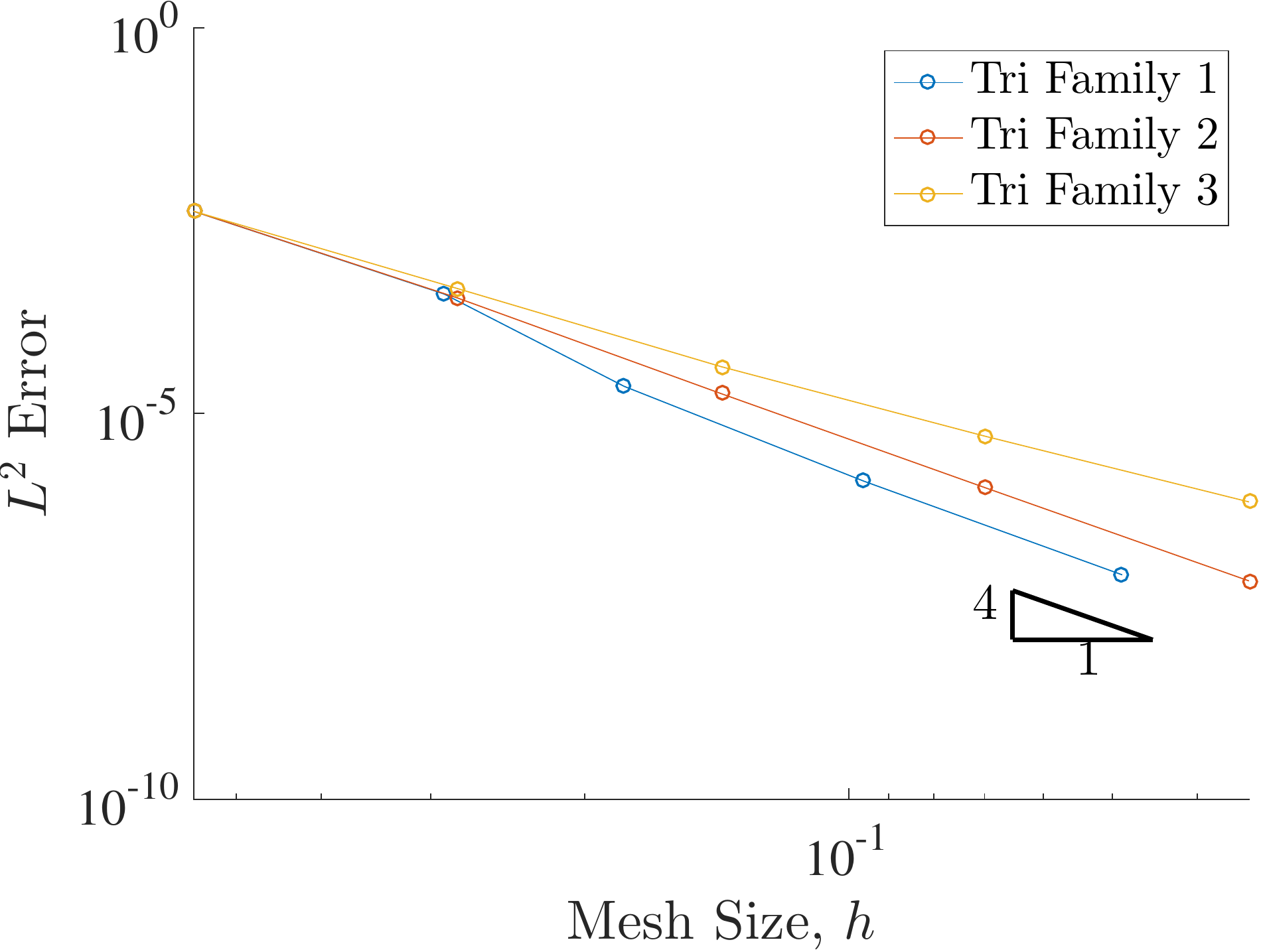}
        \caption{}
    \end{subfigure}
       \begin{subfigure}[t]{0.49\textwidth}
        \centering
            \includegraphics[width=.9\linewidth]{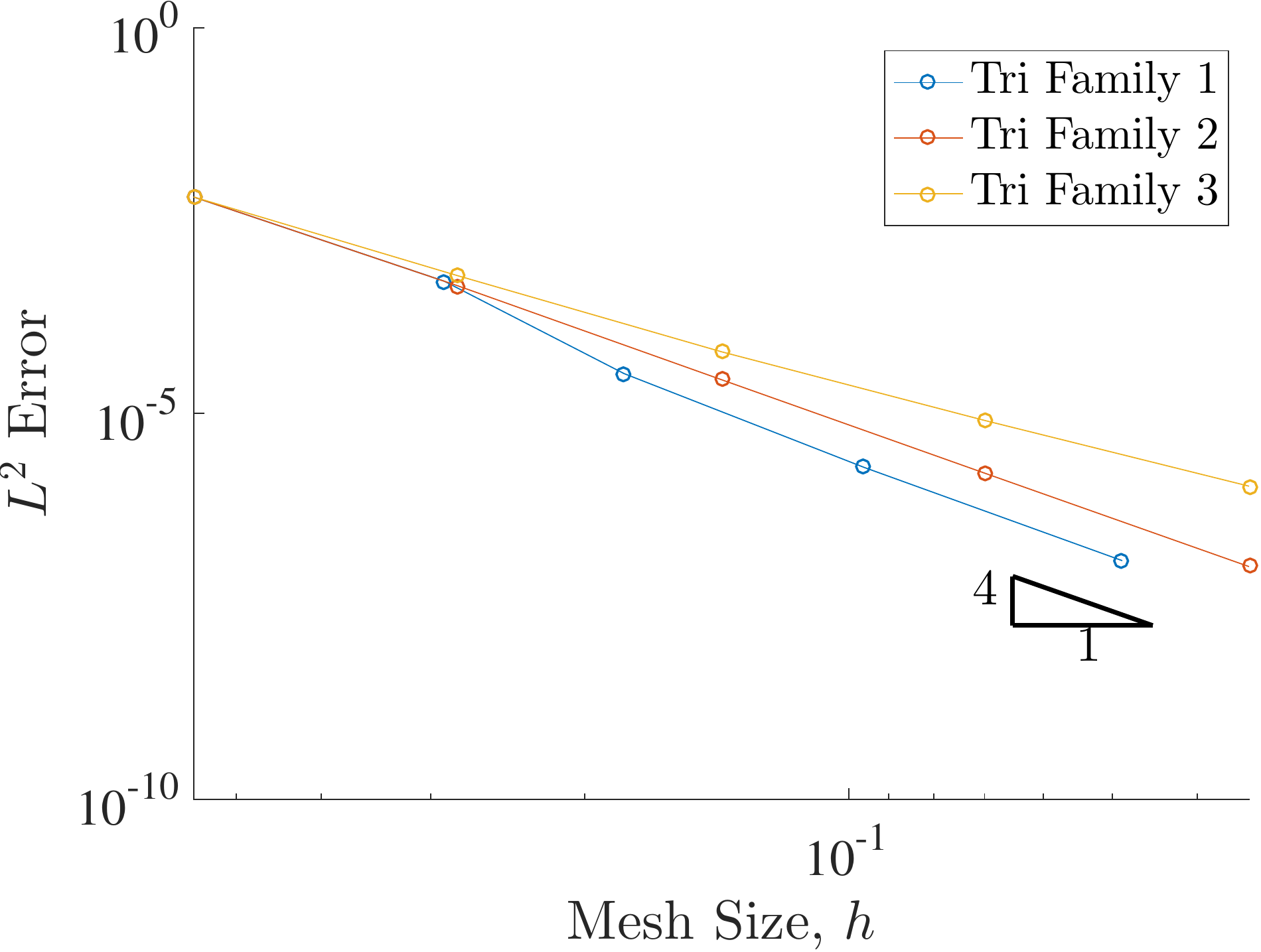}
        \caption{}
    \end{subfigure}

\caption{ Convergence plots for the triangular meshes  for the rectangular plate manufactured solution. (a) $L^2$ norm of the error for the $L^2$ projection problem. (b) $L^2$ norm of the error for Poisson's problem. }
\label{plate_tri_norm}
\end{figure}

Then, for both the quadrilateral and triangular mesh, we create the three families of refined meshes as follows.
The first family of meshes, shown in Table \ref{plate_family1}, is created by simple uniform subdivision of the original mesh.  
To create the second family of meshes, shown in Table \ref{plate_family2}, we first perform uniform subdivision on the original \textbf{\textit{linear}} mesh. 
We then create the $m^{th}$ curvilinear mesh in the family by degree elevation and again perturbing the interior control points using the prescription given by Eq. \eqref{pert1}.
The final family of meshes, shown in Table \ref{plate_family3}, is created analogously to the second family, but the perturbation distance is instead given by the equation:
\begin{equation}
d\bP_{\bi} =
\of{-1}^{i_2}\dfrac{8a}{\of{4m}^3 } \dfrac{a - \absof{\of{\bP_\bi}_{x_1}}}{a}
\label{pert2}
\end{equation}
%
%
%
%
%
%
%
\begin{table}[b!]
  \begin{center}
    \caption{Mesh Family 1 for the Rectangular Plate Manufactured Solution.}
    \begin{tabular}{|  >{\centering\arraybackslash} m{.3cm} |  >{\centering\arraybackslash} m{7.1cm} |  >{\centering\arraybackslash} m{7.1cm} |  }
      \hline
       $m$ & Quadrilaterals & Triangles  \\ 
      \hline
      1 &
      \vspace{5pt}
      \includegraphics[width=6.9cm,height=6.9cm,keepaspectratio]{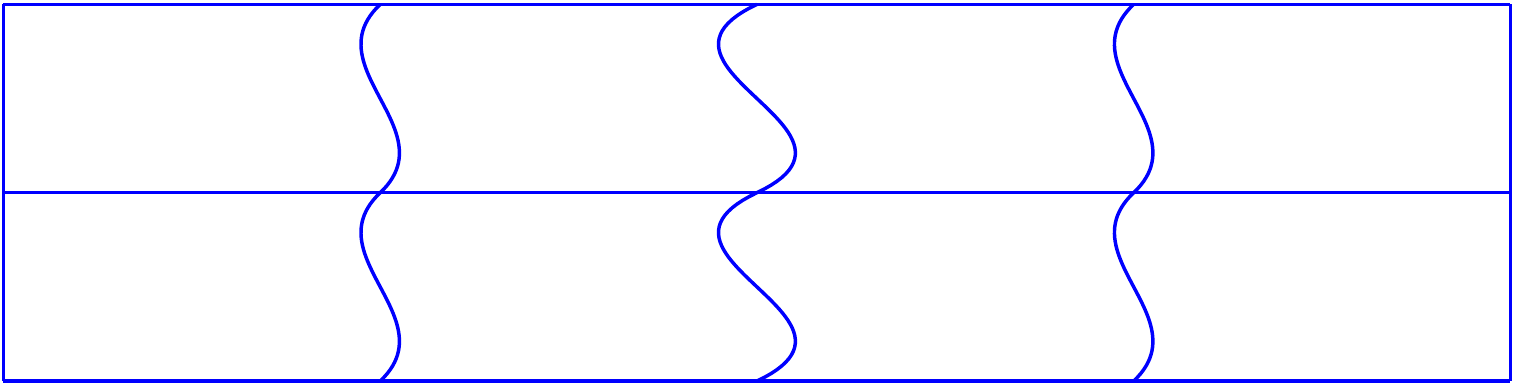} &
       \vspace{5pt}
      \includegraphics[width=6.9cm,height=6.9cm,keepaspectratio]{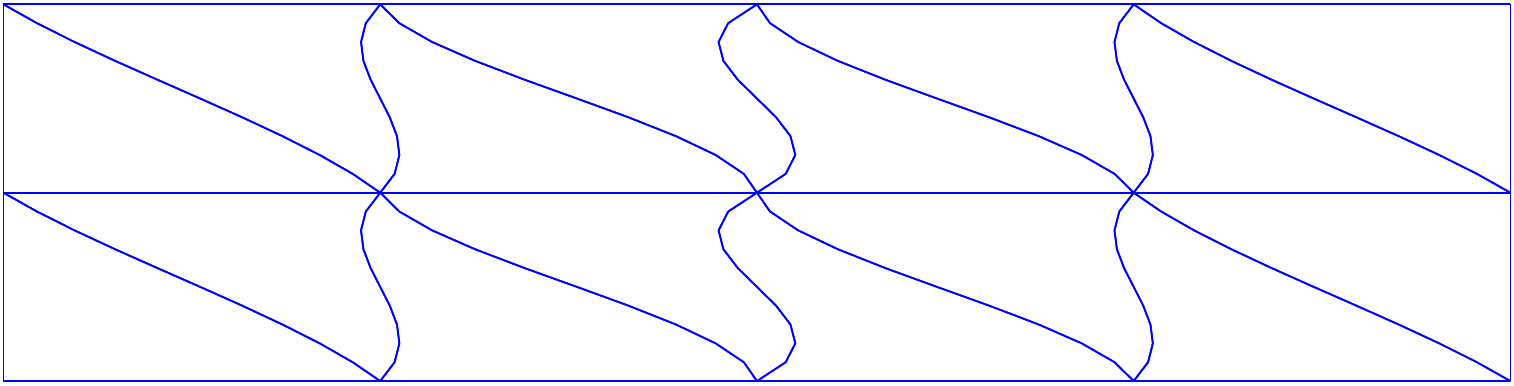}\\
      \hline
            2 &
             \vspace{5pt}
      \includegraphics[width=6.9cm,height=6.9cm,keepaspectratio]{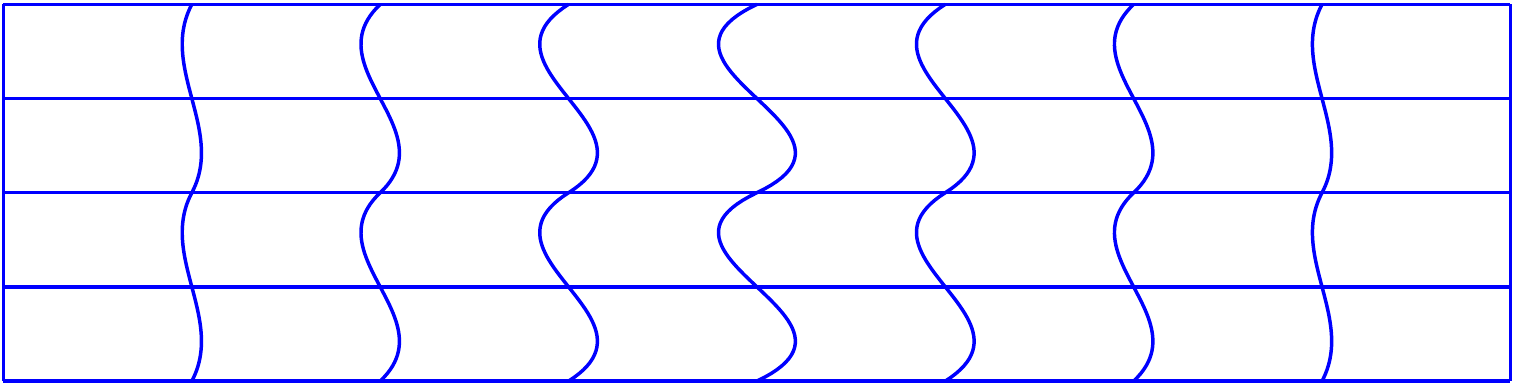} &
             \vspace{5pt}
      \includegraphics[width=6.9cm,height=6.9cm,keepaspectratio]{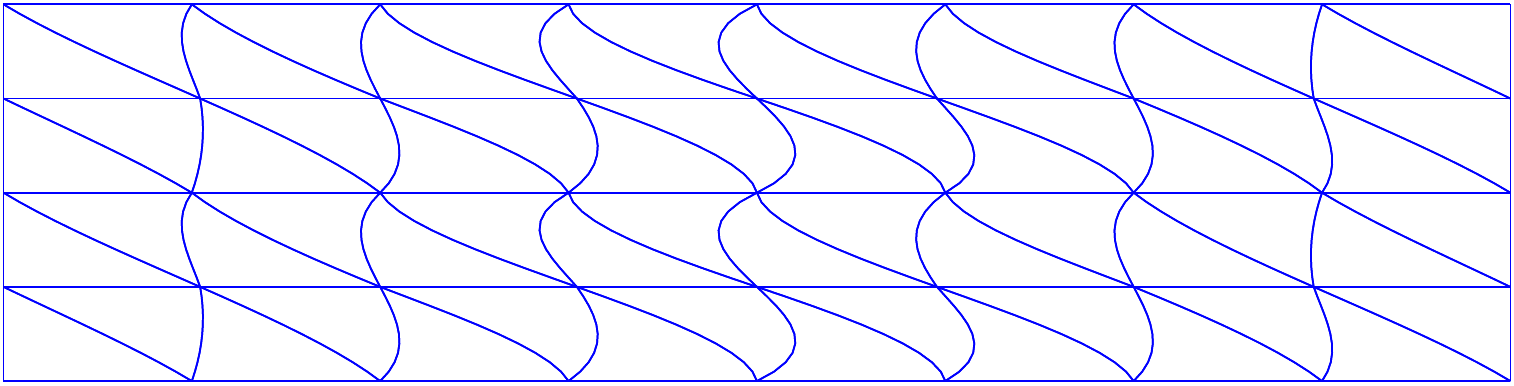}\\
	\hline
	      3 &
	       \vspace{5pt}
      \includegraphics[width=6.9cm,height=6.9cm,keepaspectratio]{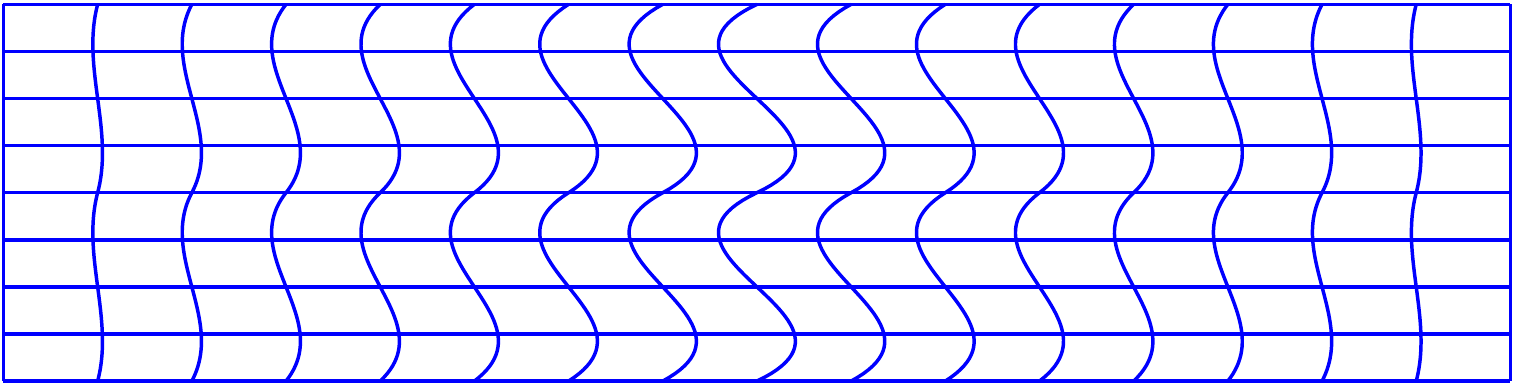} &
             \vspace{5pt}
      \includegraphics[width=6.9cm,height=6.9cm,keepaspectratio]{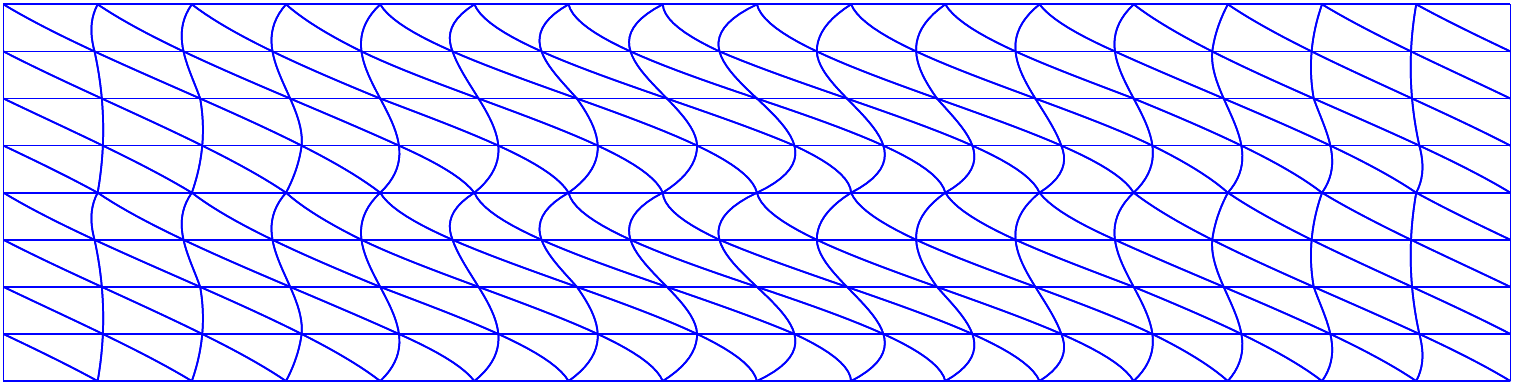}\\
      \hline
      \end{tabular}
    \label{plate_family1}
  \end{center}
\end{table}
%
%
%
\begin{table}[t]
  \begin{center}
    \caption{Mesh Family 2 for the Rectangular Plate Manufactured Solution.}
    \begin{tabular}{|  >{\centering\arraybackslash} m{.3cm} |  >{\centering\arraybackslash} m{7.1cm} |  >{\centering\arraybackslash} m{7.1cm} |  }
      \hline
       $m$ & Quadrilaterals & Triangles  \\ 
      \hline
       1 &
      \vspace{5pt}
      \includegraphics[width=6.9cm,height=6.9cm,keepaspectratio]{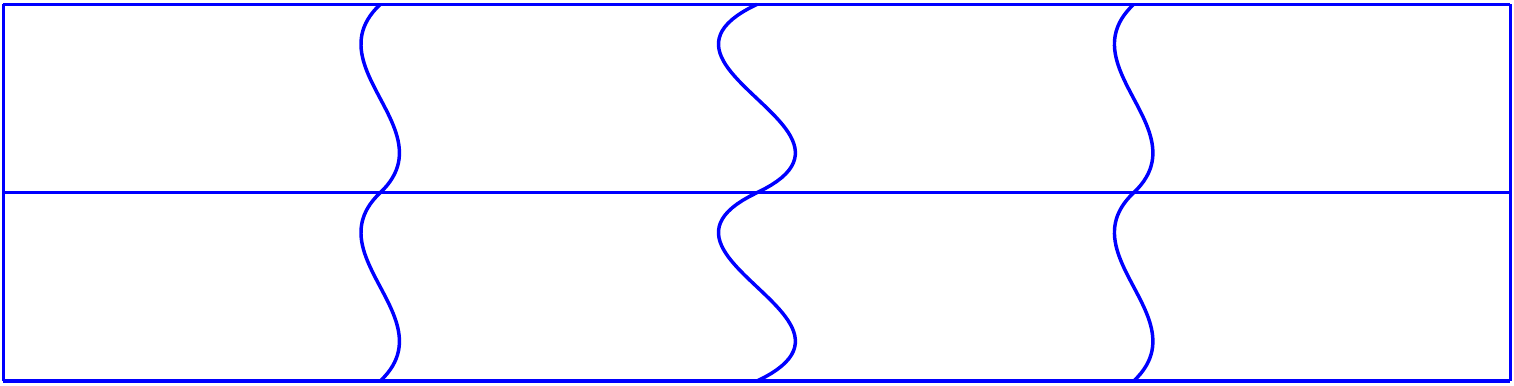} &
       \vspace{5pt}
      \includegraphics[width=6.9cm,height=6.9cm,keepaspectratio]{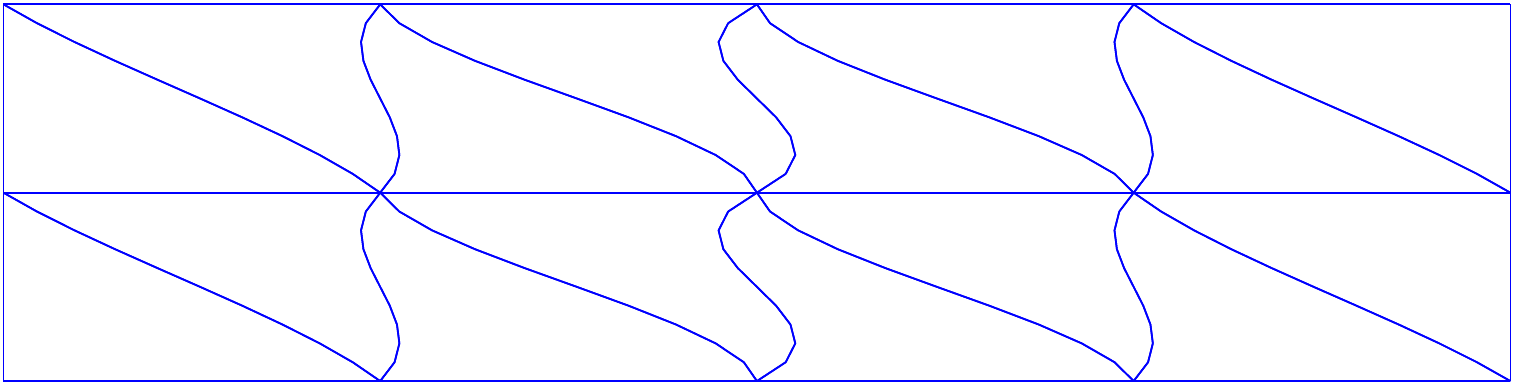}\\
      \hline
       2 &
             \vspace{5pt}
      \includegraphics[width=6.9cm,height=6.9cm,keepaspectratio]{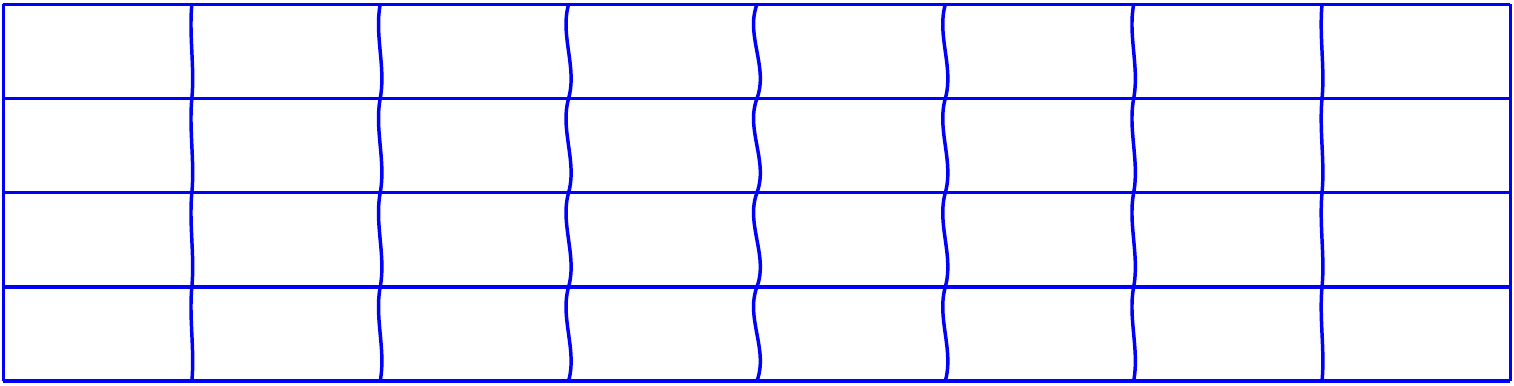} &
             \vspace{5pt}
      \includegraphics[width=6.9cm,height=6.9cm,keepaspectratio]{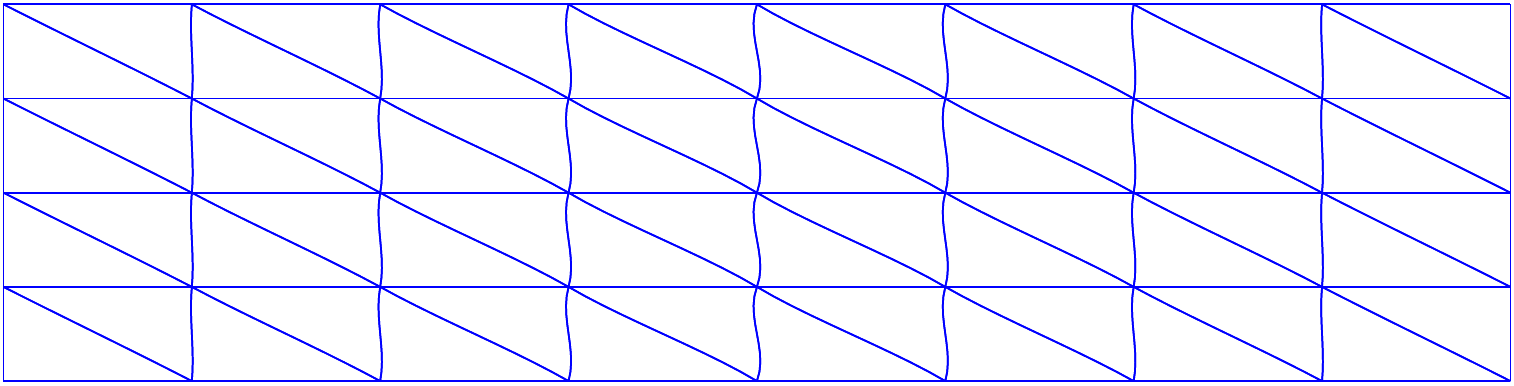}\\
	\hline
	 3 &
	       \vspace{5pt}
      \includegraphics[width=6.9cm,height=6.9cm,keepaspectratio]{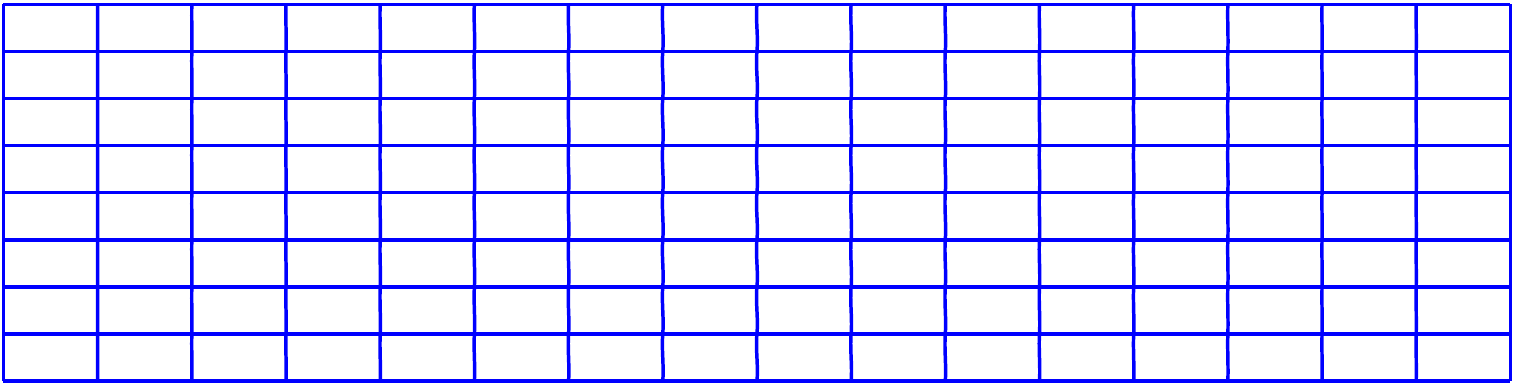} &
             \vspace{5pt}
      \includegraphics[width=6.9cm,height=6.9cm,keepaspectratio]{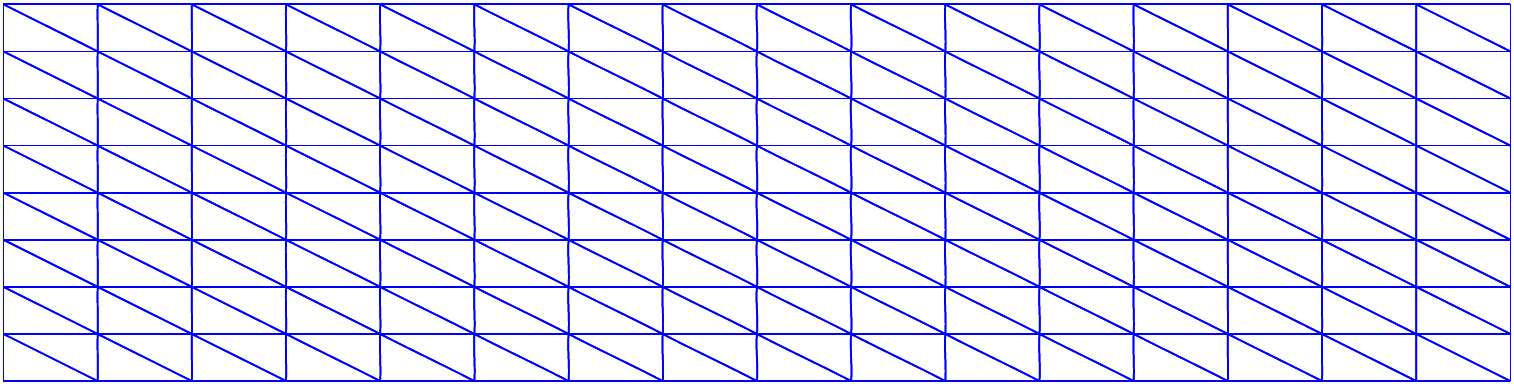}\\
      \hline
      \end{tabular}
    \label{plate_family2}
  \end{center}
\end{table}
%
%
%
\begin{table}[h!]
  \begin{center}
    \caption{Mesh Family 3 for the Rectangular Plate Manufactured Solution.} 
    \begin{tabular}{|  >{\centering\arraybackslash} m{.3cm} |  >{\centering\arraybackslash} m{7.1cm} |  >{\centering\arraybackslash} m{7.1cm} |  }
      \hline
       $m$ & Quadrilaterals & Triangles  \\ 
      \hline
       1 &
      \vspace{5pt}
      \includegraphics[width=6.9cm,height=6.9cm,keepaspectratio]{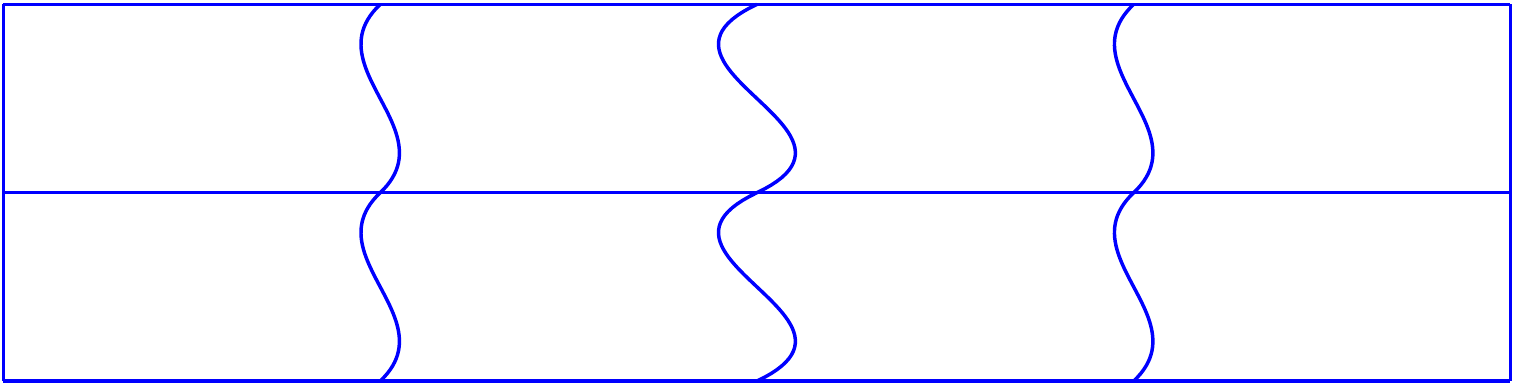} &
       \vspace{5pt}
      \includegraphics[width=6.9cm,height=6.9cm,keepaspectratio]{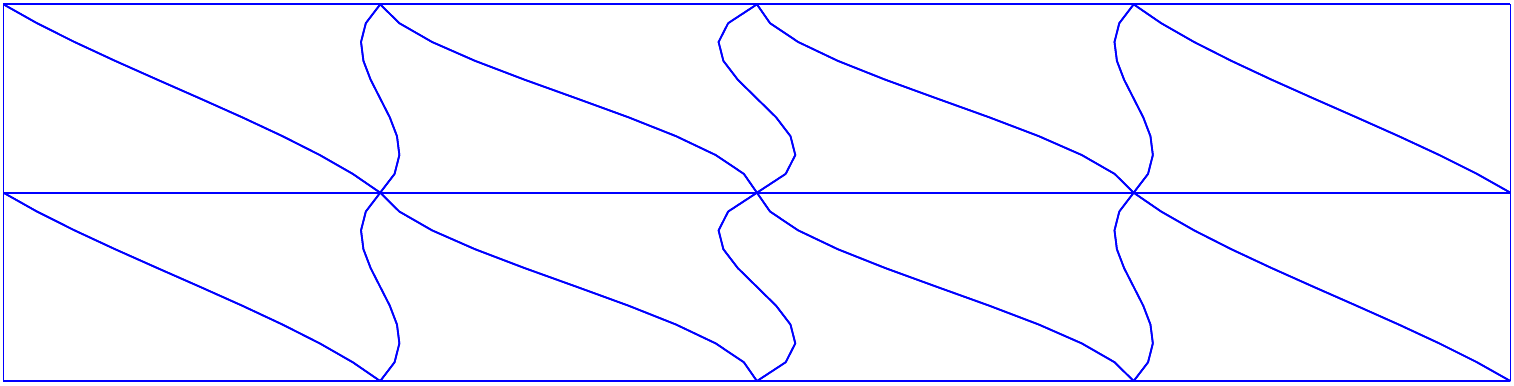}\\
      \hline
       2 &
             \vspace{5pt}
      \includegraphics[width=6.9cm,height=6.9cm,keepaspectratio]{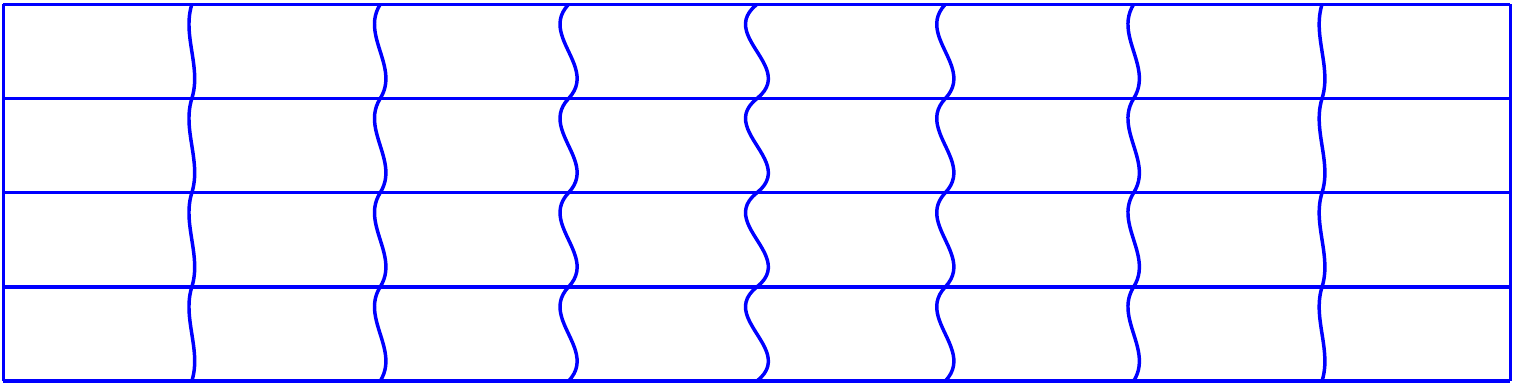} &
             \vspace{5pt}
      \includegraphics[width=6.9cm,height=6.9cm,keepaspectratio]{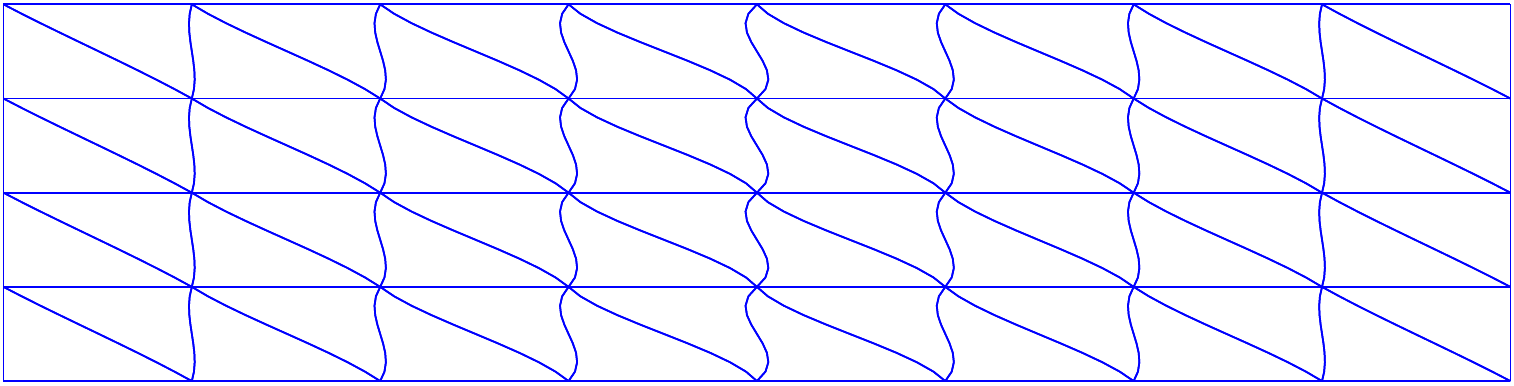}\\
	\hline
	 3 &
	       \vspace{5pt}
      \includegraphics[width=6.9cm,height=6.9cm,keepaspectratio]{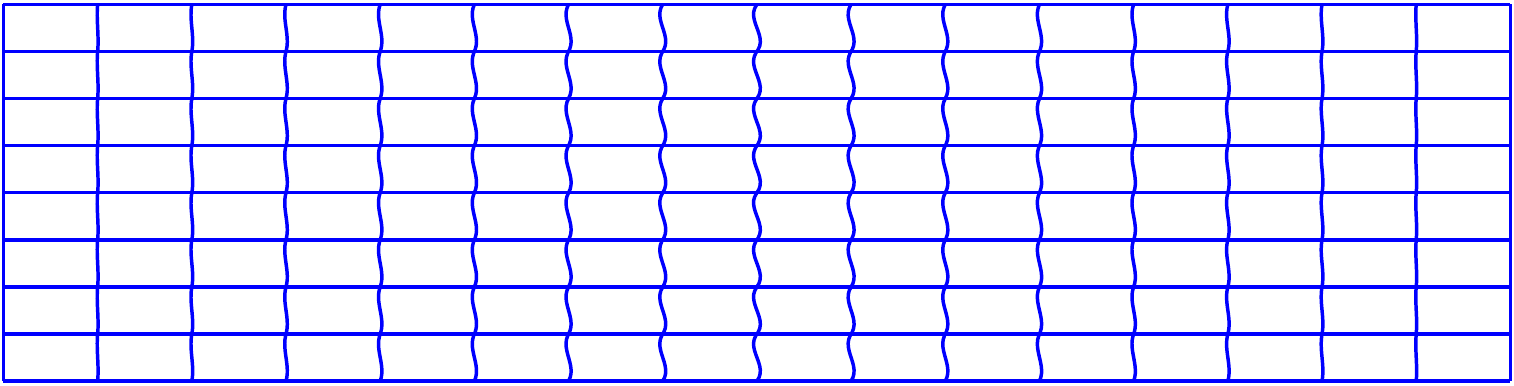} &
             \vspace{5pt}
      \includegraphics[width=6.9cm,height=6.9cm,keepaspectratio]{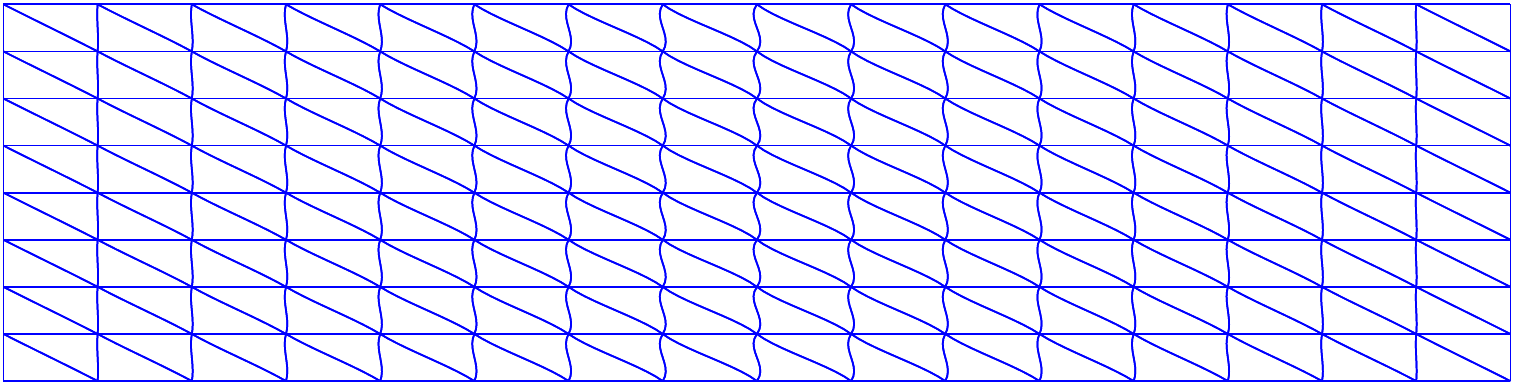}\\
      \hline
      \end{tabular}
    \label{plate_family3}
  \end{center}
\end{table}
With these three families of meshes established for both the quadrilateral and triangular case, we use the method of manufactured solutions to  study approximation error in each family of meshes. 
However, when solving partial differential equations (PDEs) using finite elements, error can be introduced not only by the element shape, but also by the choice of finite element method (e.g. Galerkin's method).
As such, over each family of meshes, we solve two problems, an $L^2$ projection problem and the Poisson problem.
We solve the $L^2$ projection over the mesh so that we may isolate the effect of element shape on approximation error.
We then consider the Poisson problem so as to consider an example with practical engineering applications.
Given the domain $\Omega$, let $\mathcal{f}$ denote a forcing function and let  $\mathcal{h}$ denote a flux across the boundary $\Gamma_{\mathcal{h}}$. Then, letting $\mathcal{S}_h$ and $\mathcal{V}_h$ denote the spaces of trial solutions (satisfying some prescribed Dirichlet boundary conditions) and test functions (satisfying homogeneous Dirichlet boundary conditions) respectively, the Poisson problem consists of finding a discrete solution $u_h \in \mathcal{S}_h$ such that for all $v_h \in \mathcal{V}_h$ :
\begin{equation}
\int\limits_{\Omega} \bnabla v_h \cdot \bnabla u_h d\Omega = 
\int\limits_{\Omega}  v_h f d\Omega +
\int\limits_{\Gamma_h}  v_h h d\Gamma
\end{equation}
For both cases, we attempt to approximate the manufactured solution
\begin{equation*}
u\of{\bx} = (x_1-a)(x_2-b)\cos\of{\dfrac{x_2}{2a} {\pi} } \sin\of{\dfrac{x_2}{b}\pi}
\end{equation*}
wherein $a$ and $b$ are the half-width and half-height of the plate centered at the origin.
To study the approximation error, we examine the convergence rate of the error $||u-u_h||$  in  the $L^2$ norm for both the $L^2$ projection and Poisson problems over each family of meshes. 
Fig. \ref{plate_quad_norm} shows error convergence plots for the three families of quadrilateral meshes, and Fig. \ref{plate_tri_norm} shows convergence plots for the three families of triangular meshes.
%
%
%
%
%
%
%
%
%
%
%
%
%
%
%
\begin{figure}[t] 
\centering
\begin{subfigure}[t]{0.48\textwidth}
        \centering
            \includegraphics[width=.95\linewidth]{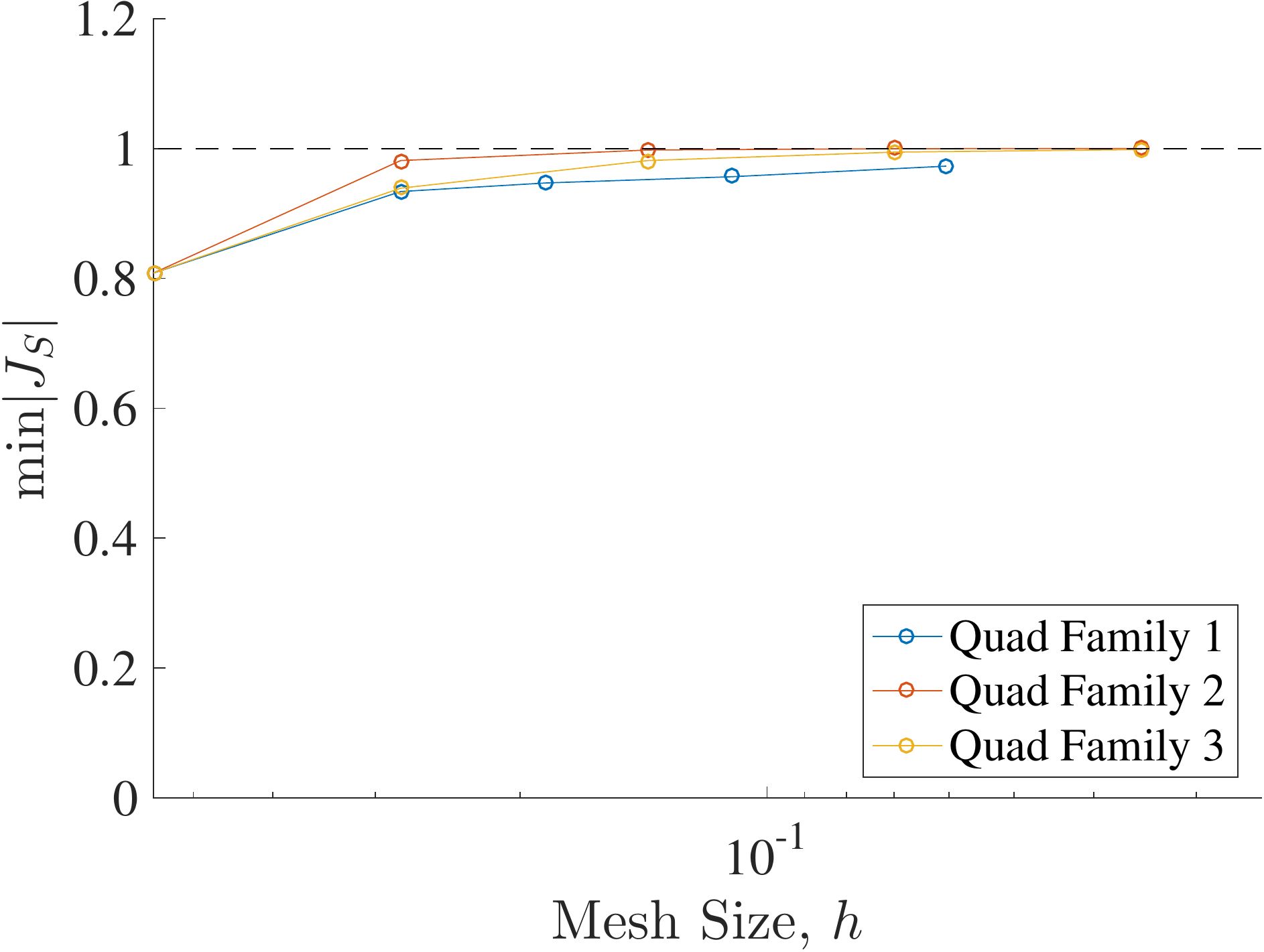}
        \caption{}
    \end{subfigure}
 \begin{subfigure}[t]{0.48\textwidth}
        \centering
            \includegraphics[width=.95\linewidth]{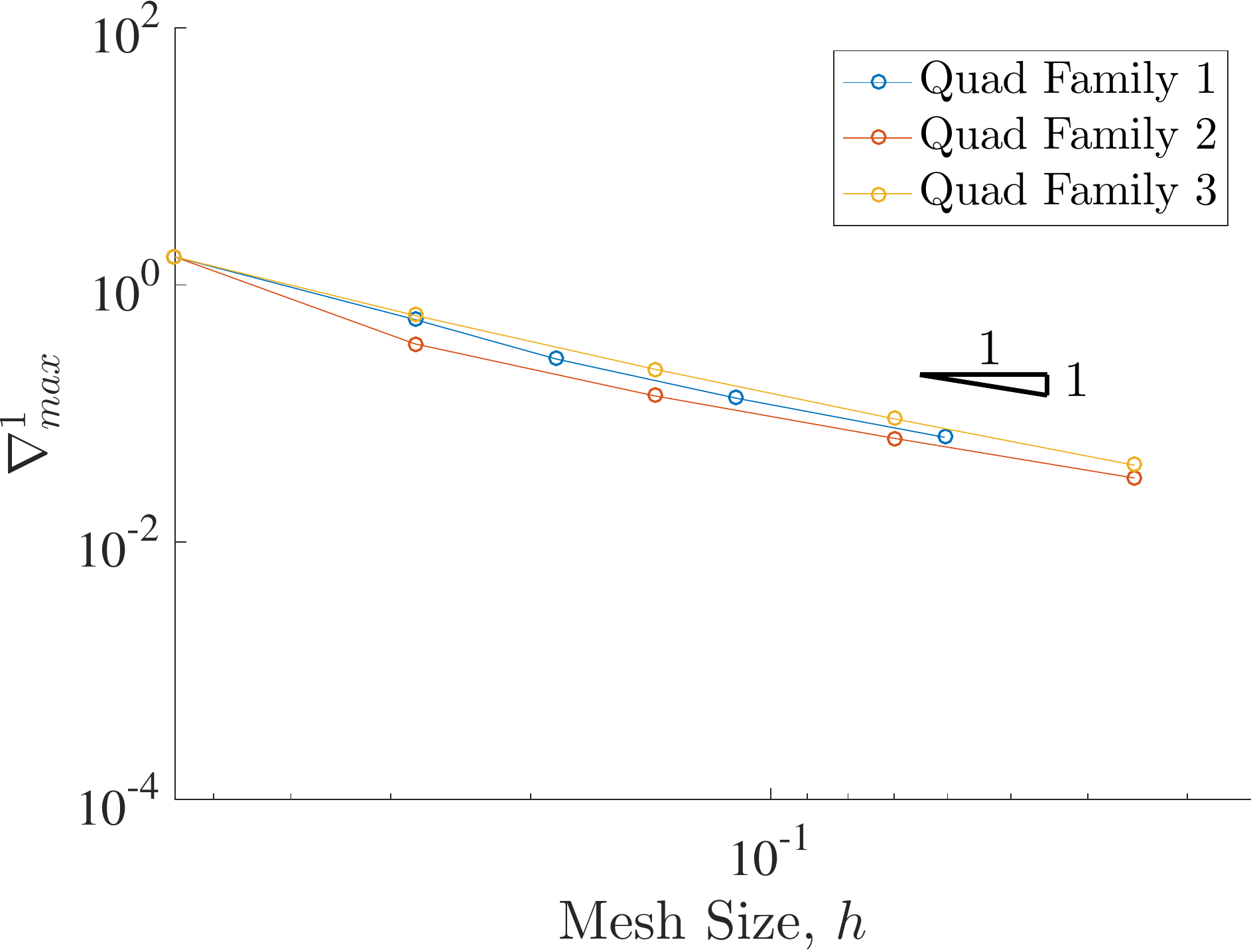}
        \caption{}
    \end{subfigure}
    \begin{subfigure}[t]{0.48\textwidth}
        \centering
            \includegraphics[width=.95\linewidth]{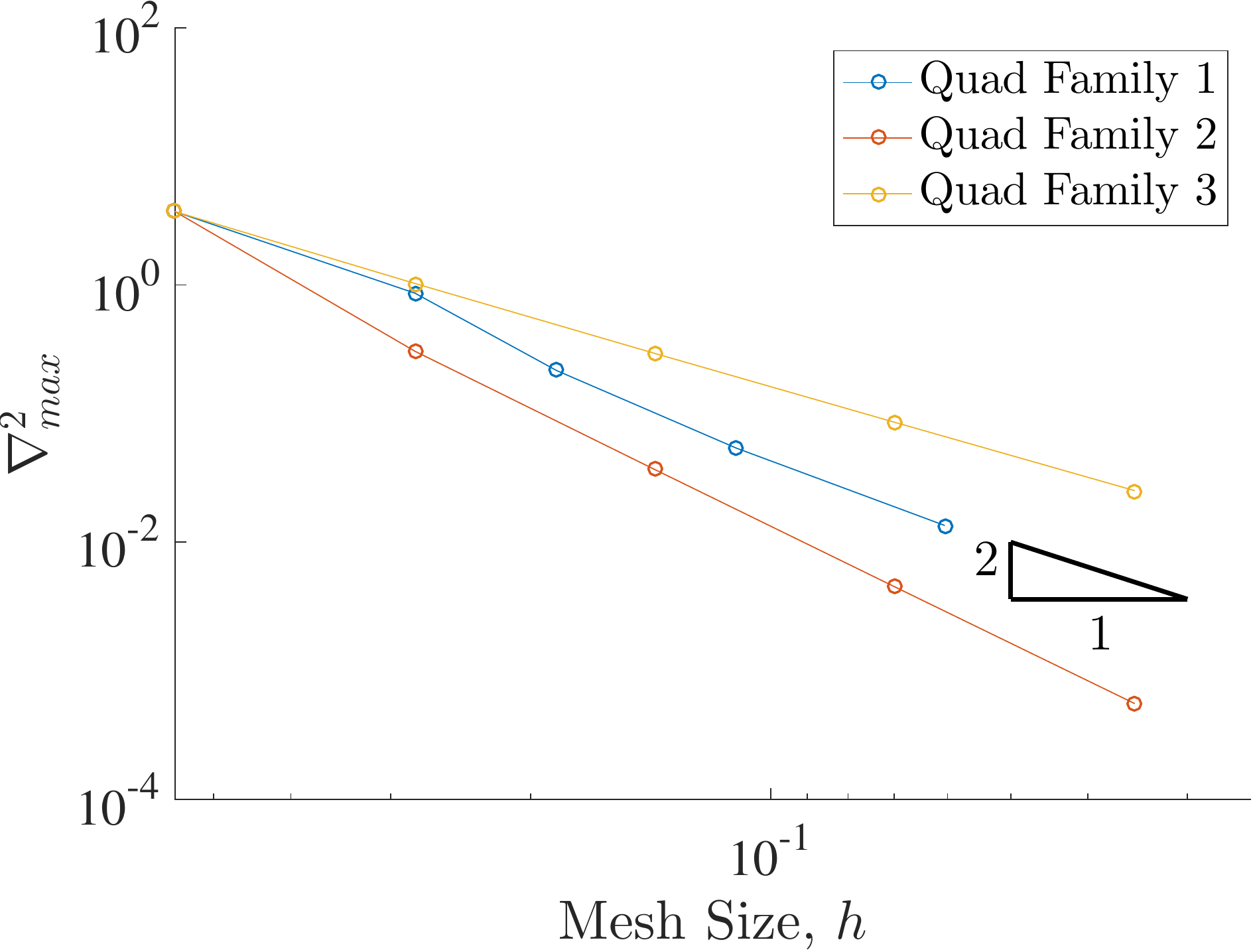}
        \caption{}
    \end{subfigure}
\begin{subfigure}[t]{0.48\textwidth}
        \centering
            \includegraphics[width=.95\linewidth]{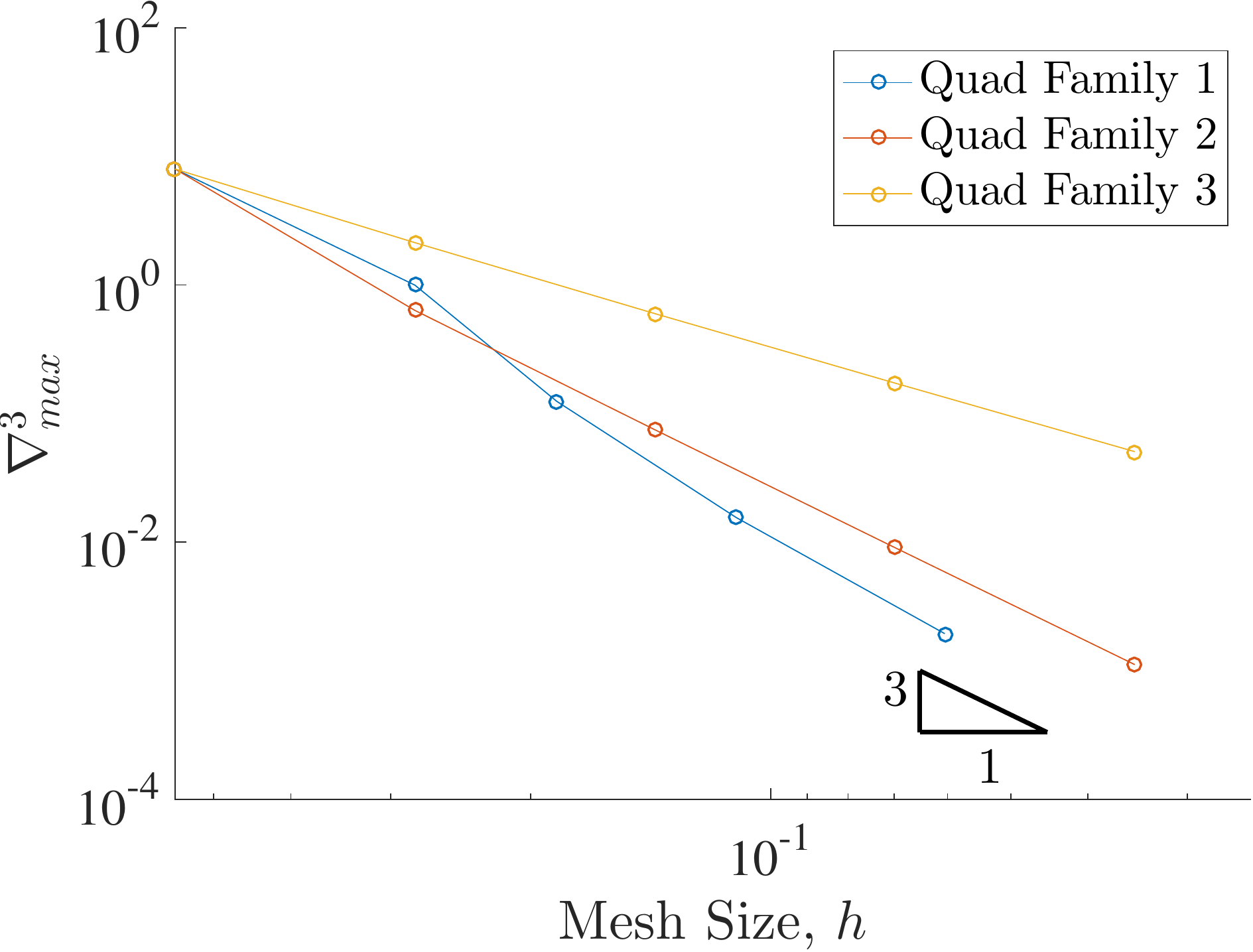}
        \caption{}
    \end{subfigure}
\caption{ Mesh distortion metrics for the quadrilateral meshes of the plate. (a) Minimum scaled Jacobian. (b) Lowest upper bound on the magnitude of the first derivatives. (c) Lowest upper bound on the magnitude of the second derivatives. (d) Lowest upper bound on the magnitude of the third derivatives.}
\label{mms_plate_quad_metrics}
\end{figure}
%
%
%
%
\begin{figure}[t] 
\centering
\begin{subfigure}[t]{0.48\textwidth}
        \centering
            \includegraphics[width=.95\linewidth]{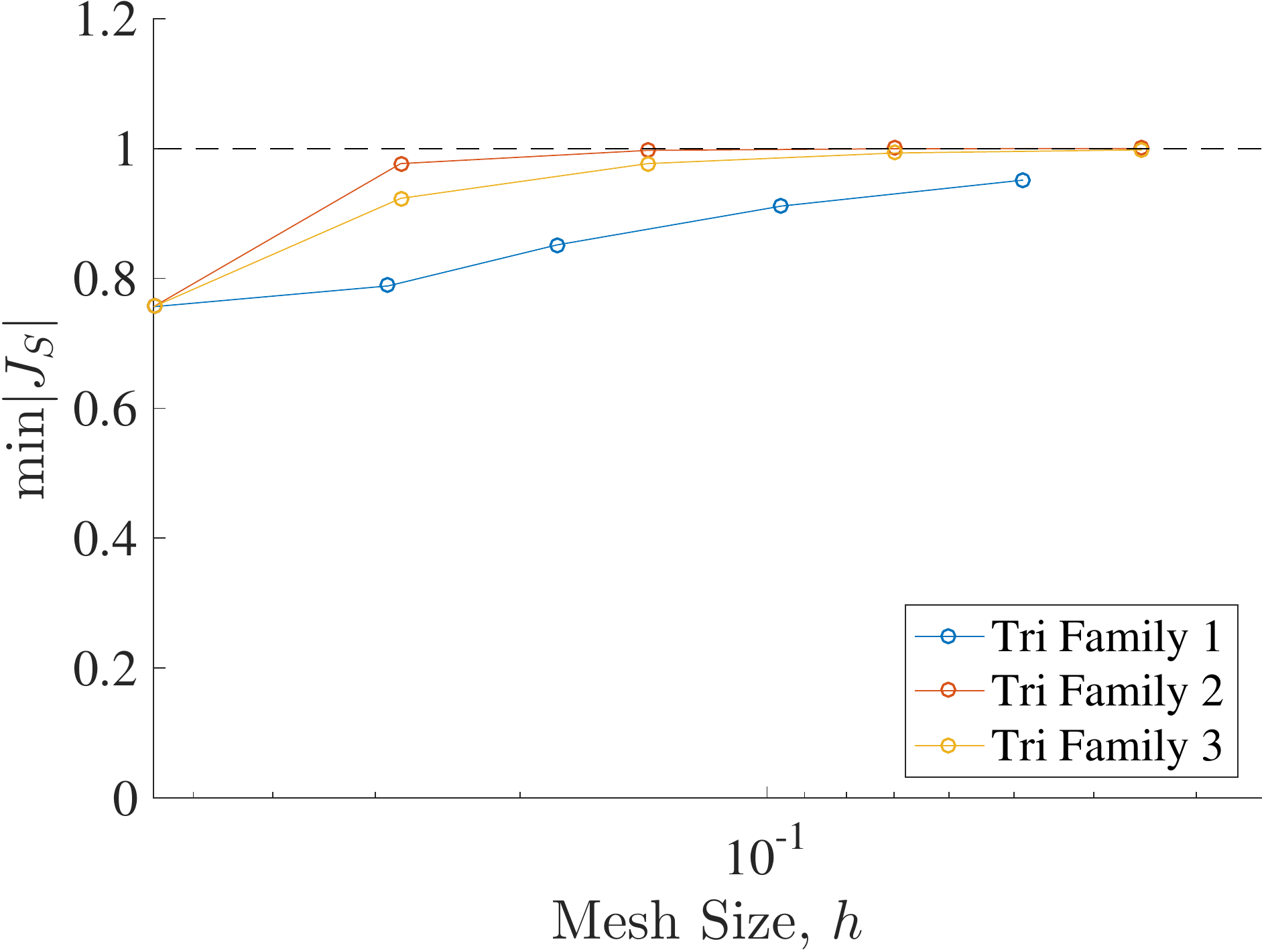}
        \caption{}
    \end{subfigure}
 \begin{subfigure}[t]{0.48\textwidth}
        \centering
            \includegraphics[width=.95\linewidth]{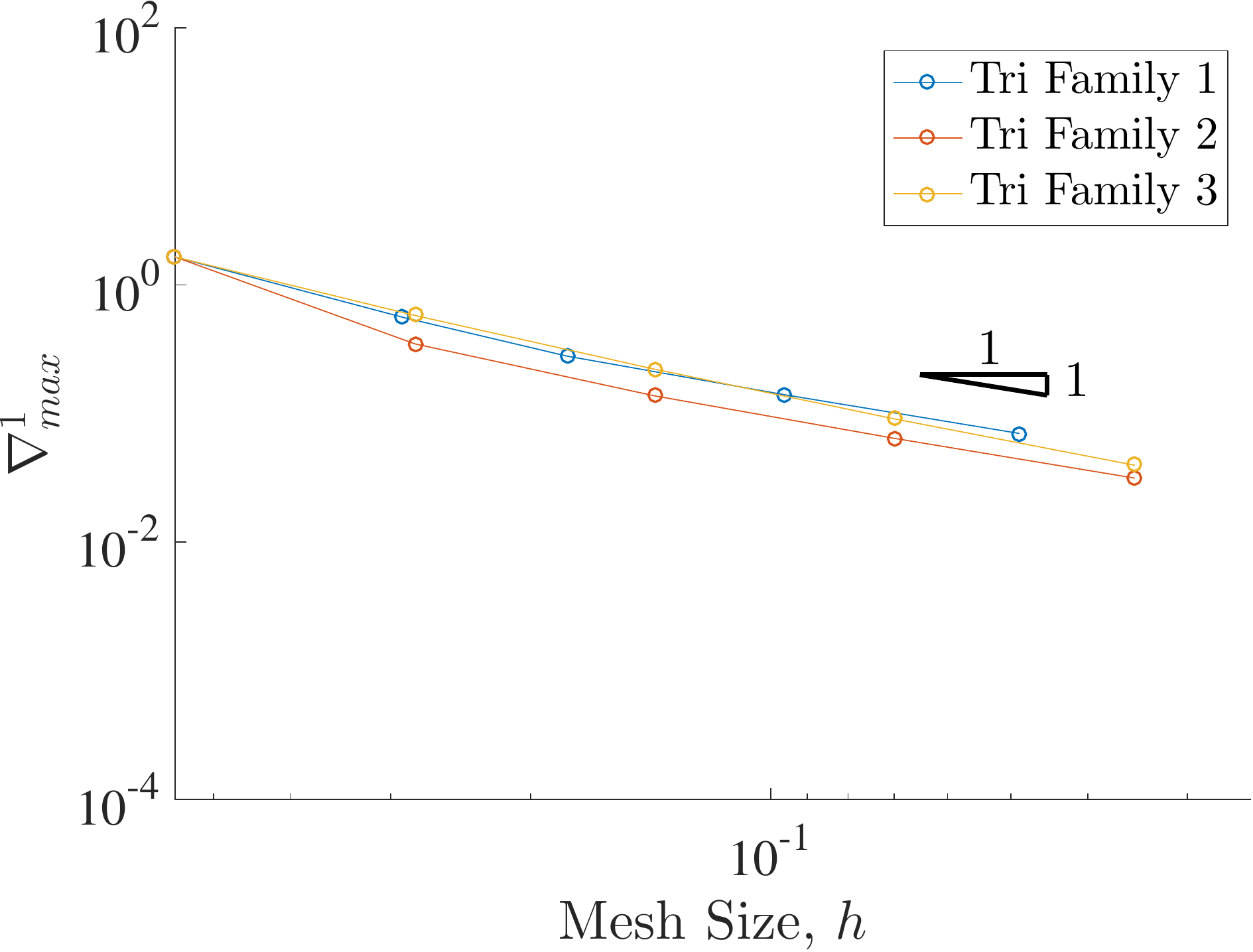}
        \caption{}
    \end{subfigure}
    \begin{subfigure}[t]{0.48\textwidth}
        \centering
            \includegraphics[width=.95\linewidth]{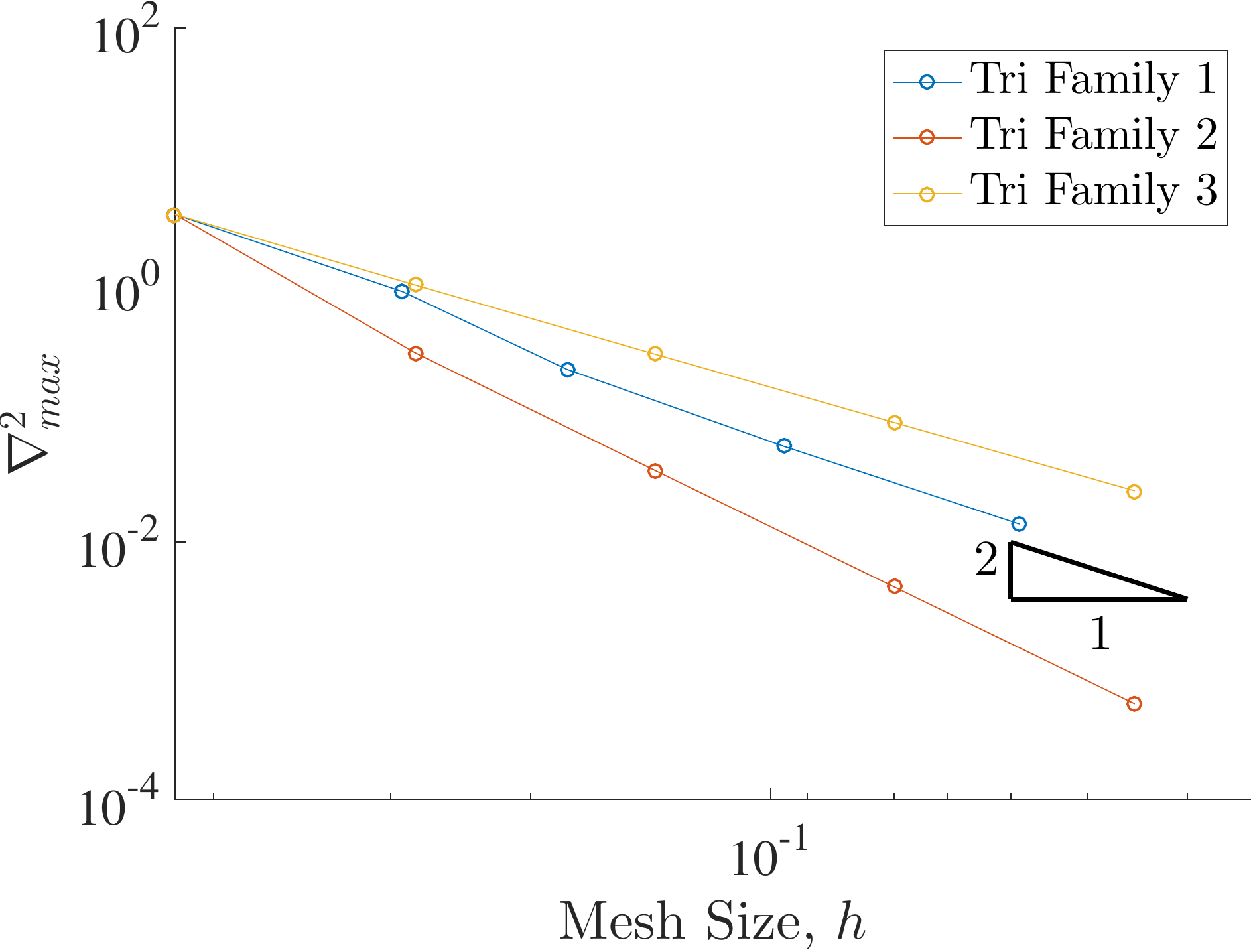}
        \caption{}
    \end{subfigure}
\begin{subfigure}[t]{0.48\textwidth}
        \centering
            \includegraphics[width=.95\linewidth]{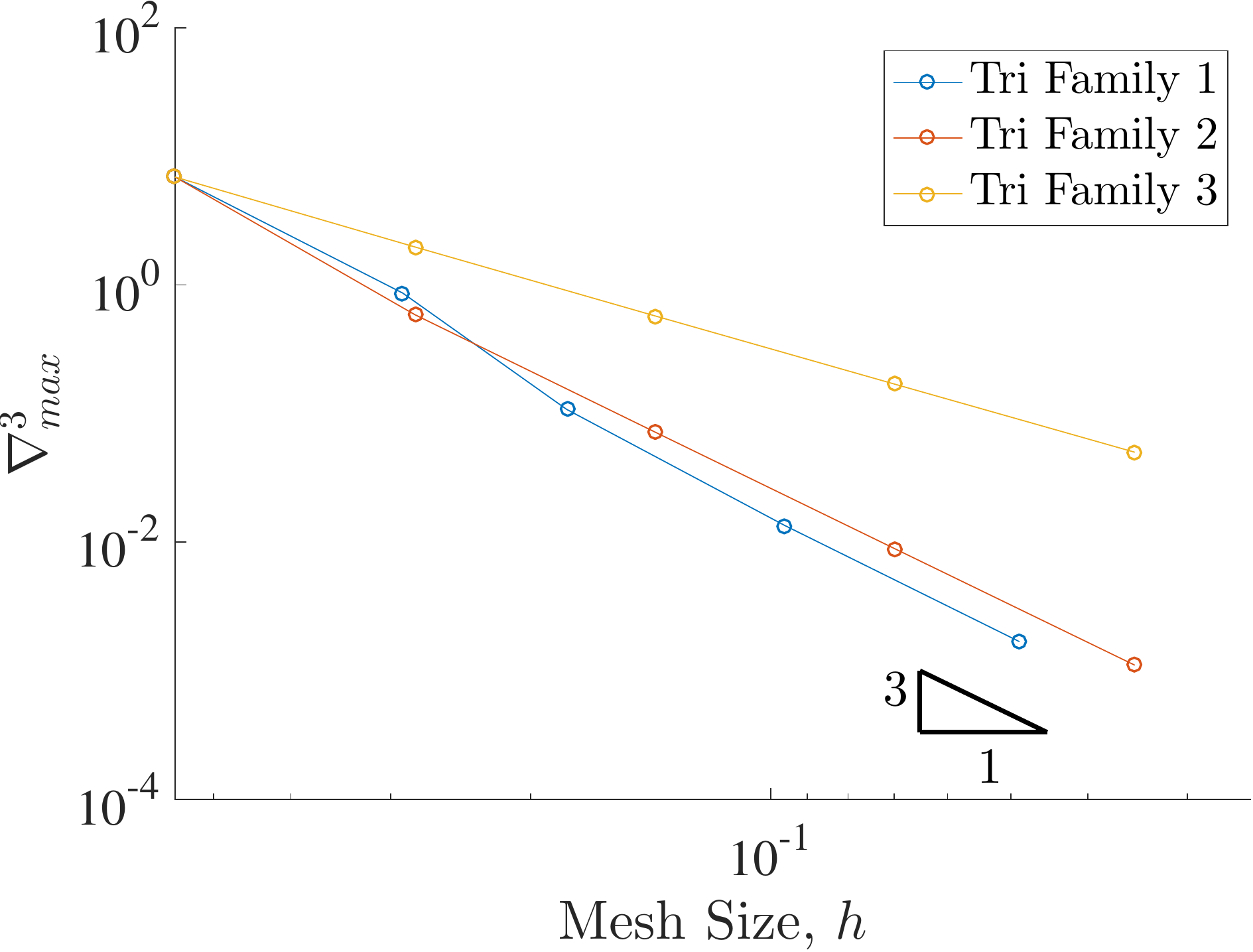}
        \caption{}
    \end{subfigure}
\caption{ Mesh distortion metrics for the triangular meshes of the plate. (a) Minimum scaled Jacobian. (b) Lowest upper bound on the magnitude of the first derivatives. (c) Lowest upper bound on the magnitude of the second derivatives. (d) Lowest upper bound on the magnitude of the third derivatives.}
\label{mms_plate_tri_metrics}
\end{figure}

From the convergence plots, we see that in all cases the solution error for problems solved over the  first and second family of meshes are converging as expected.
However, the error for the problems over the third family of meshes is converging at a less than optimal rate for both the quadrilateral and triangular meshes. 
To gain insight into this, we look to the distortion metrics for each mesh family. 
Fig. \ref{mms_plate_quad_metrics} shows distortion metrics for each family of quadrilateral meshes, and Fig. \ref{mms_plate_tri_metrics} shows distortion metrics for each family of triangular meshes.

From Fig. \ref{mms_plate_quad_metrics}a, we see that the minimum scaled Jacobian is bounded from below for every quadrilateral mesh family, and that in each case $J_S \rightarrow 1$ under mesh refinement. 
Furthermore, we note that $J_S$ is larger for the third family of meshes (irregular refinements) than it is for the first family (uniform refinements).
Similar behavior is also observed for the scaled Jacobian of the triangular meshes, shown in Fig. \ref{mms_plate_tri_metrics}a.
At first blush, these observations seems contradictory, as in both cases, the first family converges as expected, while the third does not. 

The cause of the slowed convergence rates can be explained by instead looking at the norms of the higher-order derivatives, shown in Fig. \ref{mms_plate_quad_metrics}b-d for the quadrilateral case, and Fig. \ref{mms_plate_tri_metrics}b-d for the triangular case.
For each plot, we show the lowest upper bound on the derivatives $D_{\bxi}^{\balpha}\xproj$ of order $\absof{\balpha} = k$ across the entire mesh. That is, for each $k \leq  p$, we plot the value of $\bnabla^{k}_{max}$, where:
\begin{equation}
\bnabla^{k}_{max} = \max_{e=1,...,nel}\max_{\absof{\balpha}=k}\sup_{\bxi \in \Oref} \absof{ D_{\bxi}^{\balpha}\xproj}
\end{equation}
From Fig. \ref{mms_plate_quad_metrics}d, it is readily seen that the cause of the slowed convergence for the third family of quadrilateral meshes is the fact that $||\bnabla^3 \xphys ||_ {L^\infty(\Oref)}$ is converging at approximately $O(h^2)$, whereas Cond. (I.2) requires that it converge at $O(h^3)$. 

These results serve to highlight how sensitive convergence rates for higher-order elements can be, and to motivate utility of the distortion metrics developed in this work. 
Indeed, from the Jacobian metrics shown in Fig. \ref{mms_plate_quad_metrics}a and Fig. \ref{mms_plate_tri_metrics}a, as well as visual inspection of the meshes, one might be tempted to draw the conclusion that all three families should preserve optimal convergence rates, even though we have observed that this is clearly not the case. 

\subsection{Manufactured Solution on a Plate with a Hole}
\begin{table}[b!]
  \begin{center}
    \caption{Meshes and families of weighting functions for the plate with a hole manufactured solution.} 
    \begin{tabular}{|  >{\centering\arraybackslash} m{.3cm} |  >{\centering\arraybackslash} m{4cm} |  >{\centering\arraybackslash} m{4cm} |    
    >{\centering\arraybackslash} m{4cm} | }
      \hline
       $m$ & Meshes & Weight Family 1 & Weight Family 2  \\ 
      \hline
       1 &
      \vspace{5pt}
      \includegraphics[width=4cm,height=4cm,keepaspectratio]{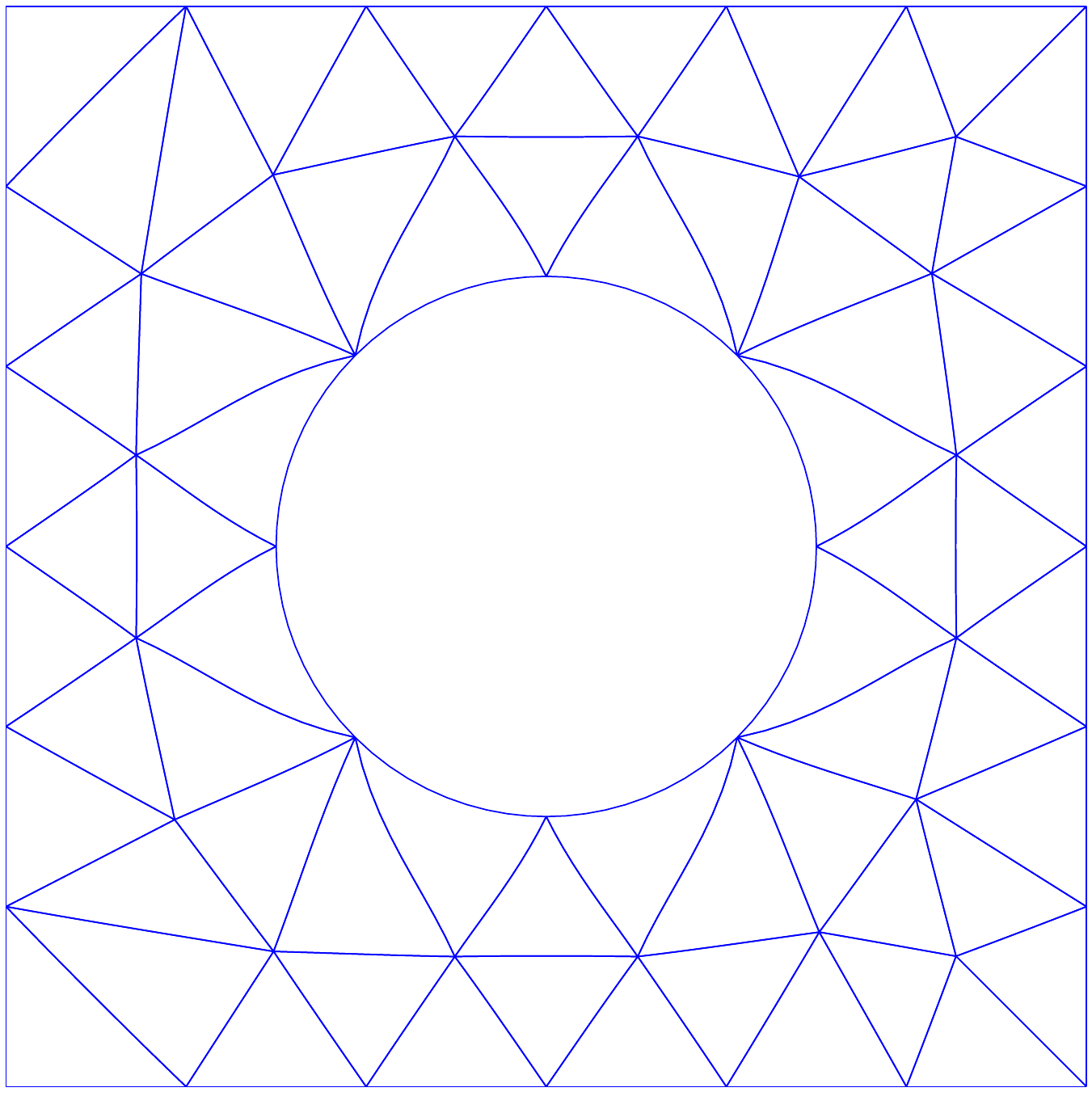} &
      \vspace{5pt}
      \includegraphics[width=4cm,height=4cm,keepaspectratio]{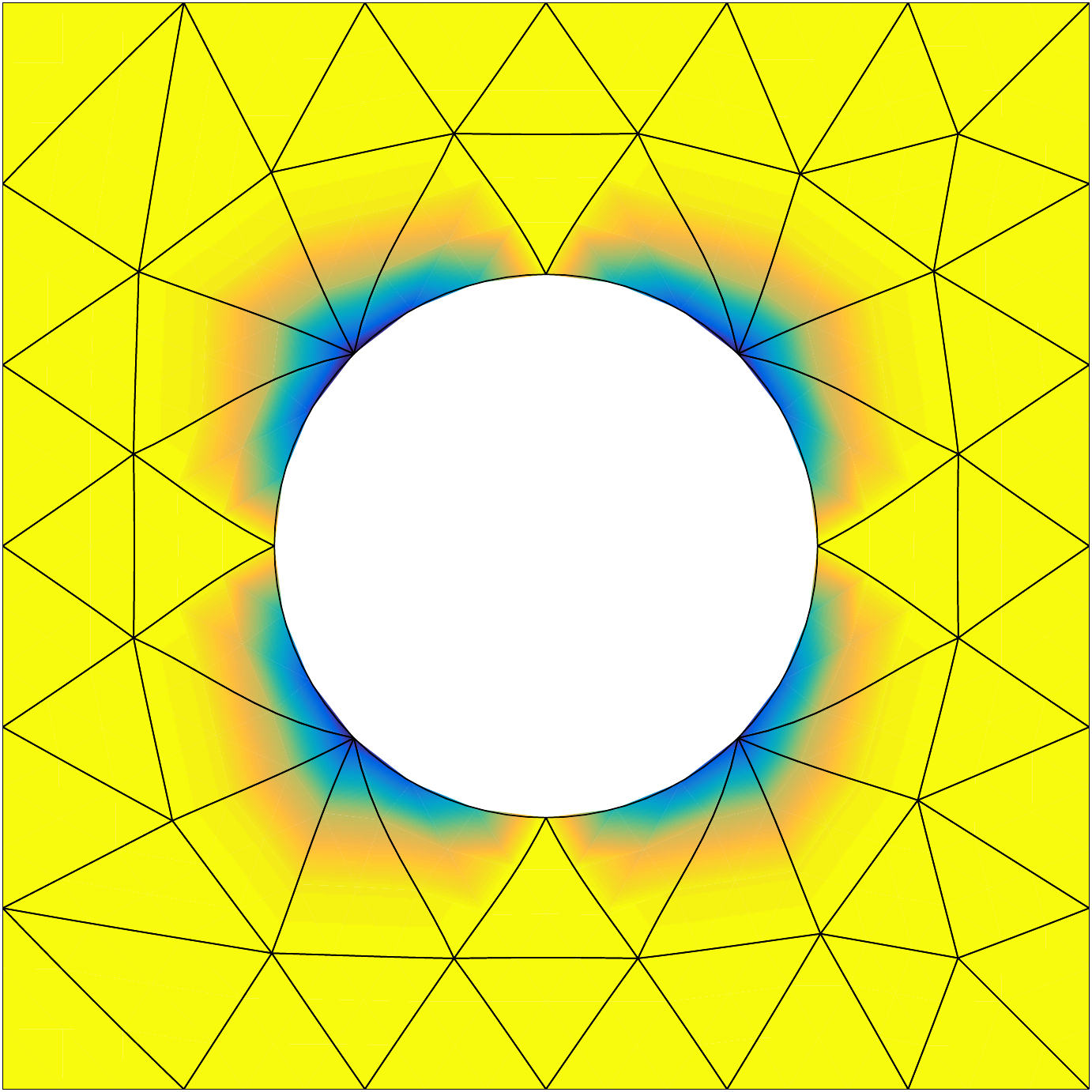} &
       \vspace{5pt}
      \includegraphics[width=4cm,height=4cm,keepaspectratio]{mms_nosmoothW_weights.pdf}\\
      \hline
             2 &
      \vspace{5pt}
      \includegraphics[width=4cm,height=4cm,keepaspectratio]{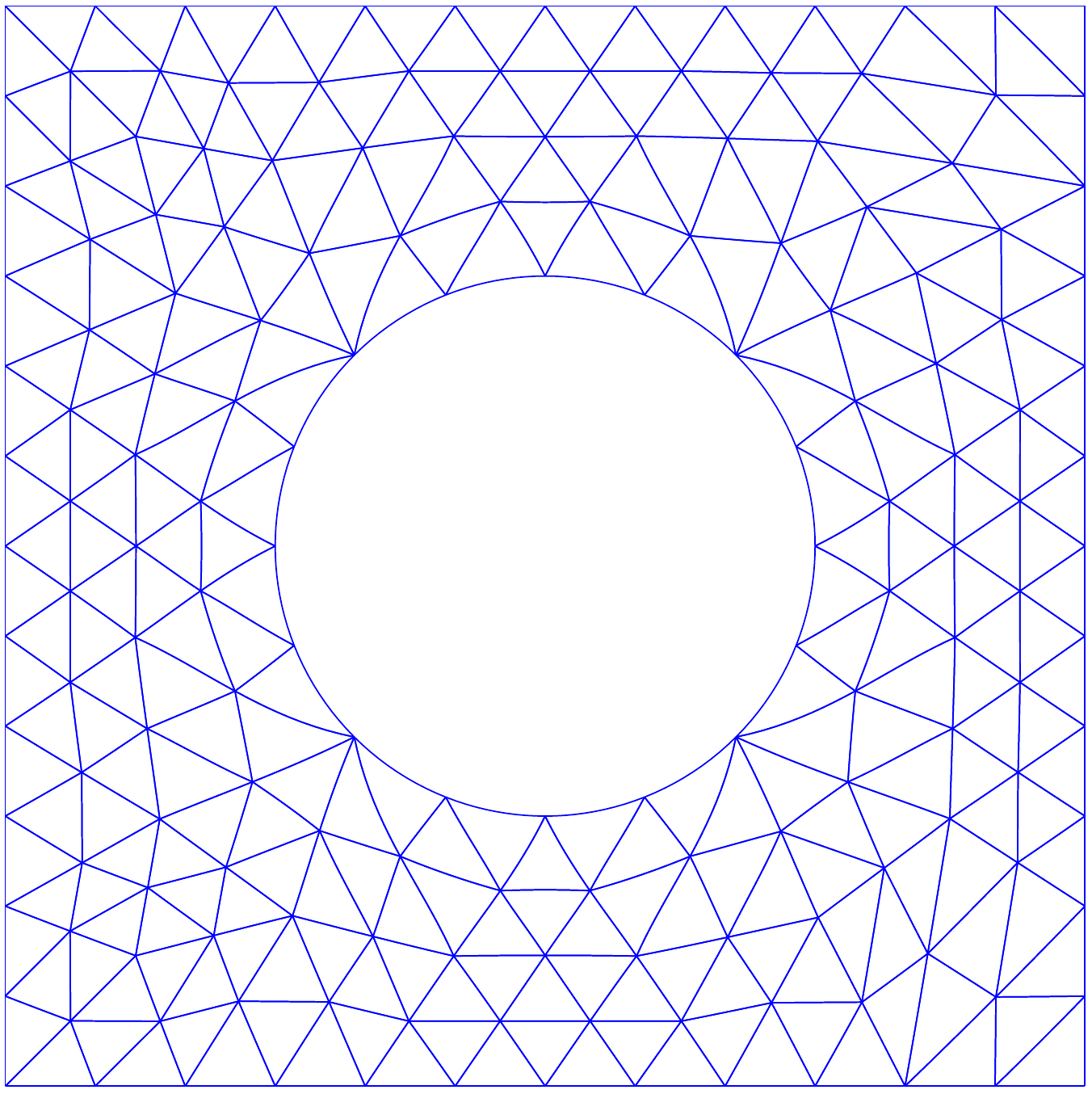} &
      \vspace{5pt}
      \includegraphics[width=4cm,height=4cm,keepaspectratio]{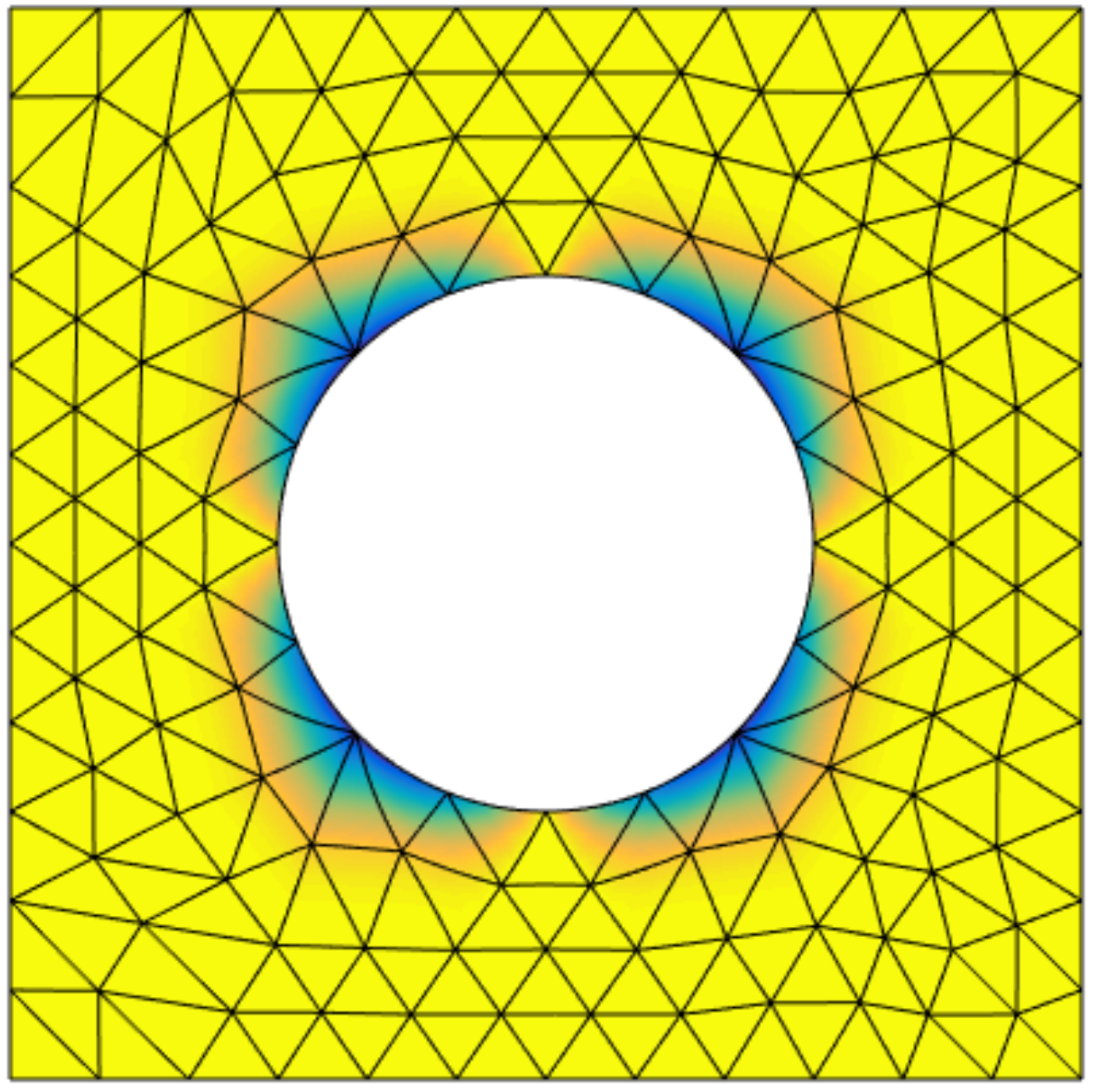} &
       \vspace{5pt}
      \includegraphics[width=4cm,height=4cm,keepaspectratio]{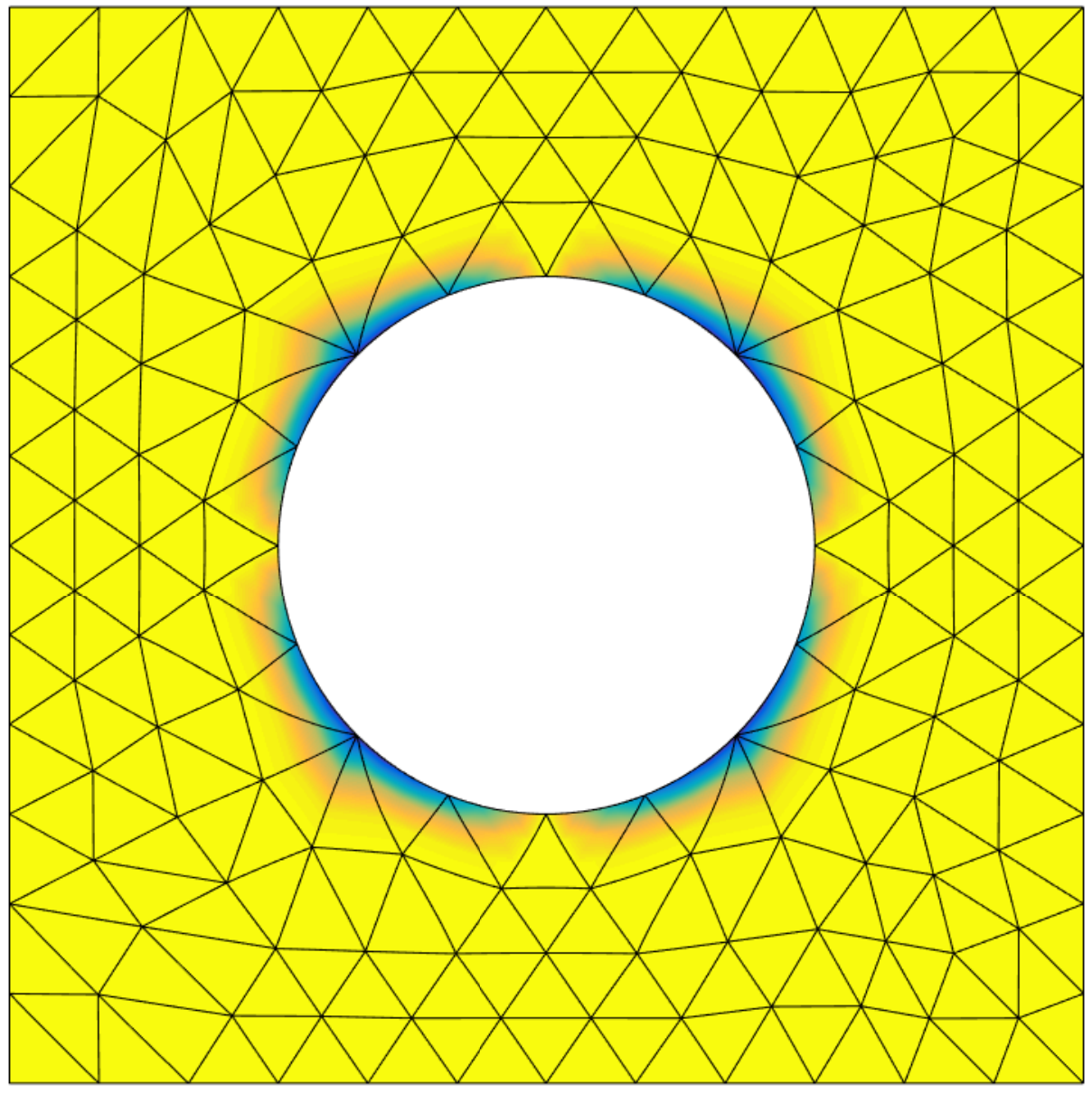}\\
      \hline
    3 &
      \vspace{5pt}
      \includegraphics[width=4cm,height=4cm,keepaspectratio]{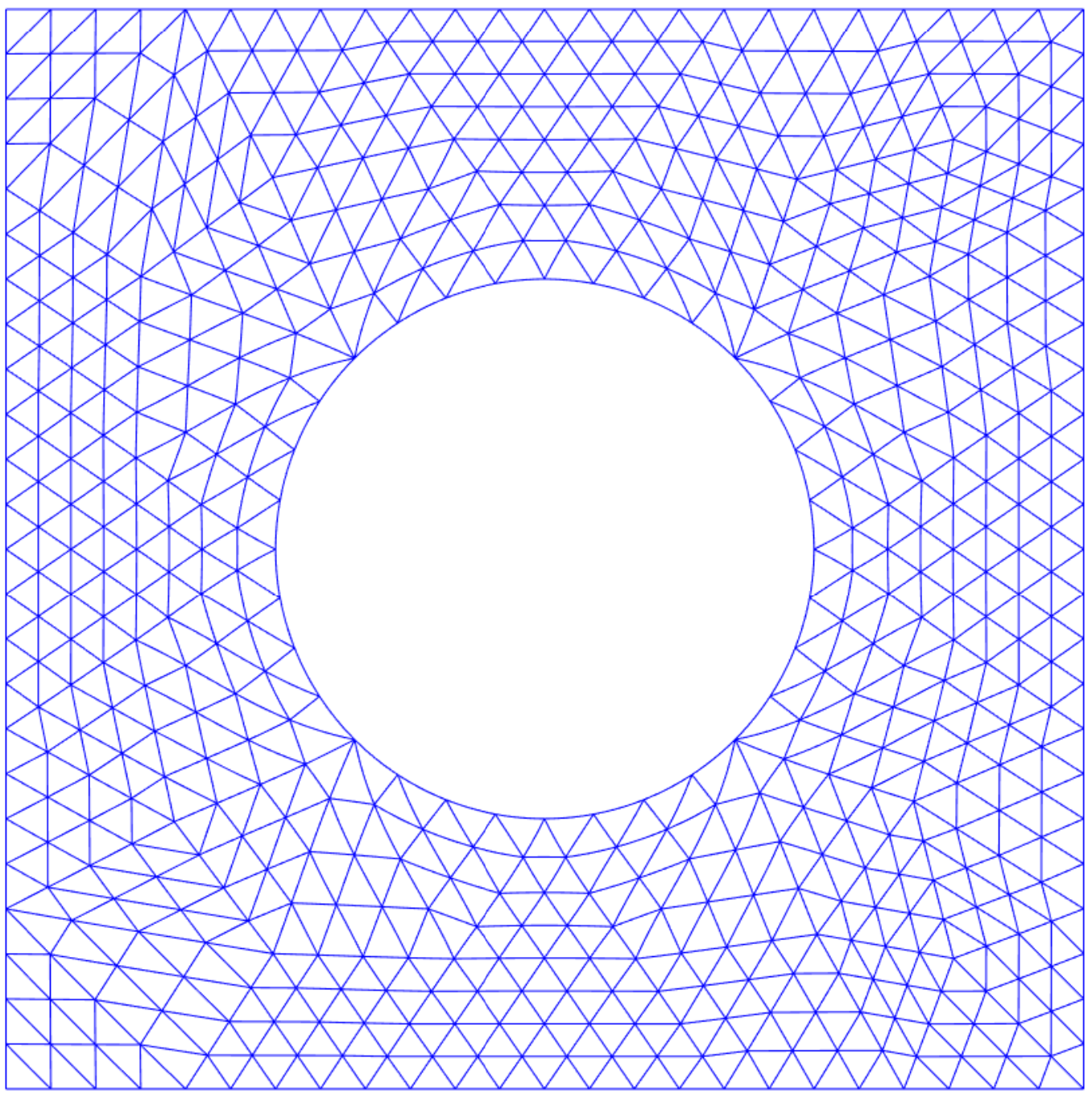} &
      \vspace{5pt}
      \includegraphics[width=4cm,height=4cm,keepaspectratio]{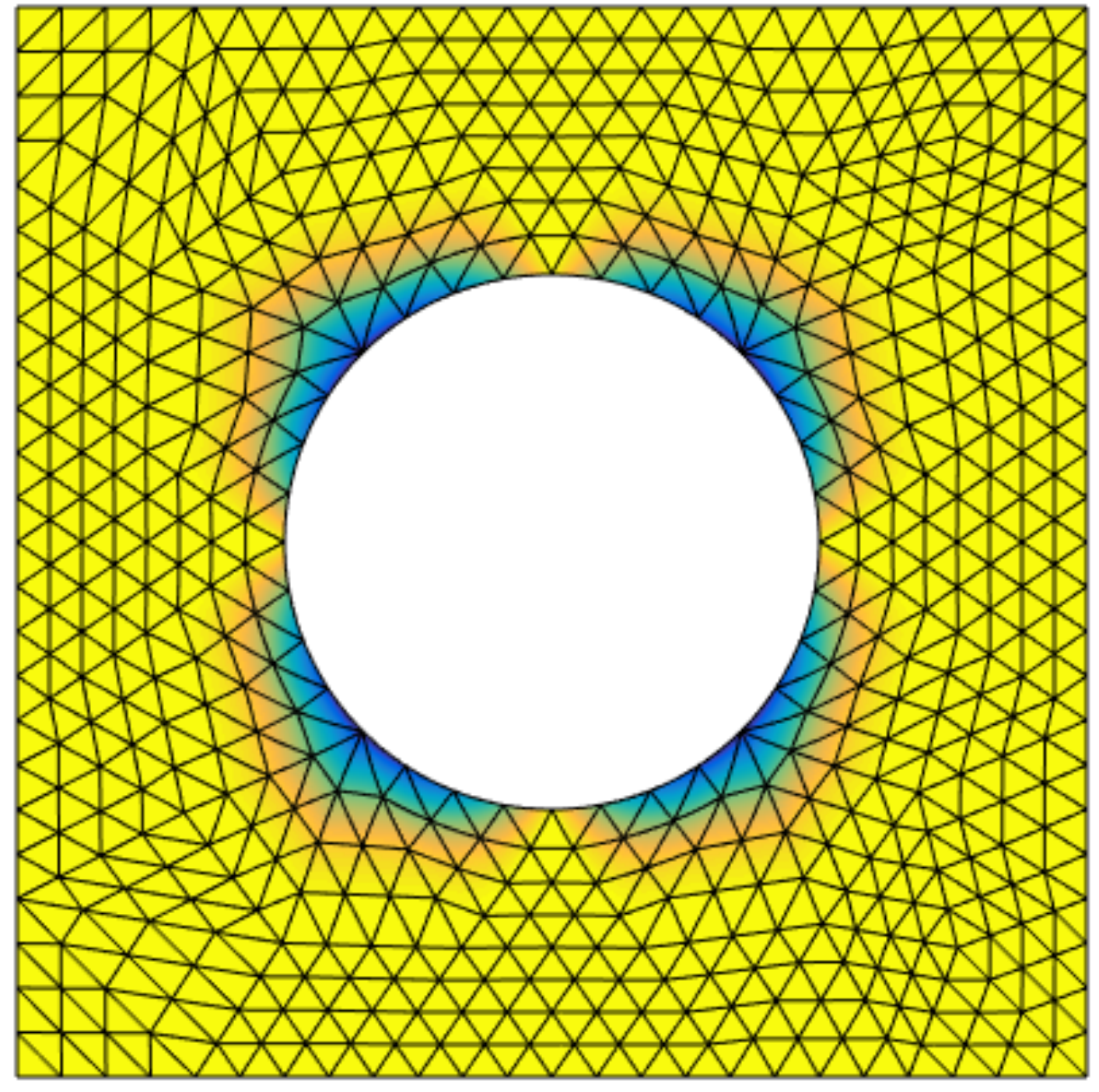} &
       \vspace{5pt}
      \includegraphics[width=4cm,height=4cm,keepaspectratio]{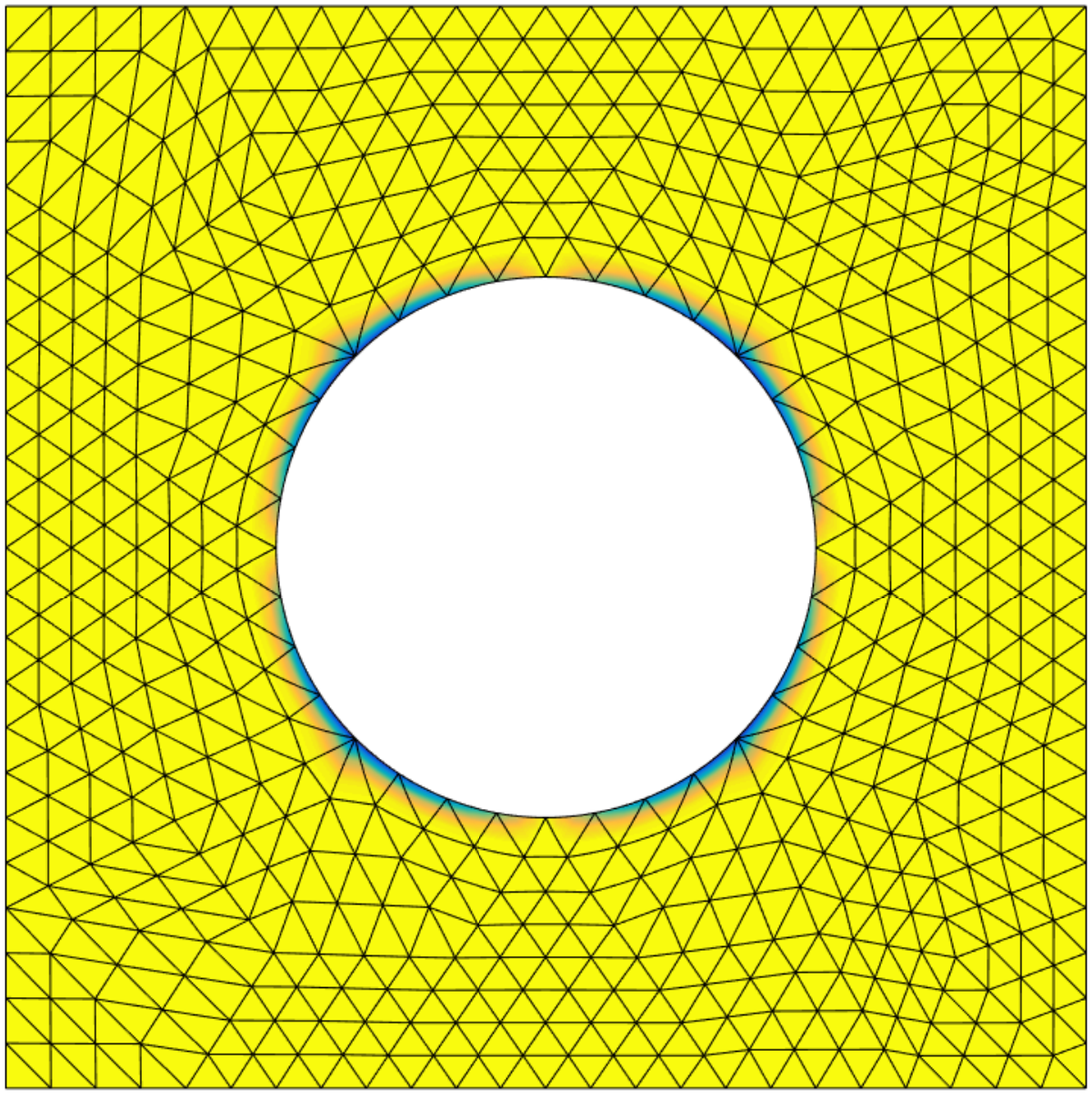}\\
      \hline
      \end{tabular}
    \label{weight_families}
  \end{center}
\end{table}

In the previous example, we examined the effect of element shape distortion on approximation error by perturbing control points. 
In this example, we examine the effect of weighting function distortion on the approximation error by perturbing control weights.
We consider a mesh of a plate with a hole, composed of cubic rational \BB triangles, shown in Table \ref{weight_families}.
We note that  since the hole in the plate is circular, we must use rational elements to capture the geometry exactly, and as a result, the control weights corresponding to the points on the boundary will be non unity.  
However, it remains to set the control weights for the interior points in the mesh. 
We consider two possible methods of setting control weights for a series of meshes. 
The first method is to simply set the weights on the circular boundary, to the appropriate values, and set all other weights to one. 
Then, under mesh refinement, we perform uniform subdivision on both the control points, and control weights 
As a result, the weighting function remains the same under mesh refinement. 
The second option we consider here is to perform uniform subdivision on the control points, but not the control weights. 
Instead, only control weights corresponding to points on the boundary are updated during each refinement step, and all other control weights are set to one.
Both of these families of refined weighting functions are shown in Table \ref{weight_families}.

\begin{figure}[t] 
\centering
        \centering
            \includegraphics[width=.5\linewidth]{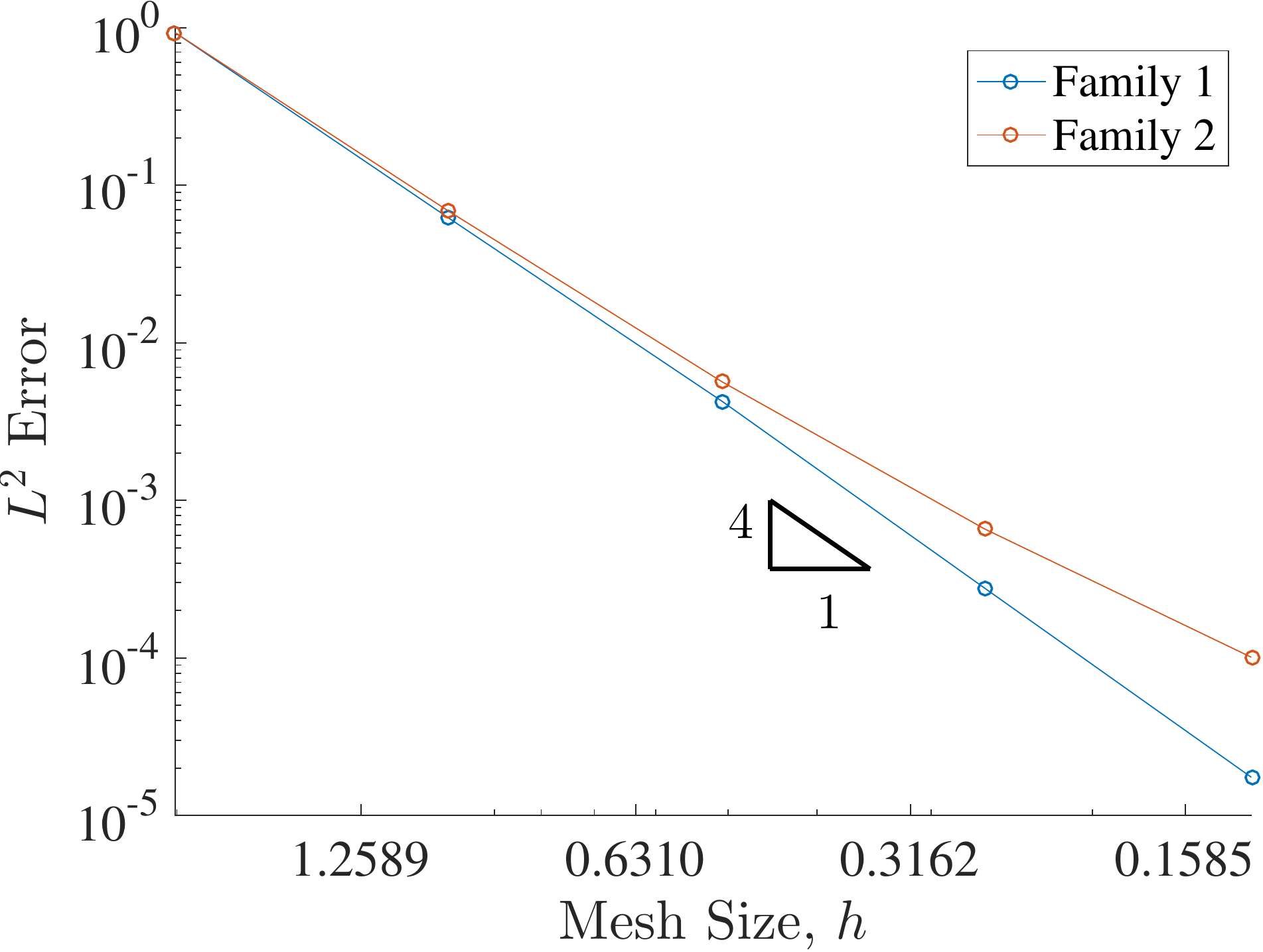}
\caption{Convergence plots for the plate with a hole manufactured solution.}
\label{plate_hole_norm}
\end{figure}

With these two families established, we study their approximation accuracy using the method of manufactured solutions. 
As before, we solve Poisson's problem, with the manufactured solution:
\begin{equation}
	u\of{x_1,x_2} = \of{x_1-a}\of{x_1+a}\of{x_2-a}\of{x_2+a}\of{r-\sqrt{x_1^2+x_2^2}}	
\end{equation}
wherein $a$ is the half-width of the square plate and $r$ is the radius of the hole, and the plate is centered at the origin.
Figure \ref{plate_hole_norm} shows the convergence rates of the $L^2$ norm of the approximation error over both mesh families

%
%
%
%
%
\begin{figure}[t!]
\centering
\begin{subfigure}[t]{0.49\textwidth}
        \centering
            \includegraphics[width=.95\linewidth]{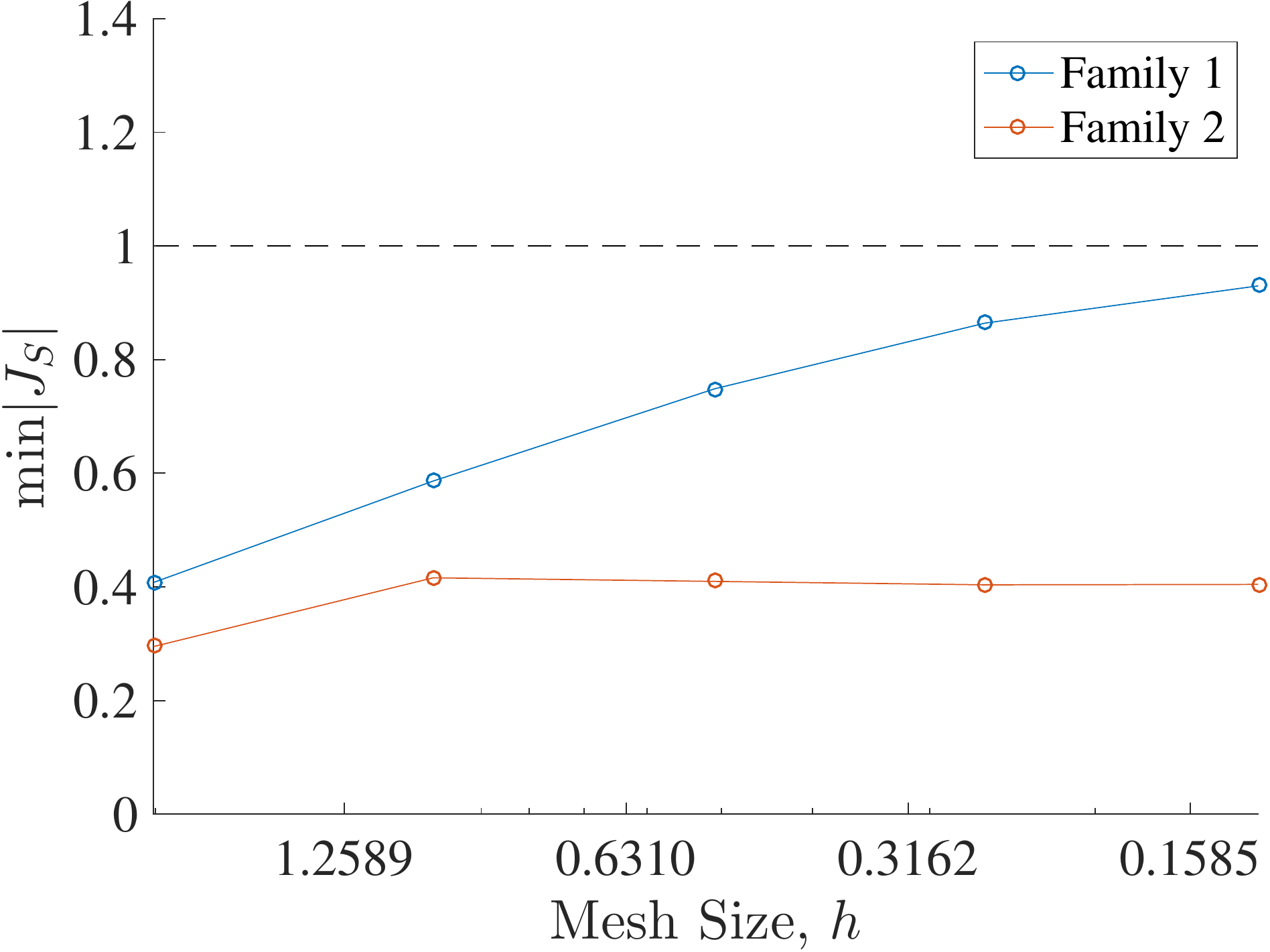}
        \caption{}
    \end{subfigure}
 \begin{subfigure}[t]{0.49\textwidth}
        \centering
            \includegraphics[width=.95\linewidth]{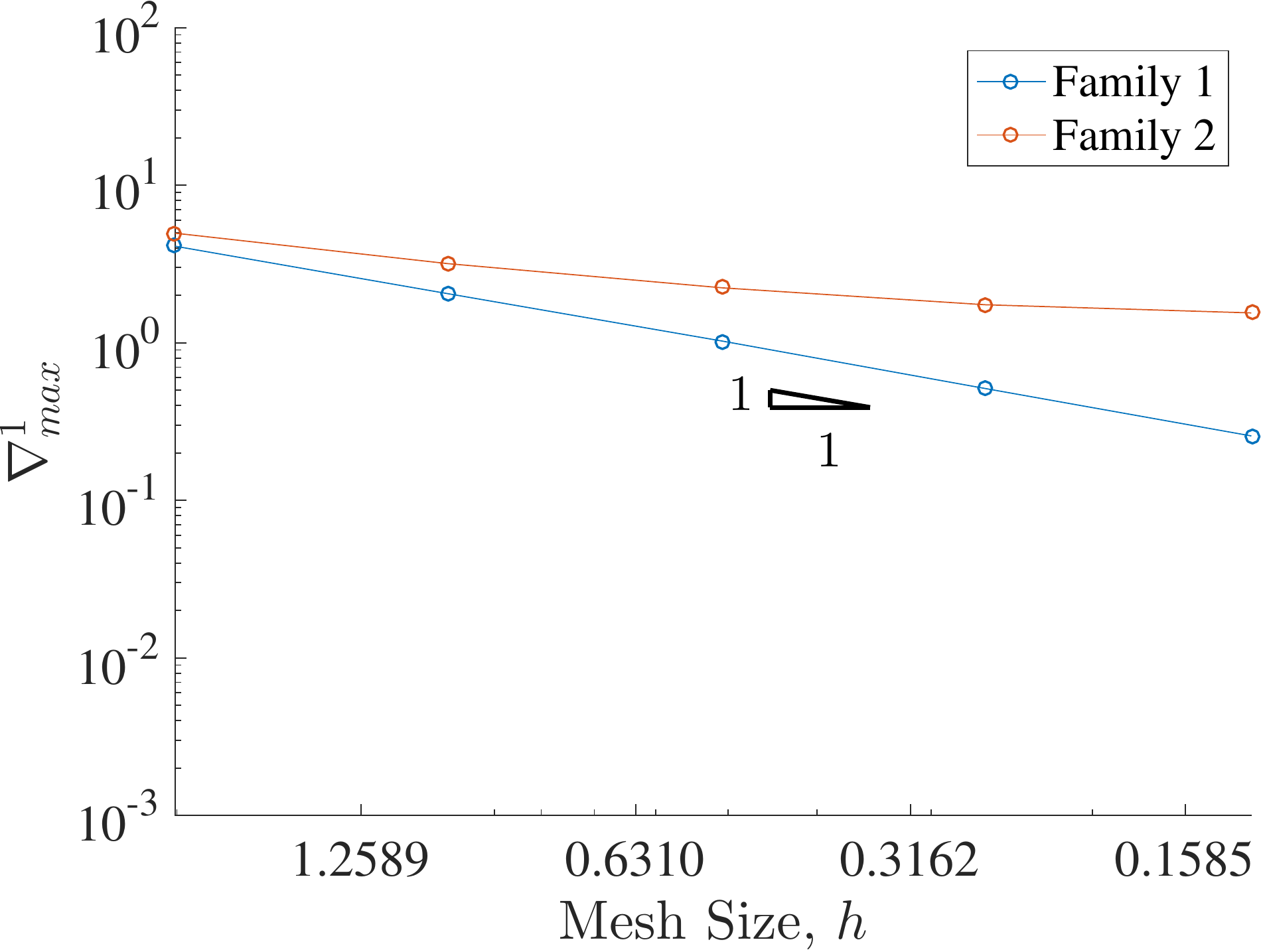}
        \caption{}
    \end{subfigure}
    \begin{subfigure}[t]{0.49\textwidth}
        \centering
            \includegraphics[width=.95\linewidth]{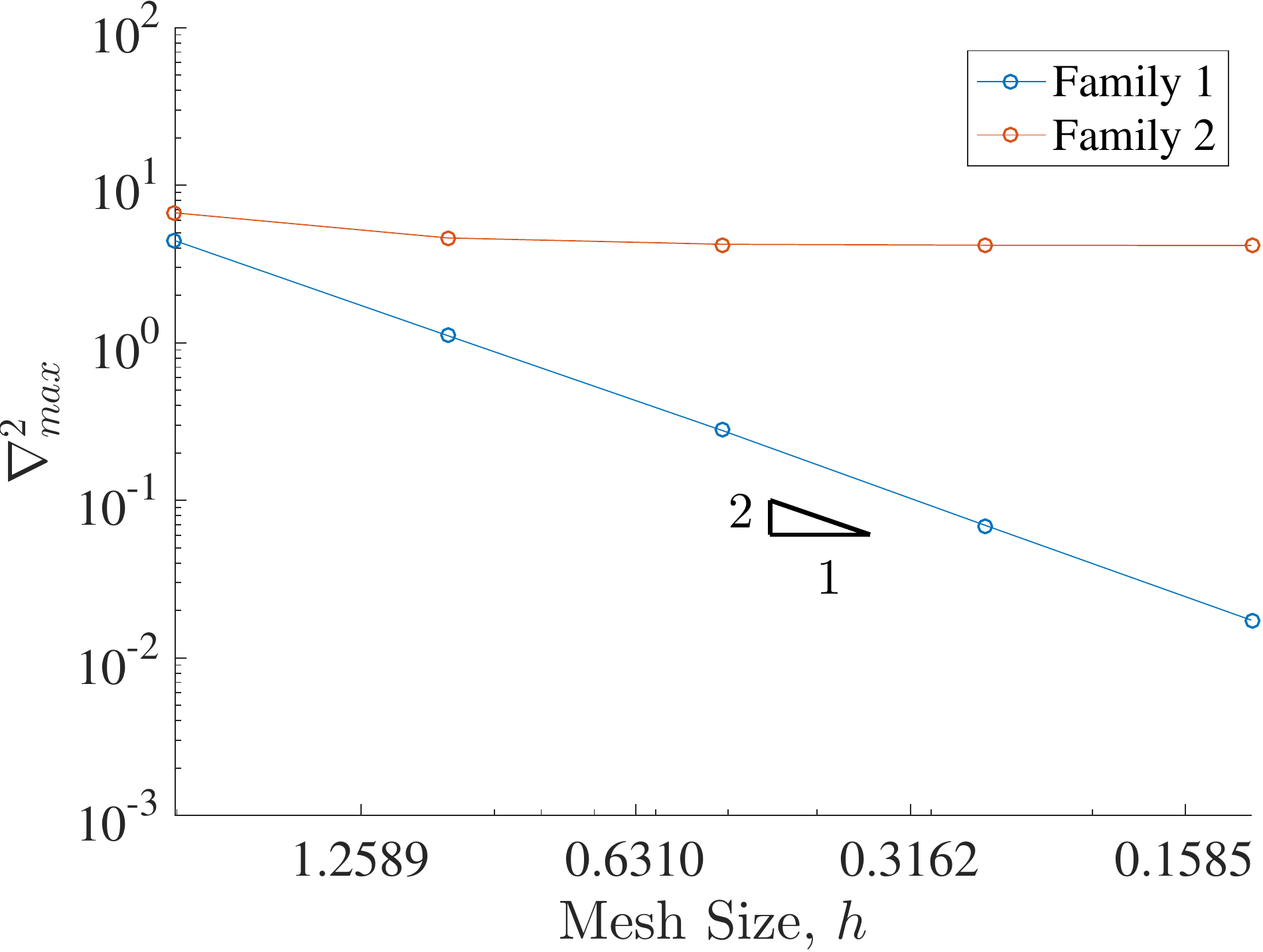}
        \caption{}
    \end{subfigure}
\begin{subfigure}[t]{0.49\textwidth}
        \centering
            \includegraphics[width=.95\linewidth]{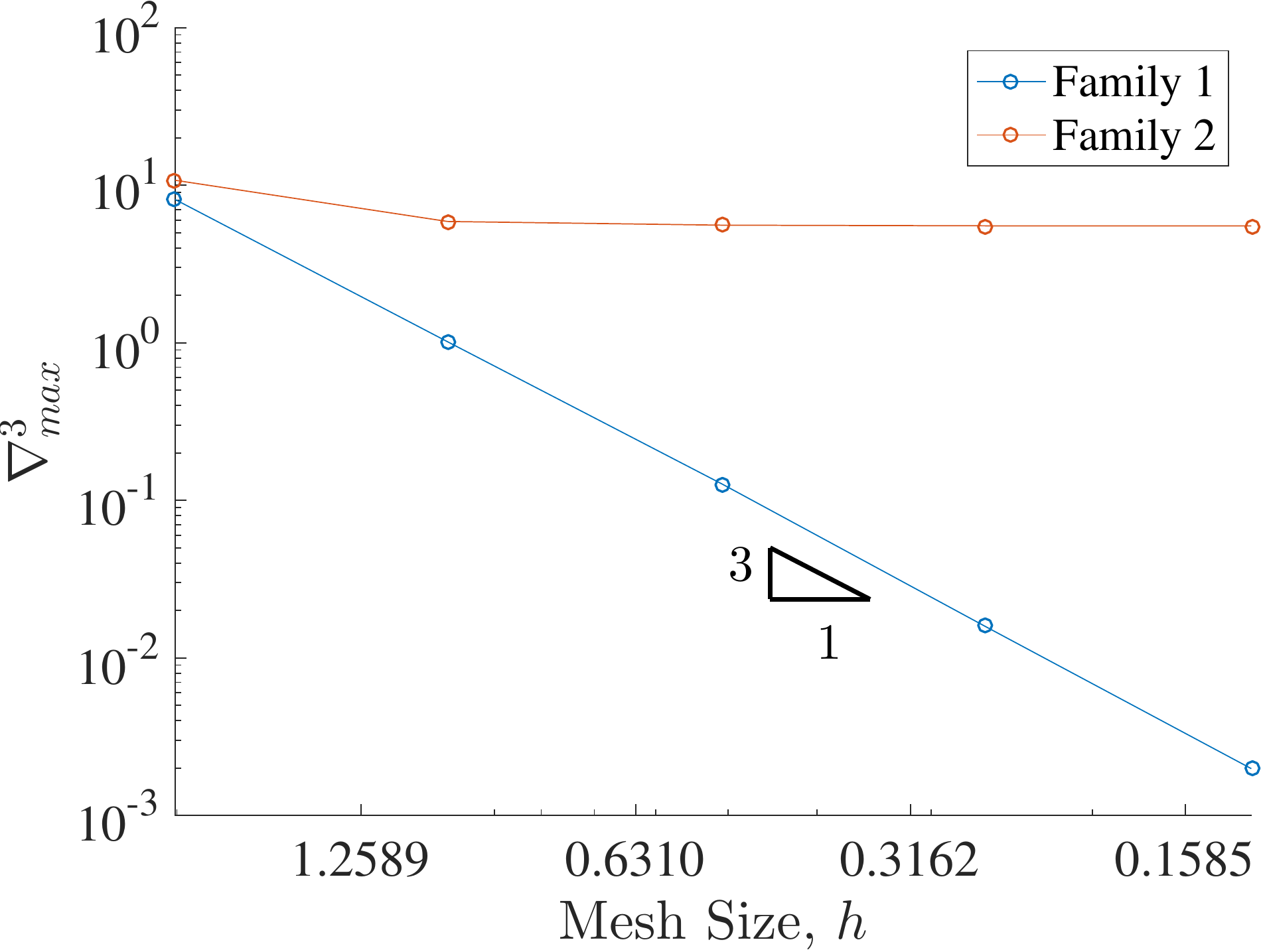}
        \caption{}
    \end{subfigure}
\caption{ Mesh distortion metrics for the triangular meshes of the plate. (a) Minimum scaled Jacobian. (b) Lowest upper bound on the magnitude of the first derivatives. (c) Lowest upper bound on the magnitude of the second derivatives. (d) Lowest upper bound on the magnitude of the third derivatives.}
\label{plate_hole_derivs}
\end{figure}

From our results, we see that the first family of meshes converges as expected, while the second does not.
Again, we examine the distortion metrics, shown in Fig. \ref{plate_hole_derivs}, to gain insight into the cause of the stalled convergence for the second family of meshes. 
We immediately see that the higher-order derivatives of the element-wise parametric mapping are not decaying to zero for the second family of meshes. 
The reason for this can be observed from the plots of the weighting functions in Table \ref{weight_families}. 
Since only the weights lying on the circular boundary are being updated, the gradient of the weighting function becomes increasingly sharp under refinement. 

\subsection{Convergence Under $p$-refinement}
Thus far, we have considered the effect of both control point distortion and control weight distortion on approximation error for cubic B\'{e}zier elements. 
We now consider the effect of mesh parameterization on convergence under $p$-refinement. 
We consider the simple case of a quarter annulus mesh with four rational bi-quadratic \BB quadrilaterals.
\begin{table}[t!]
  \begin{center}
    \caption{Families of $p$-refined meshes of the quarter annulus.} 
    \begin{tabular}{|  >{\centering\arraybackslash} m{.3cm} |  >{\centering\arraybackslash} m{3.75cm} |  >{\centering\arraybackslash} m{3.75cm} |   }
      \hline
       $p$ & Family 1 & Family 2  \\ 
      \hline
       2 &
      \vspace{5pt}
      \includegraphics[width=3.75cm,height=3.75cm,keepaspectratio]{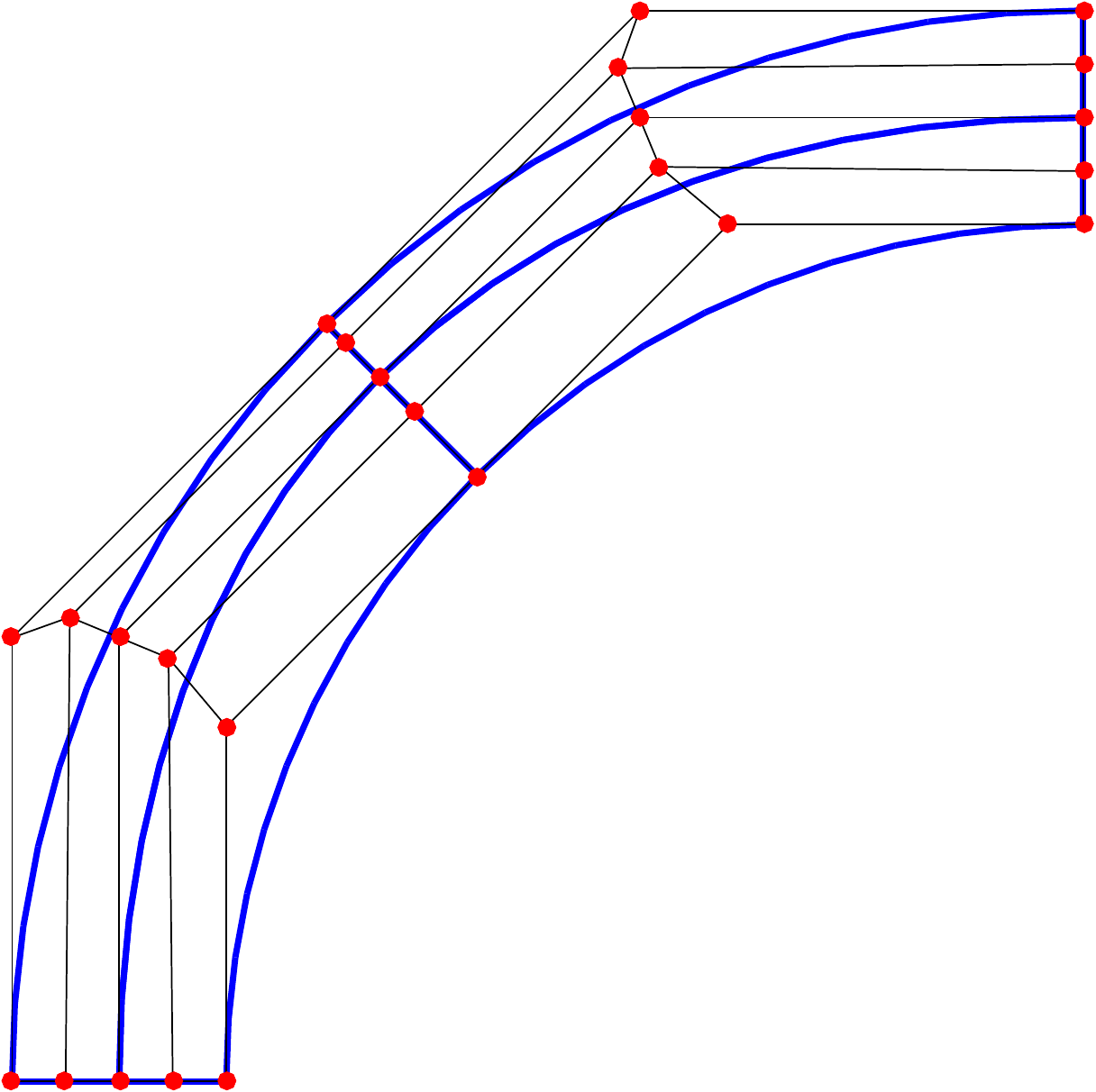} &
       \vspace{5pt}
      \includegraphics[width=3.75cm,height=3.75cm,keepaspectratio]{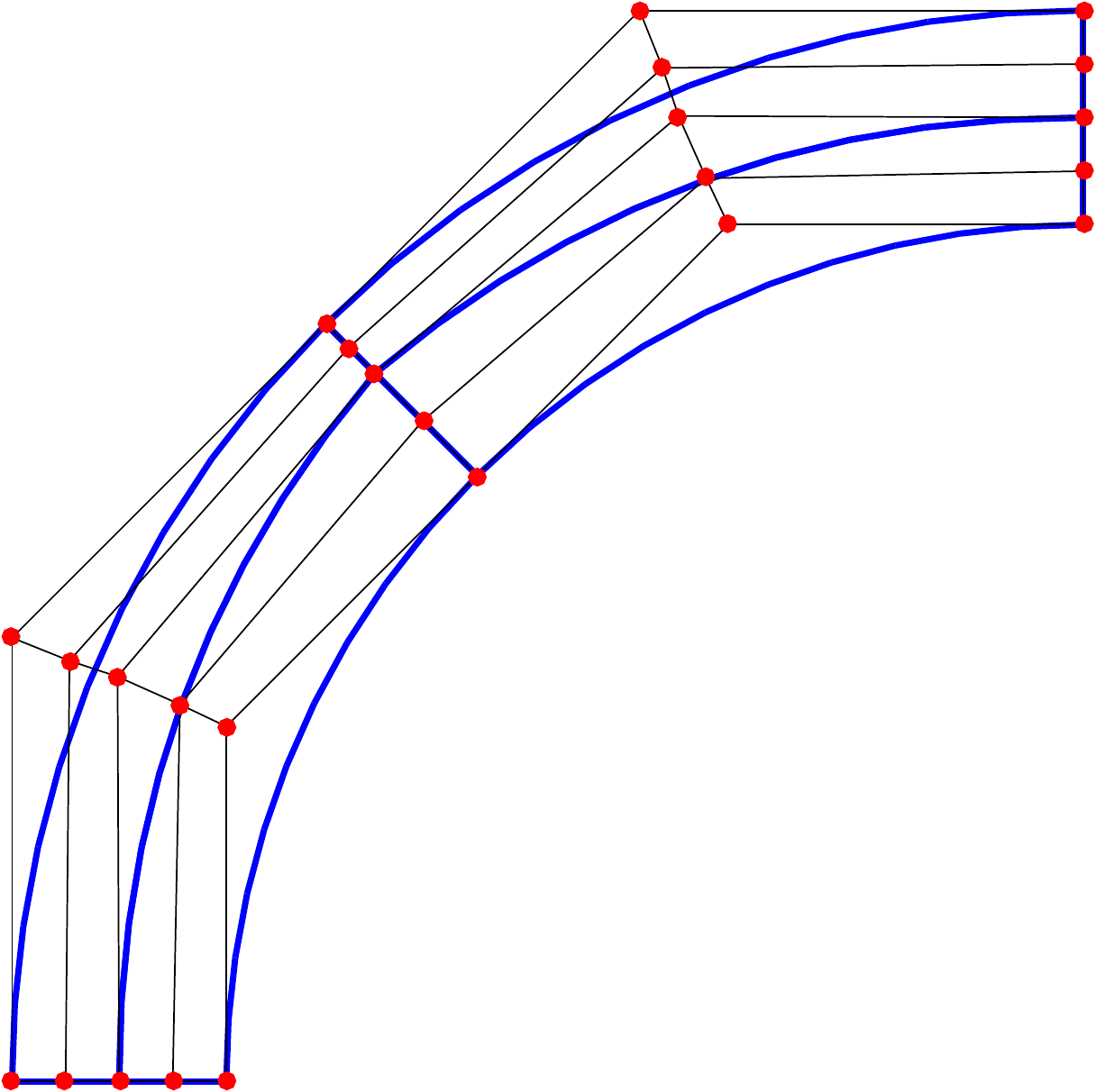}\\
      \hline
             3 &
      \vspace{5pt}
      \includegraphics[width=3.75cm,height=3.75cm,keepaspectratio]{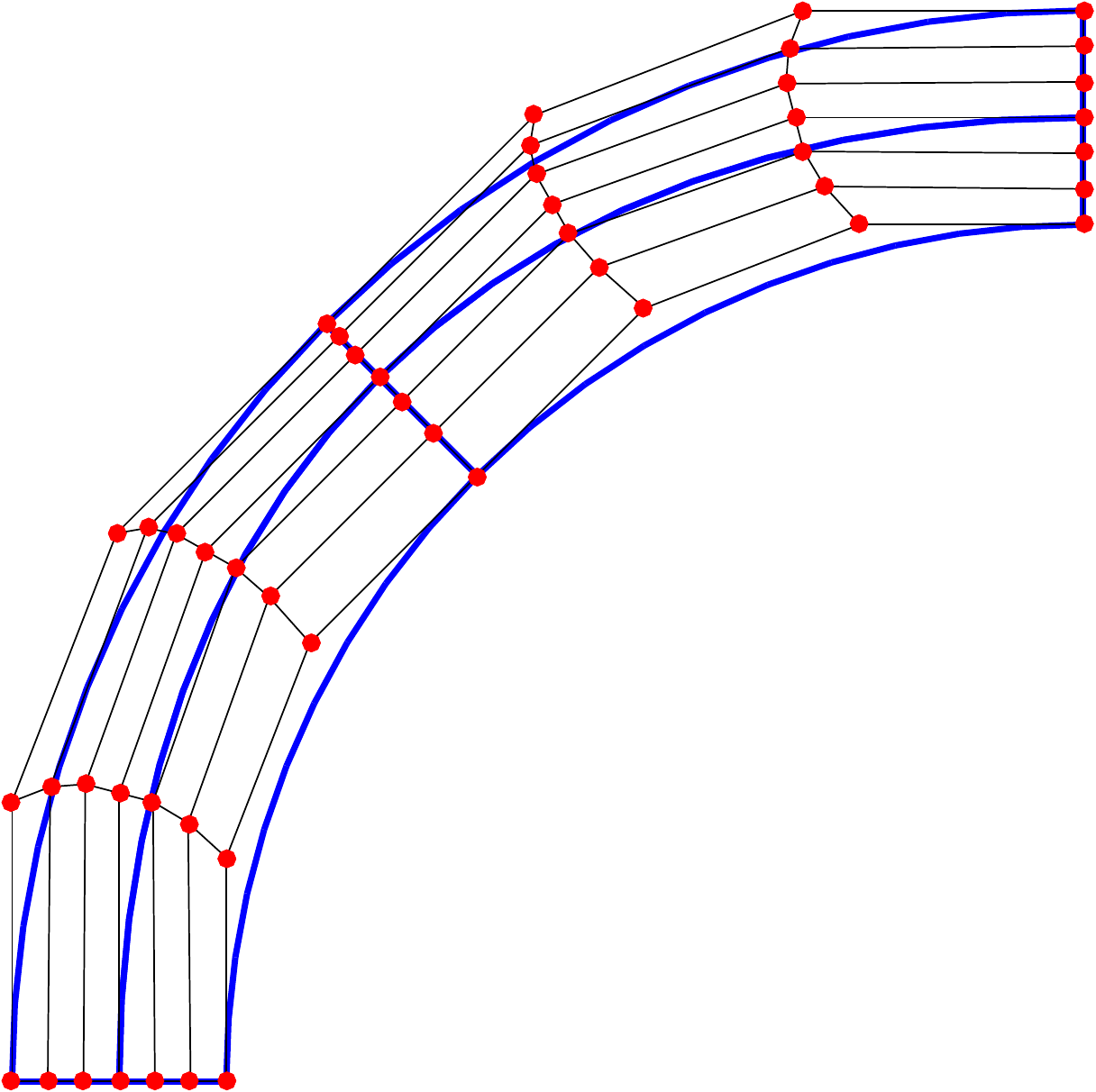} &
       \vspace{5pt}
      \includegraphics[width=3.75cm,height=3.75cm,keepaspectratio]{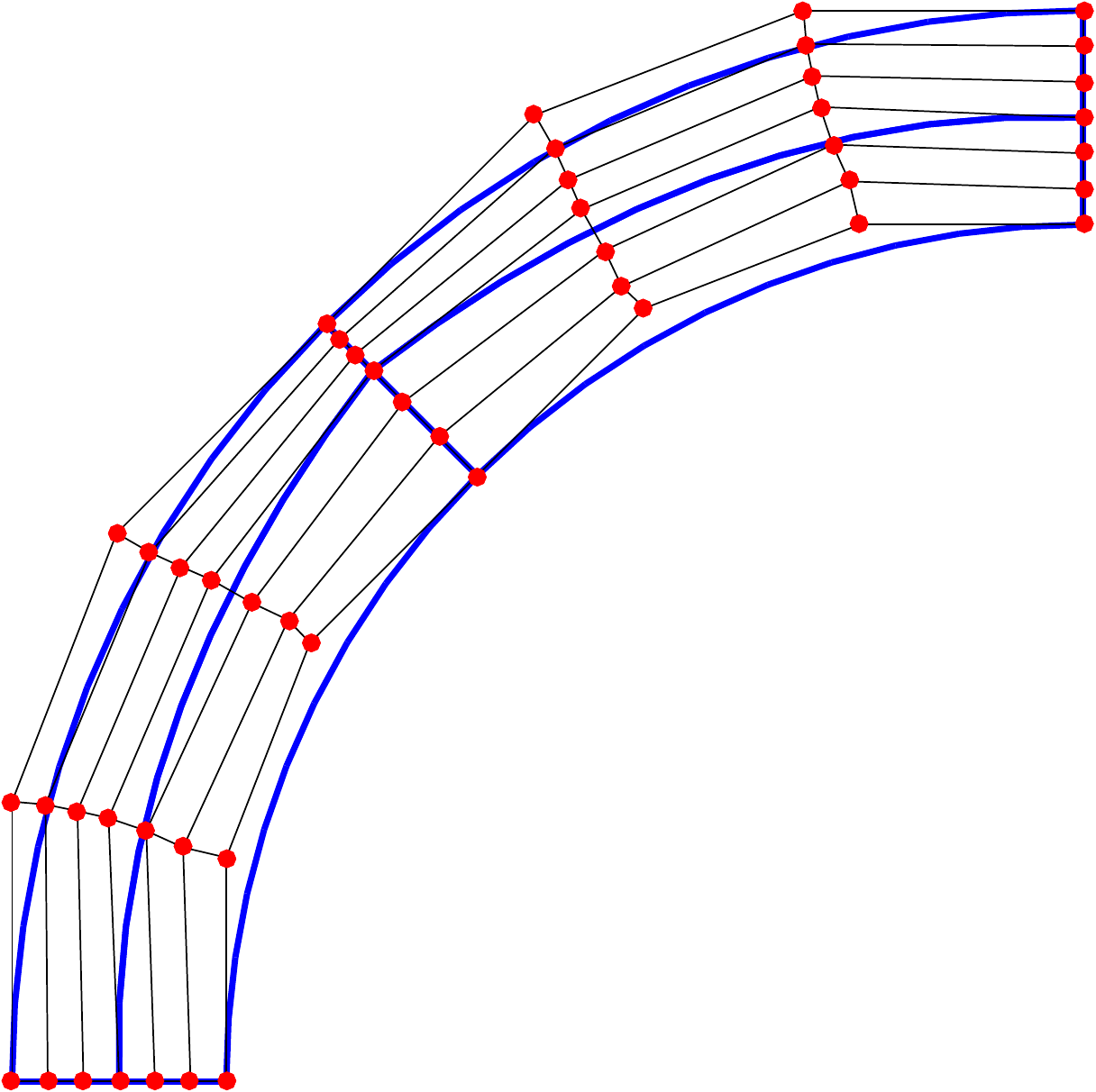}\\
      \hline
    4 &
      \vspace{5pt}
      \includegraphics[width=3.75cm,height=3.75cm,keepaspectratio]{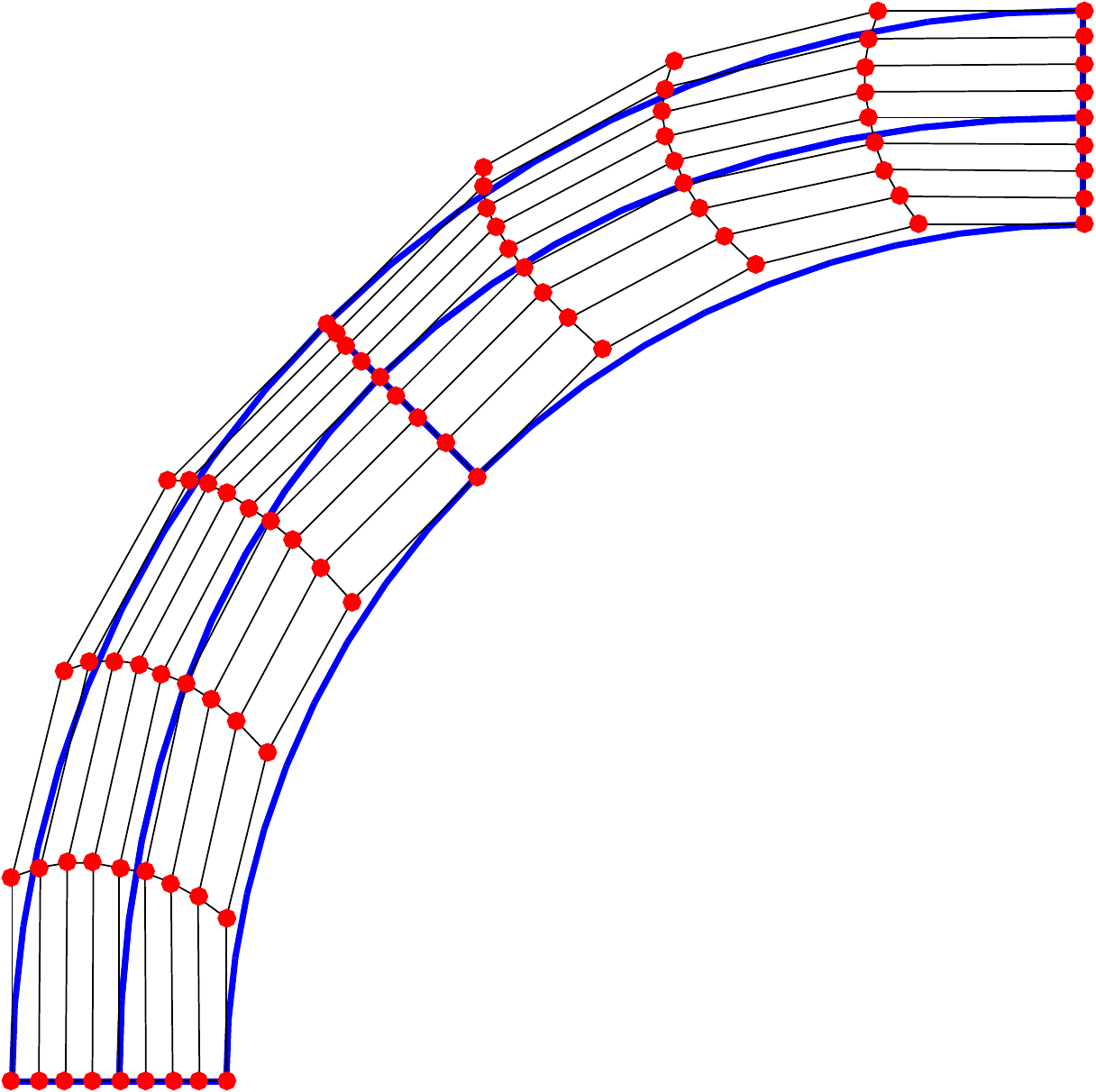} &
       \vspace{5pt}
      \includegraphics[width=3.75cm,height=3.75cm,keepaspectratio]{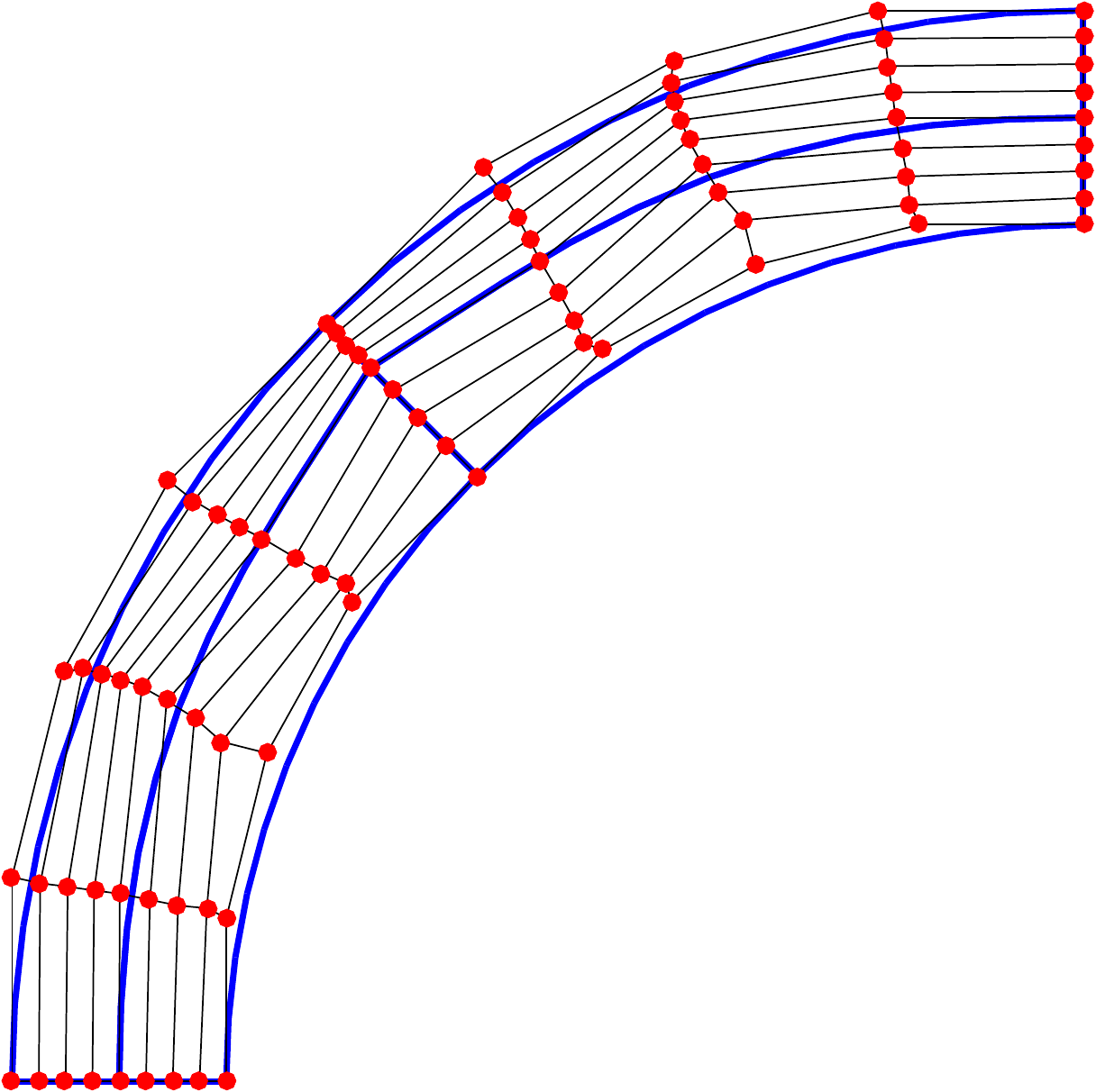}\\
      \hline
      \end{tabular}
    \label{pref_families}
  \end{center}
\end{table}

\begin{figure}[t!] 
\centering
\includegraphics[width=.5\textwidth]{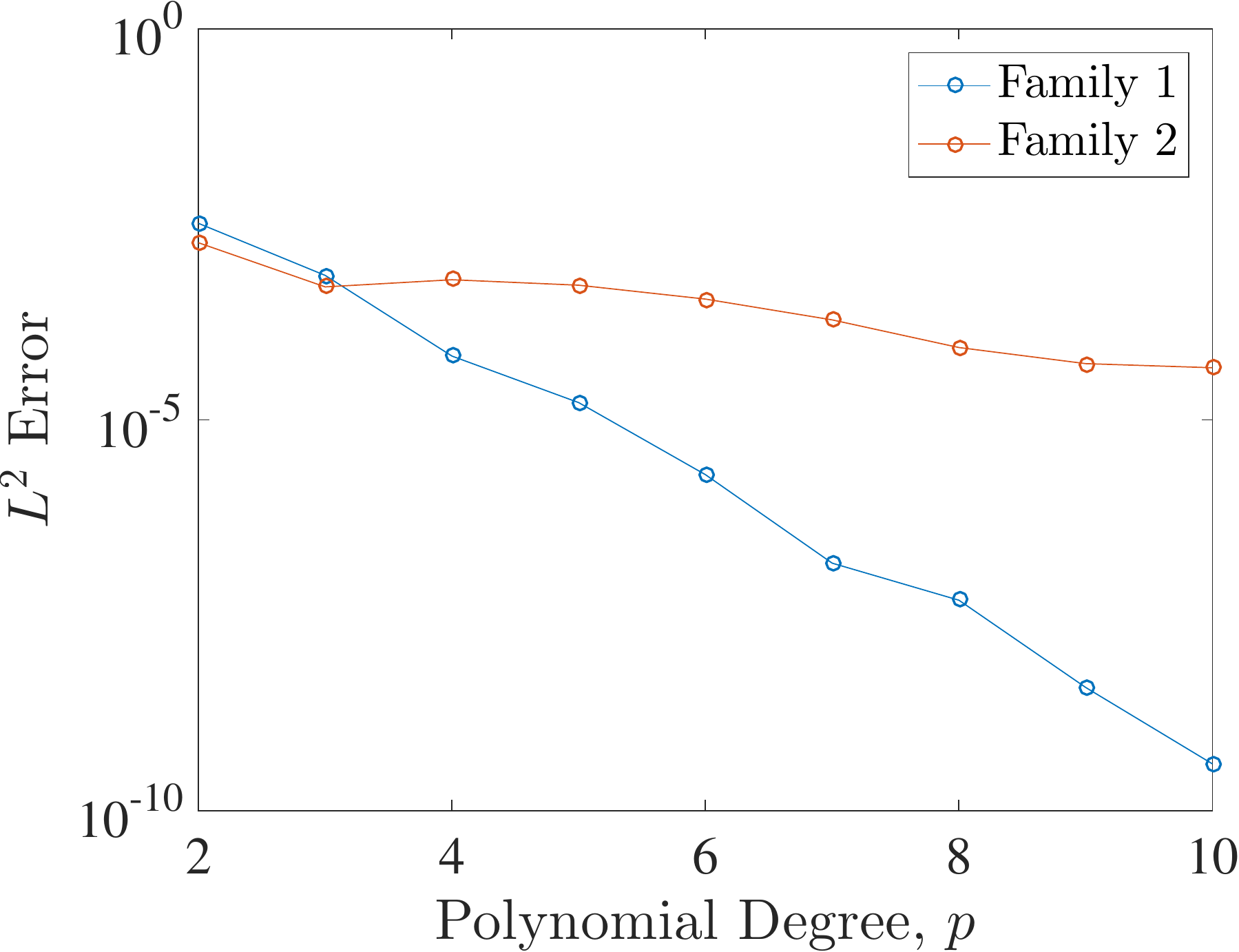}
\caption{Convergence in the $L^2$ error under $p$-refinement.}
\label{mms_pref_L2}
\end{figure}
\begin{figure}[t!] 
\centering
\includegraphics[width=.5\textwidth]{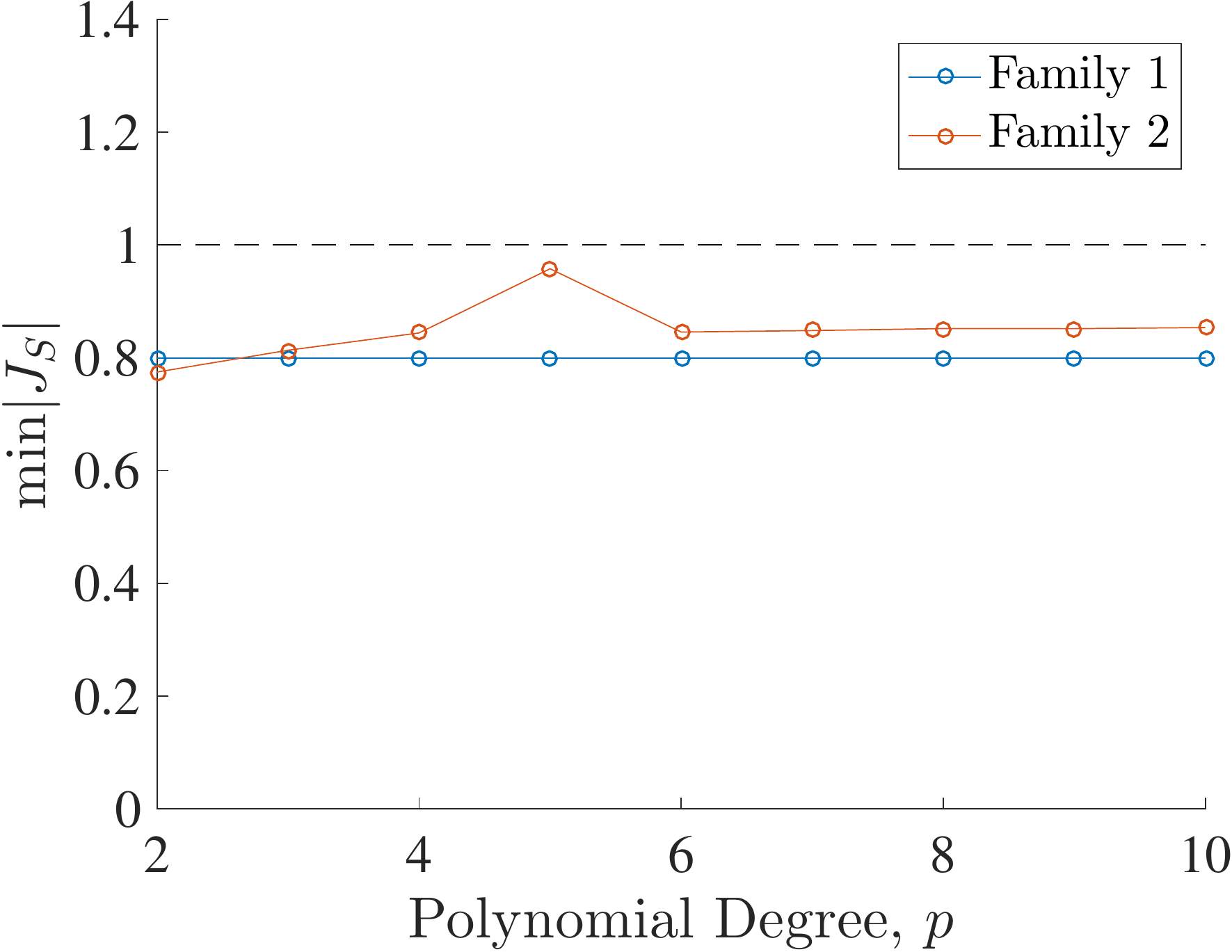}
\caption{Scaled Jacobian metrics for the $p$-refined meshes of the quarter annulus.}
\label{mms_pref_Jacobian}
\end{figure}

We consider two series of $p$-refined meshes. The first series is created by simple degree elevation of the bi-quadratic mesh. 
The second series is created using a linear elastic analogy, how higher-order meshes are typically created in practice \cite{field1988laplacian,persson_curved_2008,poya_unified_2016}. 
For each level of refinement, we degree elevate the underlying linear mesh to order $p$, and the  geometry is recovered via edge replacement \cite{engvall_isogeometric_2016}.
Then, we solve a linear elasticity problem to update the positions of the interior nodes. 
The two series of $p$-refined meshes with the B\'{e}zier control nets are shown in Table \ref{pref_families}.

\begin{figure}[t!]
\centering
\begin{subfigure}[t]{0.49\textwidth}
        \centering
            \includegraphics[width=.95\linewidth]{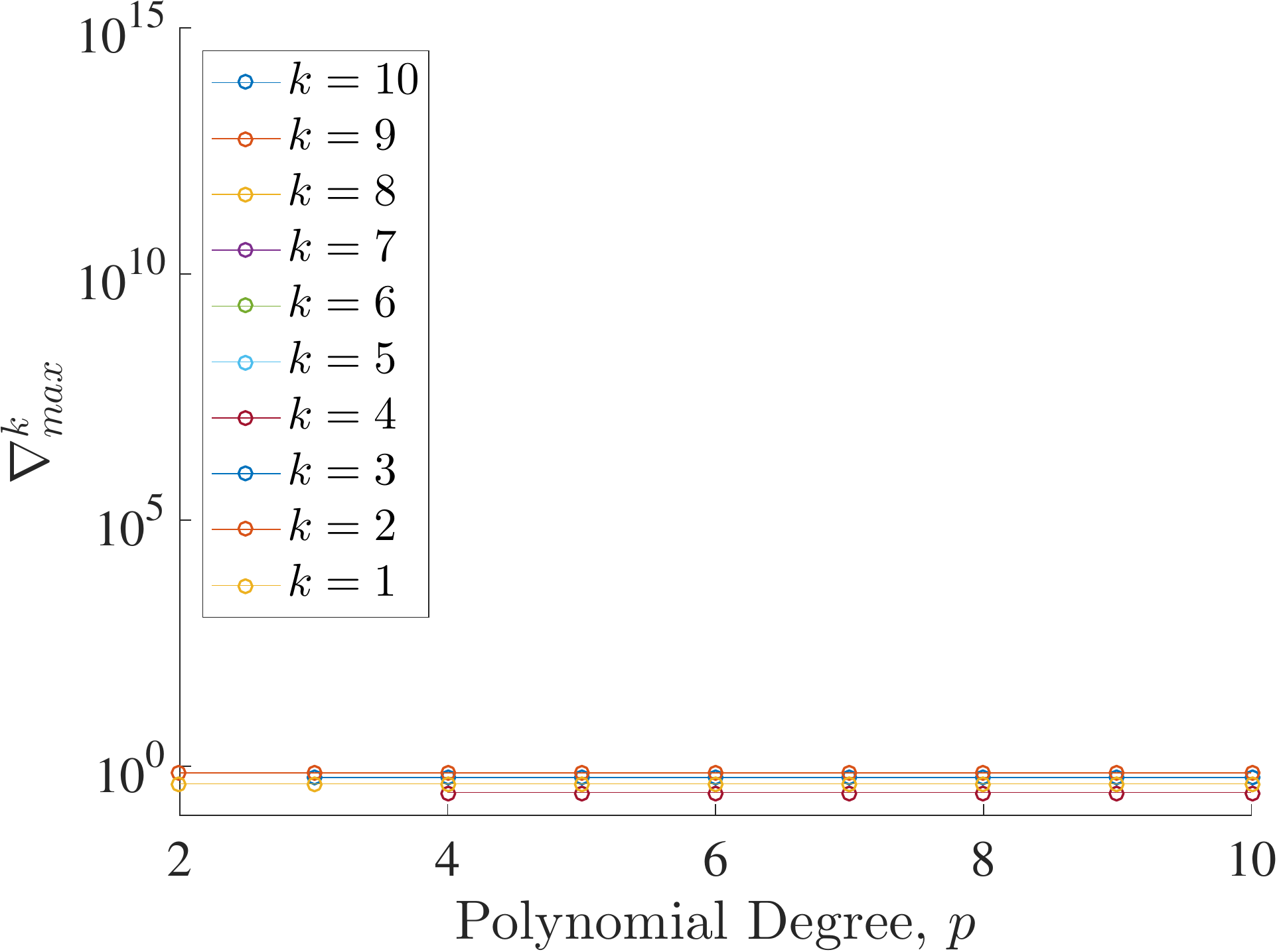}
        \caption{}
    \end{subfigure}
 \begin{subfigure}[t]{0.49\textwidth}
        \centering
            \includegraphics[width=.95\linewidth]{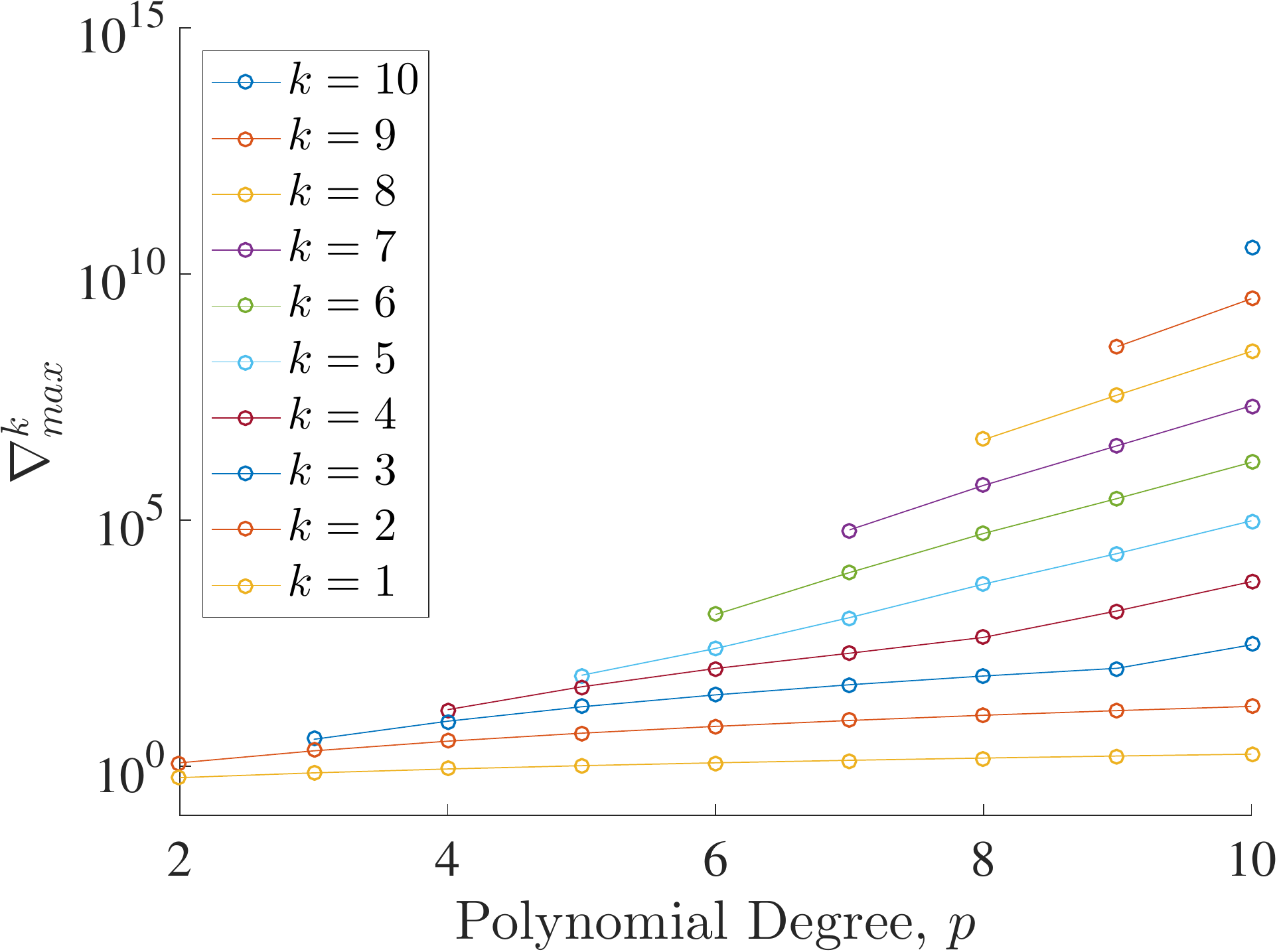}
        \caption{}
    \end{subfigure}
\caption{Lowest upper bounds on the magnitude of the higher-order derivatives of the parametric mapping for the $p$-refined meshes of the quarter annulus. The metrics for Family 1 are shown in (a) and the metrics for Family 2 are shown in (b).}
\label{QA_metrics}
\end{figure}

As before, we solve Poisson's problem, with the manufactured solution:
\begin{equation}
	u\of{x_1,x_2} = 	\dfrac{70\ln\of{\sqrt{x_1^2+x_2^2}/r_i} - 200\ln\of{\sqrt{x_1^2+x_2^2}/r_o}  }{\ln\of{r_o/r_i}}
\end{equation}
wherein $r_i$ and $r_o$ are the inner and outer radii of the quarter annulus.
The convergence plots of the $L^2$ error with respect to the polynomial degree $p$ are shown in Fig. \ref{mms_pref_L2}.
We notice immediately that the first family of meshes exhibits exponential convergence, as is expected. 
The second family, however, stagnates, even though the minimum scaled Jacobian remains well-behaved as seen in Fig. \ref{mms_pref_Jacobian}. 
The cause of this can be determined by observing the plots of the lowest upper bounds on the magnitude of the higher-order derivatives of the element-wise parametric mapping with respect to $p$, shown in Fig. \ref{QA_metrics}.
For the first family of meshes (Fig. \ref{QA_metrics}a), we see that every derivative up through order $k=10$ is bounded from above.
For the second family of meshes (Fig. \ref{QA_metrics}b), we see that  for every level of $p$ refinement, the magnitude of \textbf{\emph{every}} derivative of order $k\leq p$ increases.

\subsection{Mesh Optimization}
Thus far, we have used our mesh distortion metrics to explain sub-optimal convergence rates \textit{a posteriori}.
However, we recognize that since the mathematical theory presented in this paper relate theses mesh distortion metrics to error bounds, we should be able to use these metrics for \textit{a priori} mesh optimization.
We consider again a plate with a hole, but now with small chamfers at the corners of the plate.

We consider three families of meshes. 
In Family 1, we consider an initial coarse mesh, refined by uniform subdivision.
In Family 2, we create a new linear mesh at each refinement level $m$, enforcing a maximum edge length $h_m$ of:
\begin{equation}
h_m = \dfrac{h_1}{2^{m-1}}
\end{equation}
We then degree elevate, and recover the geometry through edge replacement. 
The mesh is then smoothed using a linear elastic analogy.
In Family 3, we begin with the same linear mesh at each refinement step as with Family 2. 
However, rather than solve a linear elasticity problem to smooth internal nodes, we instead seek to minimize the cost functional:
\begin{equation}
	Cost = \sum_{e=1}^{nel} \sum_{k=1}^{p} h^{-k}_e  \max_{\absof{\balpha}=k} \sup_{\bxi \in \Oref} \absof{ D_{\bxi}^{\balpha}\xproj}
\end{equation}
Note that the above cost functional involves a sum of order $\alpha$ scaled derivative metrics.  We further approximate the above cost functional using our computable bounds for the order $\alpha$ scaled derivative metrics from Theorem IV, resulting in the modified cost functional:
\begin{equation}
	Modified = \sum_{e=1}^{nel}  \sum_{k=1}^{p} h^{-k}_e \max_{\absof{\balpha}=k}
        \dfrac{p!}{(p - k)!} \max\limits_{\bi \in I^{p-\absof{\balpha}}}
        \absof{
          \sum_{\bj \in I^{\balpha}}
          (-1)^{{\balpha}+\bj}{{{\balpha}}\choose{\bj}}
          \widetilde{\bP}_{\bi + \bj}}
\end{equation}
The resulting three families of meshes are shown in Table  \ref{opt_families}.  Note that our new mesh optimization procedure can be interpreted as a generalization of biharmonic mesh smoothing, which minimizes a cost functional based on the second-derivatives of the parametric mapping \cite{helenbrook2003mesh}.

\begin{table}[t]
  \begin{center}
    \caption{Three mesh families for a plate with a hole and chamfered corners.} 
    \begin{tabular}{|  >{\centering\arraybackslash} m{.3cm} |  >{\centering\arraybackslash} m{4cm} |  >{\centering\arraybackslash} m{4cm} |    
    >{\centering\arraybackslash} m{4cm} | }
      \hline
       $m$ & Family 1 & Family 2 & Family 3  \\ 
      \hline
       1 &
      \vspace{5pt}
      \includegraphics[width=4cm,height=4cm,keepaspectratio]{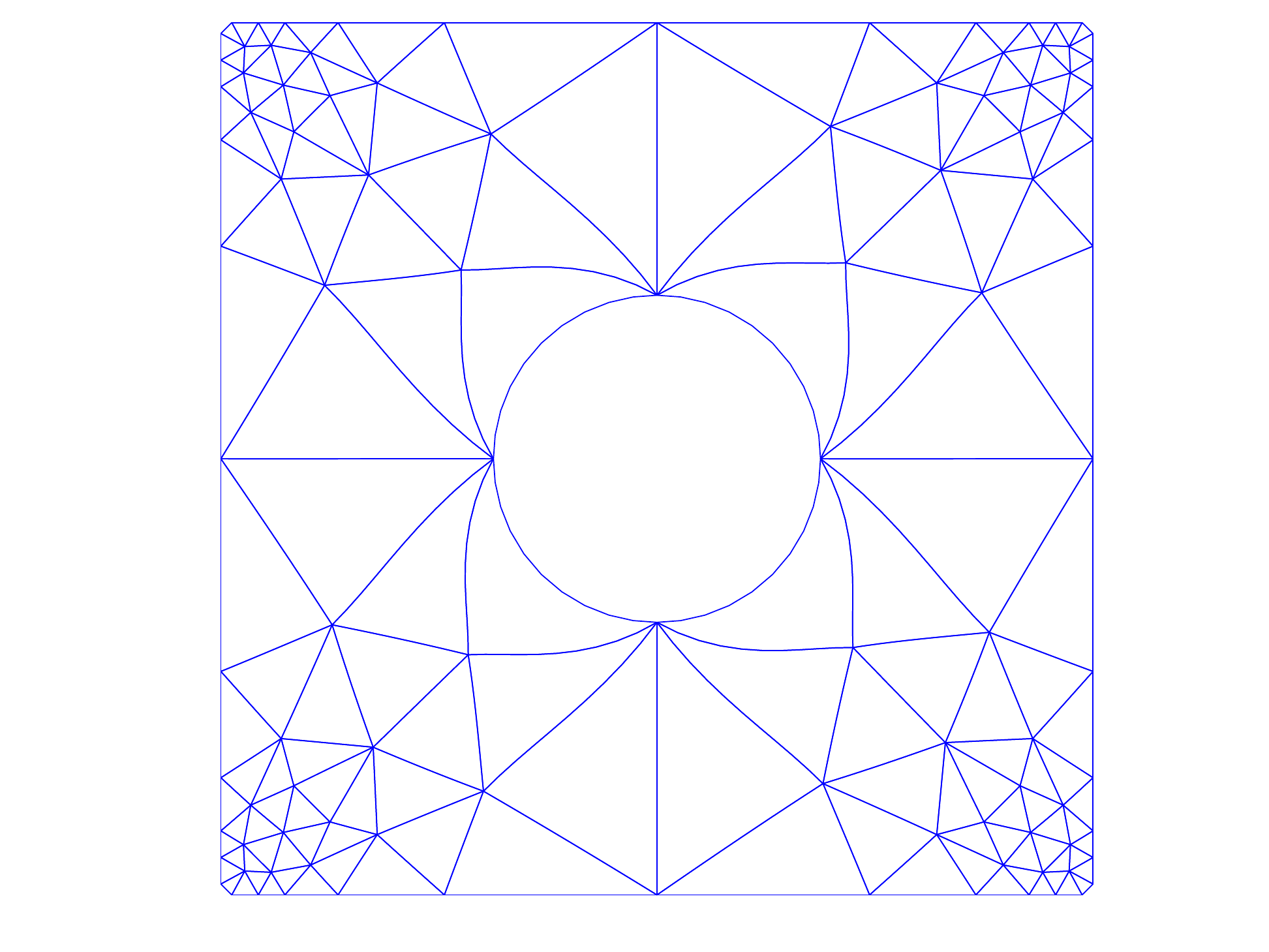} &
      \vspace{5pt}
      \includegraphics[width=4cm,height=4cm,keepaspectratio]{mms_base_mesh.pdf} &
       \vspace{5pt}
      \includegraphics[width=4cm,height=4cm,keepaspectratio]{mms_base_mesh.pdf}\\
      \hline
             2 &
      \vspace{5pt}
      \includegraphics[width=4cm,height=4cm,keepaspectratio]{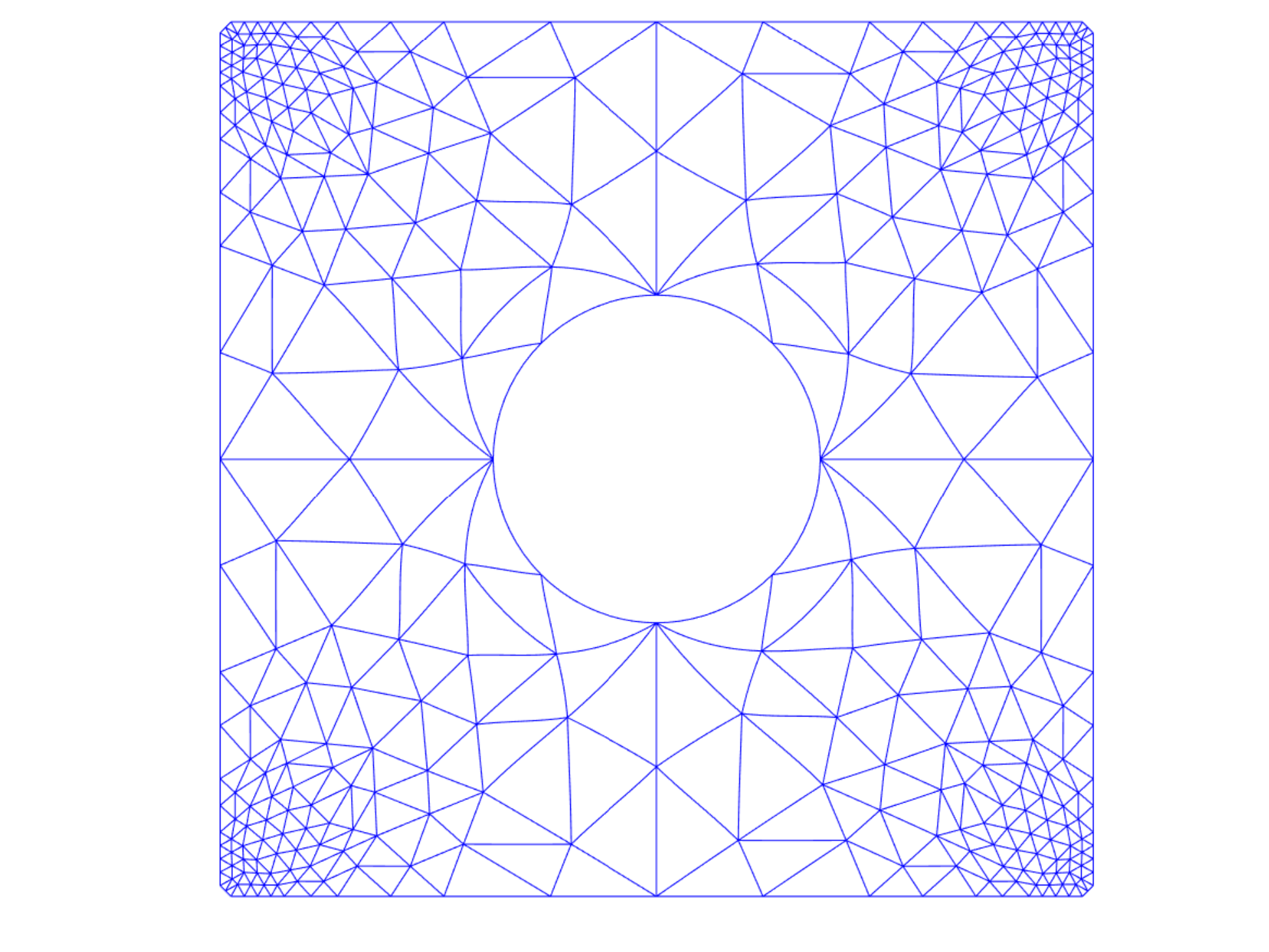} &
      \vspace{5pt}
      \includegraphics[width=4cm,height=4cm,keepaspectratio]{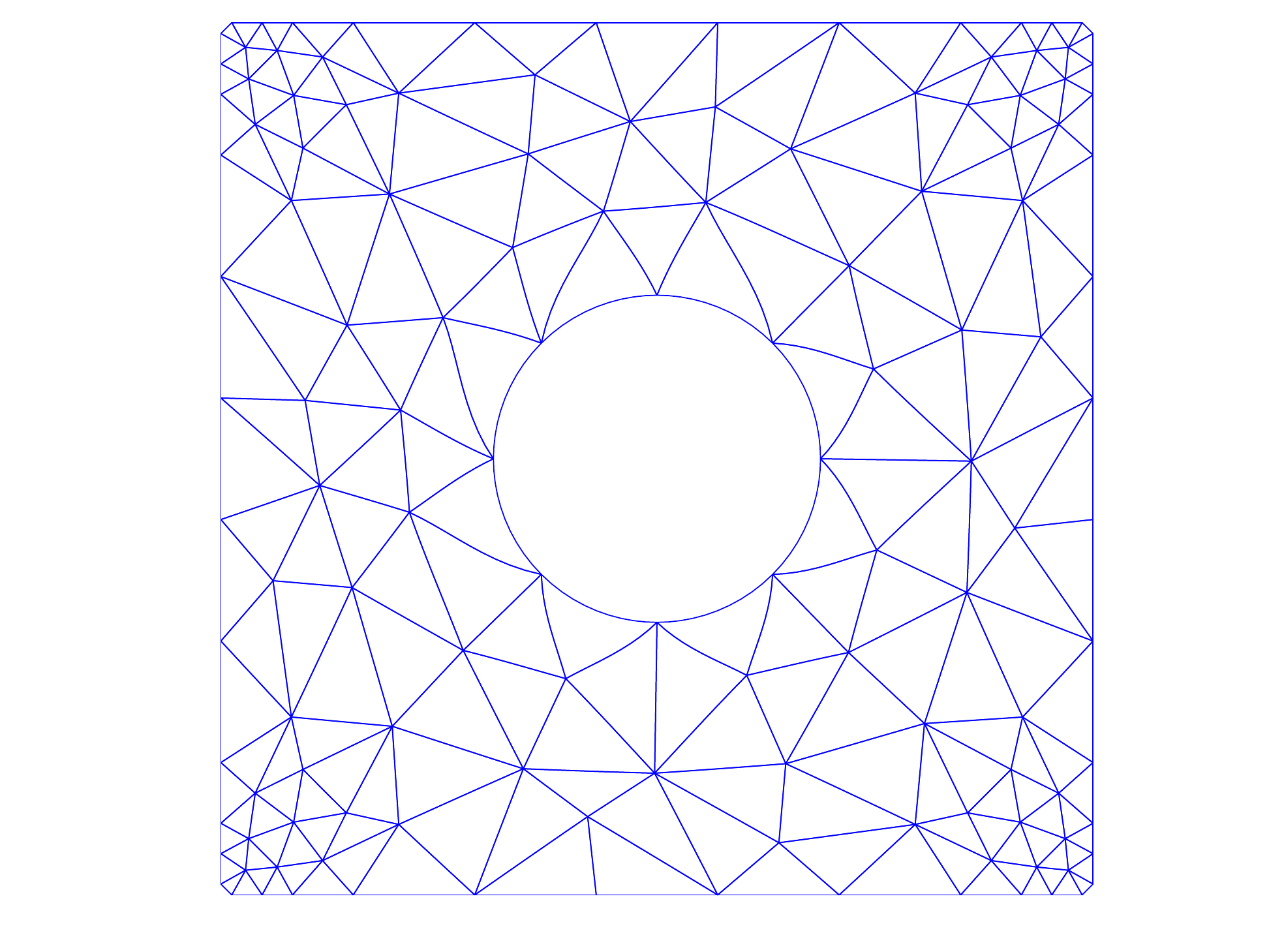} &
       \vspace{5pt}
      \includegraphics[width=4cm,height=4cm,keepaspectratio]{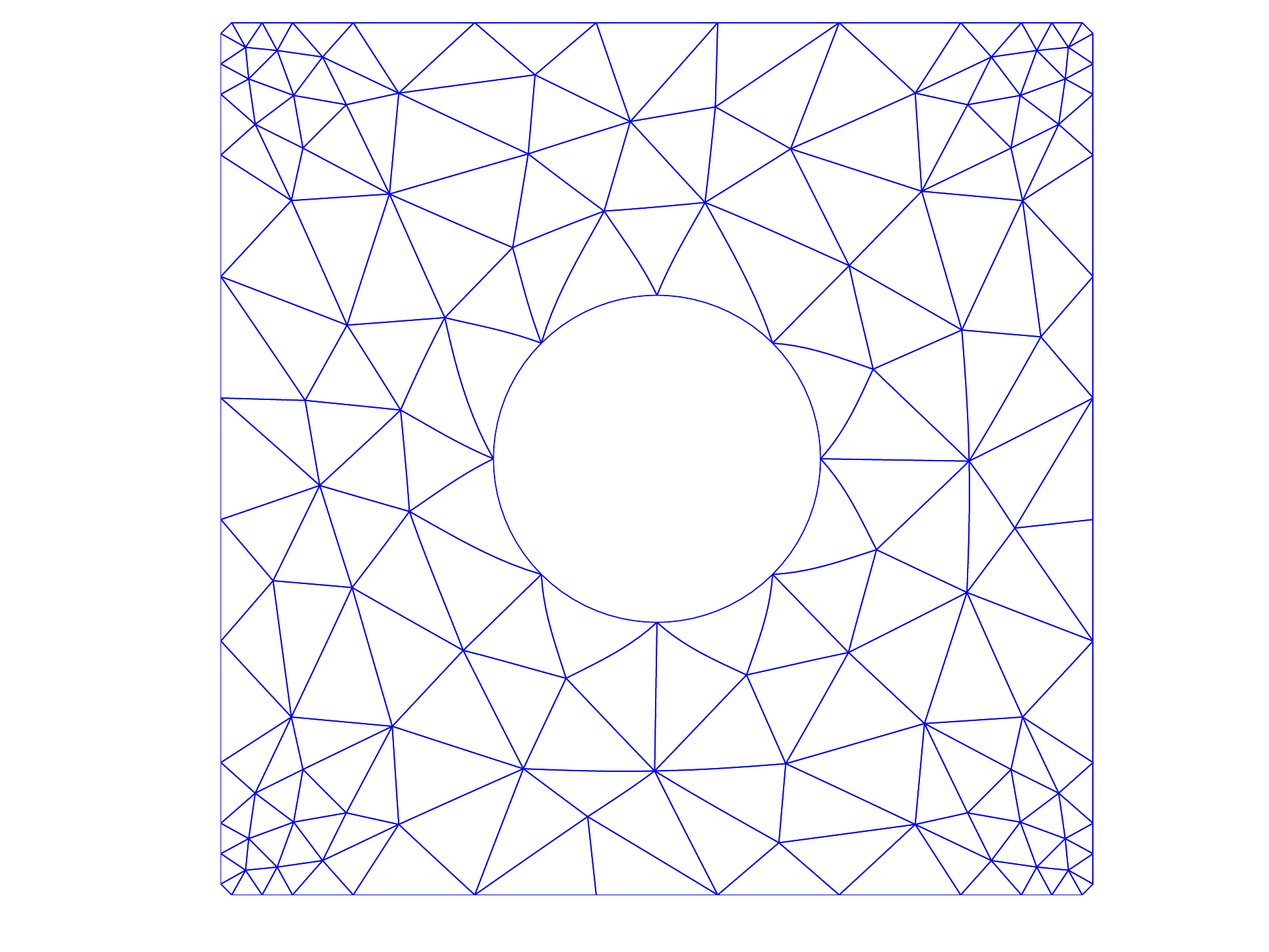}\\
      \hline
    3 &
      \vspace{5pt}
      \includegraphics[width=4cm,height=4cm,keepaspectratio]{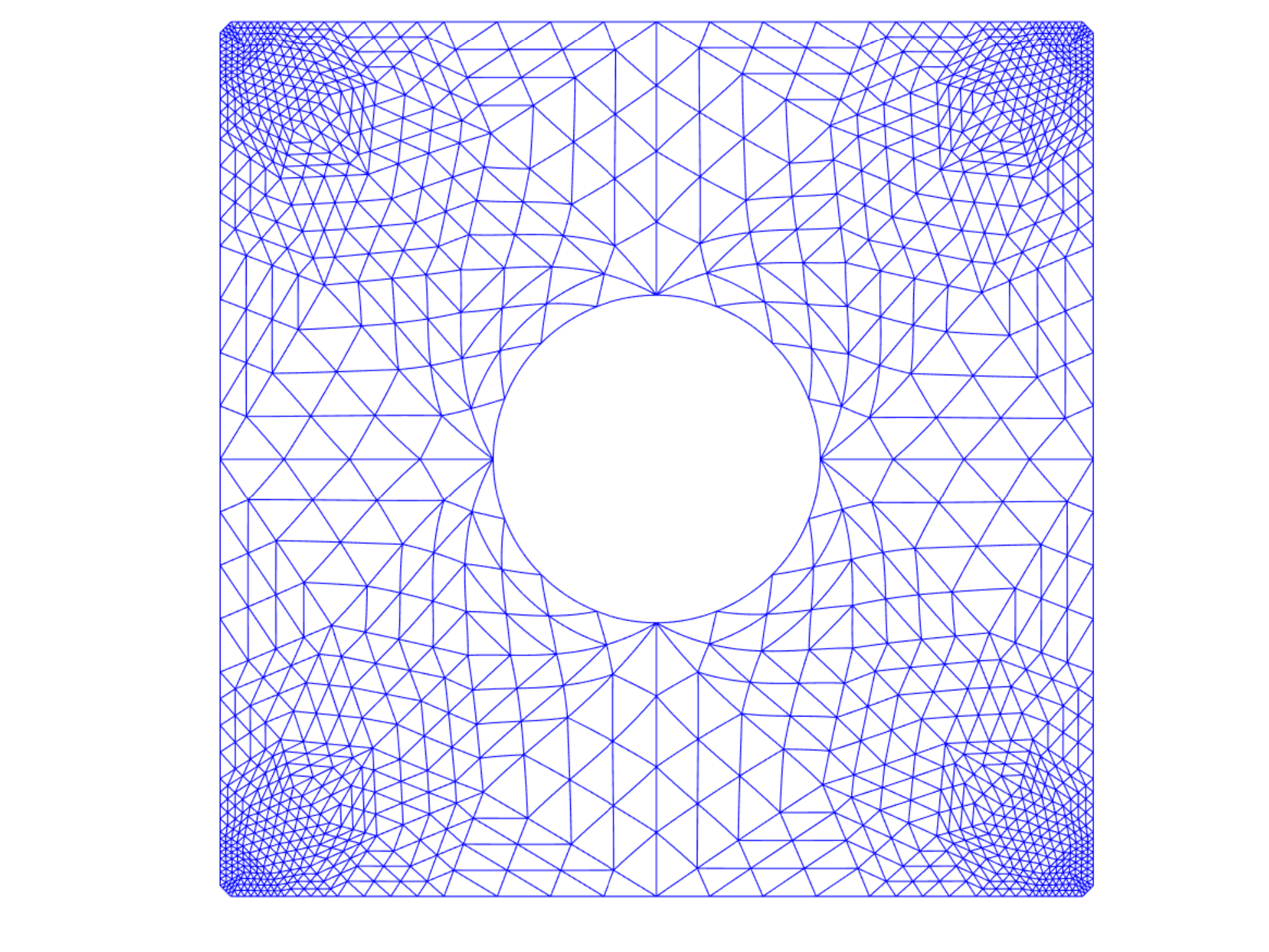} &
      \vspace{5pt}
      \includegraphics[width=4cm,height=4cm,keepaspectratio]{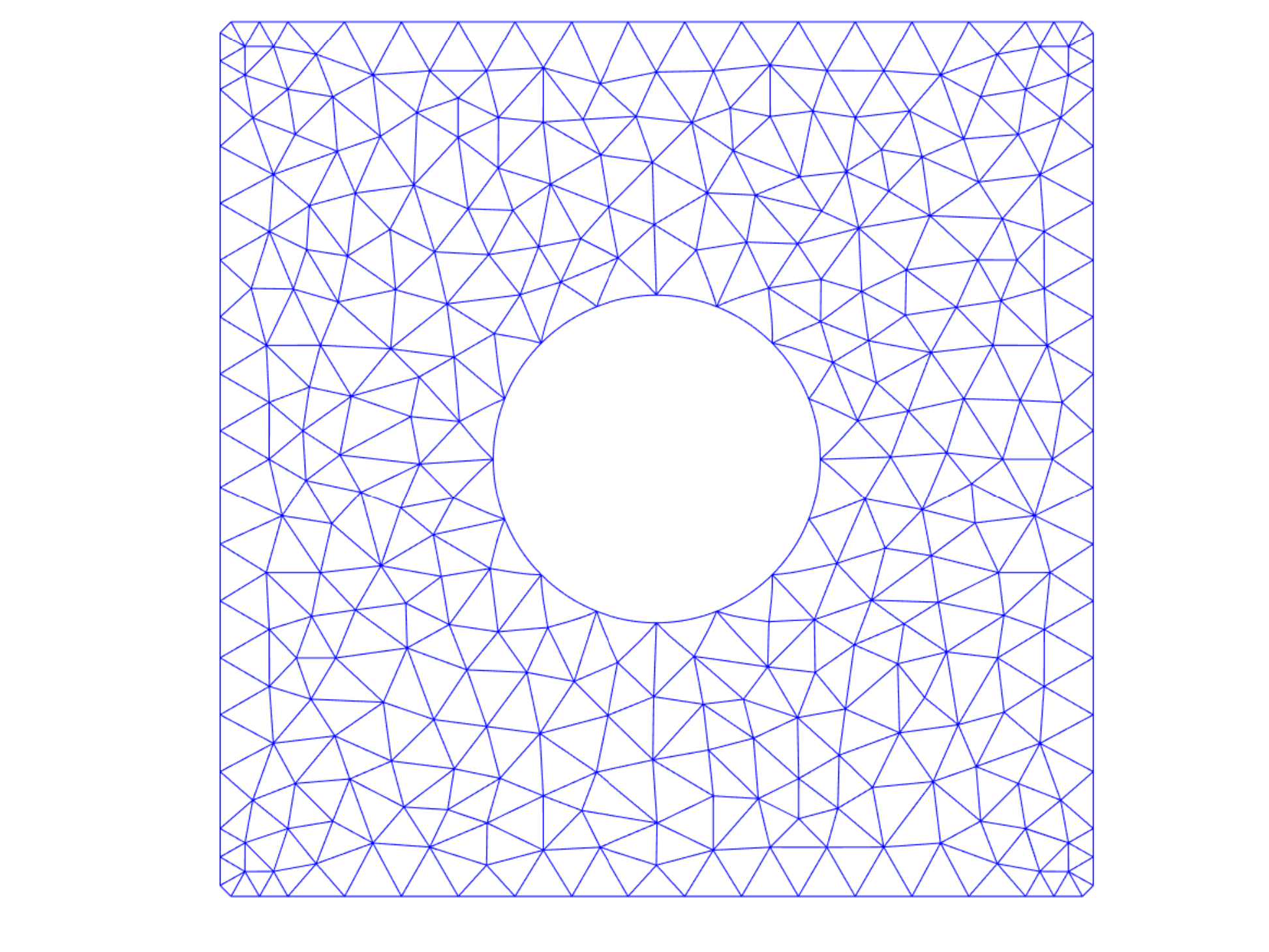} &
       \vspace{5pt}
      \includegraphics[width=4cm,height=4cm,keepaspectratio]{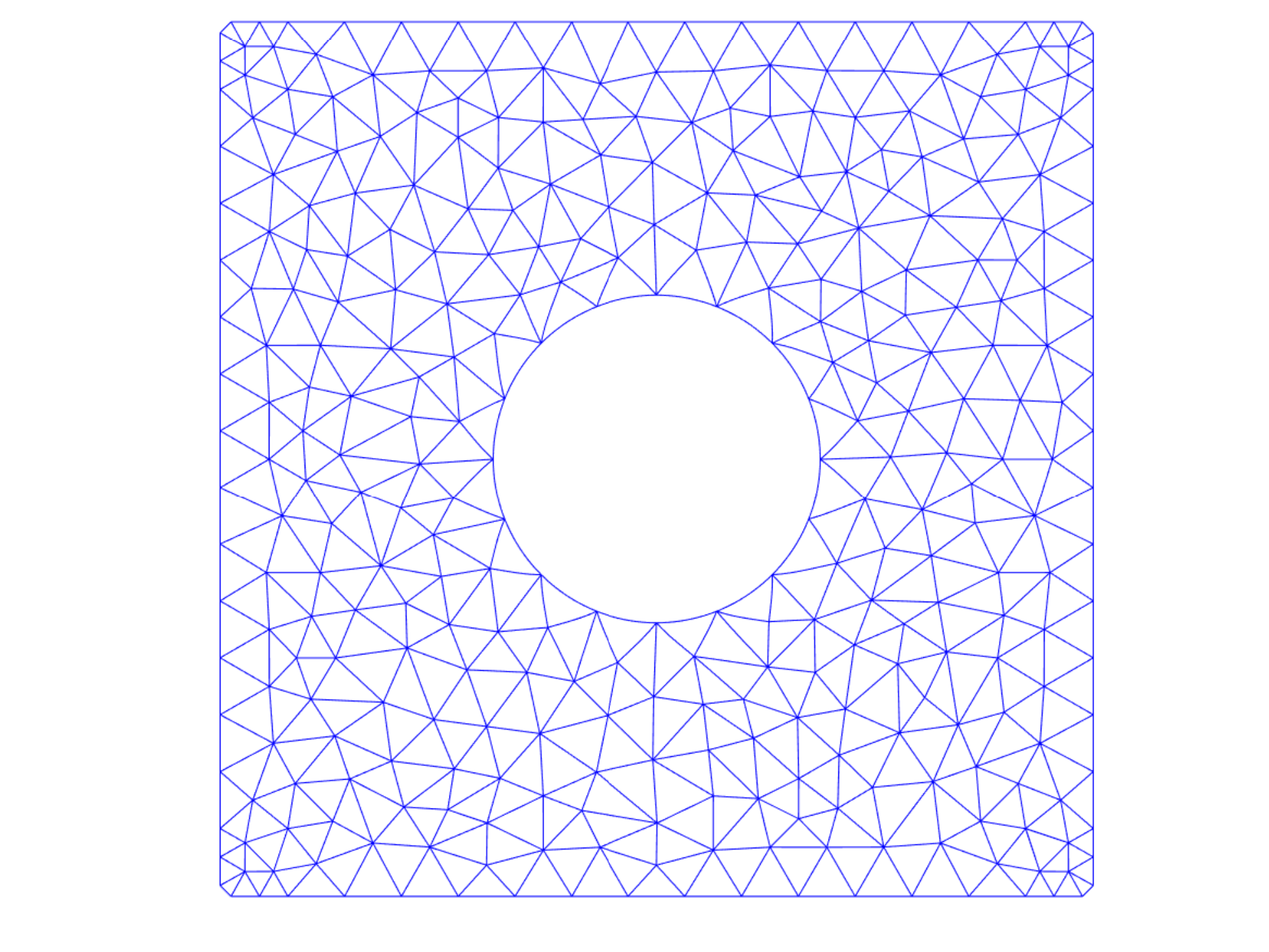}\\
      \hline
      \end{tabular}
    \label{opt_families}
  \end{center}
\end{table}

\begin{figure}[t] 
\centering
\begin{subfigure}[t]{0.49\textwidth}
        \centering
\includegraphics[width=.95\textwidth]{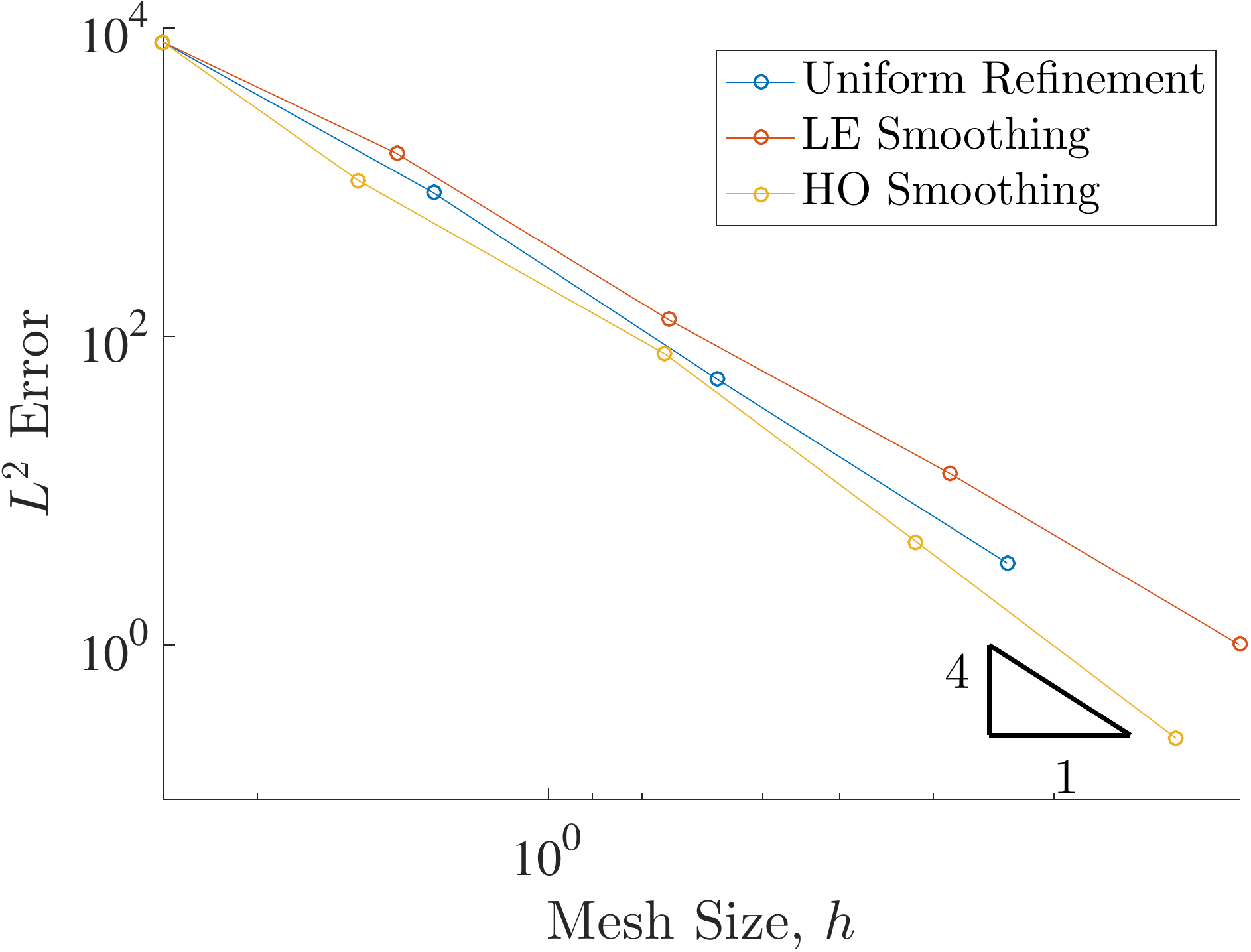}
\caption{}
\end{subfigure}
\begin{subfigure}[t]{0.49\textwidth}
        \centering
\includegraphics[width=.95\textwidth]{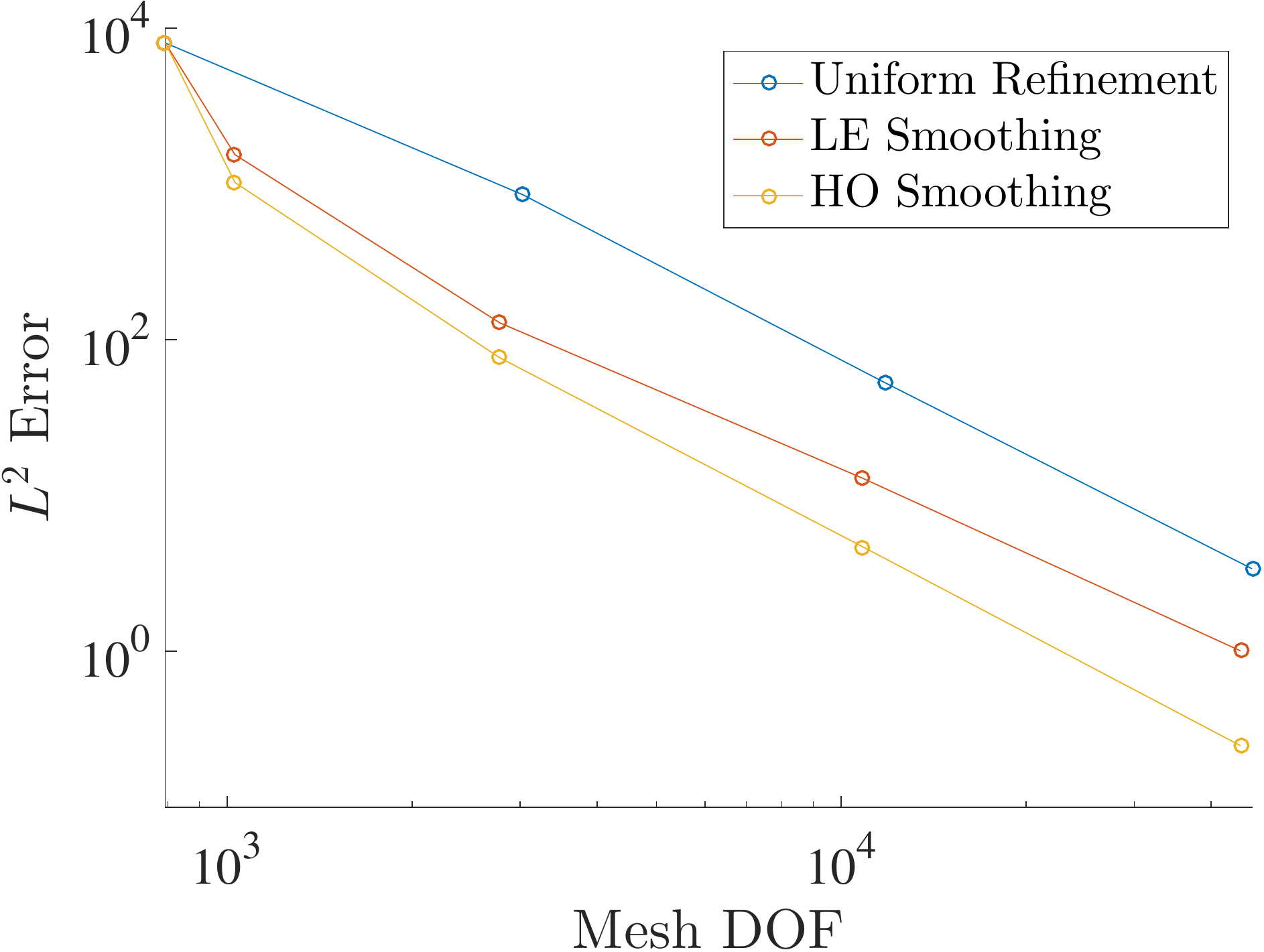}
\caption{}
\end{subfigure}
\caption{Convergence plots of the error in the $L^2$ norm. (a) Error versus mesh size. (b) Error versus nodal degrees of freedom in the mesh. Here, ``Uniform Refinement'' denotes the results associated with uniform refinement (Family 1), ``LE Smoothing'' denotes the results associated with linear elastic mesh smoothing (Family 2), and ``HO Smoothing'' denotes the results associated with our new mesh optimization procedure (Family 3).}
\label{mms_opt_L2}
\end{figure}

There is visually little difference between the meshes of Family 2, generated using linear elastic mesh smoothing, and the meshes of Family 3, generated using our new mesh optimization procedure. 
To study the effect of linear elastic mesh smoothing and our new mesh optimization procedure on the accuracy of a finite element approximation, we solve Poisson's problem over each mesh family, with the manufactured solution:
\begin{equation}
\begin{split}
	u\of{x_1,x_2} = \of{2a-c -x_1-x_2 } \of{2a-c -x_1+x_2 }\of{2a-c +x_1-x_2 }\of{2a-c +x_1+x_2 } ...\\
	...\of{x_1-a}\of{x_1+a}\of{x_2-a}\of{x_2+a}\of{r-\sqrt{x_1^2+x_2^2}}	
	\end{split}
\end{equation}
wherein $a$ is the plate half width, $r$ is the radius of the hole, and $c$ is the length of the chamfer.
We plot the convergence rate of the error in Fig. \ref{mms_opt_L2}. 
\begin{figure}[t!]
\centering
\begin{subfigure}[t]{0.49\textwidth}
        \centering
            \includegraphics[width=.95\linewidth]{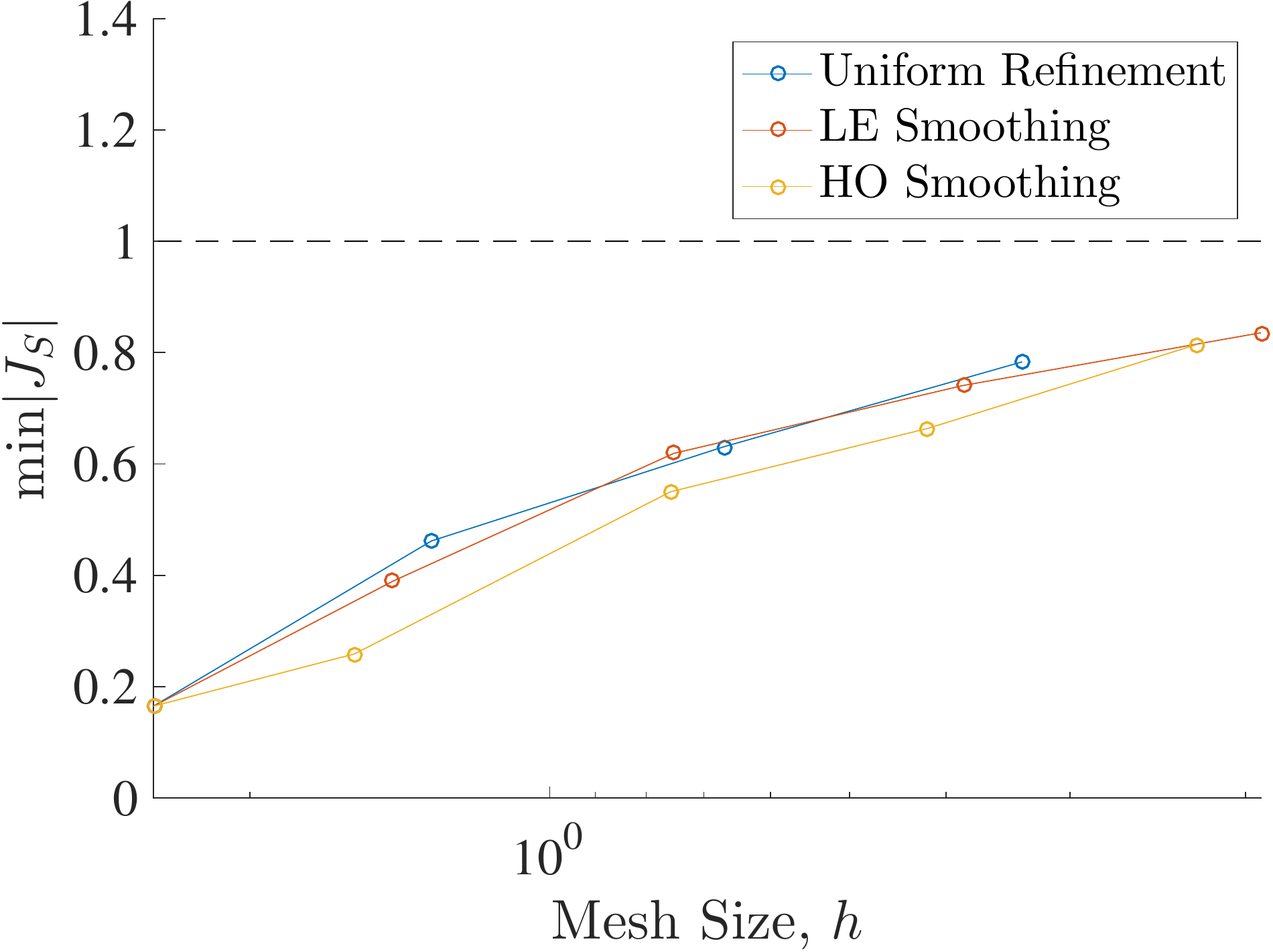}
        \caption{}
    \end{subfigure}
 \begin{subfigure}[t]{0.49\textwidth}
        \centering
            \includegraphics[width=.95\linewidth]{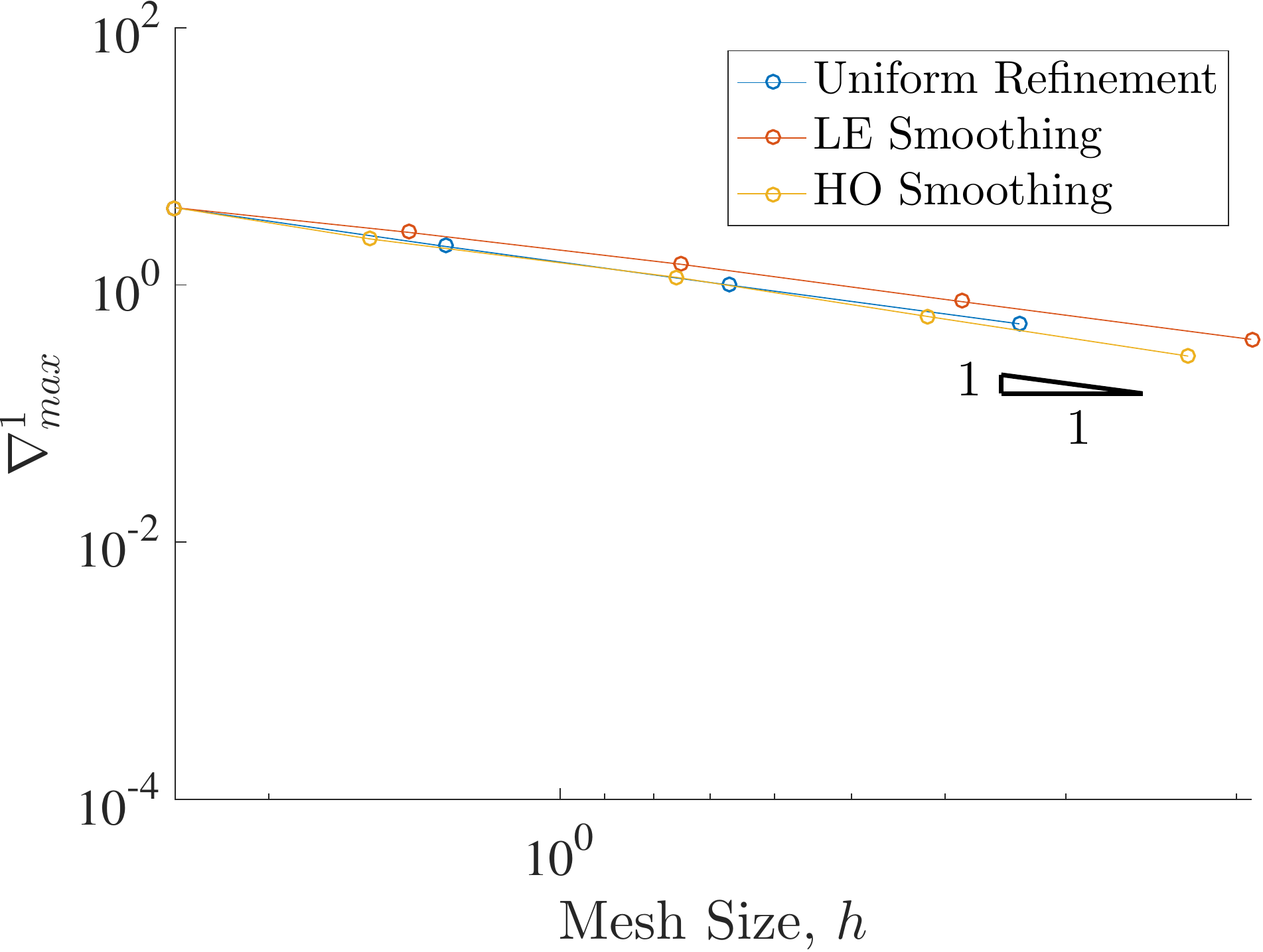}
        \caption{}
    \end{subfigure}
\begin{subfigure}[t]{0.49\textwidth}
        \centering
            \includegraphics[width=.95\linewidth]{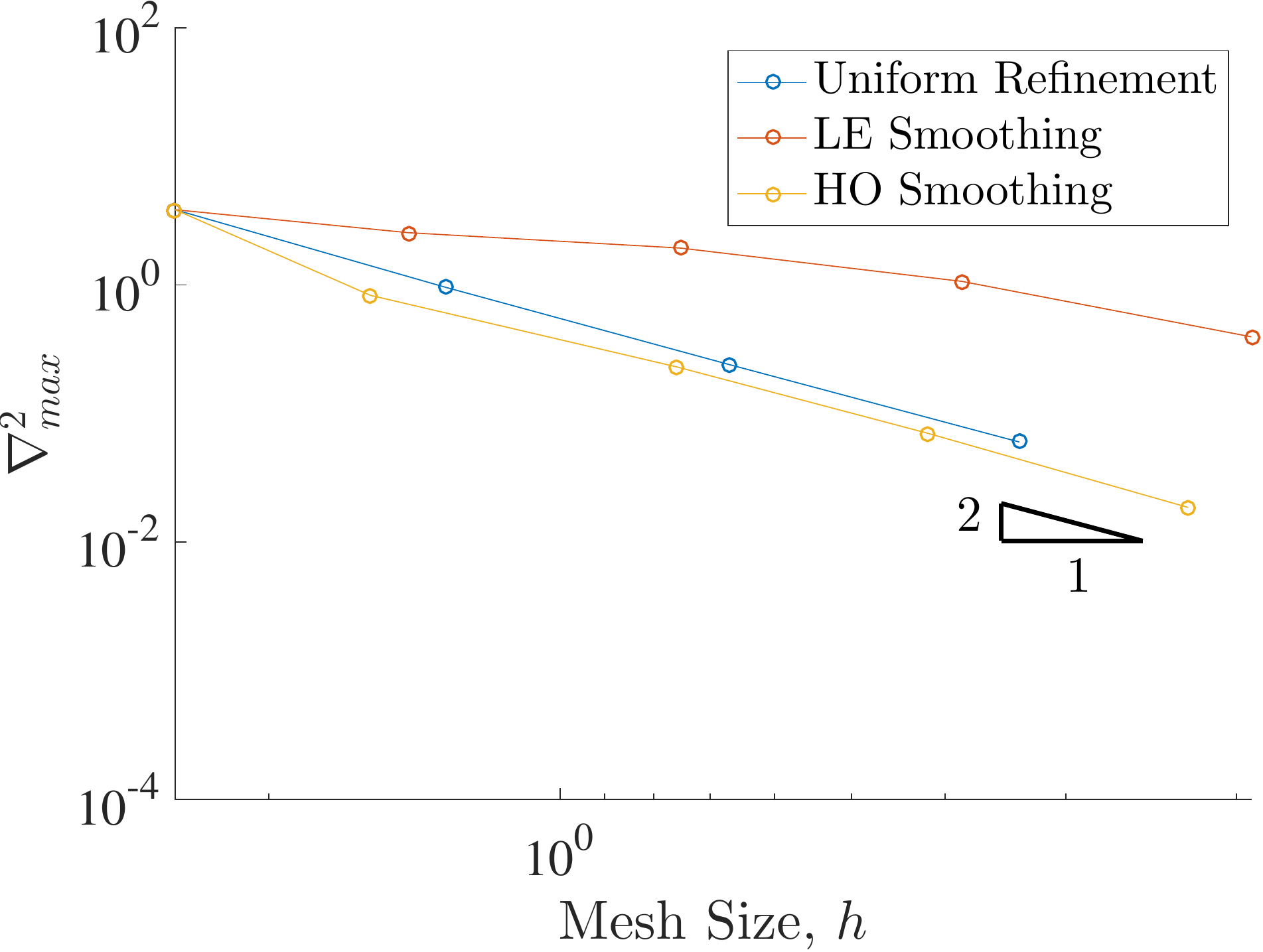}
        \caption{}
    \end{subfigure}
 \begin{subfigure}[t]{0.49\textwidth}
        \centering
            \includegraphics[width=.95\linewidth]{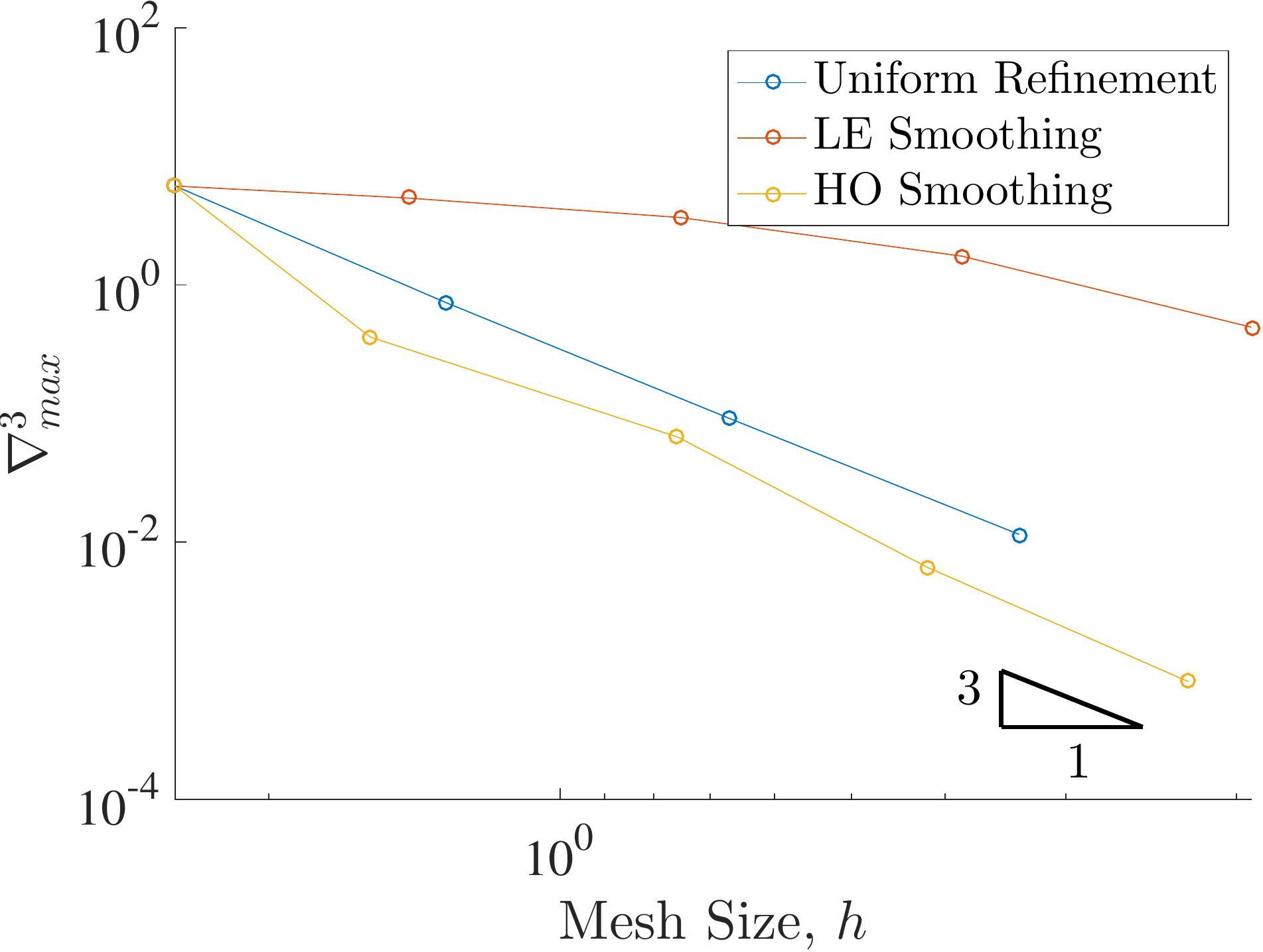}
        \caption{}
    \end{subfigure}
\caption{ Mesh distortion metrics for the triangular meshes of the plate with a hole and chamfers.  (a) Minimum scaled Jacobian. (b) Lowest upper bound on the magnitude of the first derivatives. (c) Lowest upper bound on the magnitude of the second derivatives. (d) Lowest upper bound on the magnitude of the third derivatives. Here, ``Uniform Refinement'' denotes the results associated with uniform refinement (Family 1), ``LE Smoothing'' denotes the results associated with linear elastic mesh smoothing (Family 2), and ``HO Smoothing'' denotes the results associated with our new mesh optimization procedure (Family 3).}
\label{opt_metrics}
\end{figure}
From Fig. \ref{mms_opt_L2}a, we see that both uniform subdivision (Family 1) and our new mesh optimization procedure (Family 3) outperform linear elastic smoothing (Family 2) in terms of convergence with respect to $h$. 
However, we note that the initial coarse mesh has small elements because of the chamfer. 
As such, the meshes in Family 1 become overly refined at the plate corners under uniform subdivision. 
If we instead plot solution error with respect to system degrees of freedom, as done in Fig. \ref{mms_opt_L2}b, we see that both Family 2 and Family 3 offer increased accuracy with respect to computational cost. 
Furthermore, we note that Family 3 outperforms Family 2 in every case. 
We plot element distortion metrics for the three families in Fig. \ref{opt_metrics}.
The obtained results are compelling, as they suggest that we may use our element distortion metrics to effectively optimize higher-order meshes.


\section{Conclusions} \label{conclusions}

When we began this work, it was our belief that the current element metrics used in $p$-version finite element and IGA were insufficient for generating high-quality curvilinear meshes.
To this end, we have not only developed a complete theory characterizing the effect of curvilinear mesh distortion on the approximation properties of rational \BB elements, 
but have also developed a suite of computable element distortion metrics based on this theory that are suitable for use with modern meshing algorithms.
However, while we are excited about the implications of these results, we recognize there exist a large number of future research directions.

First, we note that the error bounds presented in Theorem I hold in the limit of mesh refinement.
As a result, these can be rather loose upper bounds on the error, particularly over coarse meshes.
As a further consequence of this, the conditions presented in Theorem II are sufficient but not necessary conditions to guarantee a rational \BB finite element discretization on a curvilinear mesh will exhibit similar convergence rates as a finite element discretization over a linear mesh.
Thus, the element distortion metrics presented here are not always a good indicator of the effect of element distortion on solution accuracy.
For fluid flow simulations, boundary layer meshes consisting of highly distorted elements near walls typically yield better results per degree of freedom than their isotropic counterparts \cite{sahni2008adaptive,sahni_curved_2010}.
We are curious to see if the error bounds here can be sharpened, particularly if something is known \emph{a priori} about the PDE to be solved.

Second, while our computable distortion metrics appear to have significant promise in the context of mesh optimization, our implementation based on these metrics is not optimized for computational efficiency, and we have not benchmarked our mesh optimization procedure on large problems.
In the future, we plan to study the performance of our new mesh optimization procedure both for complex problems of engineering interest as well as on high performance computing platforms.
We are especially interested to compare the performance of our mesh optimization procedure with state-of-the-art mesh optimization procedures, such as those in \cite{gargallo-peiro_optimization_2015,remacle_optimizing_2014,toulorge_robust_2013}.

\section{Acknowledgements}
This material is based upon work supported by the National Science Foundation Graduate Research Fellowship Program under Grant No. DGE 1144083 as well as work supported by the Air Force Office of Scientific Research under Grant No. FA9550-14-1-0113.

\bibliographystyle{plain}
\bibliography{References}

\appendix
\appendixpage
\addappheadtotoc
\section{Stencils for Derivatives of \BB Elements}\label{deriv_appendix}

This appendix includes lookup tables for several common 2D \BB elements.
Each table contains stencils for every non-zero derivative over the  respective element.
Furthermore, the nodes at which to apply the given stencil are shown on a reference element.
We do not include  explicit stencils for 3D elements, as they are hard to visualize, but stencils for any 2D or 3D simplicial or tensor product element can be derived using the equations derived previously in this paper.

\begin{table}[H]
  \begin{center}
    \caption{Partial derivative stencils for quadratic  \BB triangles.}
    \begin{tabular}{ | >{\centering\arraybackslash} m{1.5cm}  | >{\centering\arraybackslash} m{3cm} |  >{\centering\arraybackslash} m{2.3cm} ||
        >{\centering\arraybackslash} m{1.5cm}  | >{\centering\arraybackslash} m{1.5cm} |  >{\centering\arraybackslash} m{2.3cm} |
      }
      \hline
      \multicolumn{6}{|c|}{  \textbf{ First Order Derivatives, $|\balpha| = 1$}}\\ 
      \hline
      Deriv. & Stencil & Evaluation Triangle & Deriv. & Stencil & Evaluation Triangle  \\ 
      \hline
      $\deldel{\xproj}{\xi_1} $&
      \includegraphics[width=2.5cm,height=2.5cm,keepaspectratio]{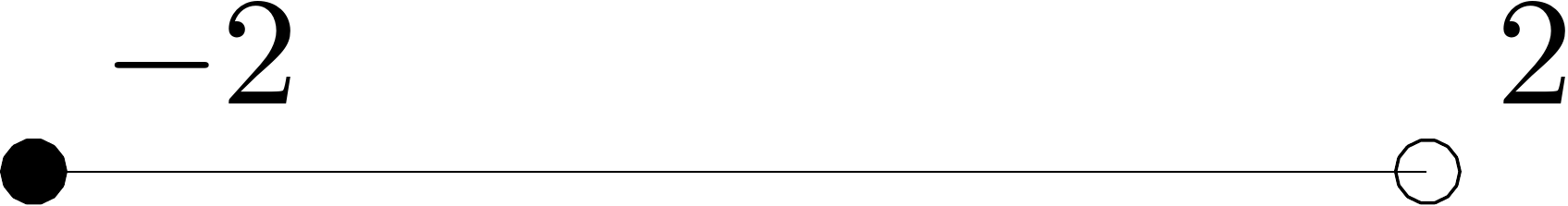} &
      \vspace{5pt}
      \includegraphics[width=2.0cm,height=2.0cm,keepaspectratio]{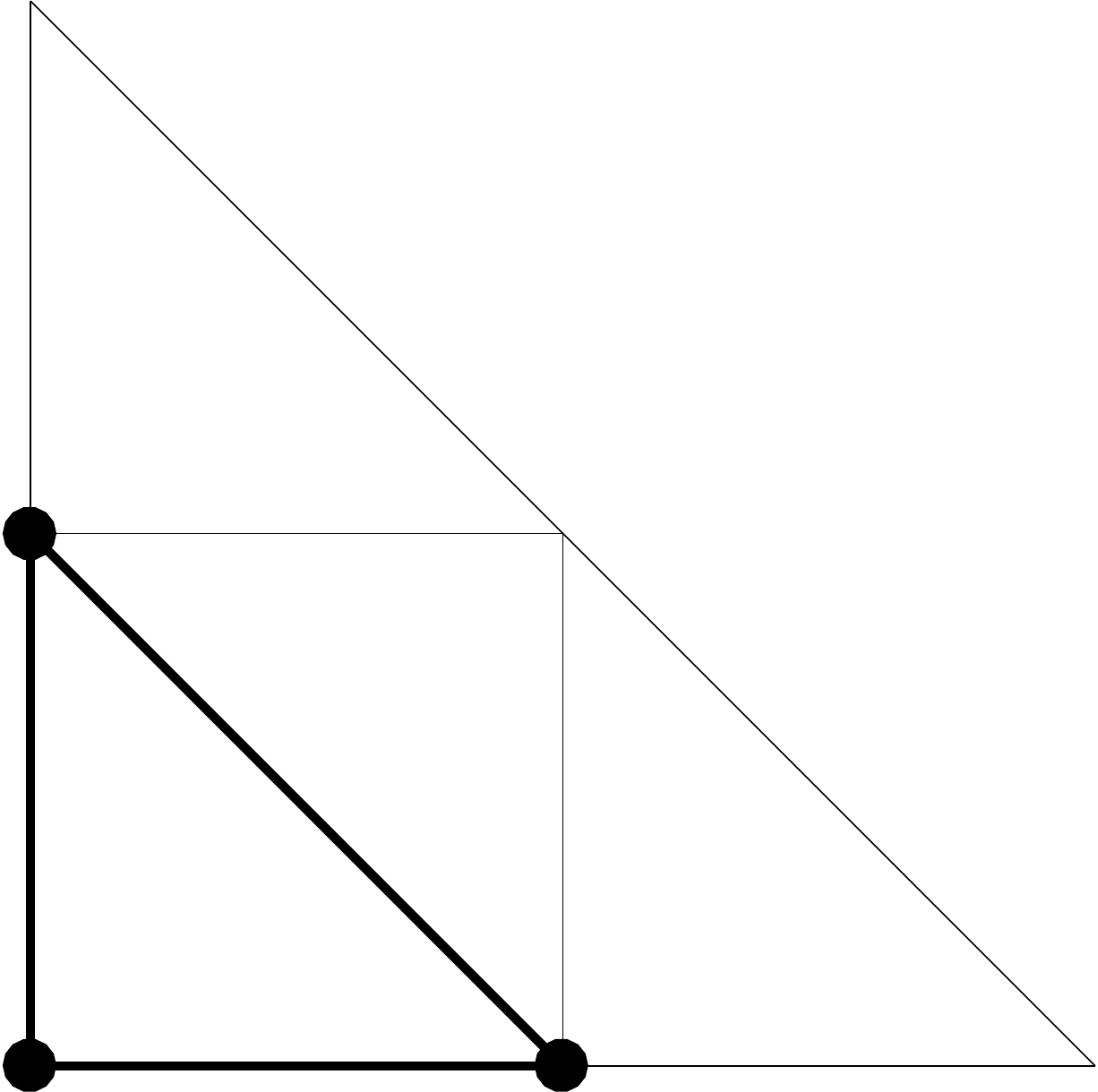} &
      $\deldel{\xproj}{\xi_2} $&
      \includegraphics[width=2.0cm,height=2.0cm,keepaspectratio]{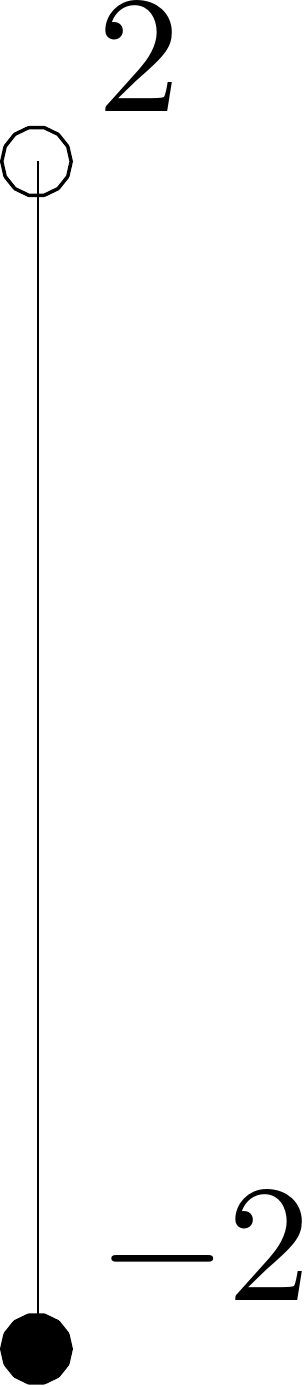} &
      \vspace{5pt}
      \includegraphics[width=2.0cm,height=2.0cm,keepaspectratio]{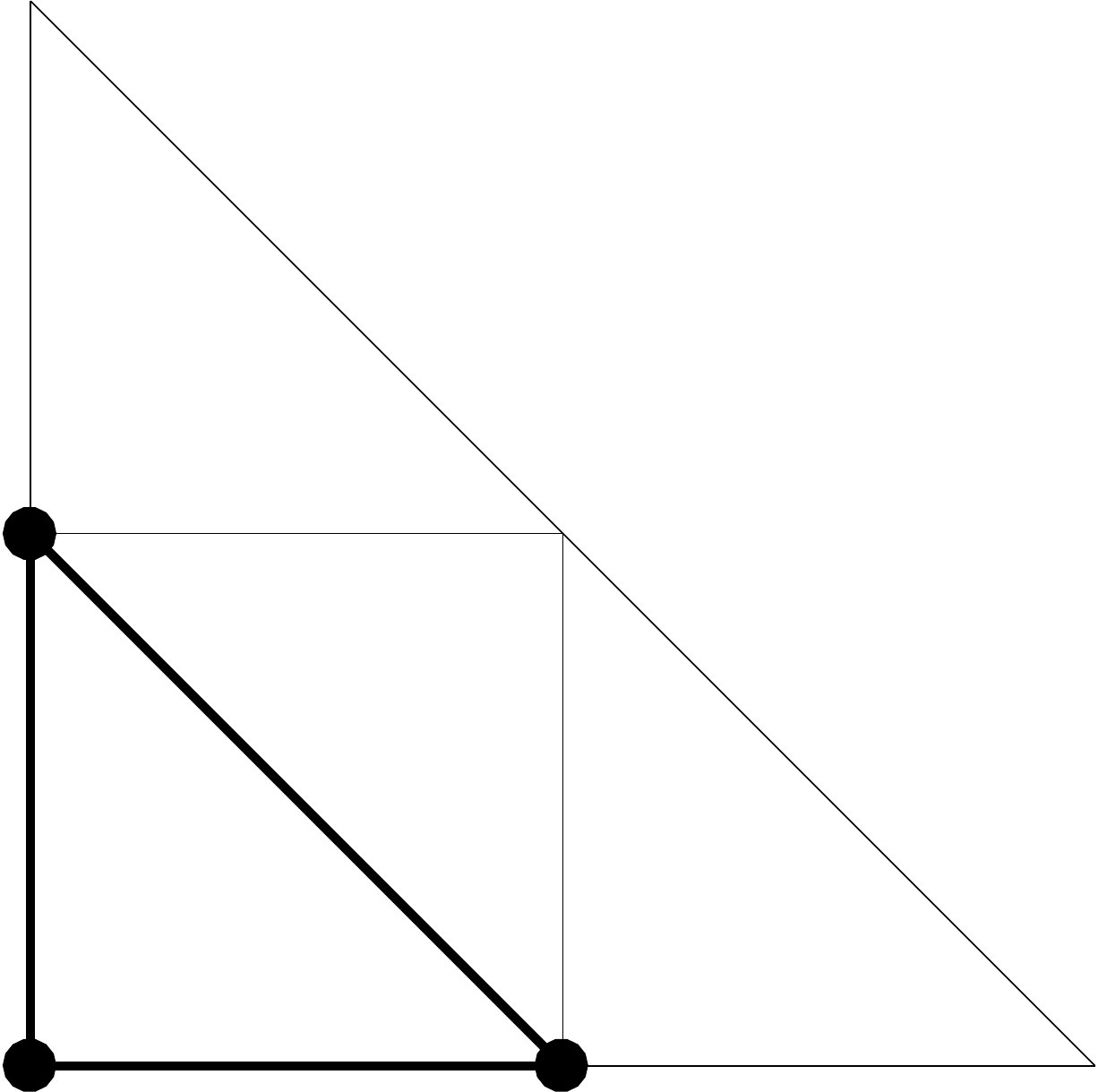} \\
      \hline
      \multicolumn{6}{c}{\vspace{8pt}}\\ 
      \hline
      \multicolumn{6}{|c|}{\textbf{Second Order Derivatives, $|\balpha| = 2$}}\\ 
      \hline
      Deriv. & Stencil & Evaluation Triangle & Deriv. & Stencil & Evaluation Triangle  \\ 
      \hline
      $\dfrac{\partial^2 \xproj}{\partial \xi_1^2} $&
      \includegraphics[width=2.5cm,height=2.5cm,keepaspectratio]{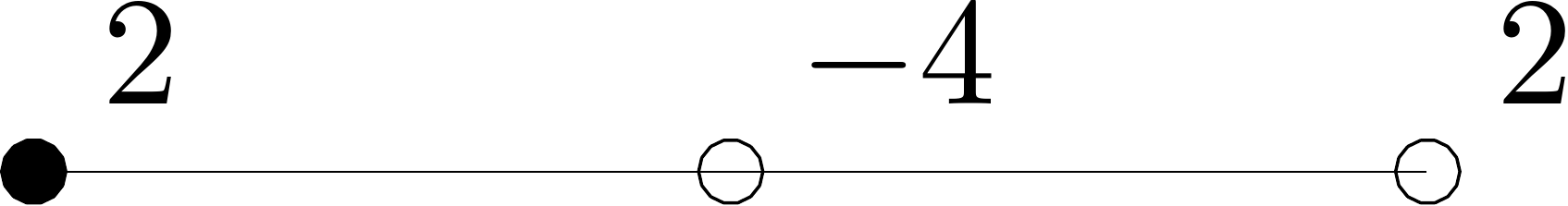} &
      \vspace{5pt}
      \includegraphics[width=2.0cm,height=2.0cm,keepaspectratio]{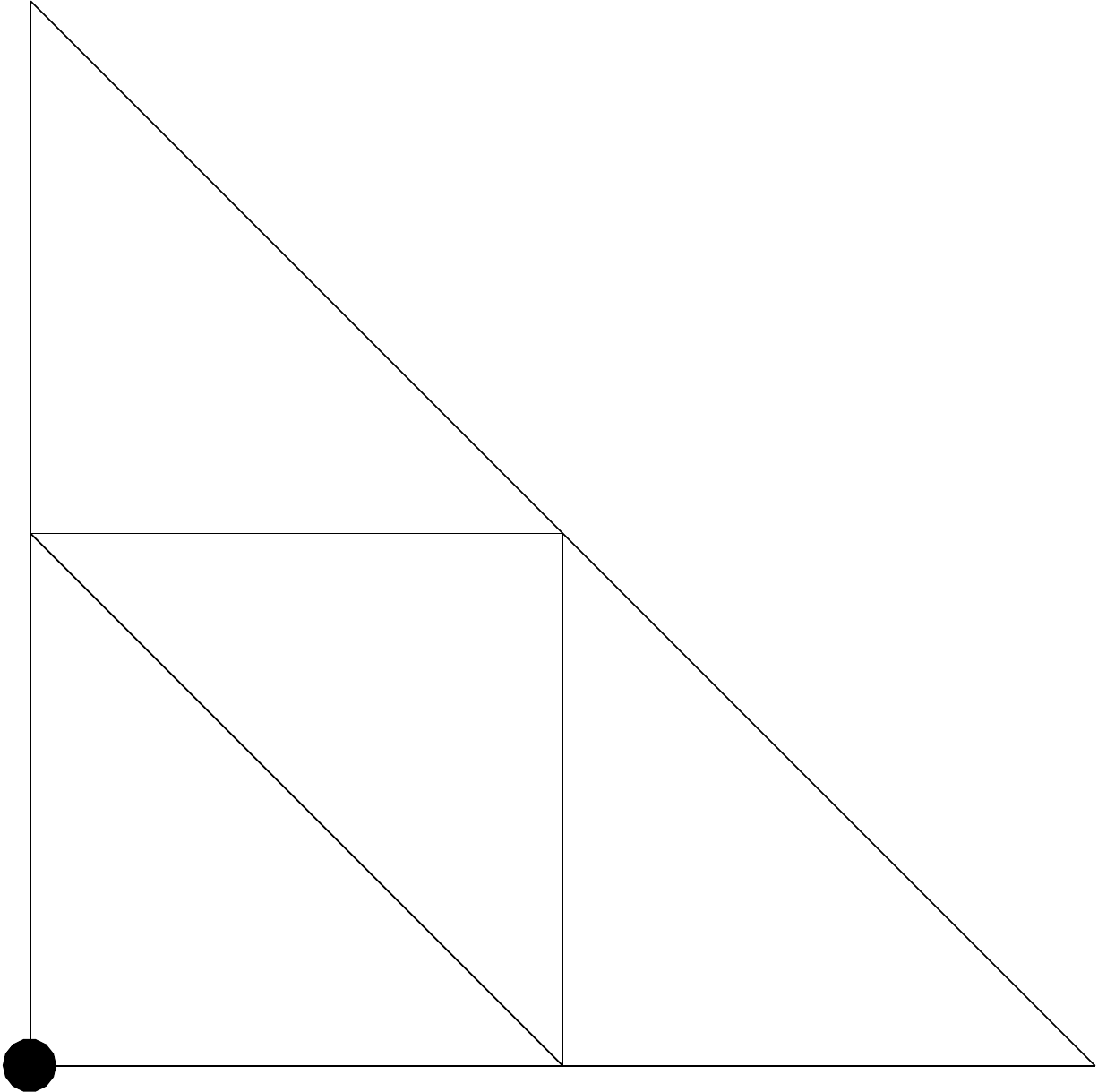} &
      $\dfrac{\partial^2 \xproj}{\partial \xi_2^2} $&
      \includegraphics[width=2.0cm,height=2.0cm,keepaspectratio]{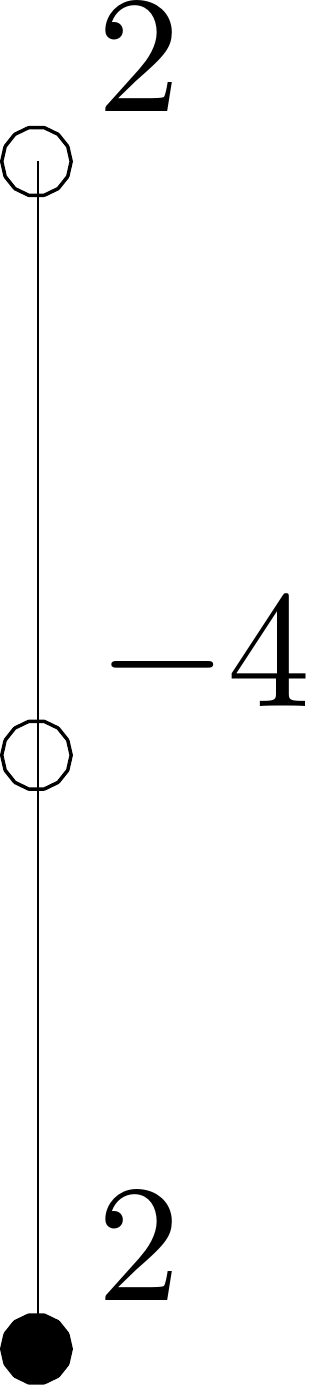} &
      \vspace{5pt}
      \includegraphics[width=2.0cm,height=2.0cm,keepaspectratio]{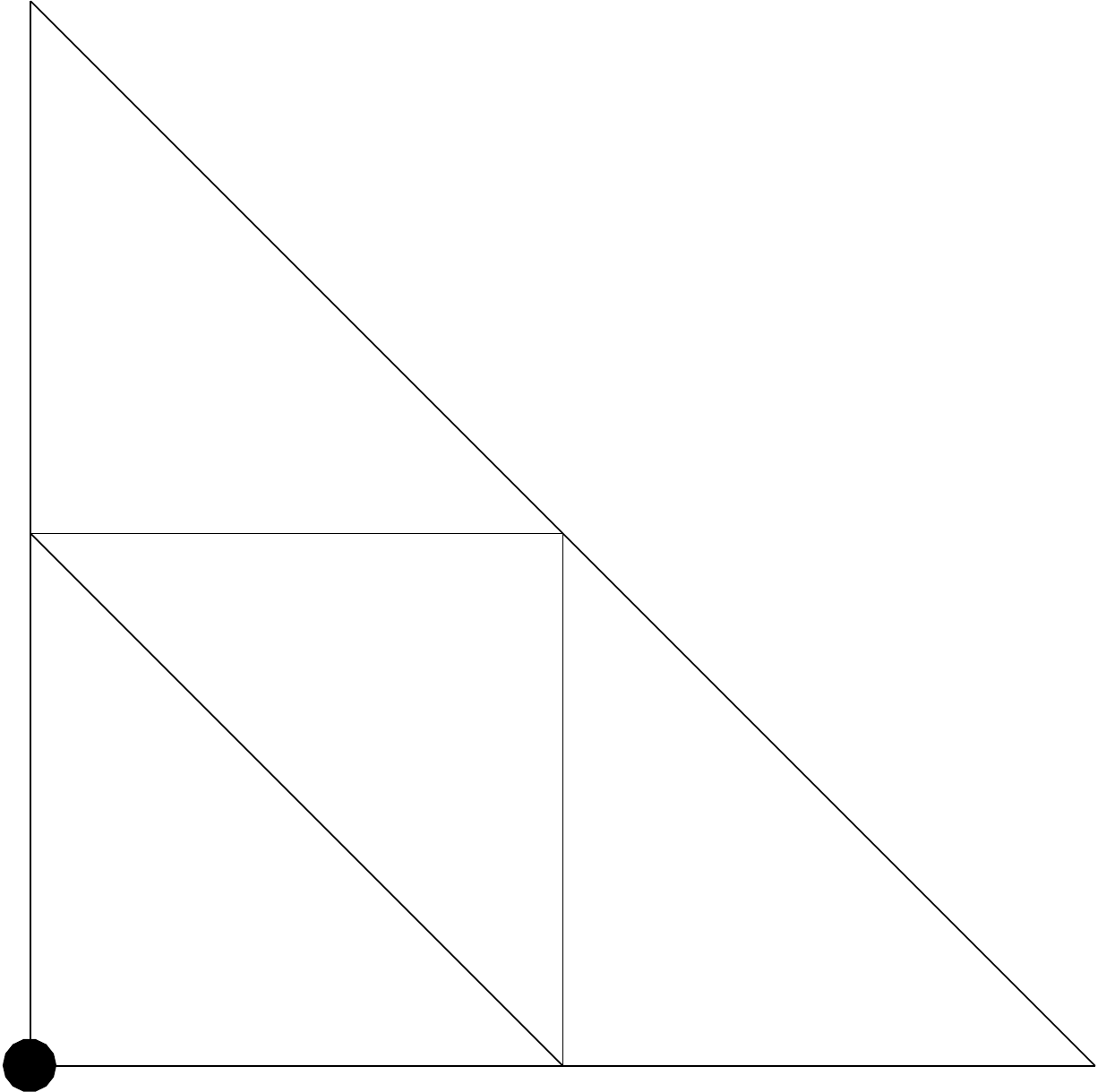} \\ 
      \hline
      $\dfrac{\partial^2 \xproj}{\partial \xi_1 \partial \xi_2} $&
      \vspace{5pt} \includegraphics[width=2.5cm,height=2.5cm,keepaspectratio]{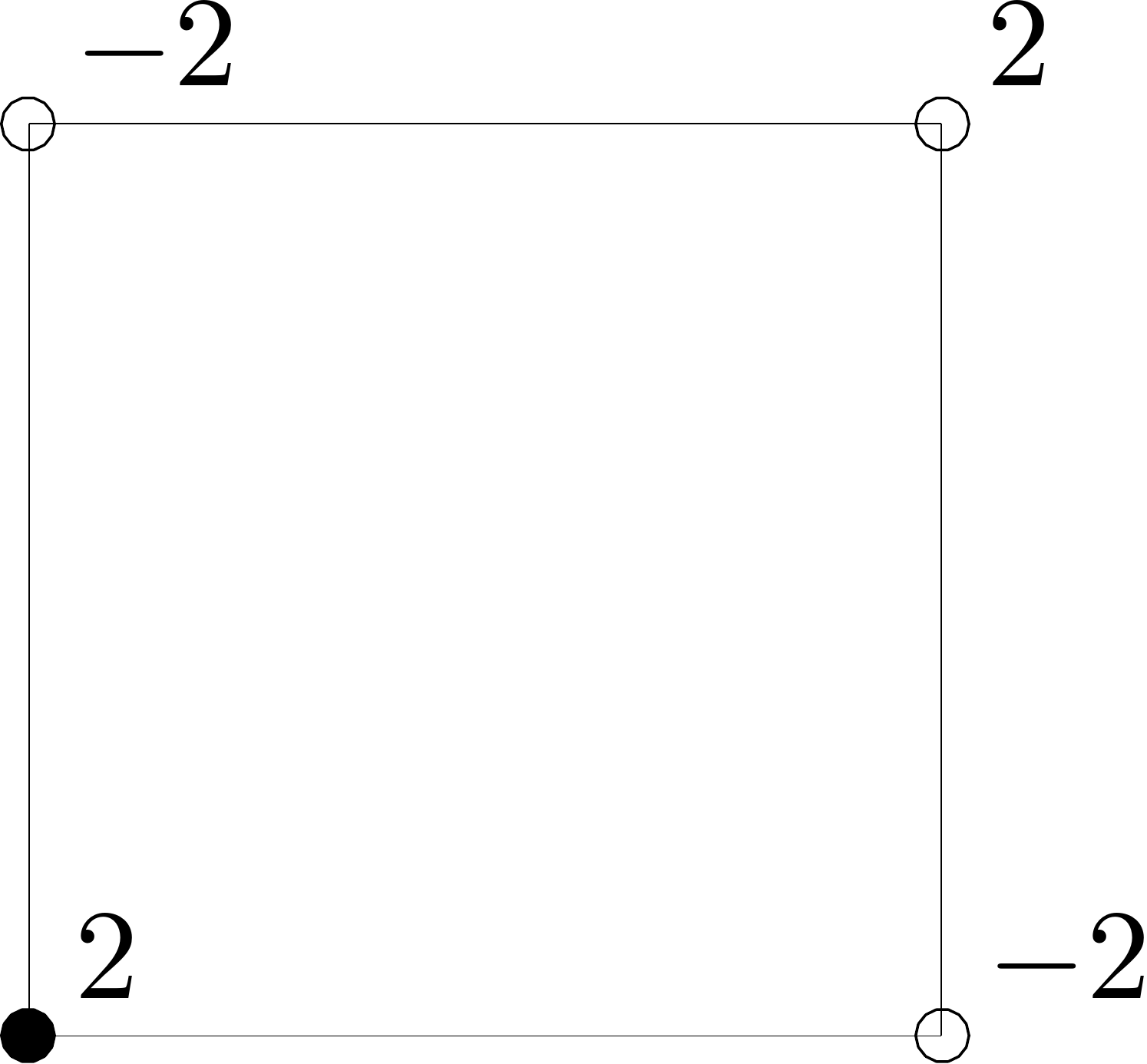} &
      \vspace{5pt}
      \includegraphics[width=2.0cm,height=2.0cm,keepaspectratio]{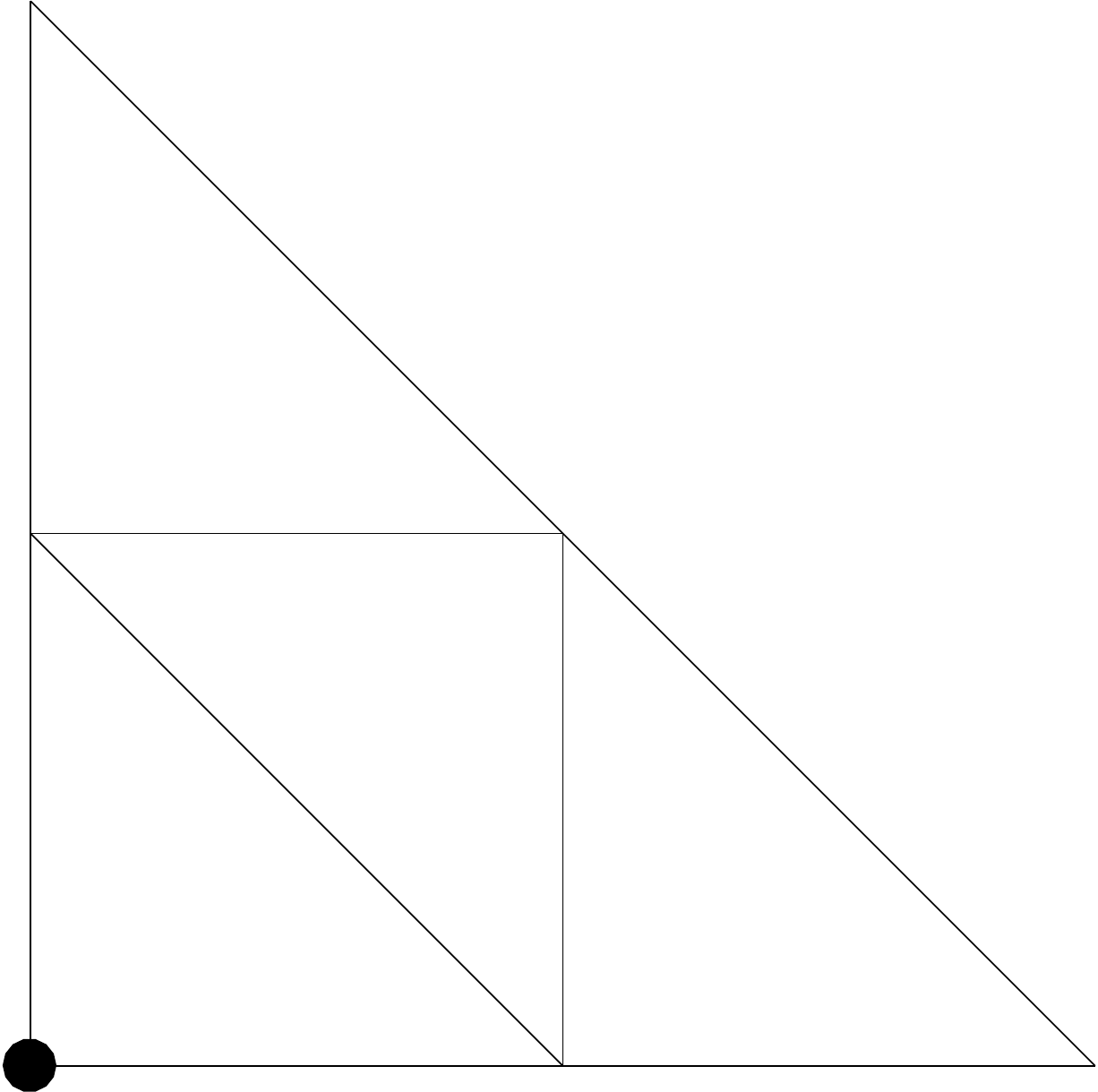} &
      \multicolumn{3}{c|}{}\\ 
      \hline
    \end{tabular}
    \label{tri2_derivsa}
  \end{center}
\end{table}

\begin{table}[H]
  \begin{center}
    \caption{Partial derivative stencils for cubic  \BB triangles.}
    \begin{tabular}{ | >{\centering\arraybackslash} m{1.5cm}  | >{\centering\arraybackslash} m{3cm} |  >{\centering\arraybackslash} m{2.3cm} ||
        >{\centering\arraybackslash} m{1.5cm}  | >{\centering\arraybackslash} m{2.0cm} |  >{\centering\arraybackslash} m{2.3cm} |
      }
      \hline
      \multicolumn{6}{|c|}{  \textbf{ First Order Derivatives, $|\balpha| = 1$}}\\ 
      \hline
      Deriv. & Stencil & Evaluation Triangle & Deriv. & Stencil & Evaluation Triangle  \\ 
      \hline
      $\deldel{\xproj}{\xi_1} $&
      \includegraphics[width=2.5cm,height=2.5cm,keepaspectratio]{dx1dy0_stencil_tri3.pdf} &
      \vspace{5pt}
      \includegraphics[width=2.0cm,height=2.0cm,keepaspectratio]{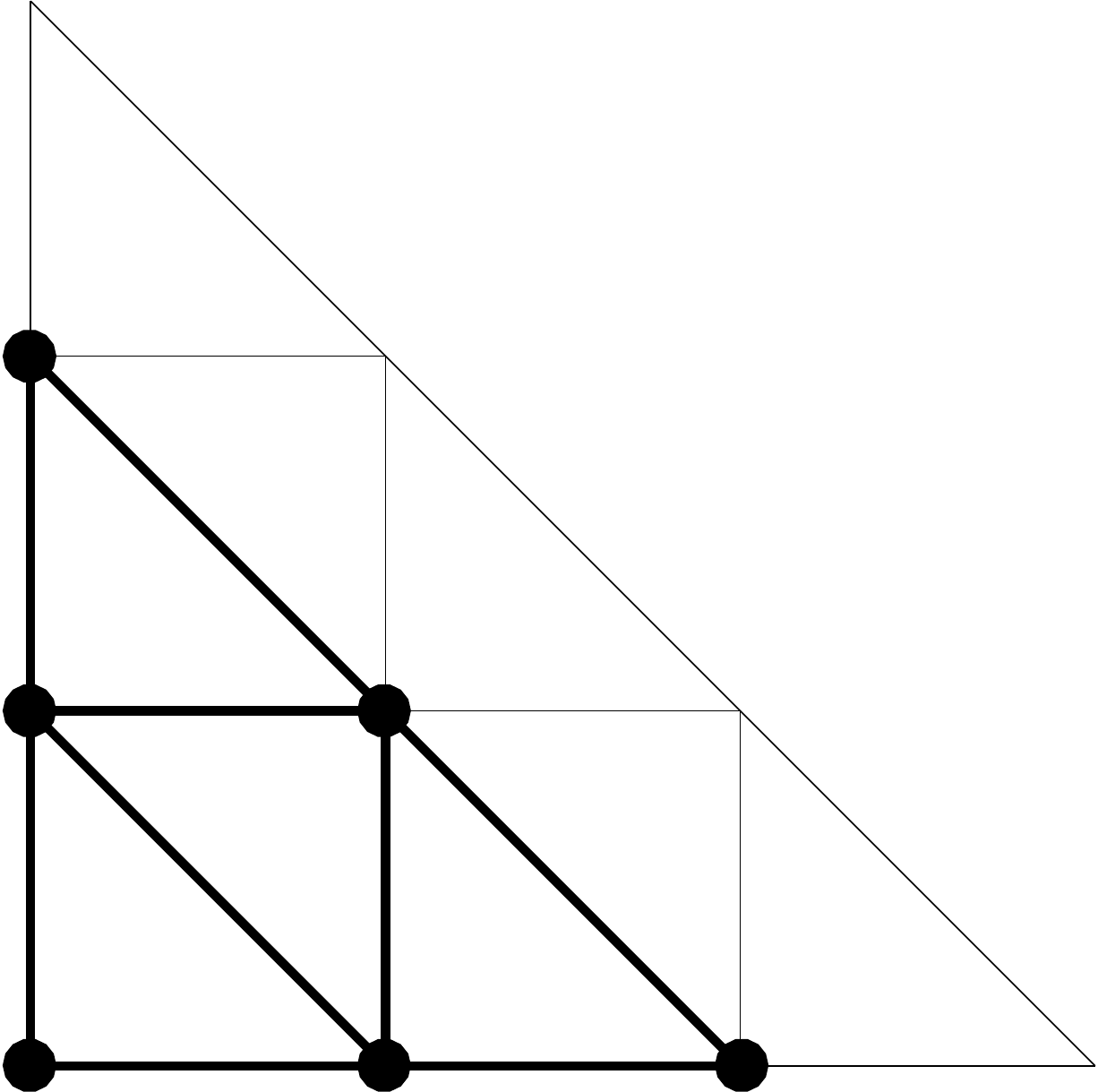} &
      $\deldel{\xproj}{\xi_2} $&
      \includegraphics[width=2.0cm,height=2.0cm,keepaspectratio]{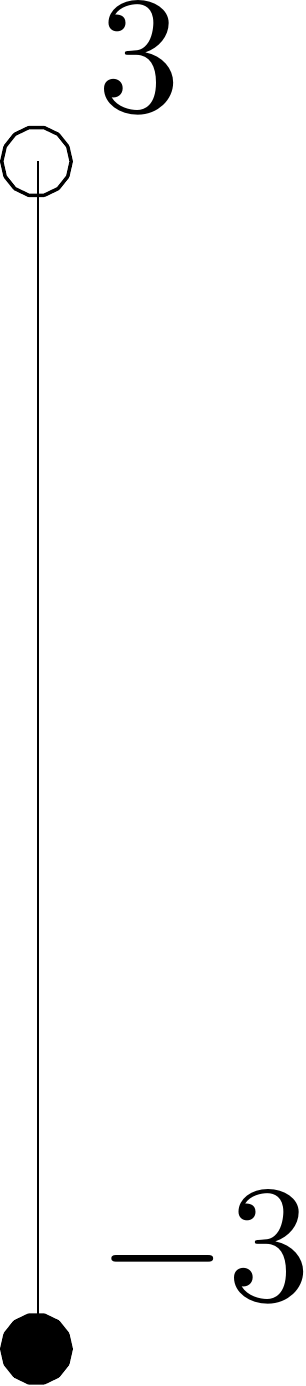} &
      \vspace{5pt}
      \includegraphics[width=2.0cm,height=2.0cm,keepaspectratio]{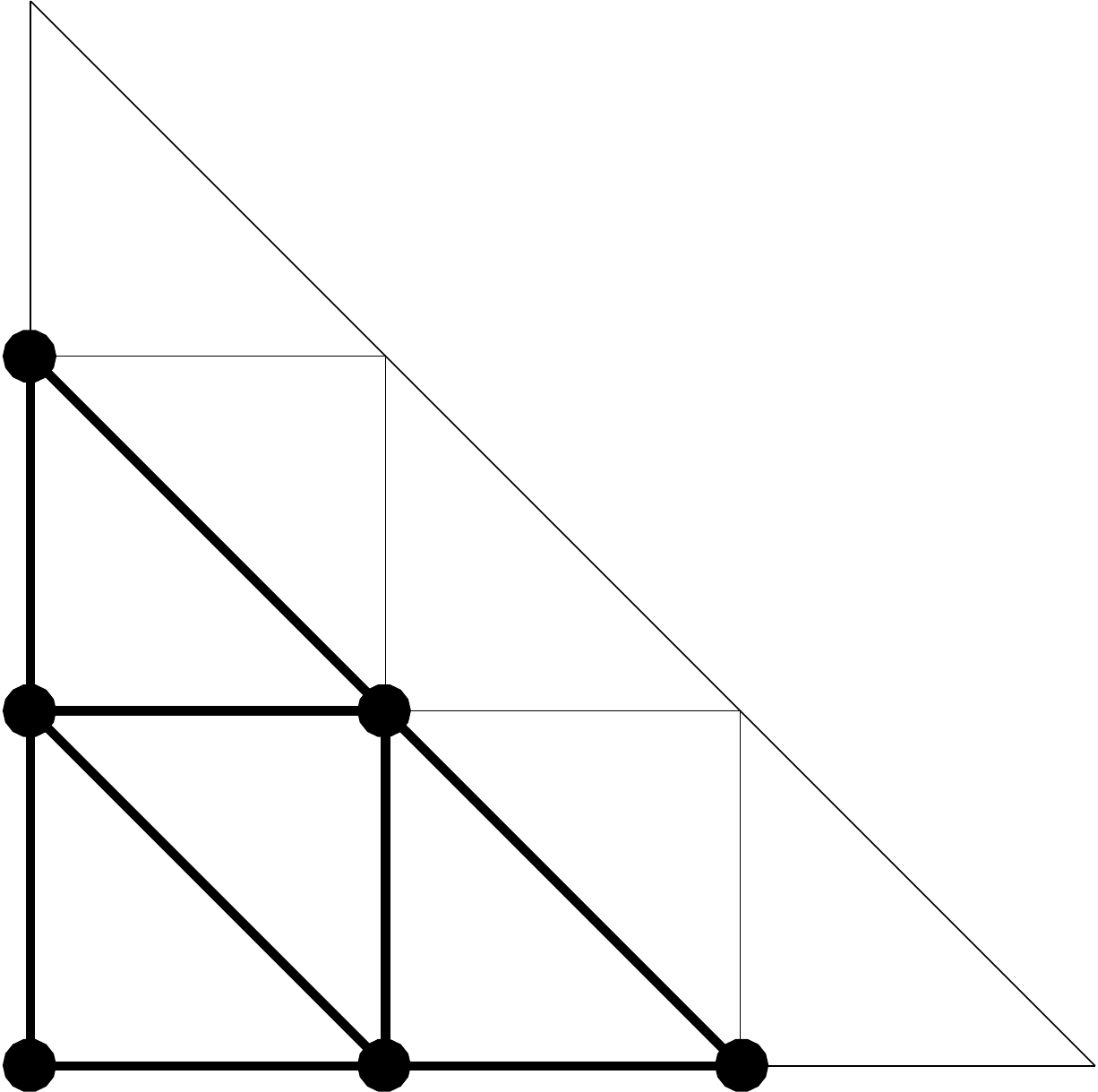} \\
      \hline
      \multicolumn{6}{c}{\vspace{8pt}}\\ 
      \hline
      \multicolumn{6}{|c|}{\textbf{Second Order Derivatives, $|\balpha| = 2$}}\\ 
      \hline
      Deriv. & Stencil & Evaluation Triangle & Deriv. & Stencil & Evaluation Triangle  \\ 
      \hline
      $\dfrac{\partial^2 \xproj}{\partial \xi_1^2} $&
      \includegraphics[width=2.5cm,height=2.5cm,keepaspectratio]{dx2dy0_stencil_tri3.pdf} &
      \vspace{5pt}
      \includegraphics[width=2.0cm,height=2.0cm,keepaspectratio]{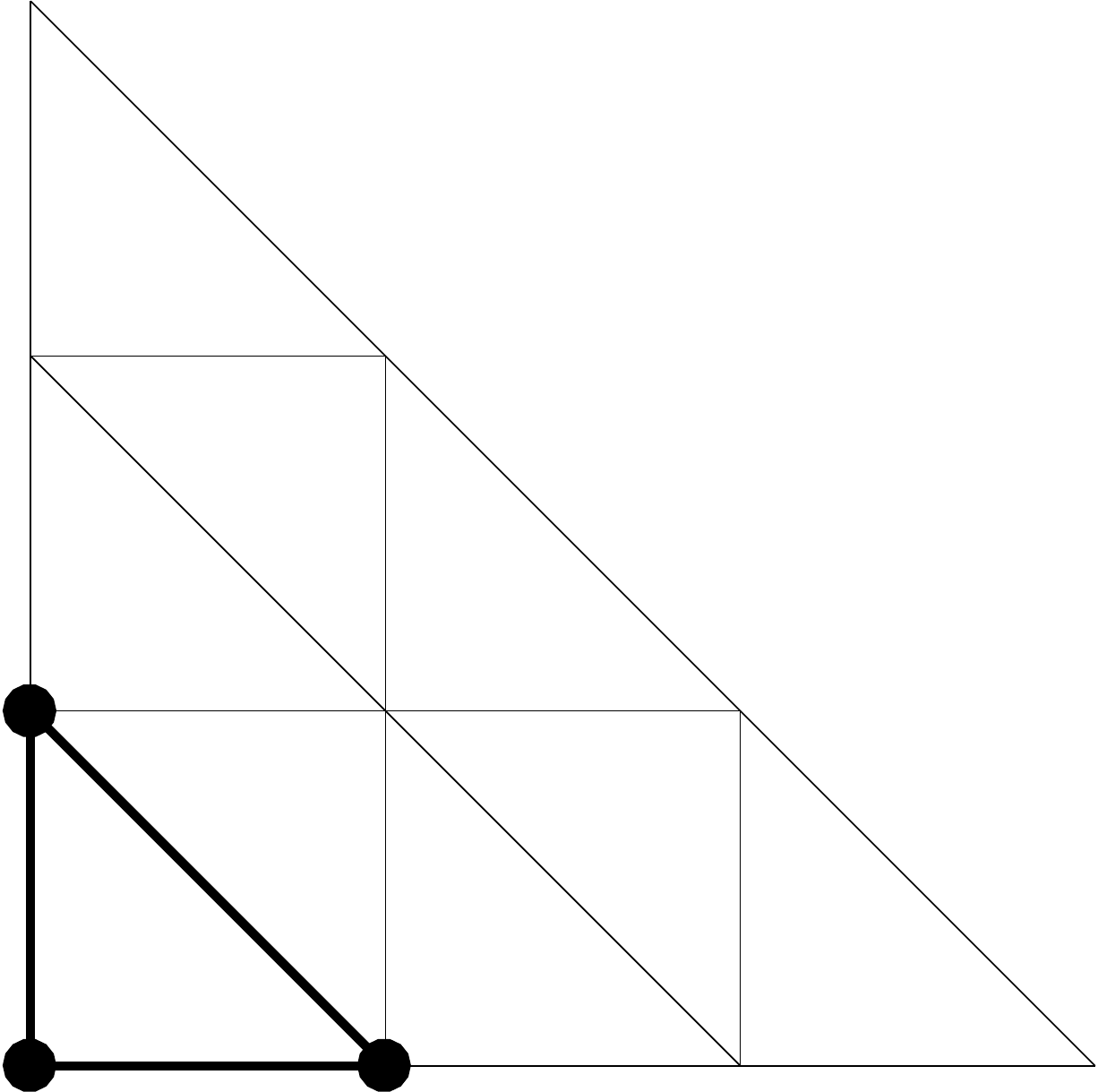} &
      $\dfrac{\partial^2 \xproj}{\partial \xi_2^2} $&
      \includegraphics[width=2.0cm,height=2.0cm,keepaspectratio]{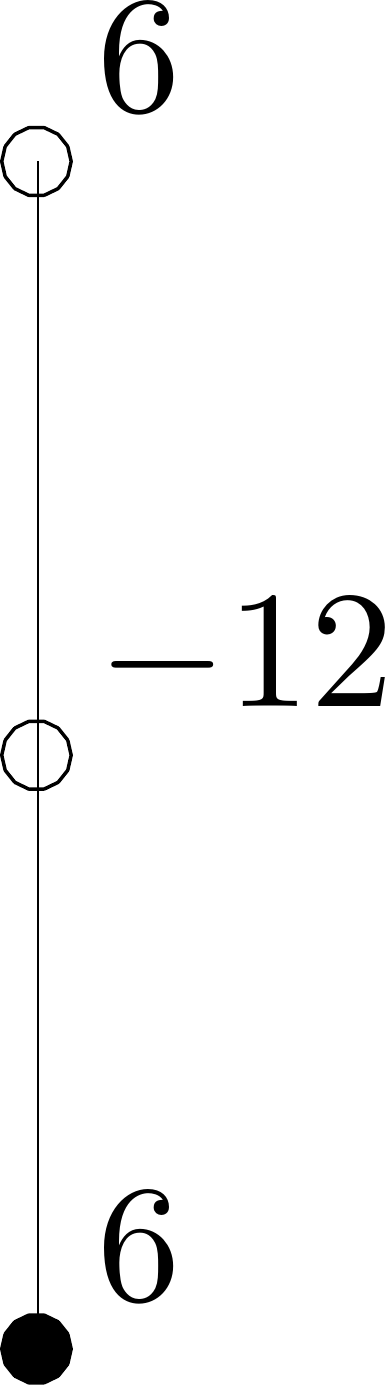} &
      \vspace{5pt}
      \includegraphics[width=2.0cm,height=2.0cm,keepaspectratio]{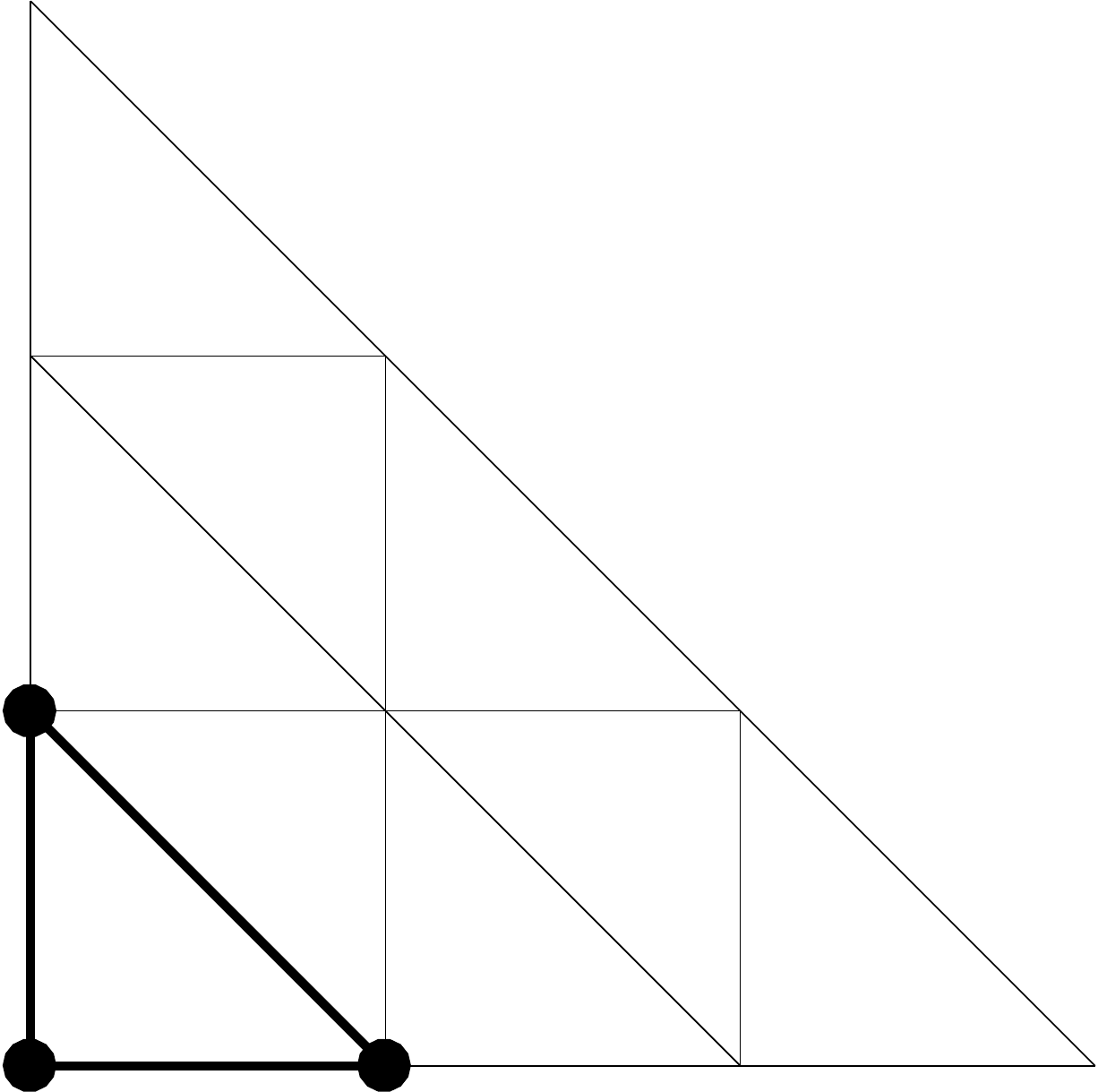} \\ 
      \hline
      $\dfrac{\partial^2 \xproj}{\partial \xi_1 \partial \xi_2} $&
      \vspace{5pt} \includegraphics[width=2.5cm,height=2.5cm,keepaspectratio]{dx1dy1_stencil_tri2.pdf} &
      \vspace{5pt}
      \includegraphics[width=2.0cm,height=2.0cm,keepaspectratio]{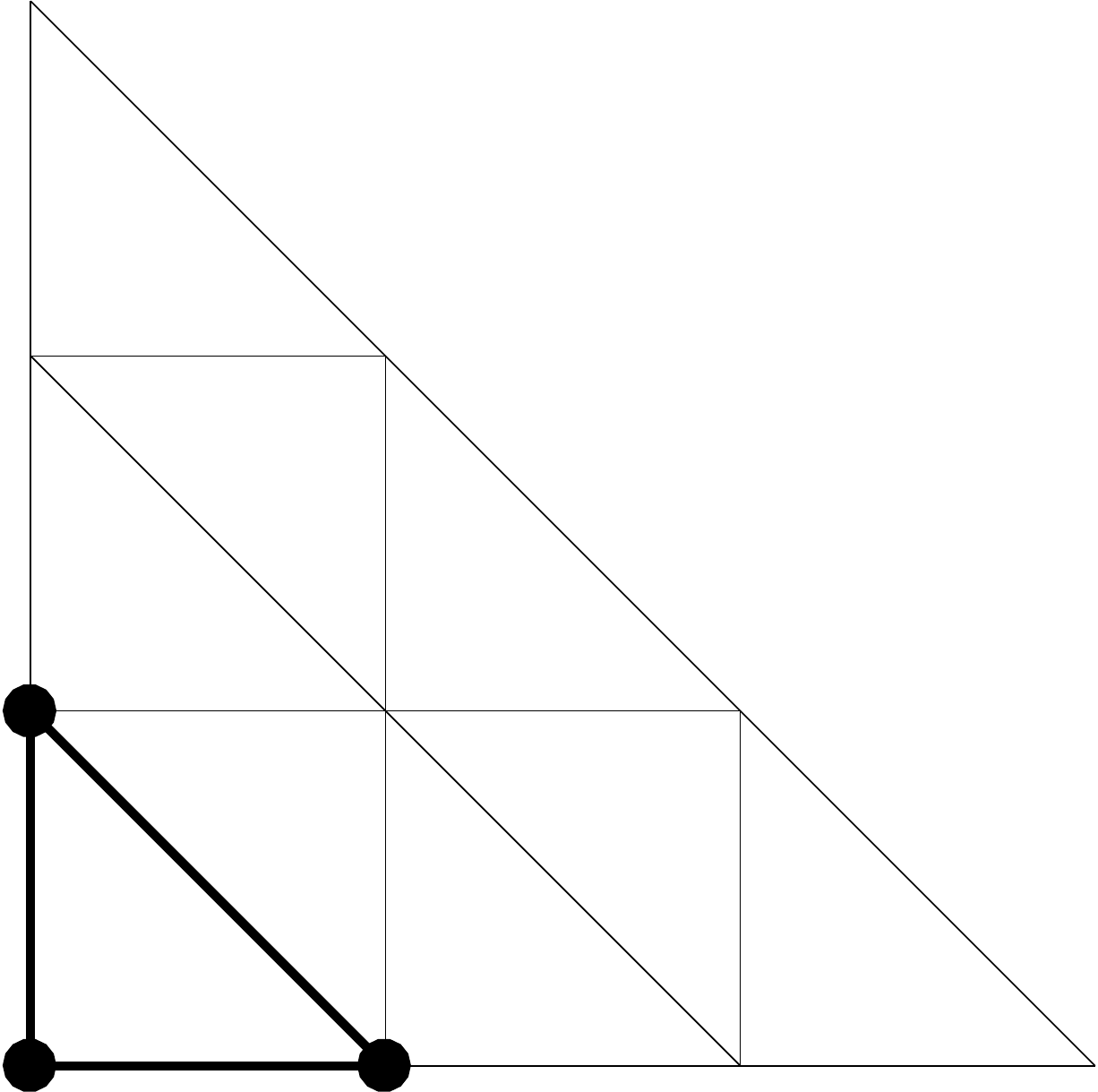} &
      \multicolumn{3}{|c|}{}\\ 
      \hline
      \multicolumn{6}{c}{\vspace{8pt}}\\ 
      \hline
      \multicolumn{6}{|c|}{\textbf{Third Order Derivatives, $|\balpha| = 3$}}\\ 
      \hline
      Deriv. & Stencil & Evaluation Triangle & Deriv. & Stencil & Evaluation Triangle  \\ 
      \hline
      $\dfrac{\partial^3 \xproj}{\partial \xi_1^3} $&
      \includegraphics[width=2.5cm,height=2.5cm,keepaspectratio]{dx3dy0_stencil_tri3.pdf} &
      \vspace{5pt}
      \includegraphics[width=2.0cm,height=2.0cm,keepaspectratio]{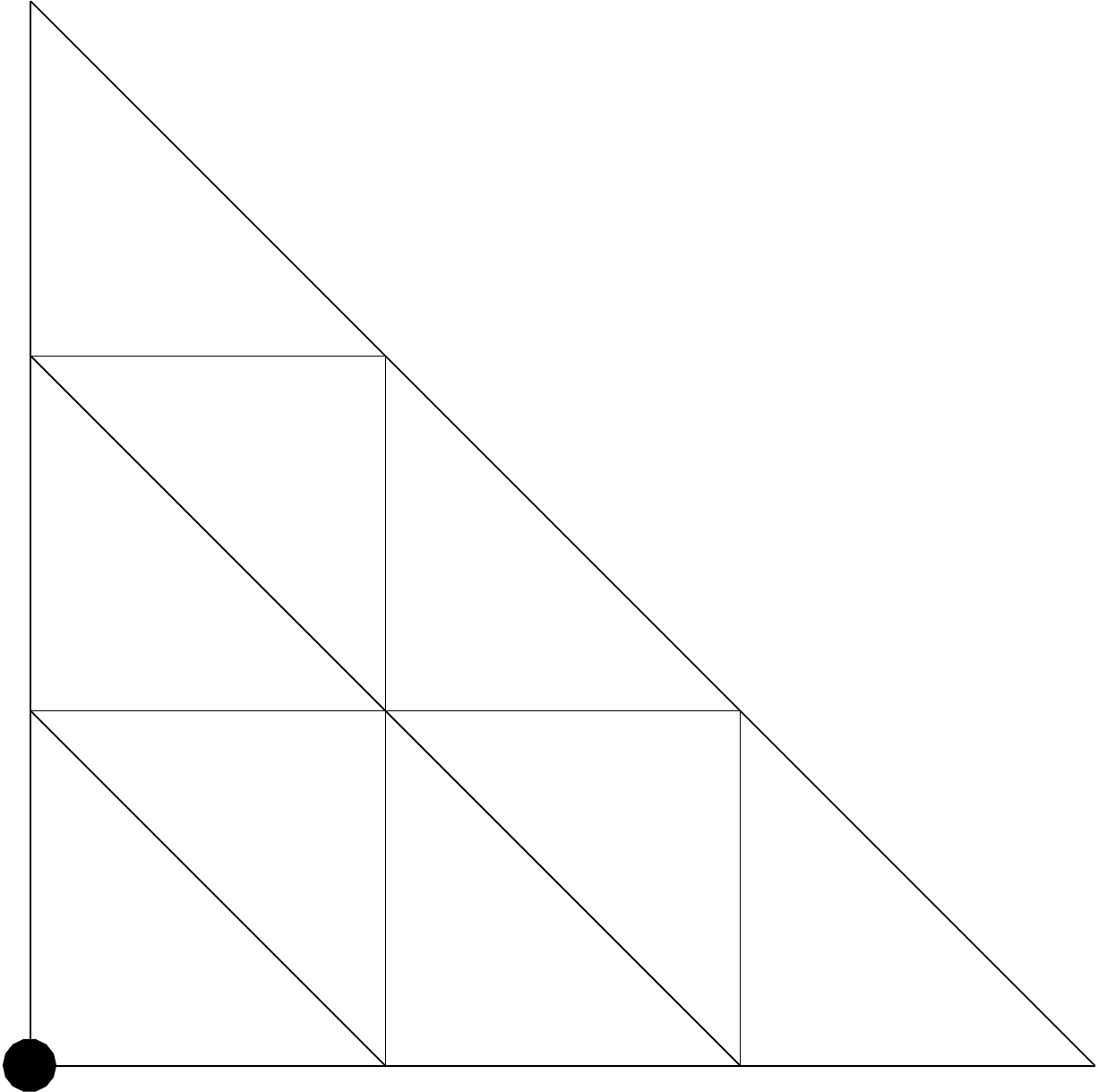} &
      $\dfrac{\partial^3 \xproj}{\partial \xi_2^3} $&
      \includegraphics[width=2.0cm,height=2.0cm,keepaspectratio]{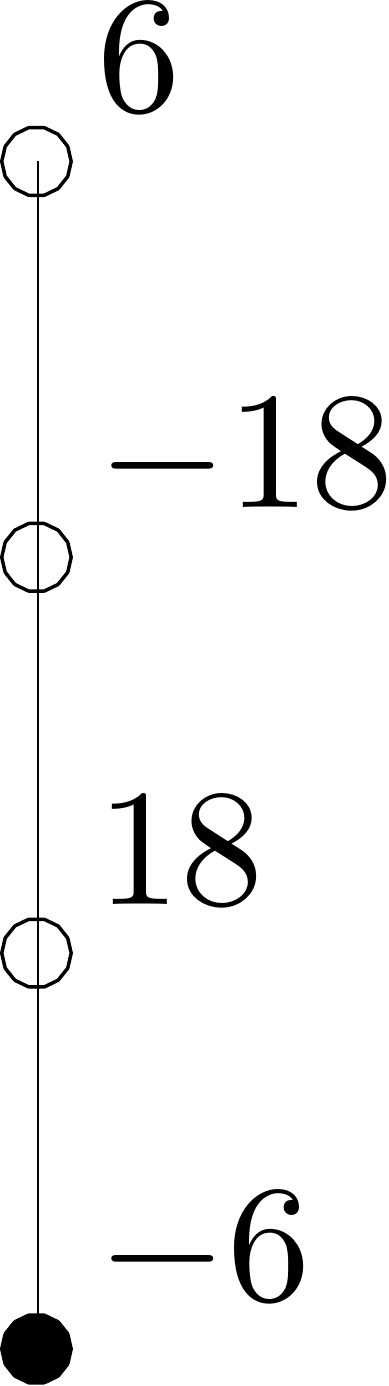} &
      \vspace{5pt}
      \includegraphics[width=2.0cm,height=2.0cm,keepaspectratio]{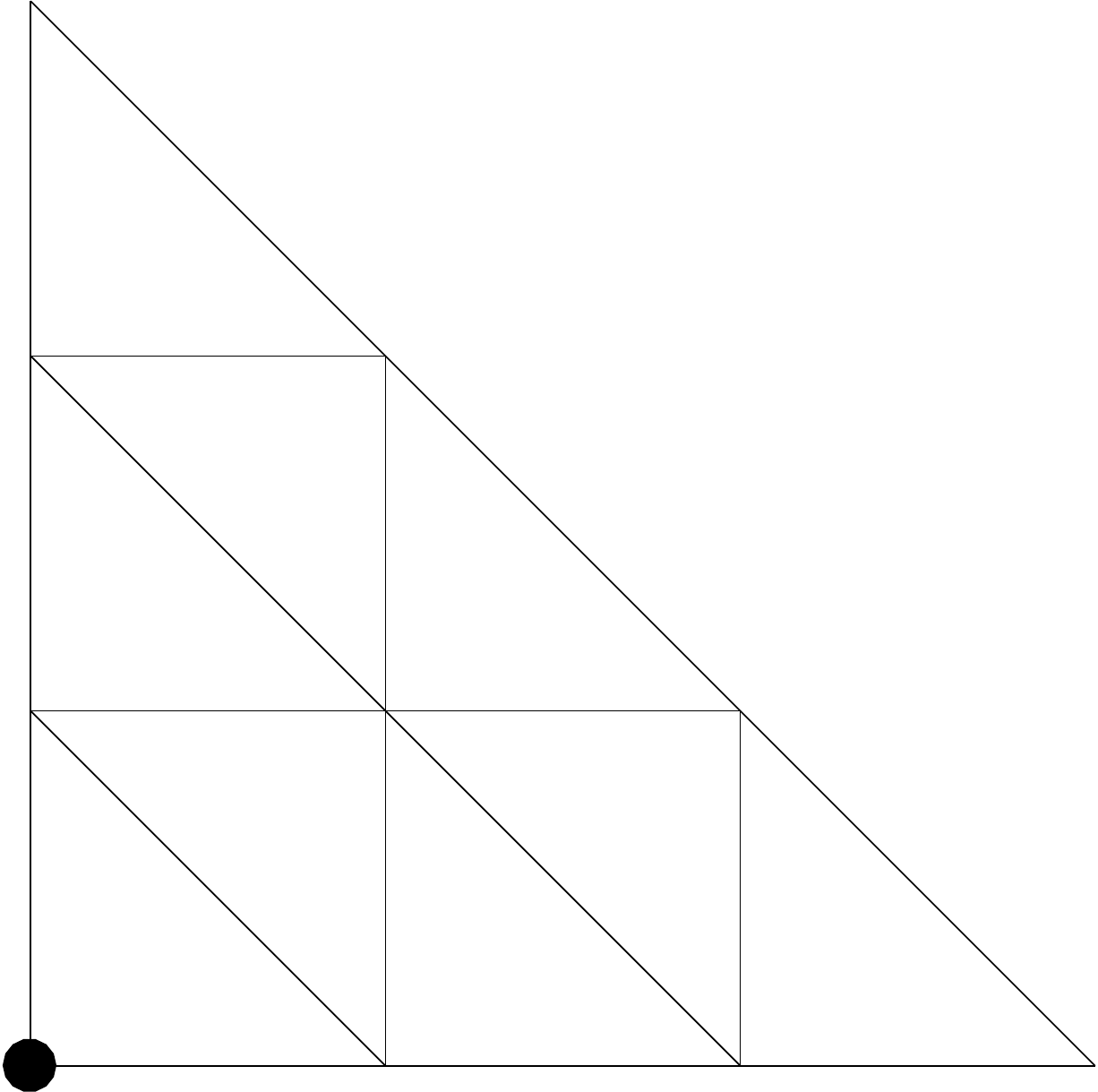} \\ 
      \hline
      $\dfrac{\partial^2 \xproj}{\partial \xi_1^2 \partial \xi_2} $&
      \vspace{5pt} \includegraphics[width=2.5cm,height=2.5cm,keepaspectratio]{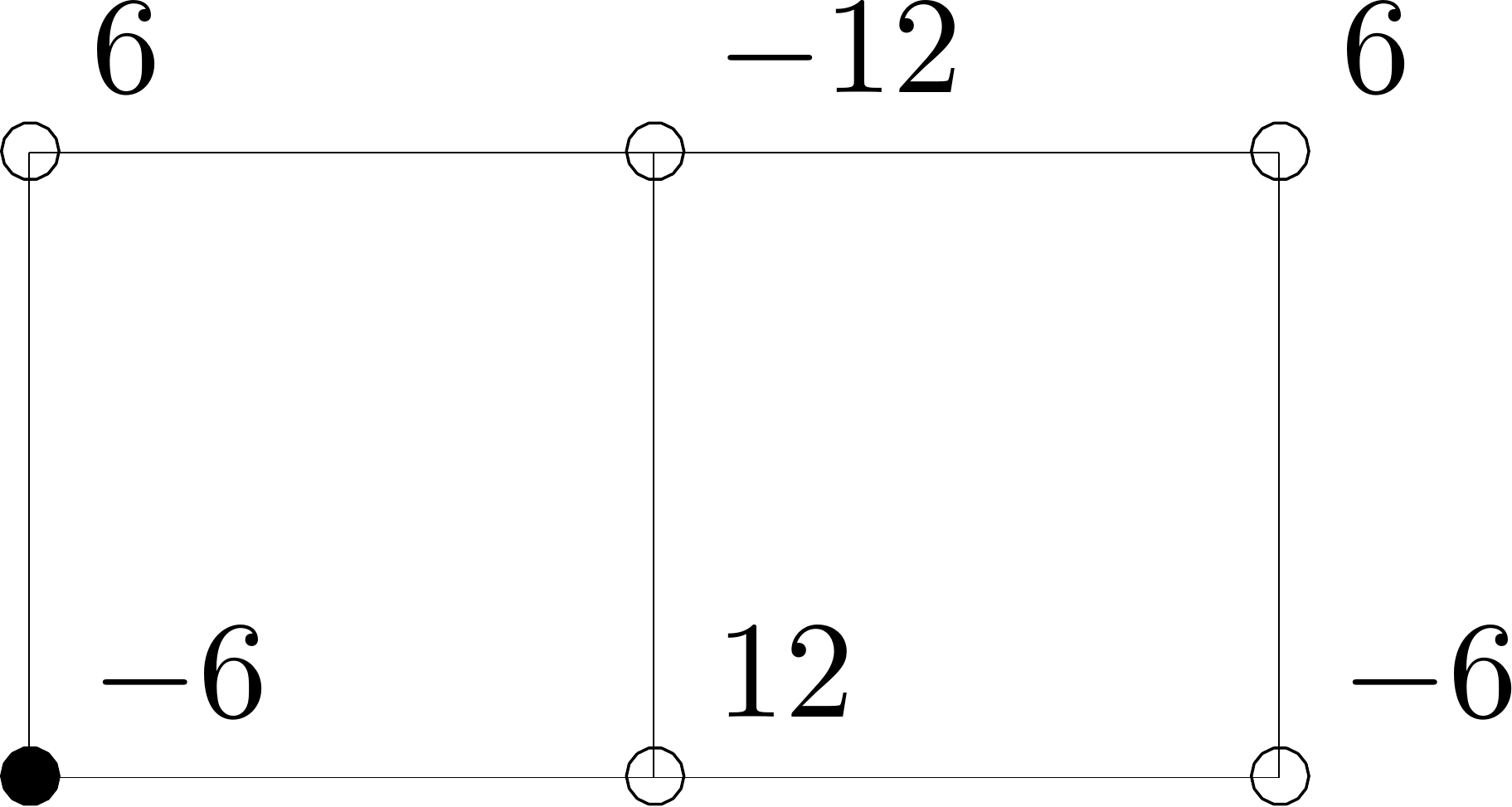} &
      \vspace{5pt}
      \includegraphics[width=2.0cm,height=2.0cm,keepaspectratio]{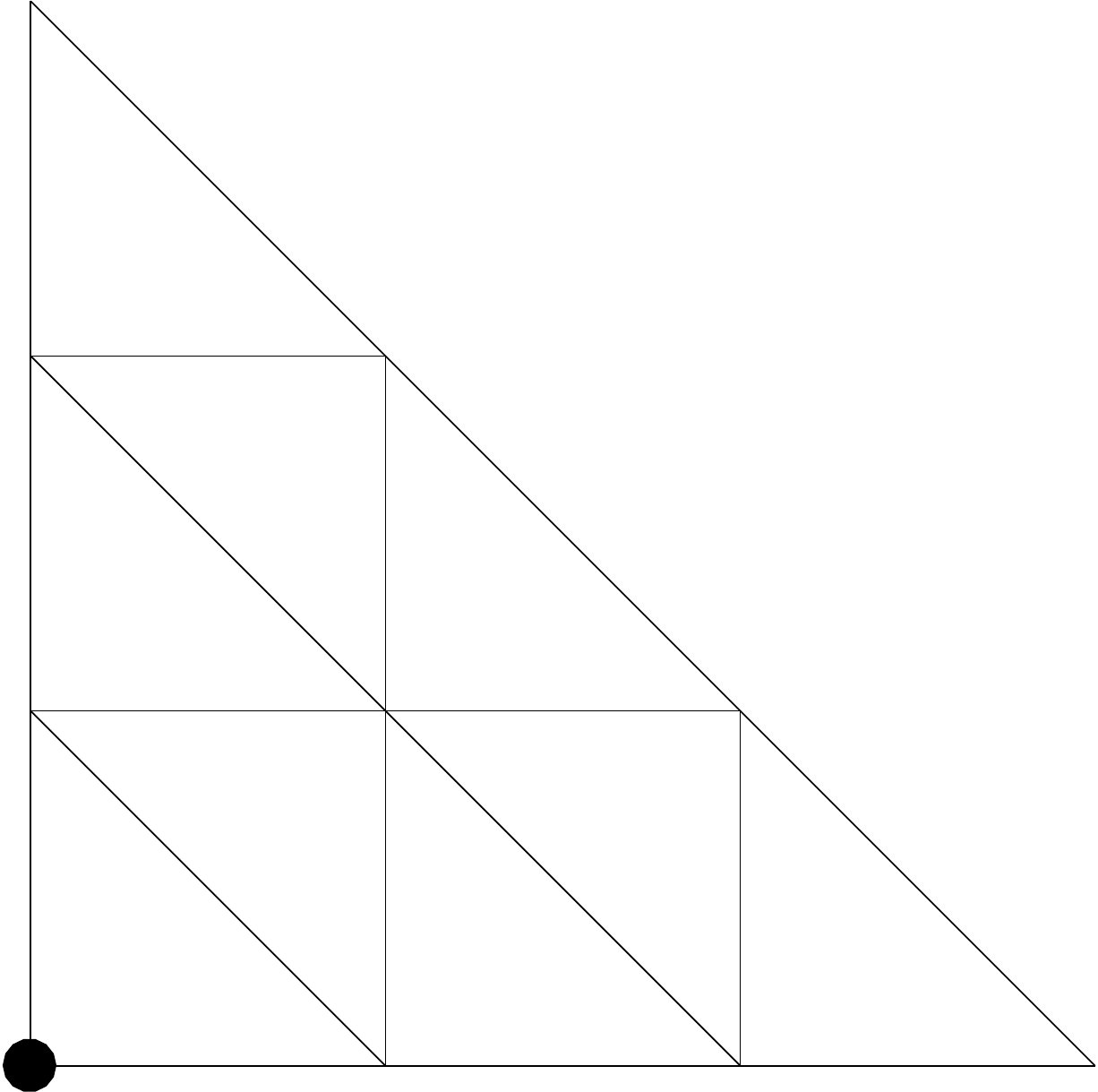} &
      $\dfrac{\partial^2 \xproj}{\partial \xi_1 \partial \xi_2^2} $&
      \vspace{5pt} \includegraphics[width=2.5cm,height=2.5cm,keepaspectratio]{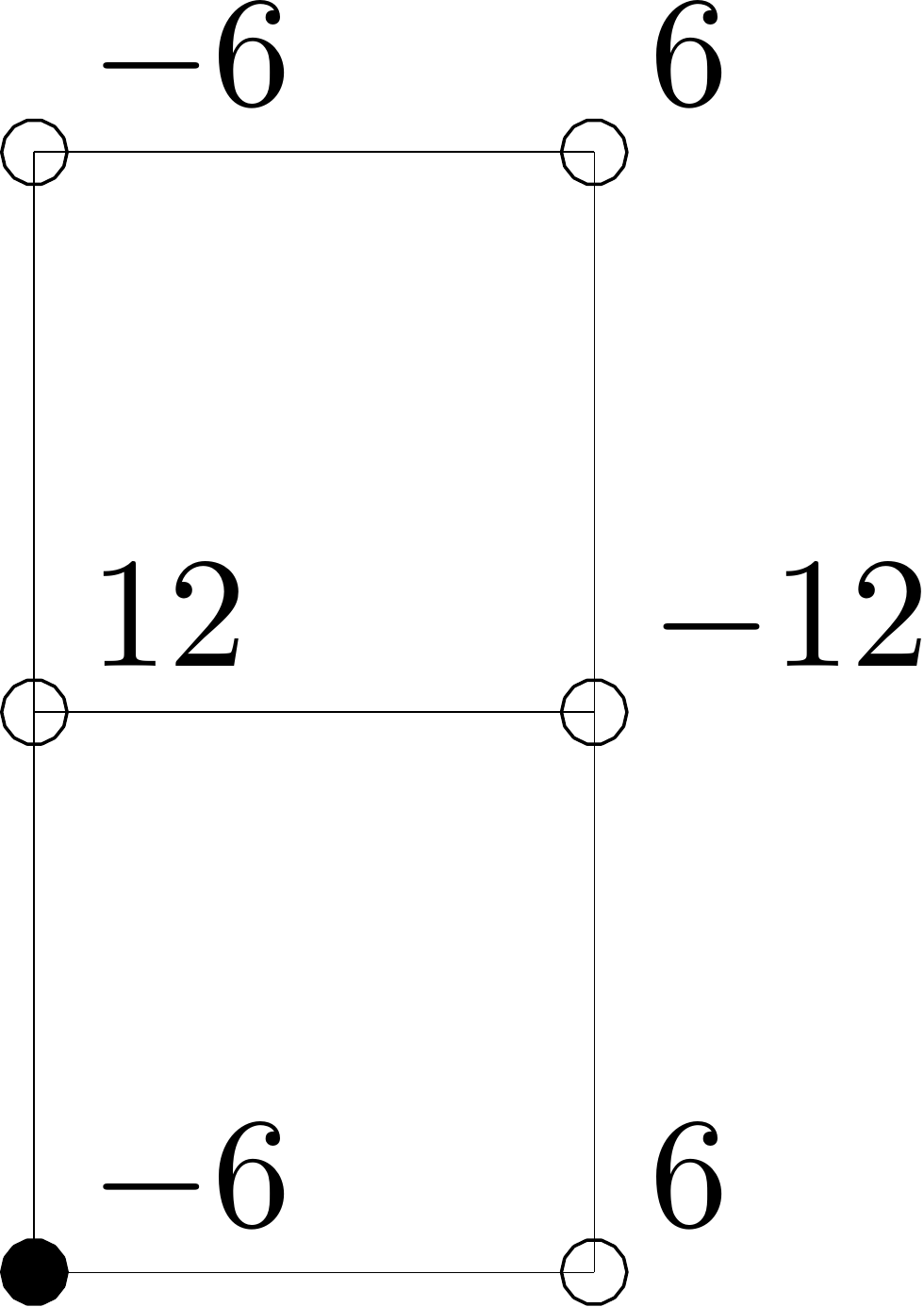} &
      \vspace{5pt}
      \includegraphics[width=2.0cm,height=2.0cm,keepaspectratio]{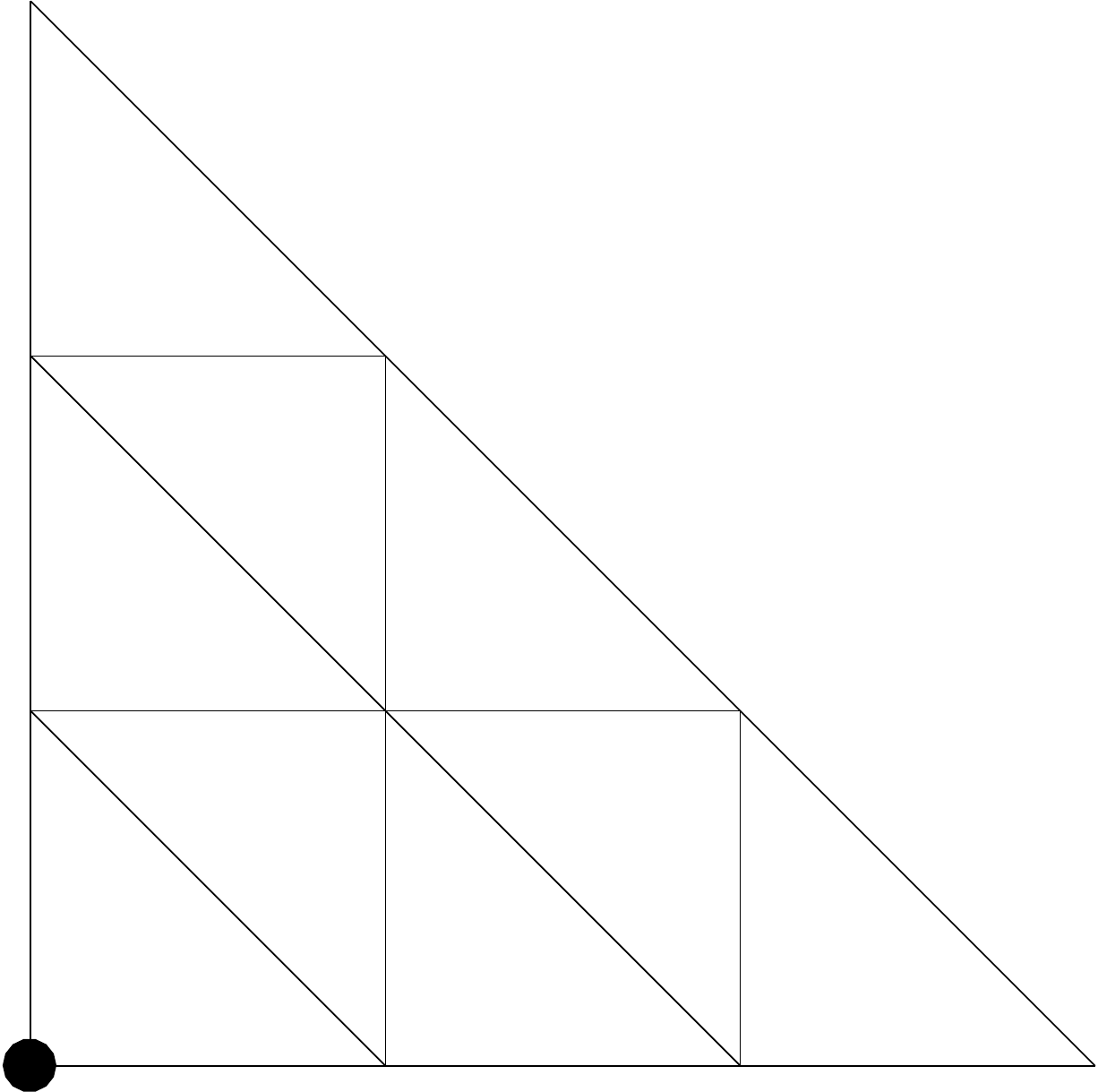}\\ 
      \hline
    \end{tabular}
    \label{tri2_derivsa}
  \end{center}
\end{table}

\begin{table}[H]
  \begin{center}
    \caption{Partial derivative stencils for bi-quadratic  \BB quadrilaterals.}
    \begin{tabular}{ | >{\centering\arraybackslash} m{1.5cm}  | >{\centering\arraybackslash} m{3cm} |  >{\centering\arraybackslash} m{2.3cm} ||
        >{\centering\arraybackslash} m{1.5cm}  | >{\centering\arraybackslash} m{2.0cm} |  >{\centering\arraybackslash} m{2.3cm} |
      }
      \hline
      \multicolumn{6}{|c|}{  \textbf{ First Order Derivatives, $|\balpha| = 1$}}\\ 
      \hline
      Deriv. & Stencil & Evaluation Quadrilateral & Deriv. & Stencil & Evaluation Quadrilateral  \\ 
      \hline
      $\deldel{\xproj}{\xi_1} $&
      \includegraphics[width=2.5cm,height=2.5cm,keepaspectratio]{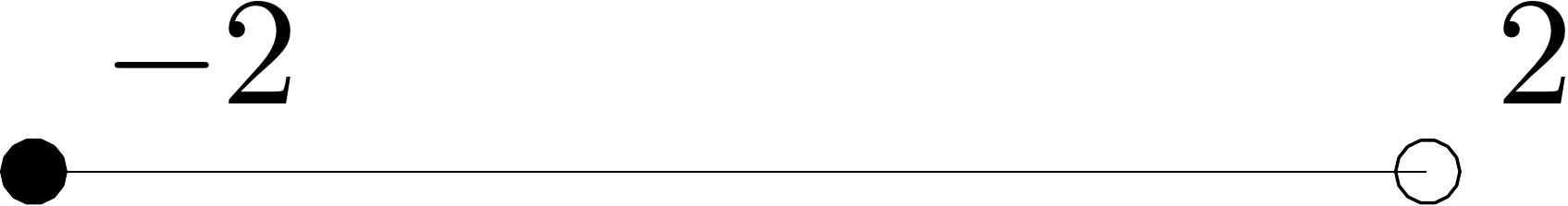} &
      \vspace{5pt}
      \includegraphics[width=2.0cm,height=2.0cm,keepaspectratio]{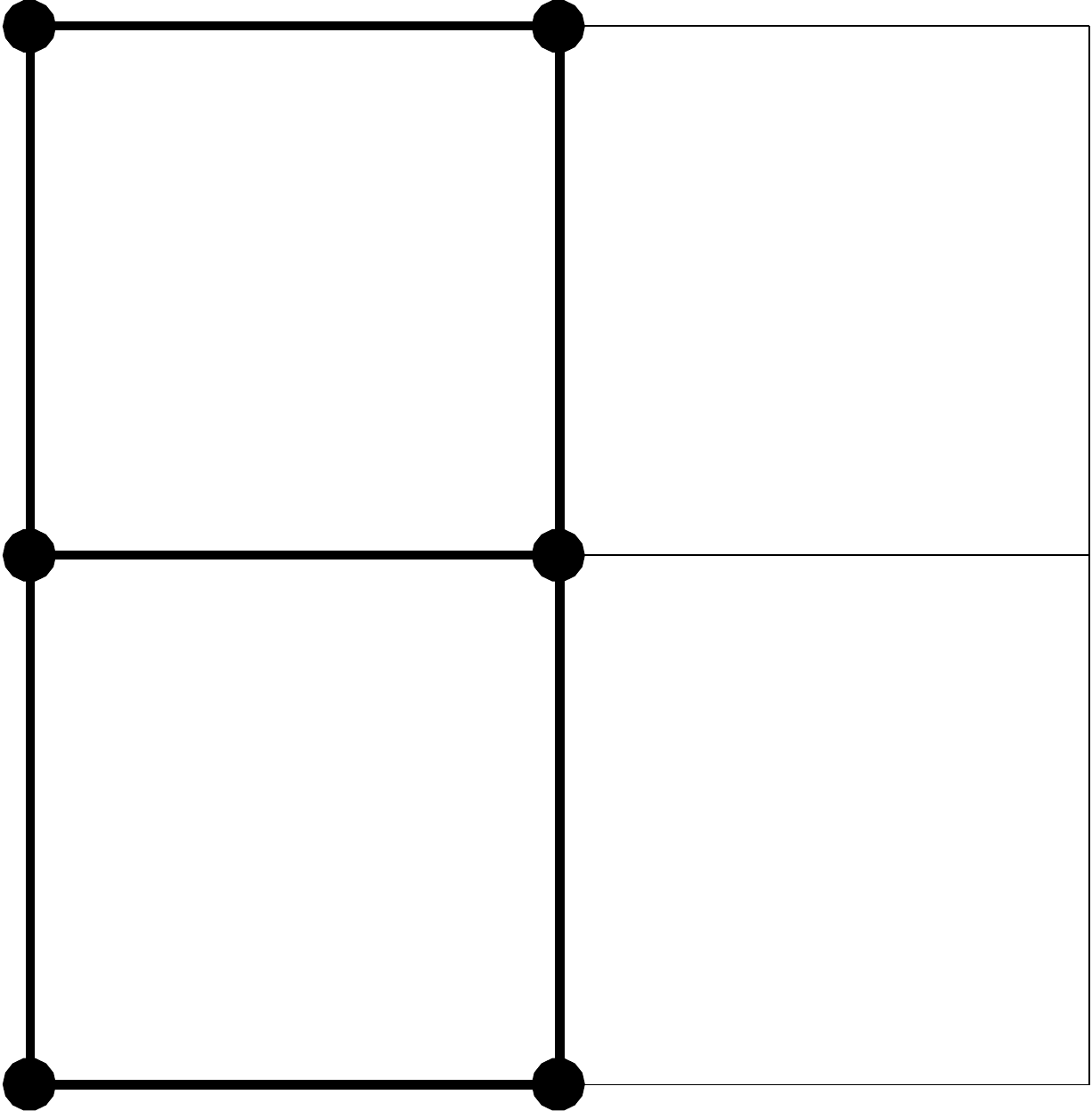} &
      $\deldel{\xproj}{\xi_2} $&
      \includegraphics[width=2.0cm,height=2.0cm,keepaspectratio]{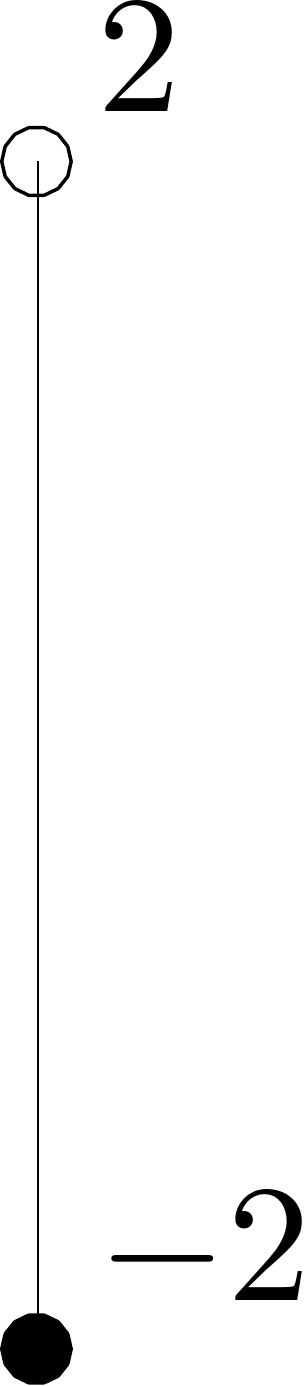} &
      \vspace{5pt}
      \includegraphics[width=2.0cm,height=2.0cm,keepaspectratio]{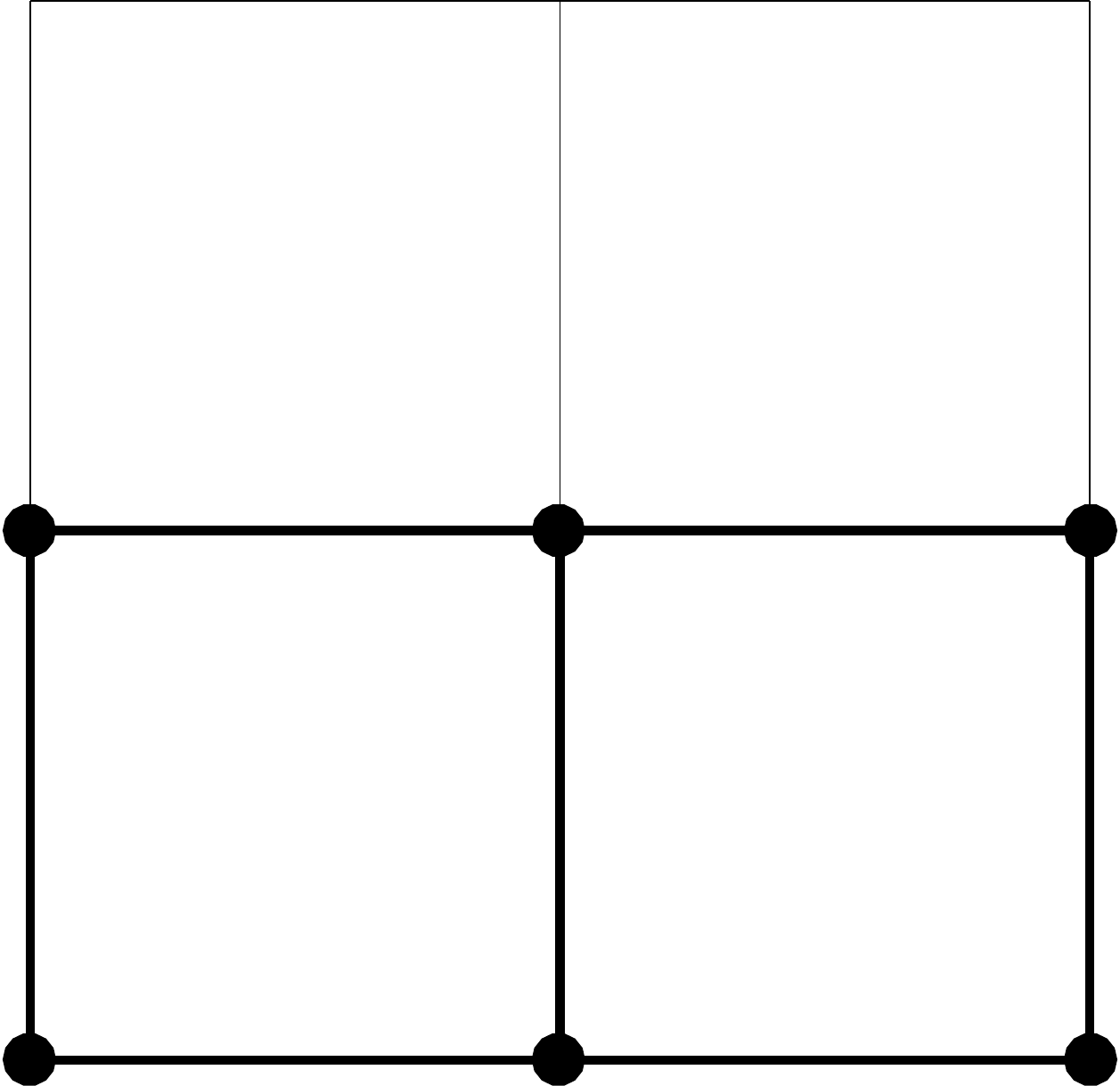} \\
      \hline
      \multicolumn{6}{c}{\vspace{8pt}}\\ 
      \hline
      \multicolumn{6}{|c|}{\textbf{Second Order Derivatives, $|\balpha| = 2$}}\\ 
      \hline
      Deriv. & Stencil & Evaluation Quadrilateral & Deriv. & Stencil & Evaluation Quadrilateral  \\ 
      \hline
      $\dfrac{\partial^2 \xproj}{\partial \xi_1^2} $&
      \includegraphics[width=2.5cm,height=2.5cm,keepaspectratio]{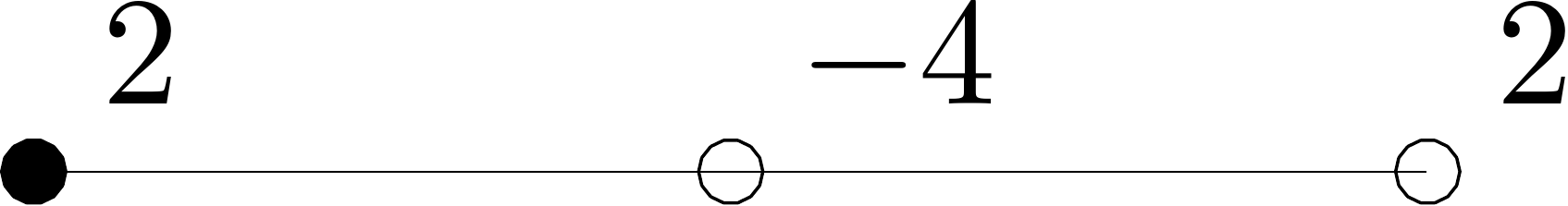} &
      \vspace{5pt}
      \includegraphics[width=2.0cm,height=2.0cm,keepaspectratio]{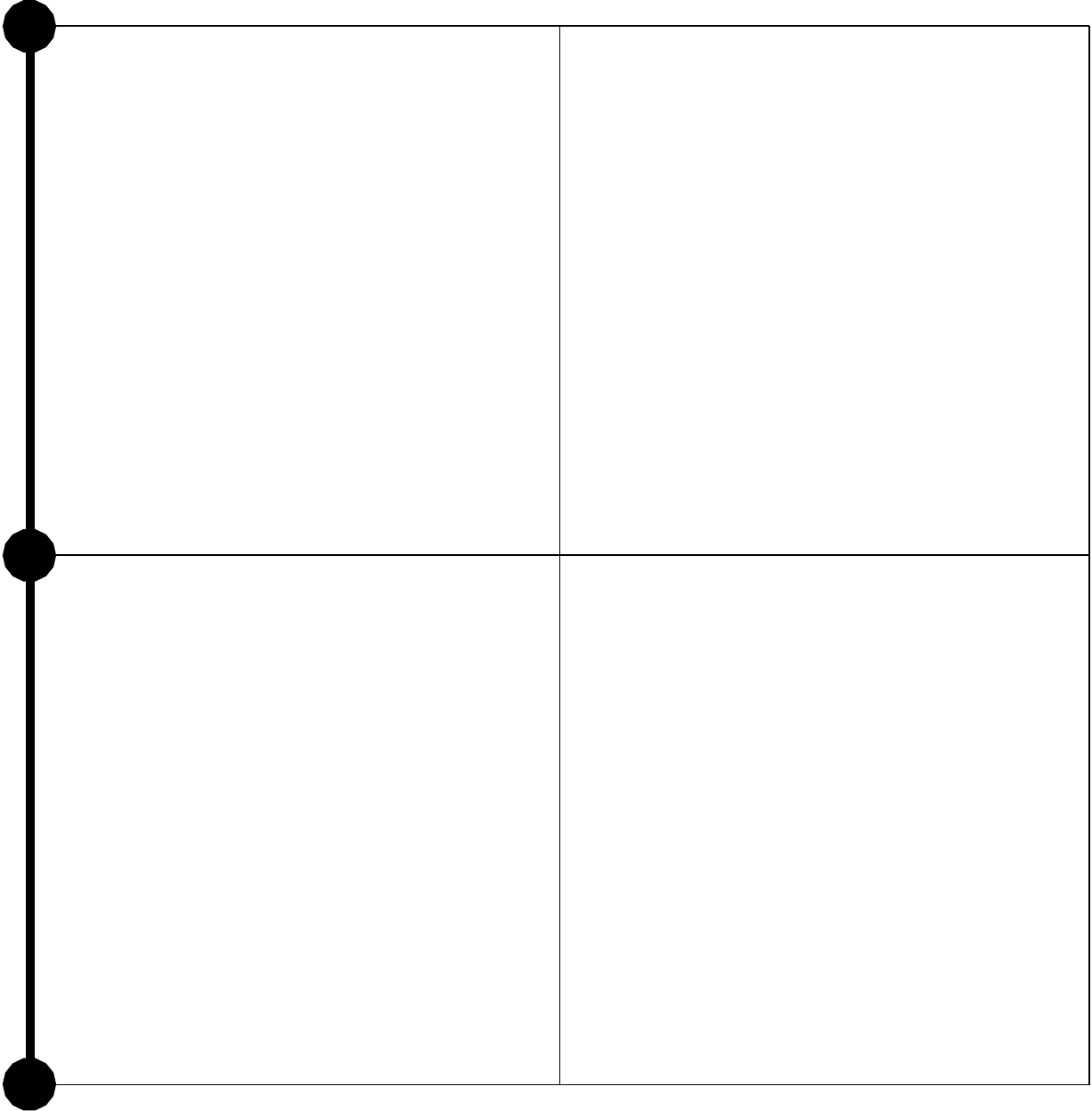} &
      $\dfrac{\partial^2 \xproj}{\partial \xi_2^2} $&
      \includegraphics[width=2.0cm,height=2.0cm,keepaspectratio]{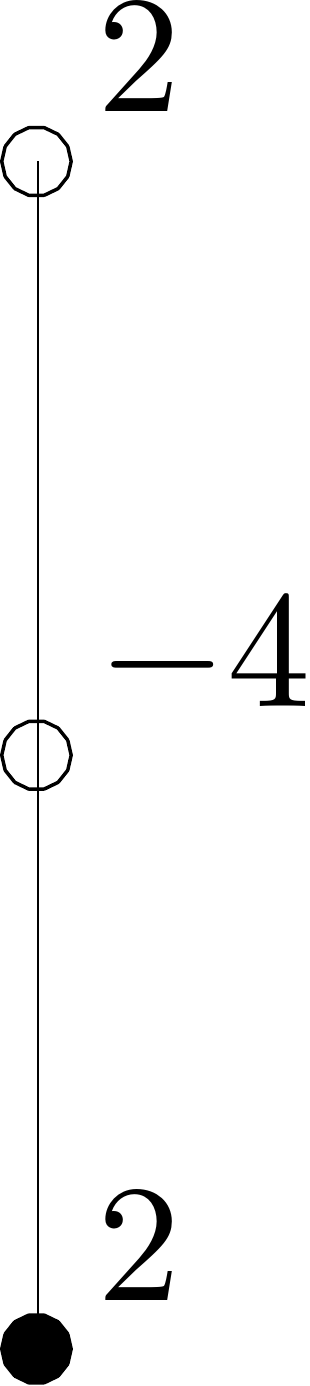} &
      \vspace{5pt}
      \includegraphics[width=2.0cm,height=2.0cm,keepaspectratio]{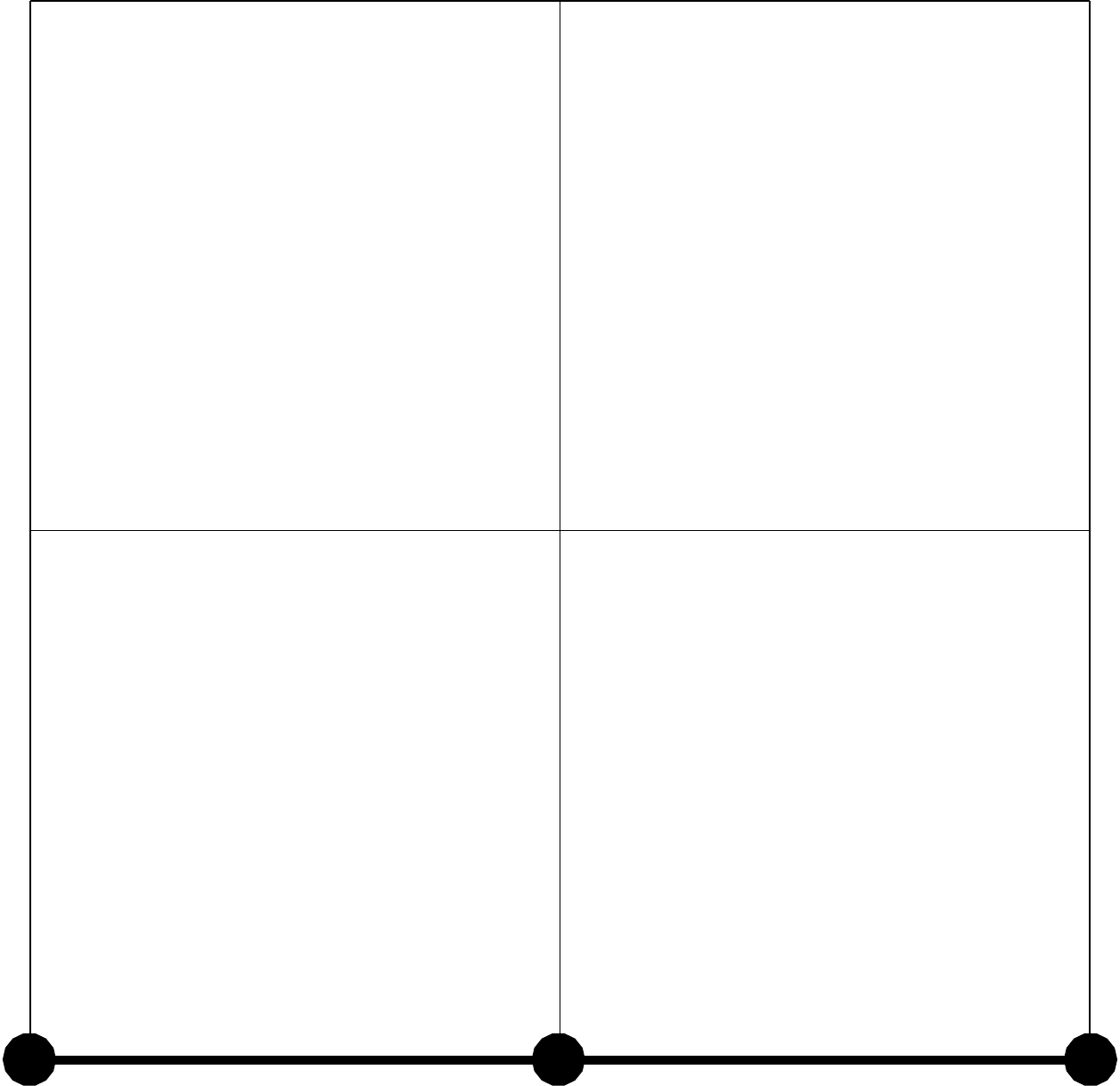} \\ 
      \hline
      $\dfrac{\partial^2 \xproj}{\partial \xi_1 \partial \xi_2} $&
      \vspace{5pt} \includegraphics[width=2.5cm,height=2.5cm,keepaspectratio]{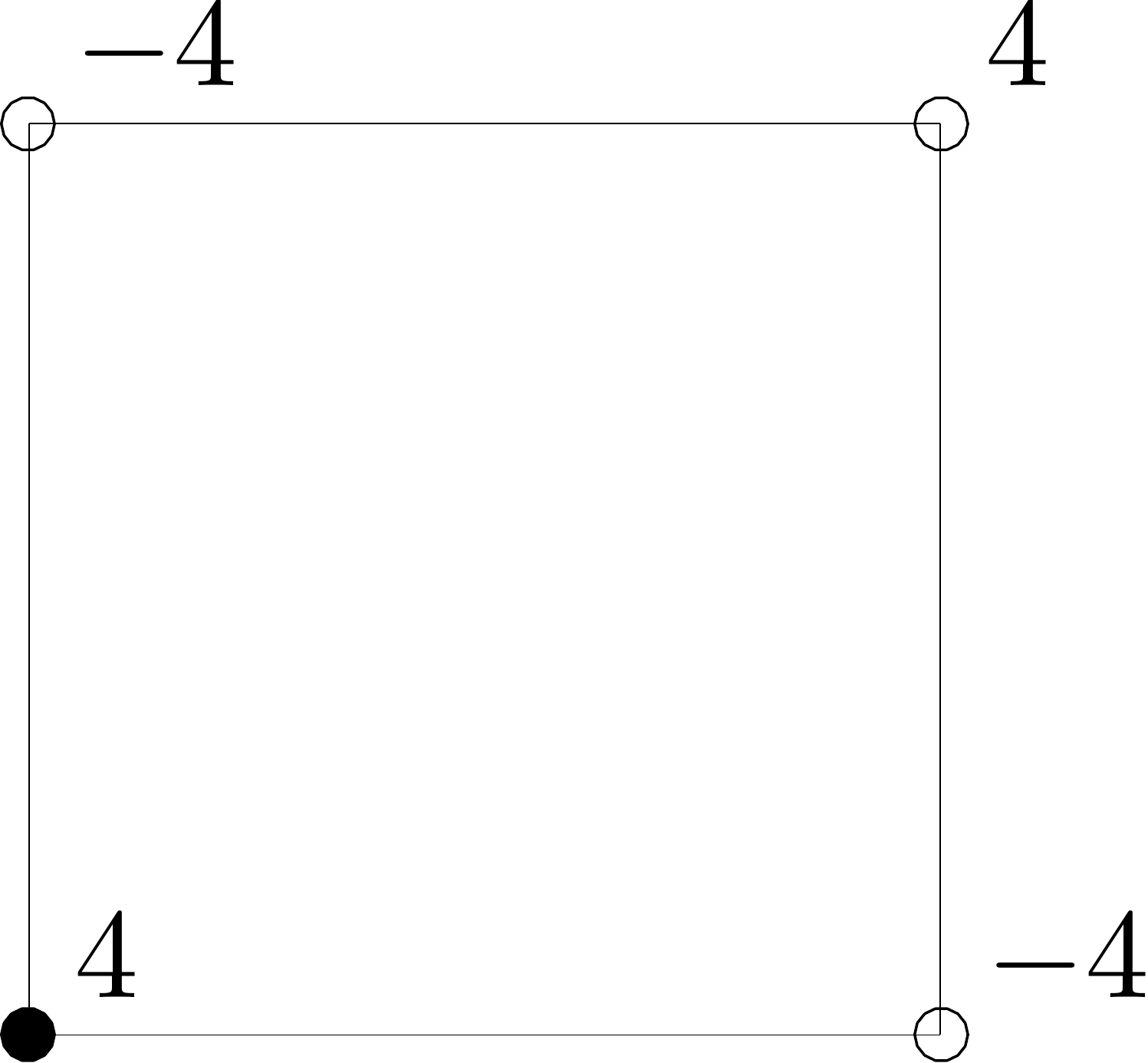} &
      \vspace{5pt}
      \includegraphics[width=2.0cm,height=2.0cm,keepaspectratio]{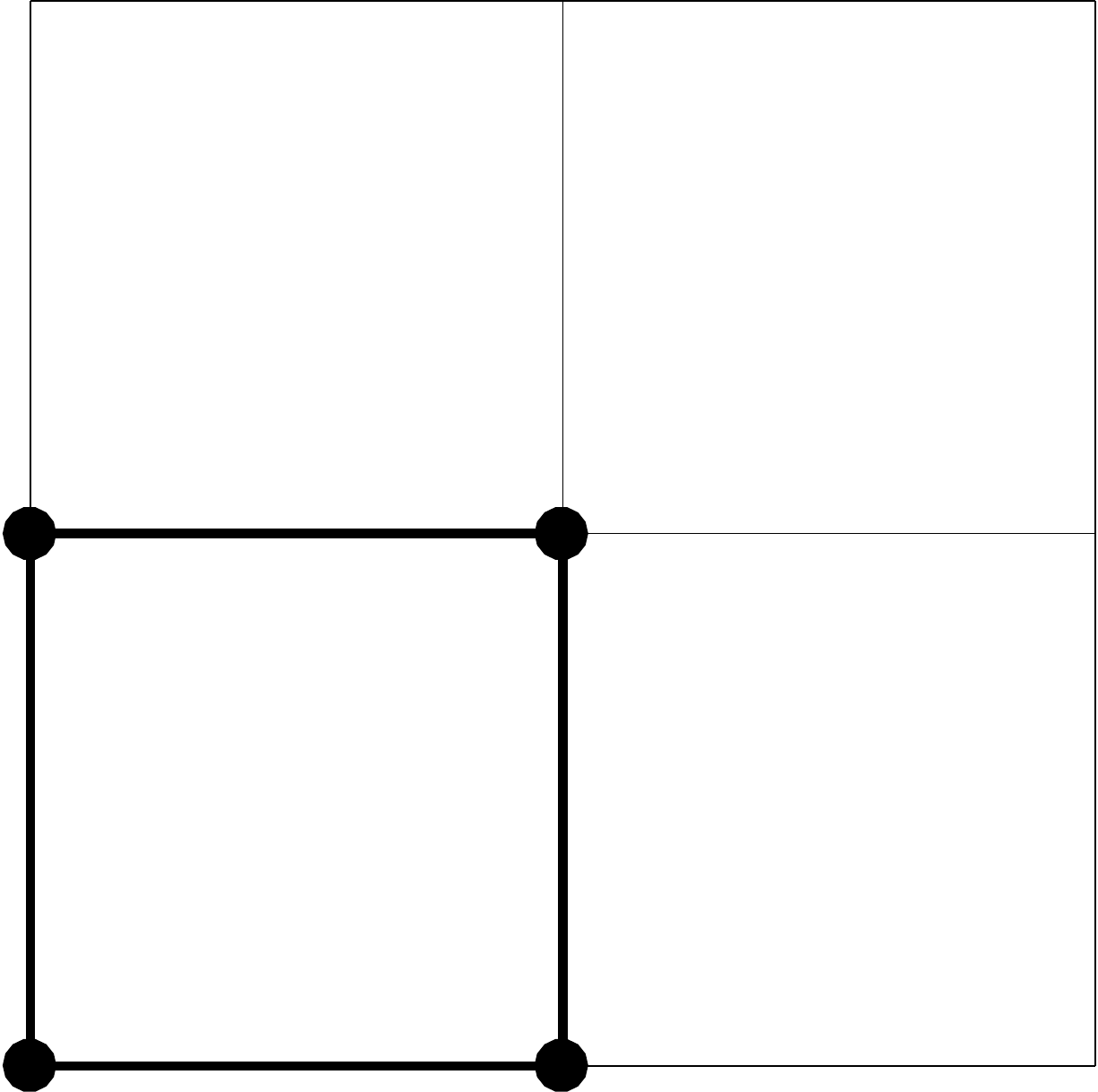} &
      \multicolumn{3}{|c|}{}\\ 
      \hline
    \end{tabular}
    \label{quad2_derivsa}
  \end{center}
\end{table}

\begin{table}[H]
  \begin{center}
    \caption{Partial derivative stencils for bi-cubic  \BB quadrilaterals.}
    \begin{tabular}{ | >{\centering\arraybackslash} m{1.5cm}  | >{\centering\arraybackslash} m{3cm} |  >{\centering\arraybackslash} m{2.3cm} ||
        >{\centering\arraybackslash} m{1.5cm}  | >{\centering\arraybackslash} m{2.0cm} |  >{\centering\arraybackslash} m{2.3cm} |
      }
      \hline
      \multicolumn{6}{|c|}{  \textbf{ First Order Derivatives, $|\balpha| = 1$}}\\ 
      \hline
      Deriv. & Stencil & Evaluation Quadrilateral & Deriv. & Stencil & Evaluation Quadrilateral  \\ 
      \hline
      $\deldel{\xproj}{\xi_1} $&
      \includegraphics[width=2.5cm,height=2.5cm,keepaspectratio]{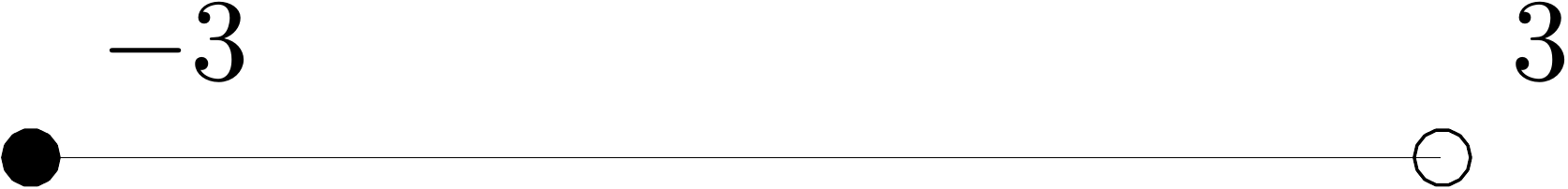} &
      \vspace{5pt}
      \includegraphics[width=2.0cm,height=2.0cm,keepaspectratio]{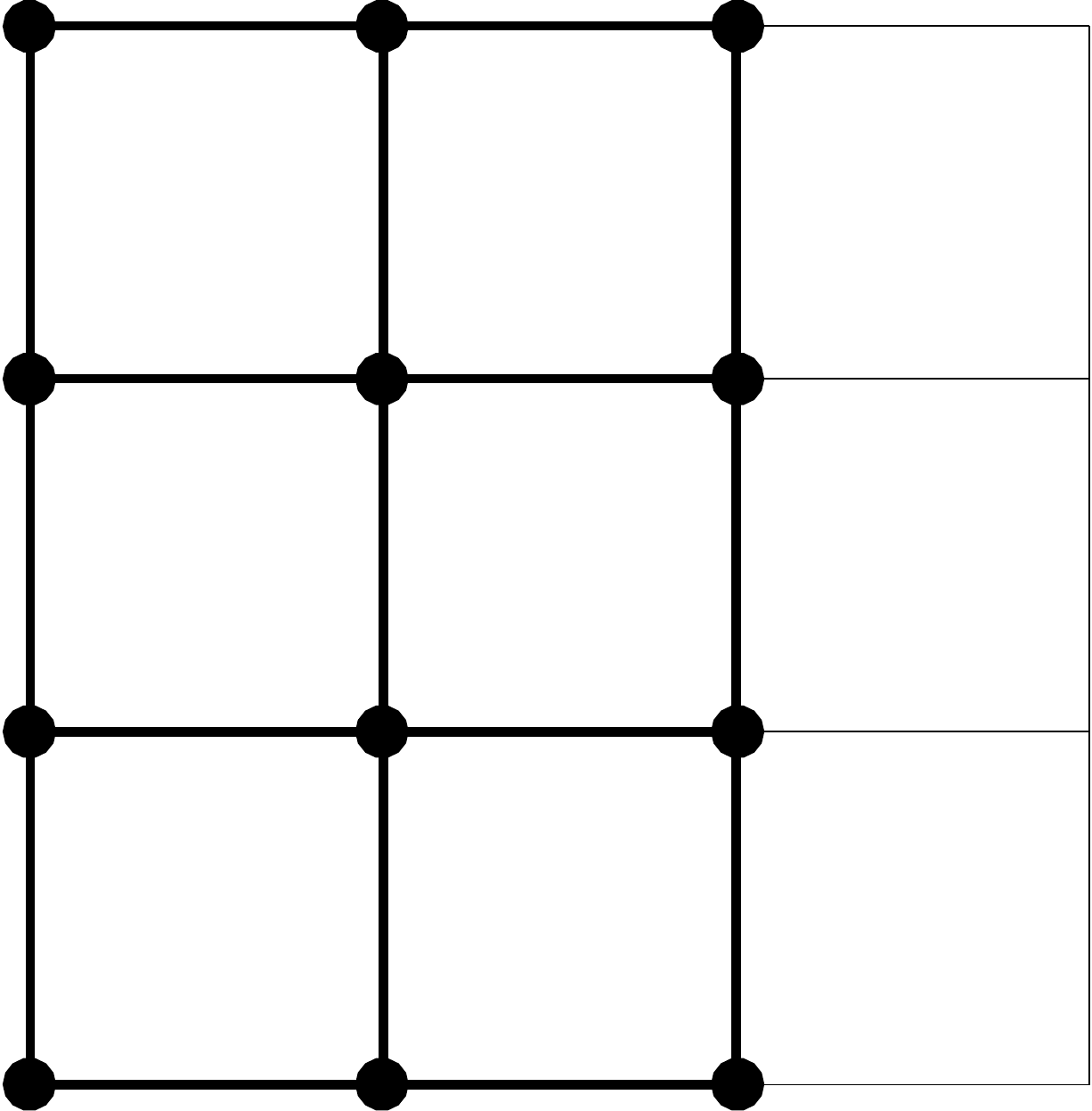} &
      $\deldel{\xproj}{\xi_2} $&
      \includegraphics[width=2.0cm,height=2.0cm,keepaspectratio]{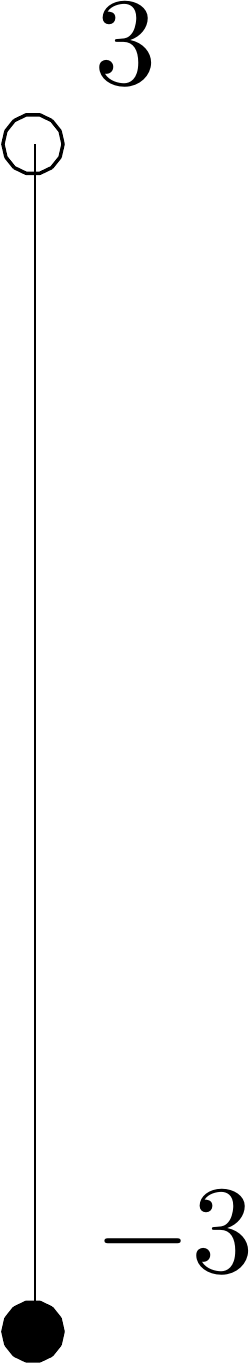} &
      \vspace{5pt}
      \includegraphics[width=2.0cm,height=2.0cm,keepaspectratio]{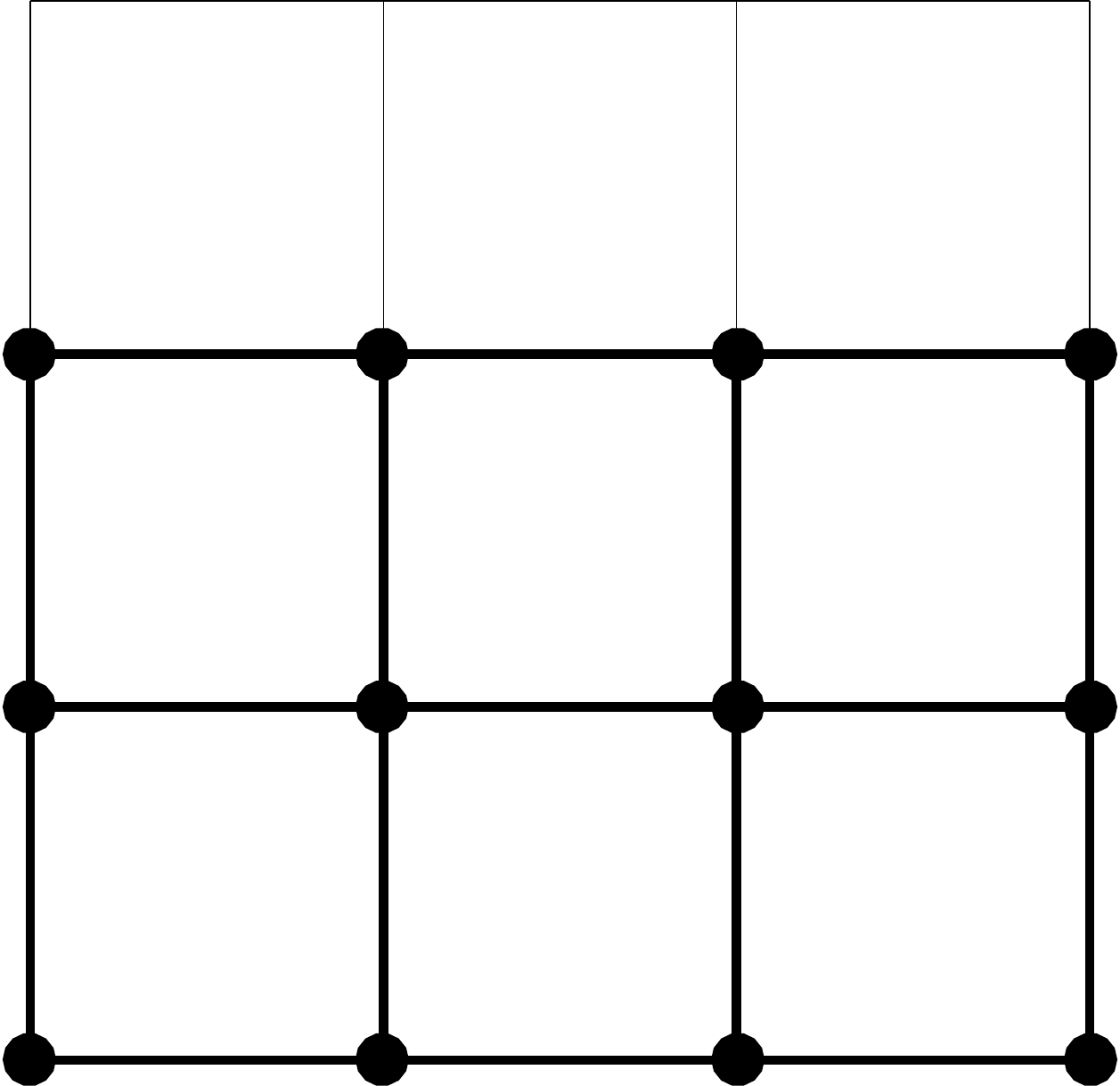} \\
      \hline
      \multicolumn{6}{c}{\vspace{8pt}}\\ 
      \hline
      \multicolumn{6}{|c|}{\textbf{Second Order Derivatives, $|\balpha| = 2$}}\\ 
      \hline
      Deriv. & Stencil & Evaluation Quadrilateral & Deriv. & Stencil & Evaluation Quadrilateral  \\ 
      \hline
      $\dfrac{\partial^2 \xproj}{\partial \xi_1^2} $&
      \includegraphics[width=2.5cm,height=2.5cm,keepaspectratio]{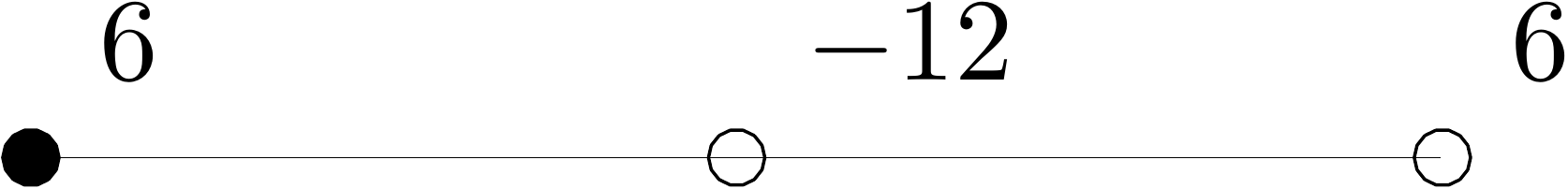} &
      \vspace{5pt}
      \includegraphics[width=2.0cm,height=2.0cm,keepaspectratio]{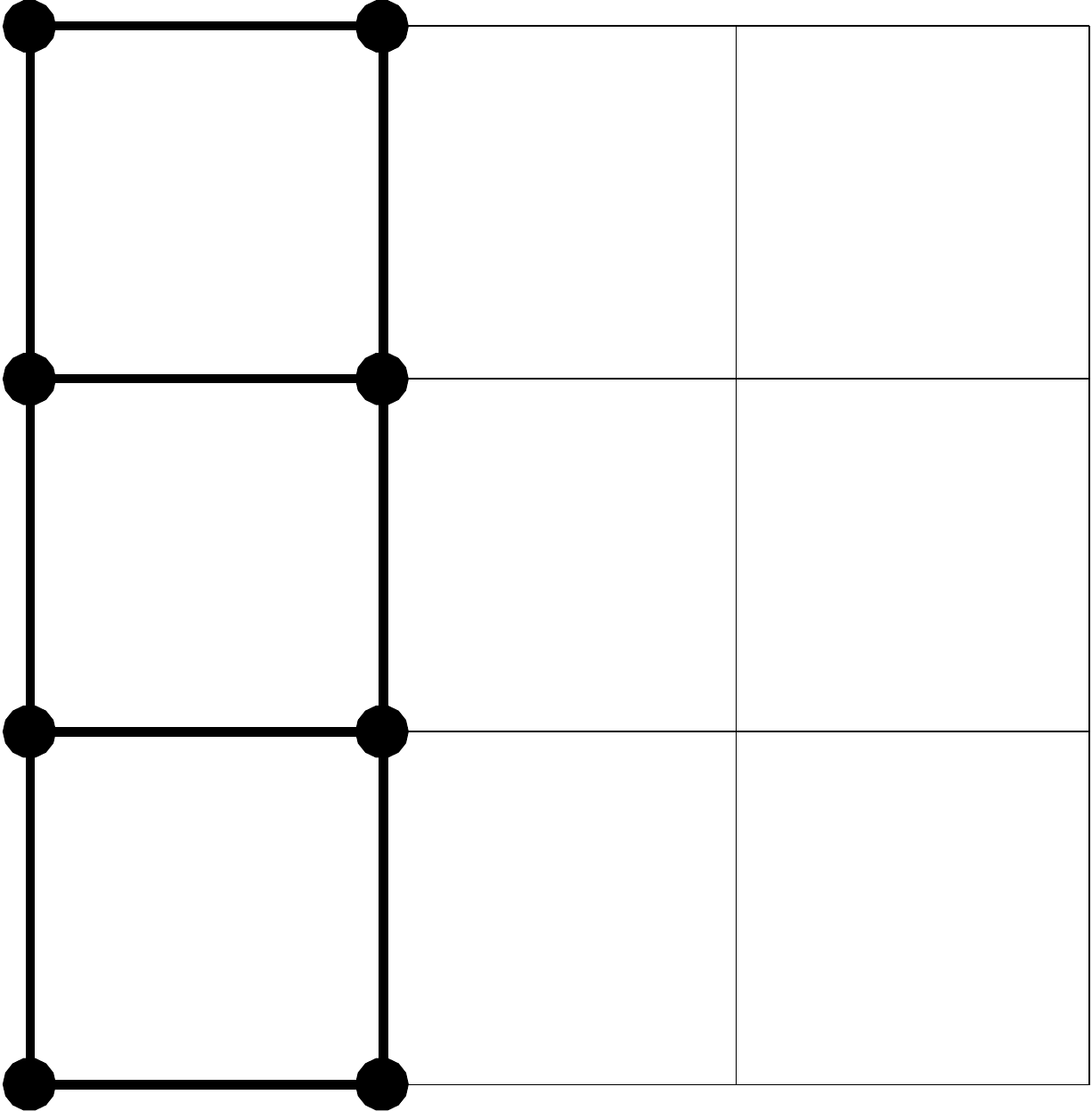} &
      $\dfrac{\partial^2 \xproj}{\partial \xi_2^2} $&
      \includegraphics[width=2.0cm,height=2.0cm,keepaspectratio]{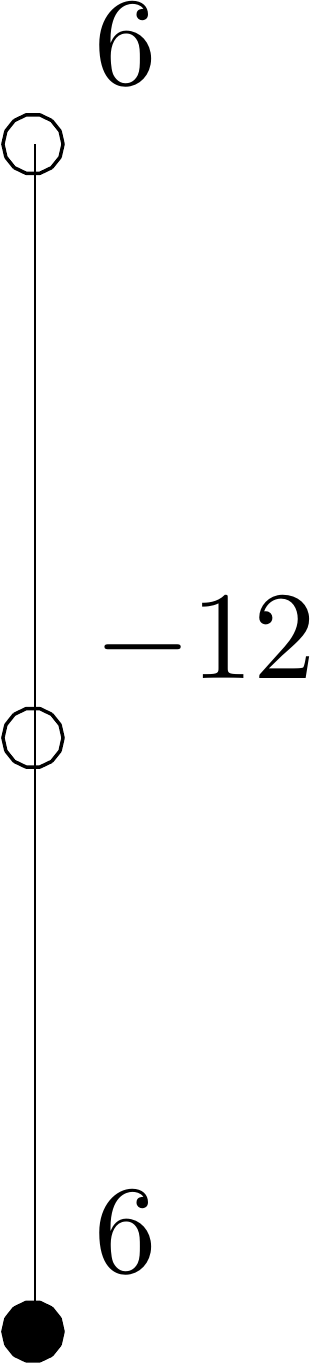} &
      \vspace{5pt}
      \includegraphics[width=2.0cm,height=2.0cm,keepaspectratio]{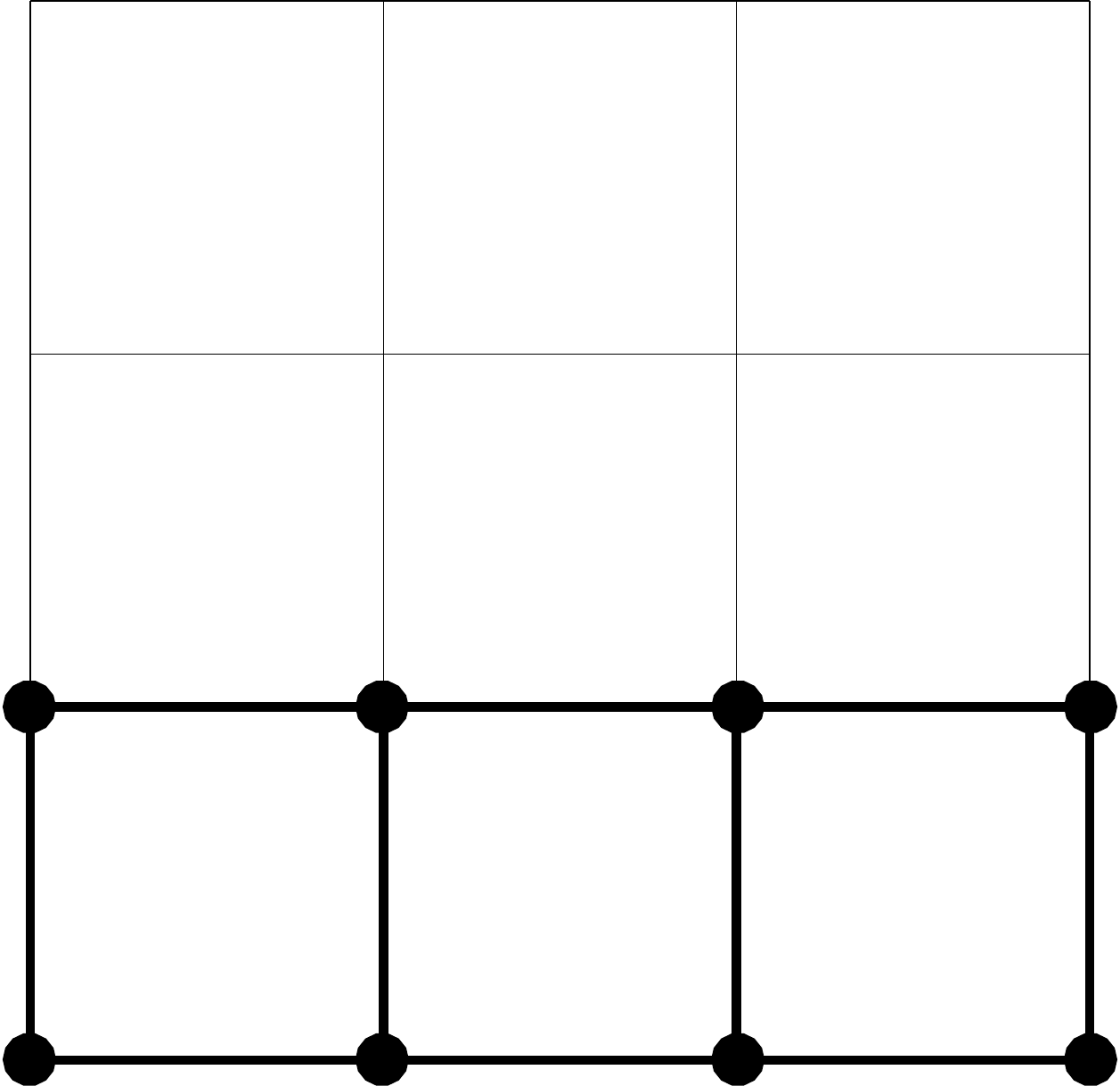} \\ 
      \hline
      $\dfrac{\partial^2 \xproj}{\partial \xi_1 \partial \xi_2} $&
      \vspace{5pt} \includegraphics[width=2.5cm,height=2.5cm,keepaspectratio]{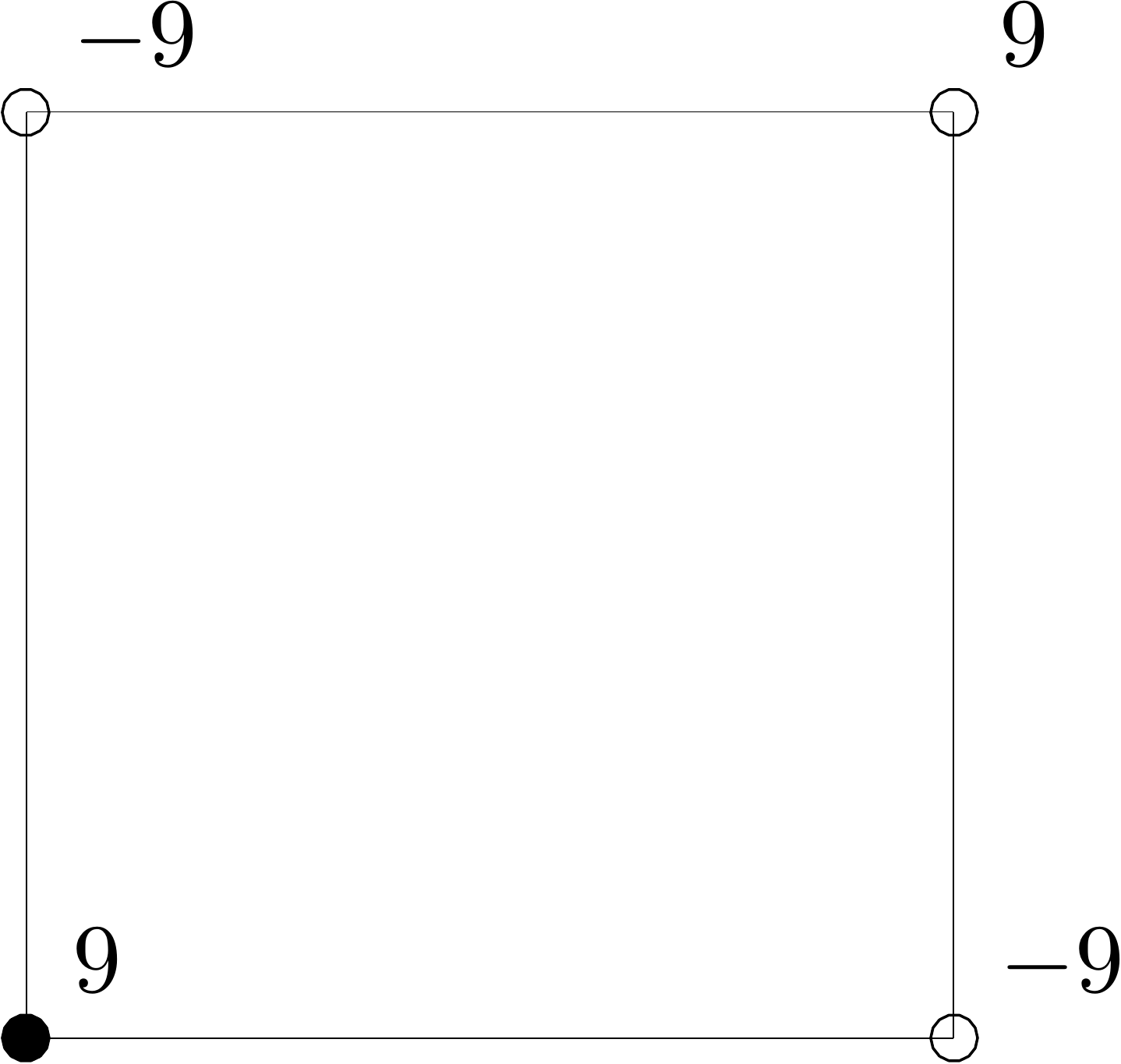} &
      \vspace{5pt}
      \includegraphics[width=2.0cm,height=2.0cm,keepaspectratio]{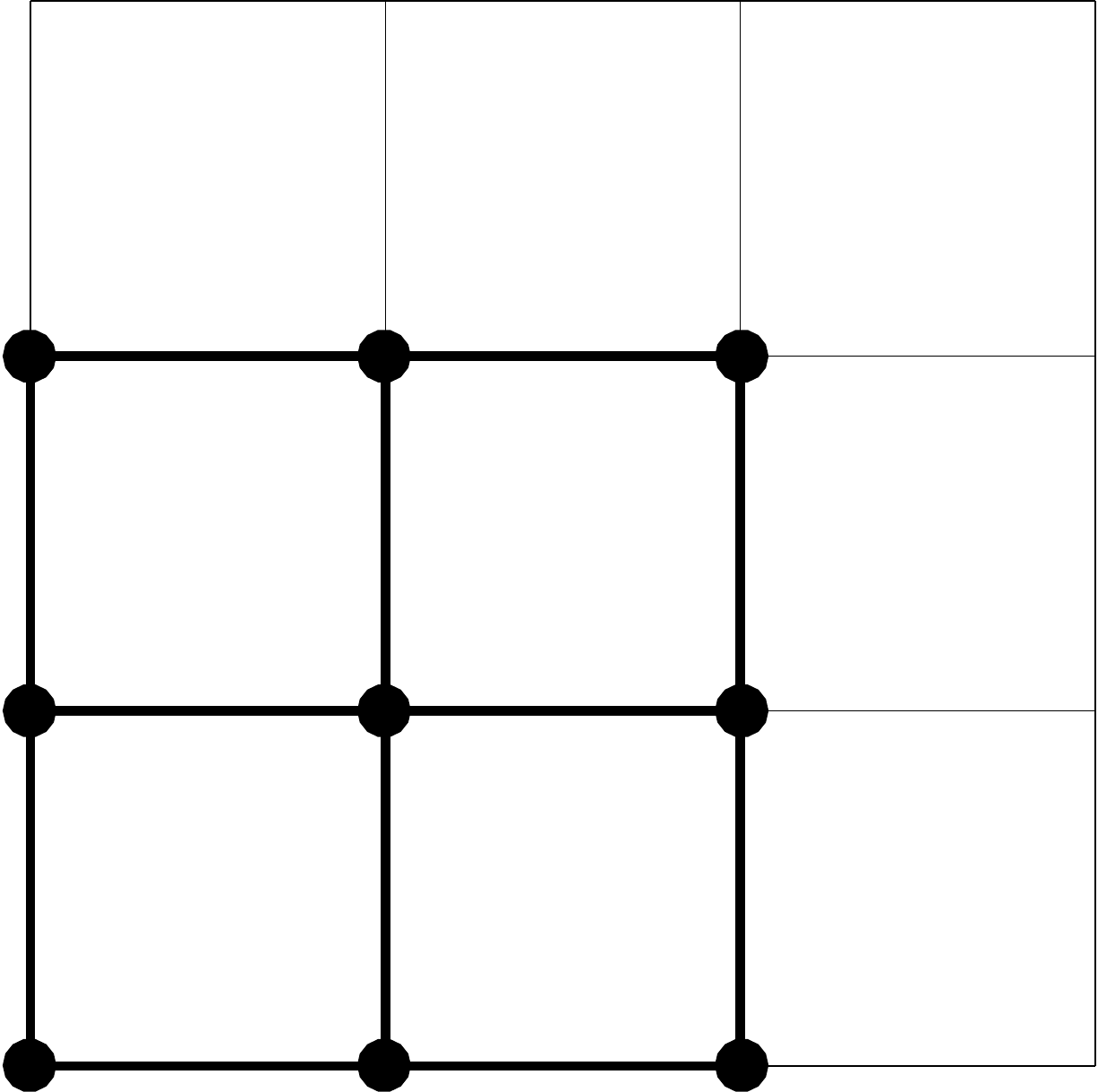} &
      \multicolumn{3}{|c|}{}\\ 
      \hline
      \multicolumn{6}{c}{\vspace{8pt}}\\ 
      \hline
      \multicolumn{6}{|c|}{\textbf{Third Order Derivatives, $|\balpha| = 3$}}\\ 
      \hline
      Deriv. & Stencil & Evaluation Quadrilateral & Deriv. & Stencil & Evaluation Quadrilateral  \\ 
      \hline
      $\dfrac{\partial^3 \xproj}{\partial \xi_1^3 } $&
      \vspace{5pt} \includegraphics[width=2.9cm,height=2.9cm,keepaspectratio]{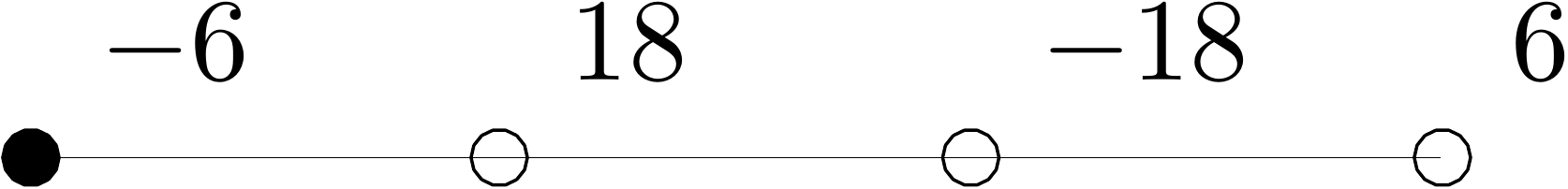} &
      \vspace{5pt}
      \includegraphics[width=2.0cm,height=2.0cm,keepaspectratio]{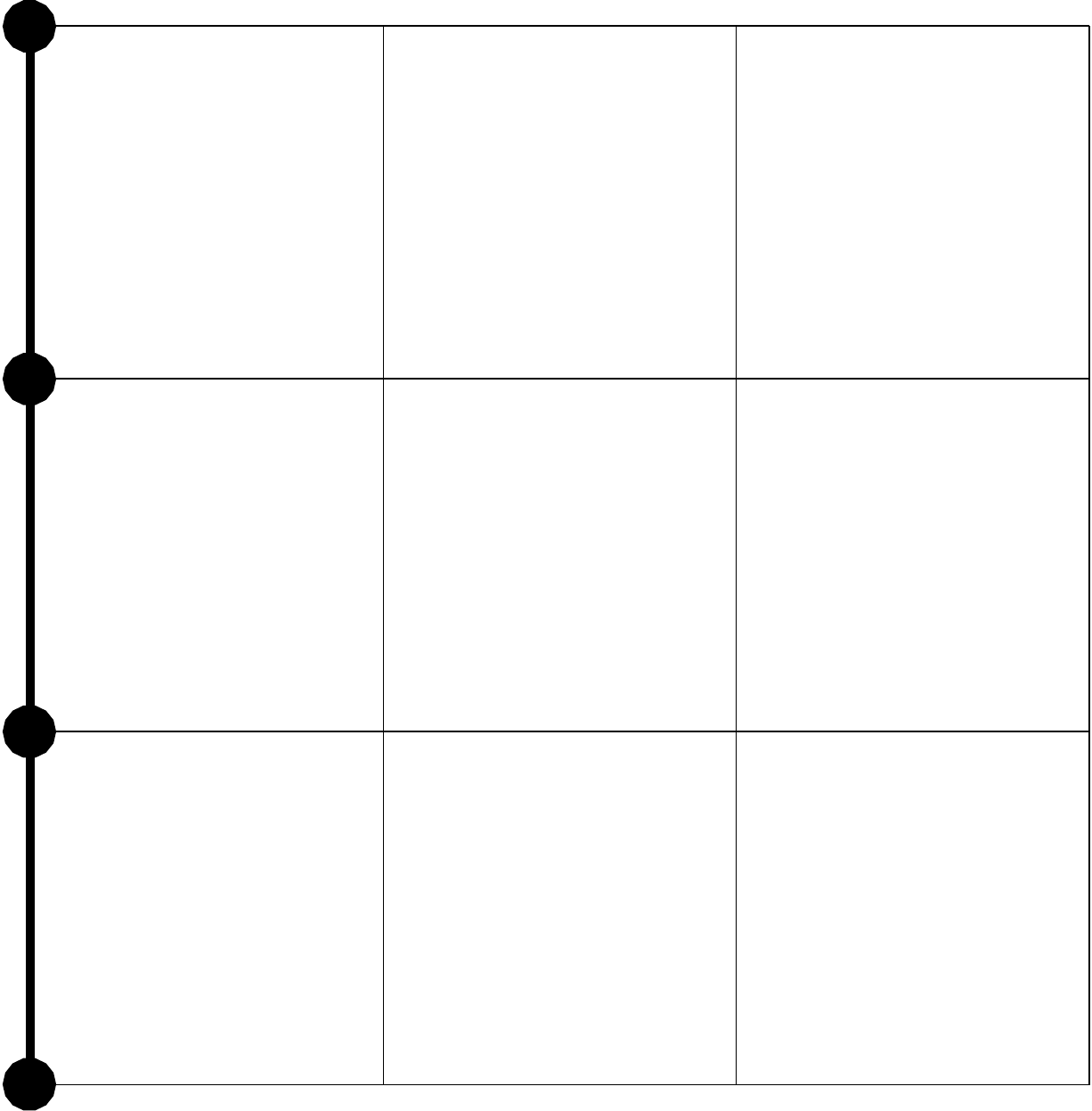} &
      $\dfrac{\partial^3 \xproj}{ \partial \xi_2^3} $&
      \vspace{5pt} \includegraphics[width=2.5cm,height=2.5cm,keepaspectratio]{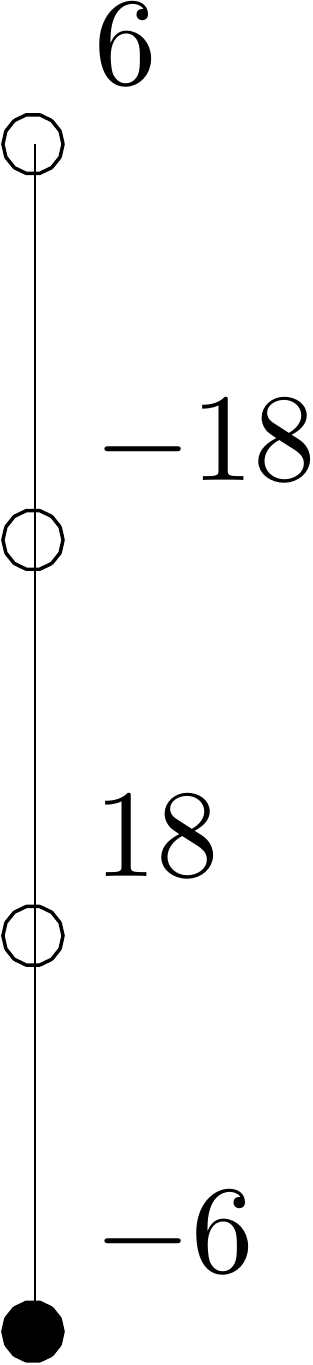} &
      \vspace{5pt}
      \includegraphics[width=2.0cm,height=2.0cm,keepaspectratio]{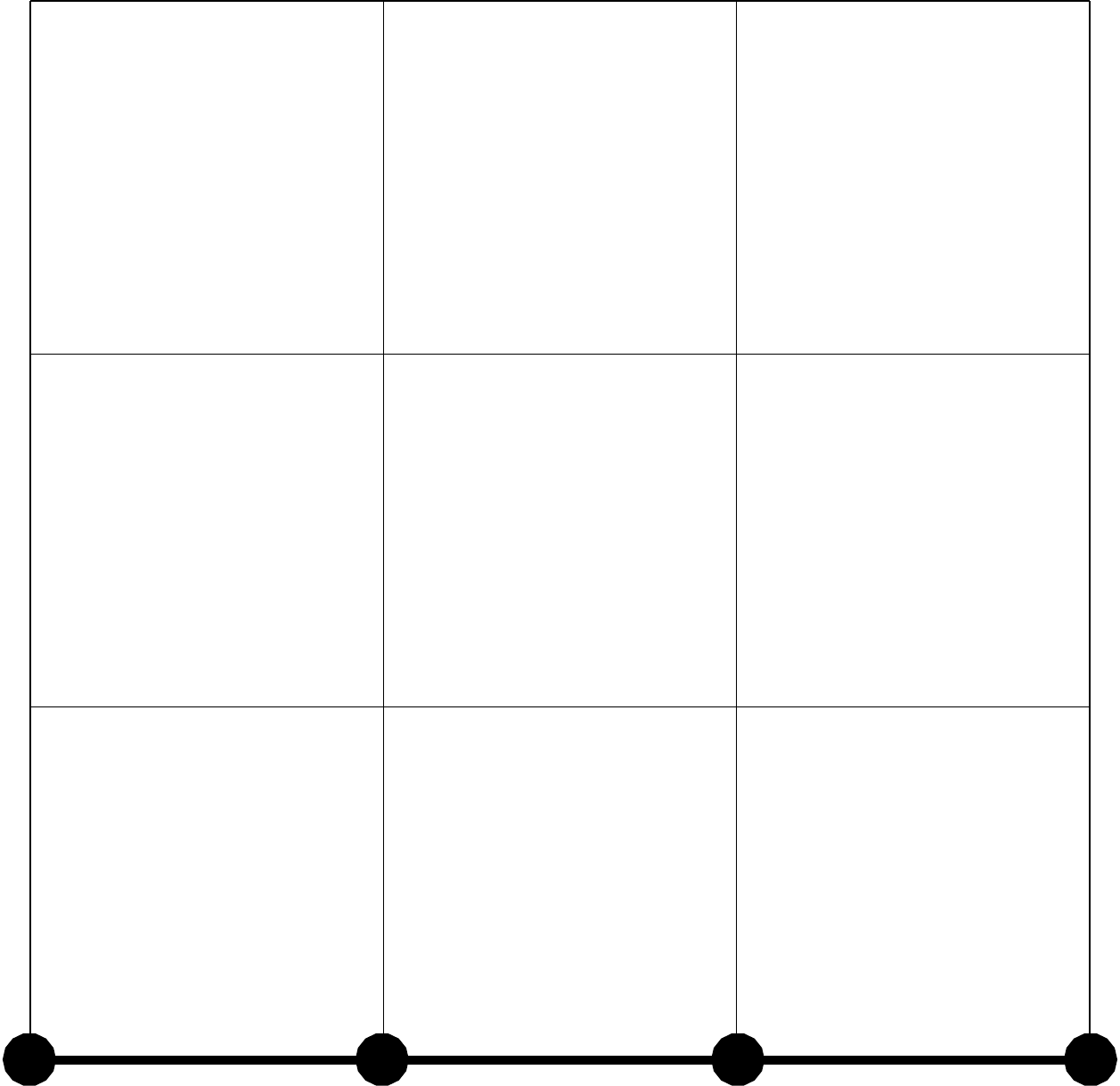}\\ 
      \hline
      $\dfrac{\partial^3 \xproj}{\partial \xi_1^2 \partial \xi_2} $&
      \vspace{5pt} \includegraphics[width=2.9cm,height=2.9cm,keepaspectratio]{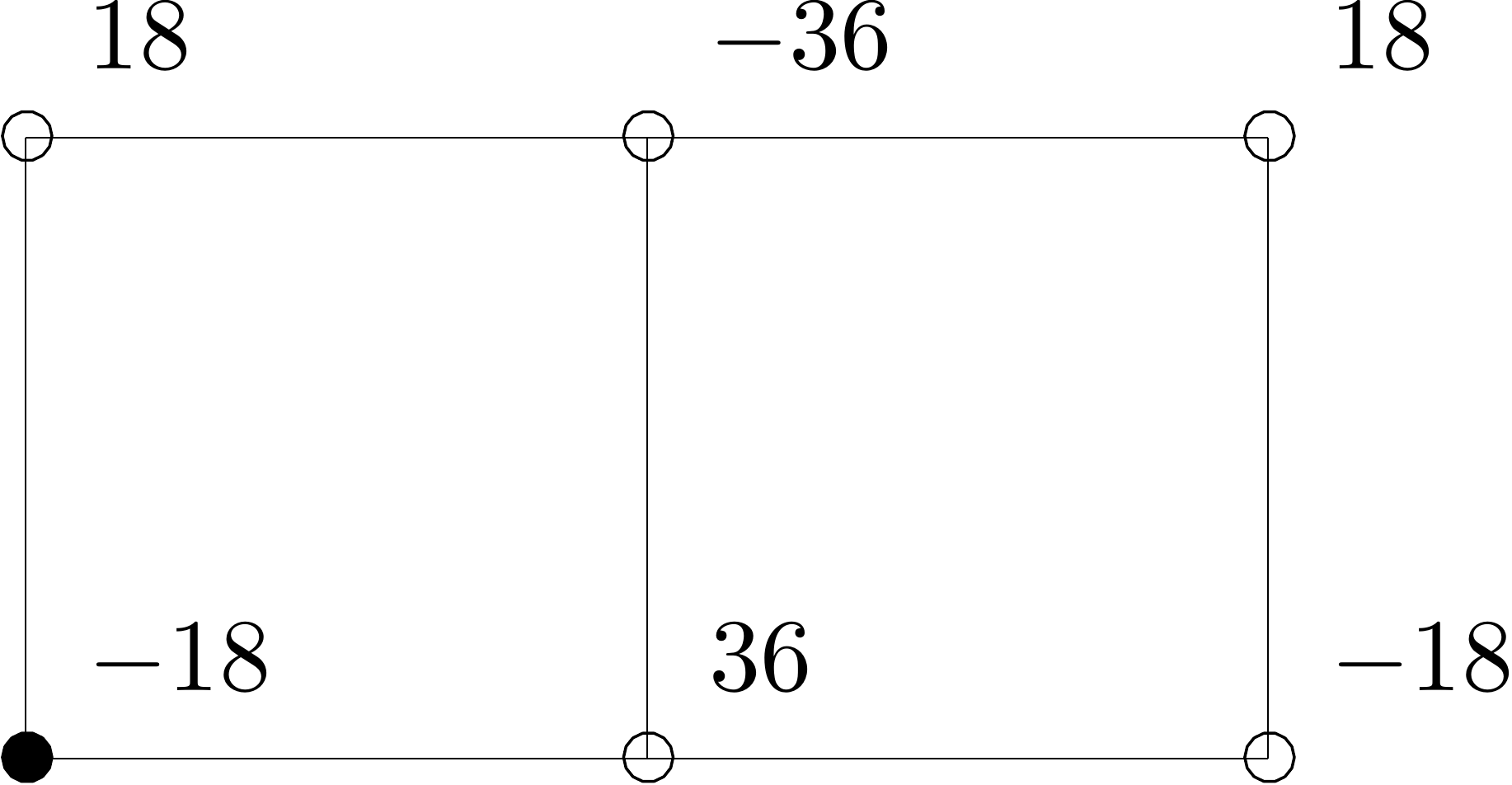} &
      \vspace{5pt}
      \includegraphics[width=2.0cm,height=2.0cm,keepaspectratio]{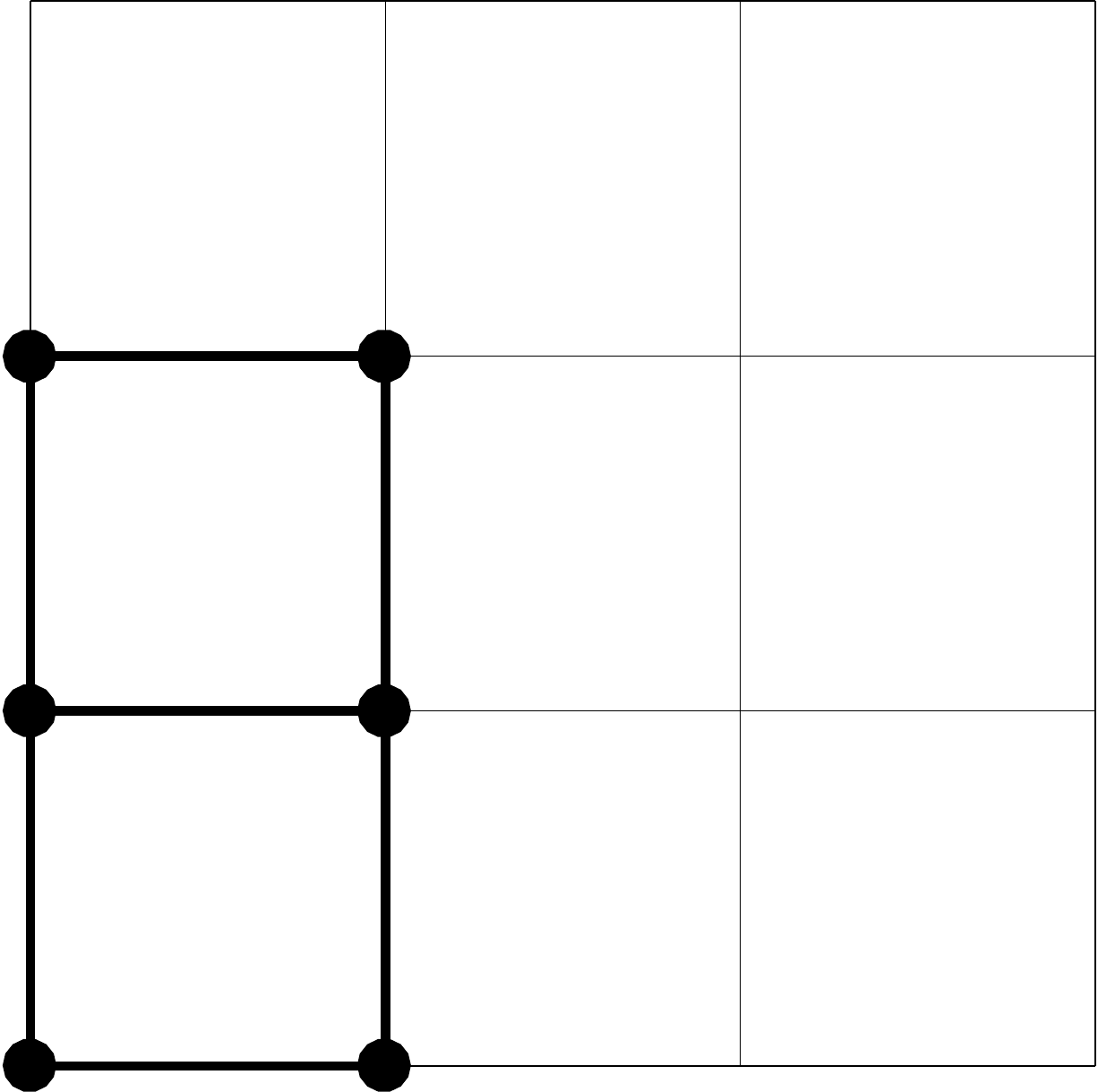} &
      $\dfrac{\partial^3 \xproj}{\partial \xi_1 \partial \xi_2^2} $&
      \vspace{5pt} \includegraphics[width=2.5cm,height=2.5cm,keepaspectratio]{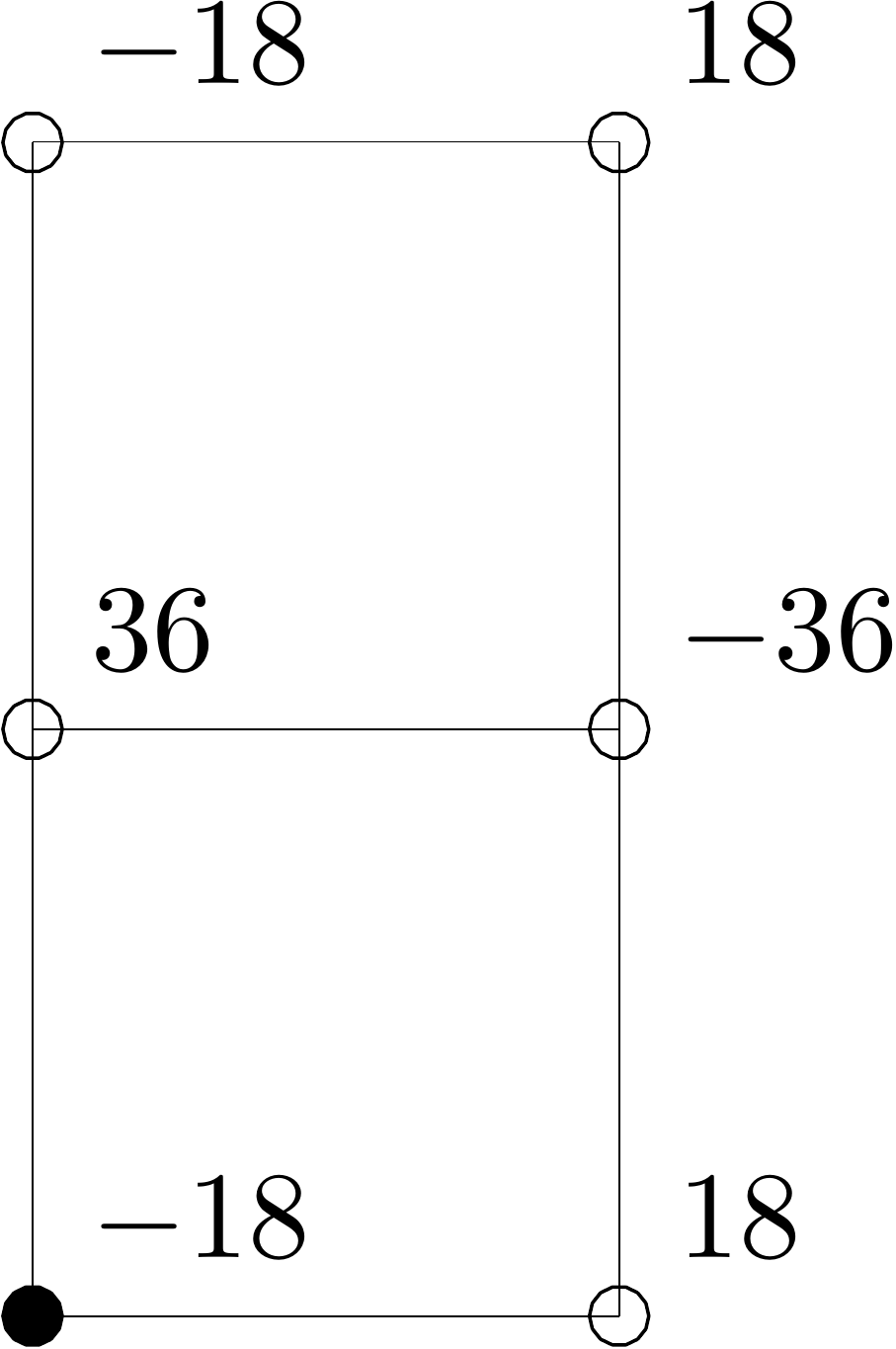} &
      \vspace{5pt}
      \includegraphics[width=2.0cm,height=2.0cm,keepaspectratio]{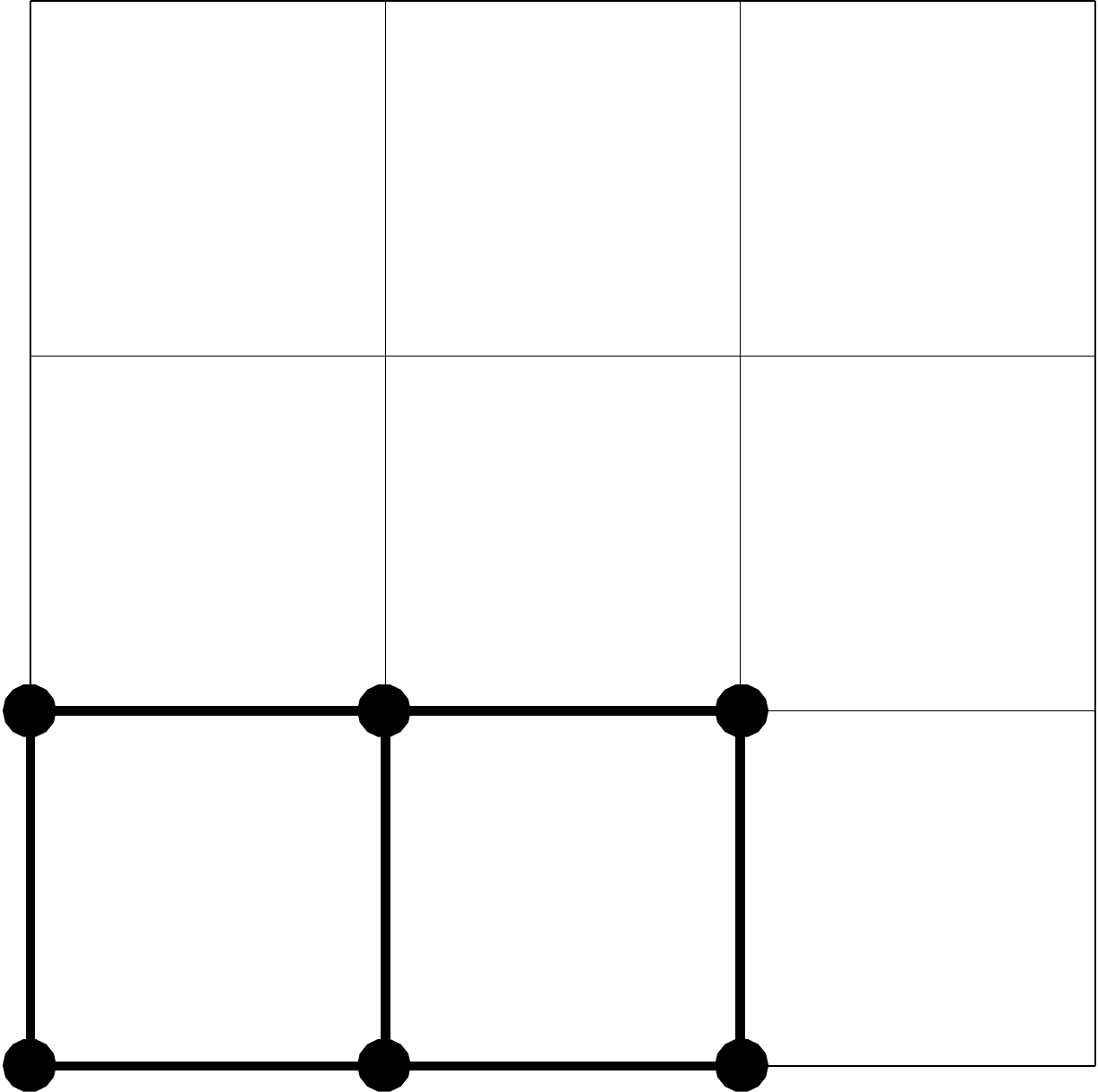}\\ 
      \hline
    \end{tabular}
    \label{quad3_derivsa}
  \end{center}
\end{table}

\end{document}